\documentclass[leqno,11pt]{amsart}
\usepackage{etoolbox}

\makeatletter
\let\old@tocline\@tocline
\let\section@tocline\@tocline

\newcommand{\subsection@dotsep}{4.5}
\newcommand{\subsubsection@dotsep}{4.5}
\patchcmd{\@tocline}
  {\hfil}
  {\nobreak
     \leaders\hbox{$\m@th
        \mkern \subsection@dotsep mu\hbox{.}\mkern \subsection@dotsep mu$}\hfill
     \nobreak}{}{}
\let\subsection@tocline\@tocline
\let\@tocline\old@tocline

\patchcmd{\@tocline}
  {\hfil}
  {\nobreak
     \leaders\hbox{$\m@th
        \mkern \subsubsection@dotsep mu\hbox{.}\mkern \subsubsection@dotsep mu$}\hfill
     \nobreak}{}{}
\let\subsubsection@tocline\@tocline
\let\@tocline\old@tocline

\let\old@l@subsection\l@subsection
\let\old@l@subsubsection\l@subsubsection

\def\@tocwriteb#1#2#3{%
  \begingroup
    \@xp\def\csname #2@tocline\endcsname##1##2##3##4##5##6{%
      \ifnum##1>\c@tocdepth
      \else \sbox\z@{##5\let\indentlabel\@tochangmeasure##6}\fi}%
    \csname l@#2\endcsname{#1{\csname#2name\endcsname}{\@secnumber}{}}%
  \endgroup
  \addcontentsline{toc}{#2}%
    {\protect#1{\csname#2name\endcsname}{\@secnumber}{#3}}}%

\newlength{\@tocsectionindent}
\newlength{\@tocsubsectionindent}
\newlength{\@tocsubsubsectionindent}
\newlength{\@tocsectionnumwidth}
\newlength{\@tocsubsectionnumwidth}
\newlength{\@tocsubsubsectionnumwidth}
\newcommand{\settocsectionnumwidth}[1]{\setlength{\@tocsectionnumwidth}{#1}}
\newcommand{\settocsubsectionnumwidth}[1]{\setlength{\@tocsubsectionnumwidth}{#1}}
\newcommand{\settocsubsubsectionnumwidth}[1]{\setlength{\@tocsubsubsectionnumwidth}{#1}}
\newcommand{\settocsectionindent}[1]{\setlength{\@tocsectionindent}{#1}}
\newcommand{\settocsubsectionindent}[1]{\setlength{\@tocsubsectionindent}{#1}}
\newcommand{\settocsubsubsectionindent}[1]{\setlength{\@tocsubsubsectionindent}{#1}}

\renewcommand{\l@section}{\section@tocline{1}{\@tocsectionvskip}{\@tocsectionindent}{\@tocsectionnumwidth}{\@tocsectionformat}}%
\renewcommand{\l@subsection}{\subsection@tocline{2}{\@tocsubsectionvskip}{\@tocsubsectionindent}{\@tocsubsectionnumwidth}{\@tocsubsectionformat}}%
\renewcommand{\l@subsubsection}{\subsubsection@tocline{3}{\@tocsubsubsectionvskip}{\@tocsubsubsectionindent}{\@tocsubsubsectionnumwidth}{\@tocsubsubsectionformat}}%
\newcommand{\@tocsectionformat}{}
\newcommand{\@tocsubsectionformat}{}
\newcommand{\@tocsubsubsectionformat}{}
\expandafter\def\csname toc@1format\endcsname{\@tocsectionformat}
\expandafter\def\csname toc@2format\endcsname{\@tocsubsectionformat}
\expandafter\def\csname toc@3format\endcsname{\@tocsubsubsectionformat}
\newcommand{\settocsectionformat}[1]{\renewcommand{\@tocsectionformat}{#1}}
\newcommand{\settocsubsectionformat}[1]{\renewcommand{\@tocsubsectionformat}{#1}}
\newcommand{\settocsubsubsectionformat}[1]{\renewcommand{\@tocsubsubsectionformat}{#1}}
\newlength{\@tocsectionvskip}
\newcommand{\settocsectionvskip}[1]{\setlength{\@tocsectionvskip}{#1}}
\newlength{\@tocsubsectionvskip}
\newcommand{\settocsubsectionvskip}[1]{\setlength{\@tocsubsectionvskip}{#1}}
\newlength{\@tocsubsubsectionvskip}
\newcommand{\settocsubsubsectionvskip}[1]{\setlength{\@tocsubsubsectionvskip}{#1}}

\renewcommand{\tocsection}[3]{%
  \@ifnotempty{#2}{\makebox[\@tocsectionnumwidth][l]{\ignorespaces#1 #2.\quad}}#3}
\patchcmd{\tocsubsection}{\indentlabel}{\makebox[\@tocsubsectionnumwidth][l]}{}{}
\patchcmd{\tocsubsubsection}{\indentlabel}{\makebox[\@tocsubsubsectionnumwidth][l]}{}{}

\newcommand{\@sectypepnumformat}{}
\renewcommand{\contentsline}[1]{%
  \expandafter\let\expandafter\@sectypepnumformat\csname @toc#1pnumformat\endcsname%
  \csname l@#1\endcsname}
\newcommand{\@tocsectionpnumformat}{}
\newcommand{\@tocsubsectionpnumformat}{}
\newcommand{\@tocsubsubsectionpnumformat}{}
\newcommand{\setsectionpnumformat}[1]{\renewcommand{\@tocsectionpnumformat}{#1}}
\newcommand{\setsubsectionpnumformat}[1]{\renewcommand{\@tocsubsectionpnumformat}{#1}}
\newcommand{\setsubsubsectionpnumformat}[1]{\renewcommand{\@tocsubsubsectionpnumformat}{#1}}
\renewcommand{\@tocpagenum}[1]{%
  \hfill {\mdseries\@sectypepnumformat #1}}

\let\oldappendix\appendix
\renewcommand{\appendix}{%
  \leavevmode\oldappendix%
  \addtocontents{toc}{%
    \protect\settowidth{\protect\@tocsectionnumwidth}{\protect\@tocsectionformat\sectionname\space}%
    \protect\addtolength{\protect\@tocsectionnumwidth}{2em}}%
}
\makeatother

\makeatletter
\settocsectionnumwidth{2em}
\settocsubsectionnumwidth{2.5em}
\settocsubsubsectionnumwidth{3em}
\settocsectionindent{1pc}%
\settocsubsectionindent{\dimexpr\@tocsectionindent+\@tocsectionnumwidth}%
\settocsubsubsectionindent{\dimexpr\@tocsubsectionindent+\@tocsubsectionnumwidth}%
\makeatother

\settocsectionvskip{2pt}
\settocsubsectionvskip{0pt}
\settocsubsubsectionvskip{0pt}

\settocsectionformat{\bfseries}
\settocsubsectionformat{\mdseries}
\settocsubsubsectionformat{\mdseries}
\setsectionpnumformat{\bfseries}
\setsubsectionpnumformat{\mdseries}
\setsubsubsectionpnumformat{\mdseries}

\let\oldtableofcontents\tableofcontents
\renewcommand{\tableofcontents}{%
  \vspace*{-5\linespacing}
  \oldtableofcontents}

\usepackage{amsmath}
\usepackage{amssymb}
\usepackage{amsthm}
\usepackage{mathtools}
\usepackage{thmtools}
\usepackage{hhline}
\usepackage{tikz-cd}
\usepackage{enumitem}
\setlist[enumerate]{label={\normalfont (\roman*)}, itemsep=1ex}

\setlist[itemize]{leftmargin=2em}
\usepackage[utf8]{inputenc}
\usepackage[T1]{fontenc}
\usepackage{lmodern}
\usepackage[babel]{microtype}
\usepackage[english]{babel}
\usepackage{relsize}
\usepackage{calc}
\usepackage{multicol}
\usepackage{float}

\usepackage[all]{xy}
\usepackage{subfig}
\usepackage{footnote}
\makesavenoteenv{tabular}
\usepackage{tikz}
\usetikzlibrary{patterns}
\usepackage{tikz-3dplot}

\usepackage{seqsplit}
\usepackage{xstring}

\usepackage{svg}

\usepackage{graphicx}
\usepackage{subcaption}

\counterwithin{figure}{section}

\usepackage{geometry}
\usepackage{xcolor} 
\colorlet{darkishRed}{red!60!black}
\colorlet{darkishBlue}{blue!60!black}
\colorlet{darkishGreen}{green!60!black}
\colorlet{darkblue}{blue!70!black}
\colorlet{darkishViolet}{violet}
\usepackage[linktoc=all]{hyperref}
\hypersetup{
	colorlinks,
    linkcolor={darkishBlue},
	citecolor={darkishGreen},
	urlcolor={darkishBlue}
}
\usepackage[nameinlink, capitalise, noabbrev]{cleveref}
\crefformat{enumi}{#2#1#3}
\crefformat{equation}{#2(#1)#3}
\crefname{equation}{}{}
\crefname{mainresult}{Theorem}{Theorems}
\let\setminus=\smallsetminus

\newcommand{\se}{\subseteq}
\newcommand{\sm}{\setminus}

\renewcommand{\leq}{\leqslant}
\renewcommand{\geq}{\geqslant}
\renewcommand{\ge}{\geq}
\renewcommand{\le}{\leq}
\newcommand{\rest}{\!\restriction}

\newcommand{\allref}[1]{\nameref{#1}~(\ref{#1})}

\newtheorem{theorem}{Theorem}[section]
\newtheorem{thm}[theorem]{Theorem}

\newtheorem{proposition}[theorem]{Proposition}

\newtheorem{corollary}[theorem]{Corollary}
\newtheorem{cor}[theorem]{Corollary}
\newtheorem{lemma}[theorem]{Lemma}
\newtheorem{lem}[theorem]{Lemma}
\newtheorem{keylemma}[theorem]{Key Lemma}
\newtheorem{setting}[theorem]{Setting}
\newtheorem{observation}[theorem]{Observation}
\newtheorem{obs}[theorem]{Observation}

\newtheorem{conjecture}[theorem]{Conjecture}

\newtheorem{claim}{Claim}[thm]

\theoremstyle{definition}

\newtheorem{eg}[theorem]{Example}
\newtheorem{definition}[theorem]{Definition}

\newtheorem{construction}[theorem]{Construction}

\theoremstyle{remark}

\newtheorem{convention}[theorem]{Convention}

\newtheorem*{claim*}{Claim}
\newtheorem{remark}[theorem]{Remark}
\newtheorem{rem}[theorem]{Remark}

\crefname{claim}{Claim}{Claims}

\newenvironment{cproof}{\noindent\textit{Proof of the claim.}}{\hfill\ensuremath{\diamondsuit}\medskip}
\usepackage{etoolbox}

\definecolor{turquoise}{RGB}{40,140,128}

\newcommand{\COMMENT}[1]{{}}
\definecolor{cMaroon}{HTML}{93152a}
\newcommand{\defnMath}[1]{{\color{darkishRed}{#1}}}
\newcommand{\defn}[1]{{\color{darkishRed}{\emph{#1}}}}
\newcommand{\casen}[1]{{\color{darkishViolet}{\emph{#1}}}}

\let\eps=\varepsilon
\let\rho=\varrho
\let\phi=\varphi

\let\preceq=\preccurlyeq
\def\N{\mathbb N}
\def\R{\mathbb R}
\def\Z{\mathbb Z}

\makeatletter

\def\calCommandfactory#1{%
  \expandafter\def\csname c#1\endcsname{\mathcal{#1}}}
\def\frakCommandfactory#1{%
  \expandafter\def\csname frak#1\endcsname{\mathfrak{#1}}}
   
\newcounter{ctr}
\loop
  \stepcounter{ctr}
  \edef\X{\@Alph\c@ctr}
  \expandafter\calCommandfactory\X
  \expandafter\frakCommandfactory\X
\ifnum\thectr<26
\repeat

\makeatletter
\patchcmd{\@tocline}
  {\hfil}
  {\leaders\hbox{\,.\,}\hfil}
  {}{}
\makeatother

\newcounter{mylabelcounter}

\makeatletter
\newcommand{\labelText}[2]{%
#1\refstepcounter{mylabelcounter}%
\immediate\write\@auxout{%
  \string\newlabel{#2}{{1}{\thepage}{{\unexpanded{#1}}}{mylabelcounter.\number\value{mylabelcounter}}{}}%
}%
}

\newlist{defenum}{enumerate}{1}
\setlist[defenum]{label={\textnormal (\upshape\thedefinition.\arabic*)}}

\newcommand{\cross}{\text{\textnormal{cr}}}

\newcommand{\sep}{\text{\textnormal{sep}}}

\newcommand{\triv}{_\text{\textnormal{t}}}

\crefname{lem}{Lemma}{Lemmas}

\usepackage{epsfig}
\usepackage{graphicx}
\usepackage{amsbsy}
\usepackage{amsmath}
\usepackage{amsfonts}
\usepackage{amssymb}
\usepackage{textcomp}
\usepackage{hyperref}
\usepackage{aliascnt}
\usepackage{mathrsfs}
\usepackage{comment}

\newcommand{\mcm}[3]{\newcommand{#1}[#2]{{\ensuremath{#3}}}} 

\mcm{\tuple}{1}{\langle #1 \rangle}
\mcm{\name}{1}{\ulcorner #1 \urcorner}
\mcm{\Nbb}{0}{\mathbb{N}}
\mcm{\Zbb}{0}{\mathbb{Z}}
\mcm{\Rbb}{0}{\mathbb{R}}
\mcm{\Cbb}{0}{\mathbb{C}}
\mcm{\Qbb}{0}{\mathbb{Q}}
\mcm{\Acal}{0}{\mathcal A}
\mcm{\Bcal}{0}{\mathcal B}
\mcm{\Ccal}{0}{\mathcal C}
\mcm{\Dcal}{0}{\mathcal D}
\mcm{\Ecal}{0}{\mathcal E}
\mcm{\Fcal}{0}{\mathcal F}
\mcm{\Gcal}{0}{\mathcal G}
\mcm{\Hcal}{0}{\mathcal H}
\mcm{\Ical}{0}{\mathcal I}
\mcm{\Jcal}{0}{\mathcal J}
\mcm{\Kcal}{0}{\mathcal K}
\mcm{\Lcal}{0}{\mathcal L}
\mcm{\Mcal}{0}{\mathcal M}
\mcm{\Ncal}{0}{\mathcal N}
\mcm{\Ocal}{0}{{\mathcal O}}
\mcm{\Pcal}{0}{{\mathcal P}}
\mcm{\Qcal}{0}{{\mathcal Q}}
\mcm{\Rcal}{0}{{\mathcal R}}
\mcm{\Scal}{0}{{\mathcal S}}
\mcm{\Tcal}{0}{{\mathcal T}}
\mcm{\Ucal}{0}{{\mathcal U}}
\mcm{\Vcal}{0}{{\mathcal V}}
\mcm{\Wcal}{0}{{\mathcal W}}
\mcm{\Xcal}{0}{{\mathcal X}}
\mcm{\Ycal}{0}{{\mathcal Y}}
\mcm{\Zcal}{0}{{\mathcal Z}}
\mcm{\Mfrak}{0}{\mathfrak M}

\mcm{\restric}{0}{\upharpoonright}
\mcm{\upset}{0}{\uparrow}
\mcm{\onto}{0}{\twoheadrightarrow}
\mcm{\smallNbb}{0}{{\small \mathbb{N}}}
\DeclareMathOperator{\preop}{op}
\mcm{\op}{0}{^{\preop}}

{\begin{array}{c}
\setlength{\unitlength}{1em}}%
{\end{array}}

\title[An approximate cactus-decomposition for vertex-connectivity]{An approximate cactus-decomposition for vertex-connectivity}

\author[Romain Bourneuf]{Romain Bourneuf${}^{*}$}

\author[Johannes Carmesin]{Johannes Carmesin${}^{**}$}

\author[Jan Kurkofka]{Jan Kurkofka${}^{\includegraphics[height=.7\baselineskip]{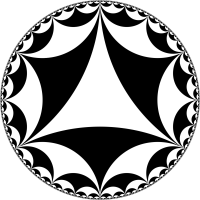}}$}

\author[Tim Planken]{Tim Planken${}^{\dagger}$}

\thanks{* Univ. Bordeaux, CNRS, Bordeaux INP, LaBRI, UMR 5800, F-33400 Talence, France.}

\thanks{** TU Freiberg, funded by the Deutsche Forschungsgemeinschaft (DFG, German Research Foundation) -- 546892829}

\thanks{$\includegraphics[height=.7\baselineskip]{Figures/MiniFarey.png}$ University of Hamburg, funded by the Deutsche Forschungsgemeinschaft (DFG, German Research Foundation) -- 566118291}

\thanks{$\dagger$ Univ. Grenoble Alpes, CNRS, Grenoble INP, G-SCOP, 38000 Grenoble, France.}

\begin{document}
\thispagestyle{empty}

\begin{abstract}
The classical cactus theorem of Dinitz, Karzanov and Lomonosov from 1976 states that the minimum edge-cuts of a {connected} graph are displayed by a cactus.
For vertex-separations, Tutte's theorem gives an exact cactus-decomposition for connectivity $2$, but exact cactus-decompositions fail for higher connectivity.
Extending earlier approximate cactus-decomposition results for connectivities~$3$ and~$4$ by Carmesin and Kurkofka [FOCS~2023] and Kurkofka and Planken [SODA~2026] to all connectivities, we prove that for every $k\in\Nbb$ there is a constant $d=d(k)$ such that every $k$-connected graph~$G$ admits a cactus-decomposition displaying $k$-separations such that every $k$-separation of $G$ has edit-distance at most $d$ to a displayed separation.

Our proof develops tools that may be useful beyond this result. The first is a three-dimensional version of the classical corner-diagram method, which we use to show that triples of crossing $k$-separations have a controlled two-dimensional behaviour. The second is a covering-space approach for studying local crossing patterns from a global perspective. 
Third, we advance the theory of flowers, which originated in matroid theory. 
\end{abstract}

\topskip0pt
\vspace*{\fill}
\maketitle

\vspace*{\fill}


\vspace{10pt}

\noindent\centerline{\textbf{AI methodology:} \cref{sec:AIstatement}.}\bigskip\bigskip

\thispagestyle{empty}

\clearpage

\begingroup
  \pagestyle{empty}
  \tableofcontents

  \clearpage
\endgroup

\section{Introduction}
\setcounter{page}{1}

\enlargethispage{2\baselineskip}

The classical cactus theorem of Dinitz, Karzanov, and Lomonosov from 1976 \cite{EdgeTutte} says that the set of minimum edge-cuts of every connected graph has the structure of a cactus\footnote{A \defn{cactus} is a connected graph in which any two cycles share at most one vertex. A cactus is \defn{strong} if its cycles are pairwise vertex-disjoint.}.
The DKL cactus representation can be computed in near-linear time \cite{HHS24b}, and hence can be used to enumerate all minimum cuts efficiently \cite{KP09,LST21}. 
It is a key ingredient in algorithms for connectivity augmentation \cite{G16,NGM97,CLP22} and for maintaining minimum cuts under edge insertions \cite{H97,DW98,GHT18}.
Vertex-connectivity encompasses edge-connectivity as a special case via the line-graph construction.
However, the direct extension of the cactus theorem to arbitrary vertex-connectivity fails. 
The connection between vertex-connectivity and tree-decompositions makes this question particularly relevant for algorithmic applications.

The quest for such an analogue goes back at least to Tutte's seminal work from 1966~\cite{TutteCon}, even predating the DKL theorem.
The search for a cactus-type representation of minimum vertex-separations\footnote{A \defn{(vertex-)separation} of a graph $G$ is a pair $(A,B)$
of subsets of $V(G)$ such that $A\cup B=V(G)$ and no edge joins
$A\setminus B$ to $B\setminus A$.
Its \defn{separator} is $A\cap B$, and its \defn{order} is $|A\cap B|$.
A \defn{$k$-separation} is a separation of order $k$.
A separation is \defn{proper} if both $A\setminus B$ and
$B\setminus A$ are nonempty, and \defn{improper} otherwise.} motivated the work of Cohen, Di Battista, Kanevsky, and Tamassia~\cite{ReinventingWheel} and subsequent structural and algorithmic investigations~\cite{KSS08,PettieTalk,StructureMinimumVertexCuts}.
The broader problem of decomposing a $k$-connected\footnote{A graph is \defn{$k$-connected} if it contains more than $k$ vertices and has no proper $(k-1)$-separation.} graph into $(k+1)$-connected regions 
has been intensively studied in the last decade \cite{confing,Grohe16,RefiningToT,ProfilesNew,AlbrechtsenRefiningToTinASS}. 
These results reveal substantially richer interactions between minimum vertex-separations than between minimum edge-cuts, making the formulation of a suitable global representation a challenge in its own right.
In this paper, we establish an approximate cactus-decomposition theorem for arbitrary vertex-connectivity.
\smallskip

\begin{figure}[ht]
    \centering
    \includegraphics[height=8.5\baselineskip]{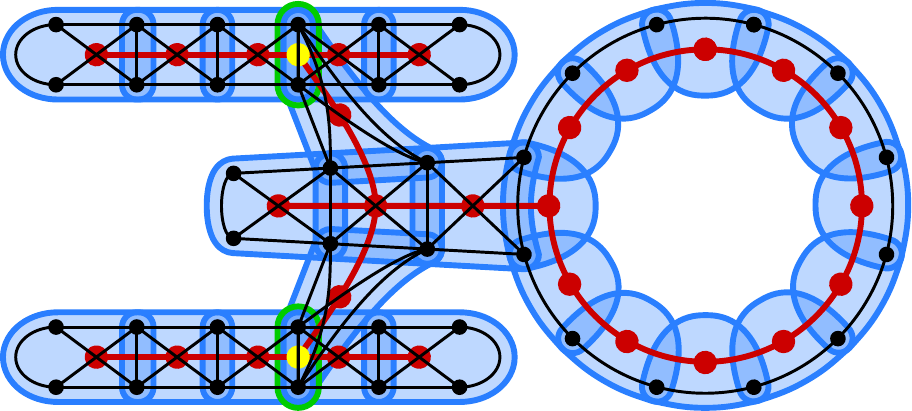}
    \caption{A $2$-connected graph (black) and a cactus-decomposition (red), with each cactus node at the centre of its corresponding blue bag and the $P$-nodes drawn in yellow. Deleting a bridge or two edges of a common cycle splits the cactus into two components; the unions of their bags form the two sides of a displayed separation. The bags of the $P$-nodes, circled in green, form correspond to $2$-separators whose deletion leaves more than two components: these may be partitioned arbitrarily into two groups, with the separator included on both sides. The cycles of the cactus encode the crossing patterns of minimum cuts, while the bridges encode the cuts that avoid every crossing.}
    \label{fig:cactus}
\end{figure}

Our starting point is the SPQR decomposition -- the algorithmic counterpart of Tutte's 2-separator theorem -- in the cactus-formulation of \cite{EP14}. 
Informally, a \emph{cactus-decomposition} of a graph $G$ is defined like a tree-decomposition, except that the decomposition tree is replaced
by a strong cactus~$C$, some of whose nodes are designated as
`$P$-nodes'. 
The decomposition displays separations of $G$ through four types
of cuts of $C$: a bridge, two edges on a common cycle, a cut
consisting of edges incident with a single $P$-node, or the empty cut\footnote{The empty cut displays $(V(G),\emptyset)$ and $(\emptyset,V(G))$.};
see \cref{fig:cactus}.
Each such cut defines a bipartition of the nodes of $C$, which
we pull back to a separation of $G$ by taking the union of the bags
on each side.
The resulting separation is \emph{displayed} by the cactus-decomposition; see \cref{subsec:overview-step2} for the formal definitions.
Cactus-decompositions extend the tree structure of nested separations in the simplest way, by allowing cycles to accommodate coherent crossings.
The SPQR theorem is implied\footnote{We can recover the SPQR tree from the decomposition cactus $C$ by contracting each cycle to a single node labelled $S$ and labelling by $R$ every node not already labelled $P$ or $S$.} by the following variant \cite{EP14}: Every $2$-connected graph has a cactus-decomposition displaying all its proper $2$-separations. 

SPQR trees have become a standard algorithmic tool, underlying classical algorithms for planarity testing and planar embeddings \cite{HT74,M03}, planar graph isomorphism \cite{HT71,HW74}, and dynamic connectivity and planarity problems \cite{ReinventingWheel,DT962,SPQRonline}.

\begin{figure}[ht]
    \centering
    \includegraphics[height=9\baselineskip]{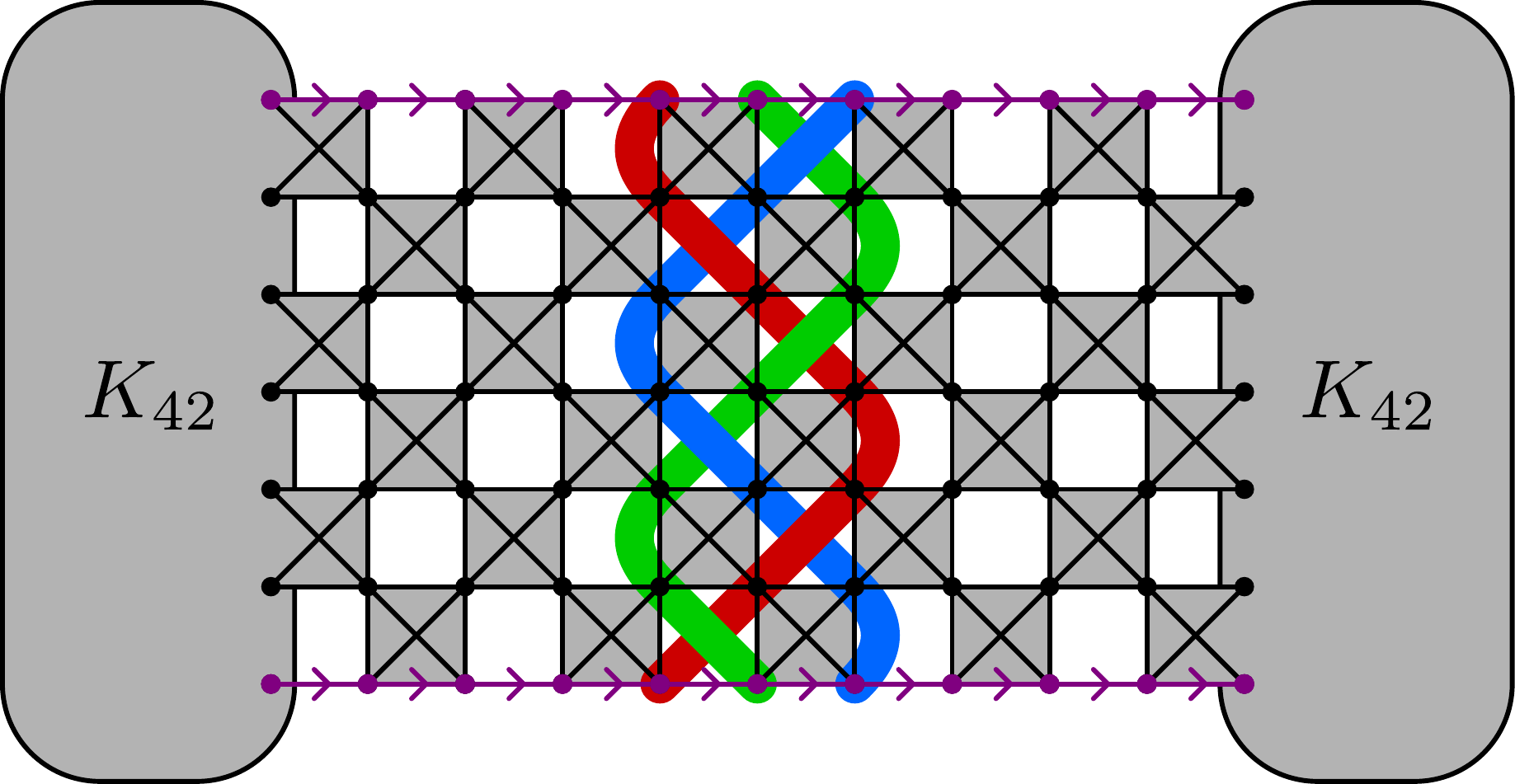}
    \caption{A $6$-connected graph obtained from a cylindrical grid by adding both diagonals in alternate faces and attaching a large clique at each end. The top and bottom rows are identified. The three coloured separations pairwise cross. They cannot be displayed at a common $P$-node, since their separators are distinct, or on a common cycle, since each separator meets all four corners of the other two. Thus no cactus-decomposition displays all the $6$-separations of this graph.}
    \label{fig:a_grid}
\end{figure}

For higher vertex-connectivity, exact cactus-decompositions are no longer possible in general: some graphs admit no cactus-decomposition simultaneously displaying all their minimum-order separations; this happens already for $3$-separations, see \cite{TriAiC} for a detailed discussion.
For general~$k$, cylindrical grids, such as the one in \cref{fig:a_grid}, provide examples of obstructions to an exact cactus-decomposition in which the necessary approximation error grows with~$k$.
Despite the linear global structure of these graphs, their $k$-separations admit too much local freedom to all be displayed by a single cactus-decomposition.

To address this for $k=3$, Carmesin and Kurkofka \cite{TriAiC} introduced \emph{tri-separations}, modified versions of $3$-separations that retain their separating properties while removing the local ambiguity that prevents an exact cactus-decomposition.
With this, they obtained a Tutte-type decomposition of $3$-connected graphs.
Kurkofka and Planken~\cite{TetraDecomp} subsequently introduced \emph{tetra-separations} and proved a Tutte-type decomposition theorem for $4$-connected graphs.
These exact results for modified notions of separation can be viewed as approximate cactus-decomposition theorems for ordinary vertex-separations: for $k=3,4$, every large enough $k$-connected graph $G$ admits a cactus-decomposition displaying $k$-separations such that every $k$-separation of $G$ has edit-distance at most $k$ to a displayed $k$-separation; here, the edit-distance between separations $(A,B)$ and $(C,D)$ is measured by $\max\{\,|A\triangle C|, |B \triangle D|\,\}$. 

Our main result generalizes this approximate cactus picture to arbitrary vertex-connectivity:

\begin{thm}[Simplified version for $k$-connected graphs]\label{k-con-intro}
For every $k\in \Nbb$, there is $d\in \Nbb$ such that every $k$-connected graph $G$ admits a cactus-decomposition displaying $k$-separations such that every $k$-separation of $G$ has edit-distance at most $d$ to a displayed separation.
\end{thm}

A graph can have quadratically many distinct separators of minimum size~\cite{K90}, and a single separator can give rise to exponentially many separations. 
For fixed $k$, our theorem gives a linear-size data structure that approximately represents all $k$-separations of a $k$-connected graph. 
These are represented implicitly by the cactus and its bags through bridges, cycles and $P$-nodes.
Compact representations of minimum vertex separations have also been studied in connection with pairwise connectivity queries~\cite{KTDC91,ReinventingWheel,StructureMinimumVertexCuts}; see \cref{subsec:related-work}.
\medskip

\emph{Unbreakability} is a standard relaxation of the notion of connectivity in parameterized algorithms, {introduced in the context} of recursive understanding \cite{KT11,CCHPP16,LRSZ18,CLPPS19}, which reduces many graph problems to suitable unbreakable instances.
In this paper, we prove a strengthening of \cref{k-con-intro} with the assumption of \lq $k$-connectivity\rq\ relaxed to \lq $(k-1)$-unbreakability\rq. This generalisation allows us to apply our result iteratively.

Formally, a graph is \defn{$(q,k)$-unbreakable} if every separation of order at most $k$ has a proper side of size at most $q$.
Thus, the $k$-connected graphs are exactly the $(0,k-1)$-unbreakable graphs with more than $k$ vertices, while allowing $q>0$ gives an approximate version of $k$-connectivity.
This relaxation allows us to look beyond the local obstructions imposed by low-degree vertices and capture the global connectivity of graphs such as square grids and cubic expanders.
Our main result holds in this more general setting:

\begin{restatable}[Main theorem]{thm}{cactusgeneral}\label{thm:main-intro}
For all $q, k\in \Nbb$, there is $d\in\Nbb$ such that every $(q,k-1)$-unbreakable graph $G$ admits a cactus-decomposition displaying $(\leq k)$-separations such that every $(\leq k)$-separation of $G$ has edit-distance at most $d$ to a displayed separation.
\end{restatable}

The approximation bound $d$ is effective and can be taken to be $2^{(q+k)^{O(1)}}$. 
We leave the dependence on $q$ and $k$ implicit in our proof, prioritising simplicity over optimizing bounds.

Iterated decompositions have proved useful in structural graph theory~\cite{TreesOfTDCs}, graph enumeration~\cite{GNR13}, and algorithms for connectivity~\cite{DKKN23,conny_aug} and graph isomorphism~\cite{DNTW10,DLNTW22}.
\cref{thm:main-intro} can be applied recursively as follows.
We first group together the bags whose nodes lie on the same cycle of the cactus, obtaining pieces with cycle-decompositions of bounded adhesion.
We retain these cyclic pieces, as well as the bags at $P$-nodes (which have size at most~$k$).
Each remaining bag can be described by a $k$-unbreakable minor of $G$, to which we apply \cref{thm:main-intro} with $k$ increased by one to reveal finer structure, see \cref{app:unbreakable-minors}.

\subsection{Examples and the need for approximation}

We now illustrate our theorem on a few families of graphs with a simple underlying structure, showing how the cactus-decomposition can recover this structure. Some of these examples also justify the need for approximation.

Already for $k=3,4$, there are $k$-connected graphs whose $k$-separations cannot all be displayed by a single cactus-decomposition.
Toroidal hexagonal and square grids provide examples for $k=3$ and $k=4$, respectively.
These examples extend to higher connectivity, for instance by adding universal vertices, or by replacing each vertex $v$ by a clique $K(v)$ and replacing each edge $uv$ with all possible edges between $K(u)$ and~$K(v)$.
However, these graphs are $(q,k)$-unbreakable for some~$q$ depending only on~$k$, so every $k$-separation has a side of bounded size.
It is therefore natural to accept such graphs as basic pieces of our decomposition, leaving them intact and focusing on sufficiently broad separations, meaning those with both sides sufficiently large.

\begin{figure}[ht]
    \centering
    \includegraphics[width=\linewidth]{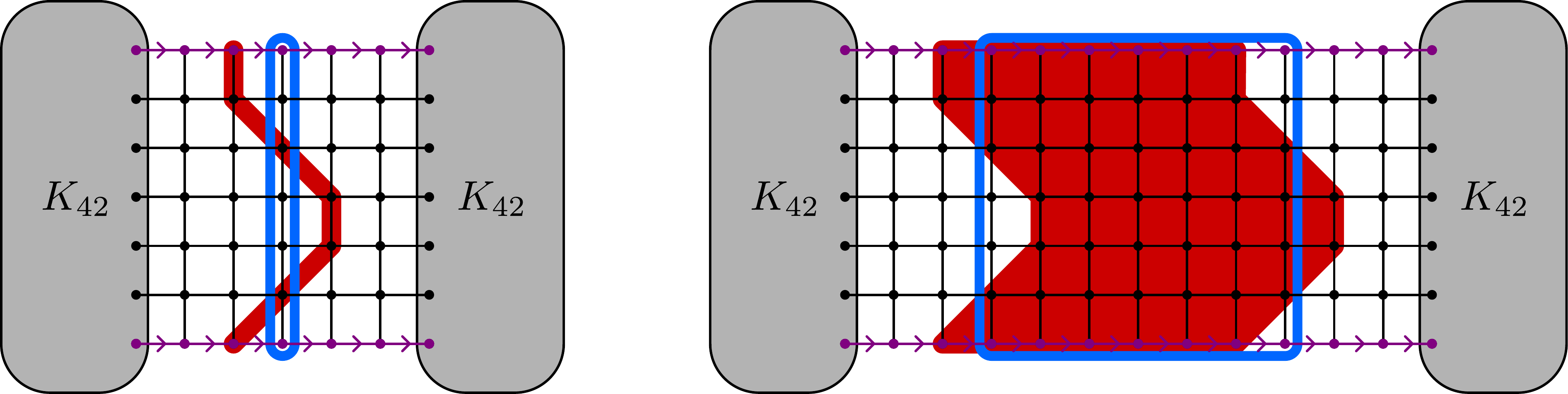}
    \caption{Both panels show a cylindrical grid of width~$6$, with the top and bottom rows identified. Left: the red $6$-separation is close to the blue separation defined by a single vertical cross-section. Right: in the $7$th power of the grid, the red vertices form a minimum separator. The corresponding separation is close to the blue separation, whose separator consists of seven consecutive cross-sections.}
    \label{fig:cylinder}
\end{figure}

To see why even these separations require approximation, consider a long cylindrical grid of width~$k$, with a large clique attached at each end; see \cref{fig:cylinder}.
This graph is $(q,k-1)$-unbreakable for some~$q$ depending only on~$k$.
Its vertical cross-sections define separations of order~$k$.
These separations admit local modifications that preserve their order but produce crossing patterns that cannot all be displayed by a single cactus-decomposition, for the same reason as in \cref{fig:a_grid}.
Since the end cliques can be arbitrarily large, all the separations involved can be arbitrarily broad.

Nevertheless, as illustrated in \cref{fig:cylinder}, every sufficiently broad $k$-separation differs from a straight one only on a number of vertices bounded in terms of~$k$, independently of the length of the cylinder.
The cactus-decomposition constructed in our proof displays precisely the straight separations, thereby recovering the intrinsic linear structure of the graph.

Taking the $r$-th power of this graph yields a $k'$-connected example with $k'=rk$, provided $r\ge k$ and the end cliques are sufficiently large.
Now, the vertices of any $r$ consecutive vertical cross-sections form a minimum separator.
The cactus-decomposition constructed in our proof displays precisely these separations, again recovering the intrinsic linear structure of the graph.
More generally, as illustrated in \cref{fig:cylinder}, a minimum separator consists of $r$ consecutive vertices on each longitudinal path, with starting positions differing by at most one between neighbouring paths.
These variable bands give rise to the same obstruction to an exact cactus representation, now among minimum-order separations of a $k'$-connected graph.
Thus, strengthening unbreakability to $k$-connectivity does not eliminate the obstruction.
Nevertheless, every minimum separation still differs from one displayed by our decomposition only on a bounded number of vertices.

One might hope to recover an exact theorem by generalising tri- and tetra-separations to a suitable notion of reduced mixed-separation, whose separator may contain both vertices and edges.
In the cylinder examples, however, the crossing patterns can involve a number of vertices quadratic in the width of the cylinder, and it is not clear how a reduction favouring edges over vertices could resolve them.

\begin{figure}[ht]
    \centering
    \captionsetup{margin=0pt}
    \begin{minipage}[t]{0.48\textwidth}
        \centering
        \includegraphics[height=11\baselineskip]{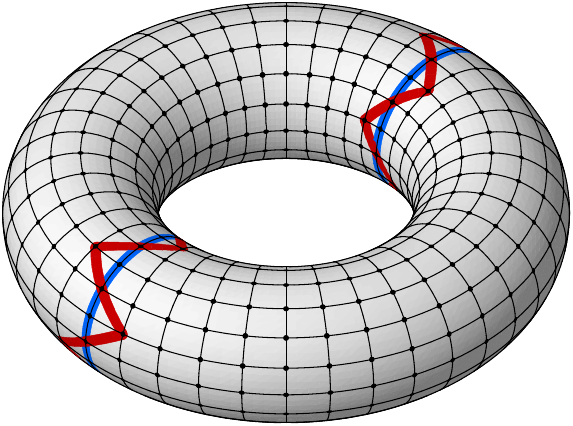}
        \caption{A toroidal grid with meridians of length~$m=16$ and longitudes of length~$40$, so $k=2m=32$. The two blue meridians together form the separator of a $k$-separation. The red cycles together form the separator of another $k$-separation, whose sides agree with those of the blue one outside the two separators.}
        \label{fig:torus}
    \end{minipage}\hfill
    \begin{minipage}[t]{0.48\textwidth}
        \centering
        \includegraphics[height=12\baselineskip]{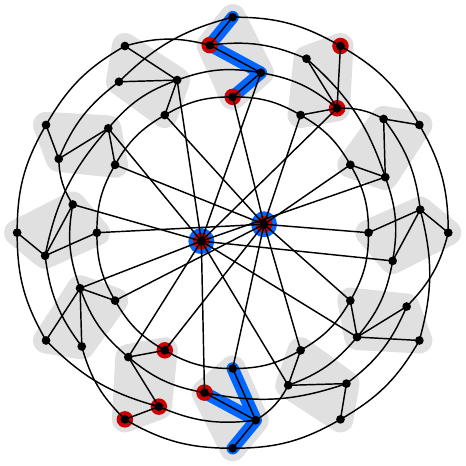}
        \caption{Twelve trees with $t=4$ vertices each, joined cyclically by perfect matchings, and $\ell=2$ hubs, so $k=2t+\ell=10$. The two blue trees and the hubs form the separator of a $k$-separation. The red vertices and the hubs form the separator of another $k$-separation, whose sides agree with those of the blue one outside the two separators.}
        \label{fig:wheel}
    \end{minipage}
\end{figure}

A toroidal grid whose longitude is much longer than its meridian provides a cyclic analogue of the cylinder; see \cref{fig:torus}.
Writing~$m$ for the length of the meridian and setting~$k=2m$, this graph is again $(q,k-1)$-unbreakable for some~$q$ depending only on~$k$.
Pairs of meridional cross-sections define separations of order~$k$, and every sufficiently broad $k$-separation differs from one of these only on a bounded number of vertices, as shown in \cref{fig:torus}.
Our construction displays precisely the separations defined by these pairs of cross-sections, thereby recovering the intrinsic cyclic structure of the graph.

Finally, consider the following generalisation of the wheels introduced by Geelen and Joeris~\cite{GeelenJoerisGrid}.
Take a long cyclic sequence of arbitrary $t$-vertex trees, join consecutive trees by arbitrary perfect matchings, and add $\ell$ hubs, each adjacent to a vertex in every tree; see \cref{fig:wheel} for an example.
For $k=2t+\ell$, these graphs are again $(q,k-1)$-unbreakable for some~$q$ depending only on~$k$.
Pairs of trees together with all the hubs define separations of order~$k$.
Every sufficiently broad $k$-separation differs from one of these only on a number of vertices bounded in terms of~$k$.
Our cactus-decomposition recovers the cyclic arrangement of the trees, with the hubs forming a common core of the displayed separators.

\subsection{A short outline of the proof of \texorpdfstring{\cref{k-con-intro}}{Theorem 1.1} (see \texorpdfstring{\cref{sec:overview}}{Section 2} for a detailed overview).}\label{subsec:brief-overview}
We construct the cactus-decomposition of \cref{k-con-intro} in two steps.
First, we build a tree-decomposition $\mathcal{T}$ whose bags, roughly speaking, either have higher connectivity or possess a cyclic {or star-like} structure.
We then refine $\mathcal{T}$ into the desired cactus-decomposition $\mathcal{C}$ by replacing each bag of the latter type with the corresponding cycle {or star}-decomposition{,  with the centres of the stars becoming $P$-nodes}.

{\bf Step 1: Constructing the tree-decomposition $\Tcal$.}
We construct our tree-decomposition from a maximal nested set of carefully selected separations and use maximality to obtain the desired approximation guarantee.
A maximality argument also underlies a modern proof of the tangle-tree theorem~\cite{CarmesinToTshort}, where the aim is to distinguish the highly connected regions of a graph~$G$, formalised by its maximal tangles.
We briefly recall the argument.
Given a graph $G$, let $S$ be the set of all separations efficiently distinguishing two tangles, and choose a maximal nested subset $T\subseteq S$.
The separations in $T$ induce a tree-decomposition $\Tcal$: since they are nested, we can successively split~$G$ along them to obtain the bags of the decomposition.
One then shows that the maximality of $T$ ensures that $\Tcal$ distinguishes all maximal tangles.

We now return to the proof of \cref{k-con-intro} and explain how we use maximality in this setting.
Given $k\in\Nbb$, choose constants $\beta$ and~$\gamma$ depending only on $k$; their precise values are not important for this overview.
The approximation bound $d$ will be larger than both.

Let $G$ be a $k$-connected graph.
We say that two $k$-separations $(A,B)$ and $(C, D)$ \defn{heavily cross} if all four sets $A\cap C$, $A\cap D$, $B\cap C$ and $B\cap D$ have size at least $\gamma$.
Intuitively, these are the separations that will be captured by the cyclic parts of the cactus, so we set them aside in the first step.
Let $S$ be the set of all $k$-separations $(A,B)$ with $|A|,|B|\geq\beta$ that are not heavily crossed.
Choose a maximal nested subset $T\subseteq S$, and let $\Tcal$ be the tree-decomposition induced by $T$.

Using corner diagram arguments, we show that every separation in $S$ has bounded edit-distance to some separation in $T$.
Thus, $\Tcal$ approximately displays all separations in $S$.
Separations~$(A,B)$ with~$|A|<\beta$ or $|B|<\beta$ are already within edit-distance $\beta$ of $(\emptyset,V(G))$ or $(V(G),\emptyset)$, respectively, and are therefore covered by our choice of $d$.
It remains to refine $\Tcal$ into a cactus-decomposition that also approximately displays the heavily crossed separations.

{\bf Step 2: Construction of the cactus-decomposition $\Ccal$.}
We say that a separation $(A,B)$ \defn{lives} in a bag of a tree-decomposition if its separator is contained in that bag.
Another corner diagram argument shows that every heavily crossed $k$-separation has bounded edit-distance to a separation that lives in a bag of $\Tcal$ and remains heavily crossed (for a slightly smaller value of $\gamma$).

We can therefore build $\Ccal$ by treating the bags of $\Tcal$ separately.
For each bag, we seek a cycle- {or star-}decomposition of its \defn{torso} that approximately displays all the heavily crossed $k$-separations living there.
Here, the torso is the graph induced by the bag, with each intersection with a neighbouring bag completed to a clique.
Constructing these {local} decompositions is the heart of the proof and relies on the following new ideas: 
\begin{enumerate}
\item\label{itm:3d-corner}
Each separation $(A,B)$ partitions the vertex set into three parts: $A\setminus B$, $A\cap B$ and~$B\setminus A$.
Three separations therefore determine $27$ regions.
We show that when three $k$-separations pairwise heavily cross, their crossing pattern is much more constrained: many of these regions must be empty, and the whole configuration admits a cyclic description.
In particular, heavily crossing separations cannot form truly three-dimensional crossing patterns, such as the one in \cref{Overview:fig:3Dcube}.

\begin{figure}[ht]
    \centering
    \includegraphics[height=8\baselineskip]{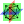}
    \caption{Three separations cross three-dimensionally in the line-graph of the cube. Their separators are indicated in red, green and blue.}
    \label{Overview:fig:3Dcube}
\end{figure}

\item\label{itm:flowers}
Families of separations with a common separator can be displayed by a star-decomposition.
For the cyclic case, we use ideas from covering spaces in topology to combine the cyclic descriptions of triples in \ref{itm:3d-corner} into a single cycle-decomposition that approximately displays the family.
Such cycle- and star-decompositions correspond to \defn{flowers} from the matroid literature \cite{{Structure3sepsMatroids},{Aikin_flowers}}; \cref{fig:flower99} illustrates the cyclic case.

\begin{figure}[ht]
    \centering
    \includegraphics[height=8\baselineskip]{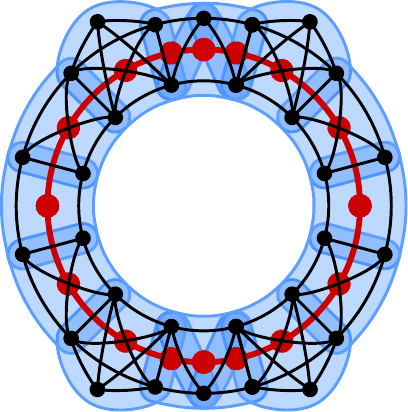}
    \caption{A cycle-decomposition (blue bags, red decomposition cycle) of a graph~$G$ (black) that is a $4$-flower}
    \label{fig:flower99}
\end{figure}

\item
{Several crossing families may occur in the same bag of $\mathcal{T}$, giving rise to several flowers.}
However, by choosing $d$ sufficiently large compared to $\beta$ and $\gamma$, we show that one flower per bag suffices: we can select a flower $\Fcal$ in each bag so that every $k$-separation living in this bag is within bounded edit-distance of a separation displayed by $\Fcal$ or $\Tcal$.
This requires a systematic analysis of flowers which, together with the canonical constructions via coverings in \ref{itm:flowers}, {extends the structural theory of flowers in graphs}, and provides, for the first time, canonical constructions of flowers.

\end{enumerate}

This completes the second step of the construction and thus our brief
proof overview.
We give a more detailed proof overview in \cref{sec:overview}, expanding on each of these three ideas.

\subsection{A brief real-world example with a cyclic structure}

\begin{figure}[ht]
  \centering
  \includegraphics[width=12.5cm]{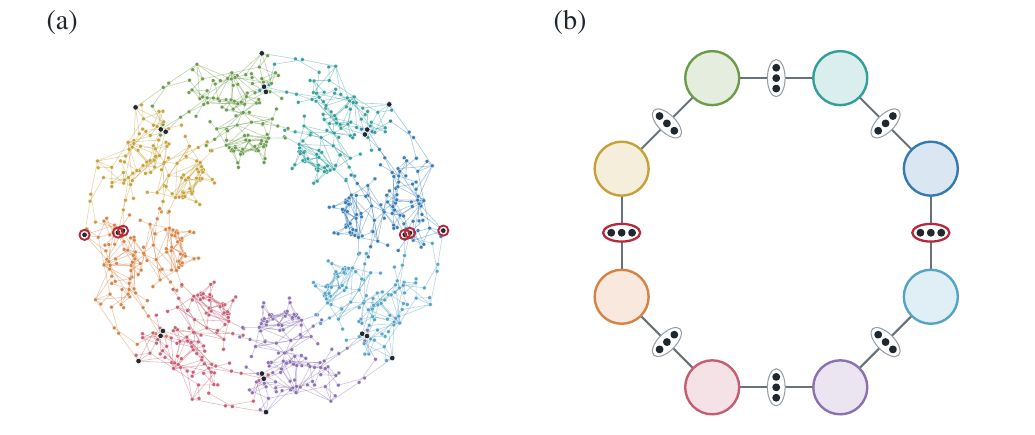}
\caption{
An eight-petal flower in a biological contact graph. Left: the largest component of the contact graph of the human nuclear-pore model PDB 7TBJ (Petrovic et al.). Vertices represent protein fragments; edges record contacts within $5\,\text{Å}$ between non-hydrogen atoms.
Vertex positions are obtained by projecting the fragment centres onto the $xy$-plane, so distinct vertices may overlap in the drawing.
In particular, the six red vertices form a separator. This is difficult to discern in the projected drawing at the displayed scale, but becomes clear at $400\%$ magnification.
Right: the eight bags form a cycle with three-vertex overlaps; matching colours identify their interiors on the left. The red $6$-separator leaves 469 vertices on each exclusive side.
\newline
}
\label{fig:npc-flower}
\end{figure}

\Cref{fig:npc-flower} shows an eight-petal flower in a contact graph
obtained from the human nuclear-pore model PDB~7TBJ of Petrovic
et al.~\cite{Petrovic2022NPC}.

Improving the constants in our main theorem would be interesting: with bounds small enough to apply here, an iterative application might eventually look like this:
the graph has a unique $6$-tangle that, for every separation $(A,B)$ of order at most five, selects a side containing a vertex from each of the eight coloured petal interiors outside $A\cap B$.  
We would start to decompose iteratively with the block-cutvertex theorem, Tutte's theorem at every 2-block and then the tri-separator theorem at each 3-connected torso, etc. 
Eventually, we would obtain a $5$-unbreakable minor carrying this $6$-tangle. Then we would apply the cactus decomposition for $k=6$, and the cactus structure it would compute would be similar to the cyclic decomposition of \cref{fig:npc-flower}.
The underlying computations have been worked out in detail, see \cite{ChatGPT2026NuclearPore}.

\subsection{Applications and perspectives} \hspace{1cm} \\
\indent \underline{Tangles and flowers.}
A \emph{tangle} describes a highly connected region of a graph by indicating, for each separation of small order, on which side it lies.
The tangle-tree theorem provides a tree-decomposition that distinguishes these regions~\cite{GMX}.
We would like to extend this viewpoint to include large flowers as highly connected regions in a broader sense.
This leads us to consider \emph{bramble orientations}, a relaxation of tangles, and to group related orientations describing the same flower.
We seek an analogue of the tangle-tree theorem that recognises these flower regions and distinguishes them from one another and from tangles.
In a follow-up paper~\cite{FollowUp}, we build on the structural results developed here to obtain an approximate tangle-tree theorem that also detects sufficiently large flowers.
\smallskip

\underline{Connectivity augmentation.}
The \textsc{Connectivity Augmentation} problem asks whether a graph can be made $k$-connected by adding at most $\ell$ edges from a prescribed set.
If the graph is already $(k-1)$-connected, this amounts to selecting edges that cross \emph{every} proper $(k-1)$-separation.
For edge-connectivity, the DKL theorem reduces increasing the connectivity by one to augmentation on a cactus.
This reduction has led to FPT and approximation algorithms for minimum-cost augmentation~\cite{MV15,Traub_2023}.

For vertex-connectivity, Carmesin and Ramanujan~\cite{conny_aug} obtained FPT algorithms parameterized by~$\ell$ for $k\leq4$, using exact cactus representations.
Korhonen and Thorup~\cite{KT26} extended this result to arbitrary $k$ with running time $2^{O(\ell\log(\ell+k))}n^{O(1)}$, using approximate structural arguments.
For $k=2$, Koana and Kumabe~\cite{KK26} recently obtained an FPT algorithm with single-exponential dependence on $\ell$ for arbitrary input graphs.
A natural question is whether our approximate cactus representation can help achieve such a dependence on $\ell$ for larger fixed $k$.
\smallskip

\underline{Canonicity and symmetry.}
The canonicity of the Tutte-decomposition underlies linear-time and logspace algorithms for planar graph isomorphism~\cite{HW74,DLNTW22}.
Retaining cyclic torsos matters for canonicity: decomposing a long cycle into bounded-size bags necessarily breaks some of its symmetries~\cite{Isomorphism}.
The decomposition in \cref{thm:main-intro} approximately represents the full family of separations of small order, which is canonical.
Its displayed separation system is therefore preserved by every automorphism up to bounded edit-distance, yielding a form of approximate canonicity.

Since canonical decompositions are invariant under graph automorphisms, they are particularly useful for studying highly symmetric graphs.
Using their canonical tetra-decomposition~\cite{TetraDecomp}, Kurkofka and Planken classified all finite connected vertex-transitive graphs of connectivity below five~\cite{TetraTrans}.
Our decomposition may help extend their classification to higher connectivity.
\smallskip

\underline{Cops and robber.}
In the classical helicopter cops and robber game, $k$ cops win precisely when the graph has treewidth less than $k$~\cite{GMcopsRobber}.
Consider the randomized variant with an invisible robber, no repeated visits by the cops, and mixed strategies\footnote{A mixed strategy is a probability distribution over deterministic strategies.} for both players.
For $k\le 4$, Planken~\cite{CopsInvisibleRobber} uses the corresponding cactus-decomposition theorems to determine the maximum winning probability of $k$ cops on $k$-connected graphs.
Our cactus representation may support similar analyses at higher connectivity.
\smallskip

\underline{Local analogues.}
Local separators disconnect bounded-radius balls without necessarily separating the whole graph.
Local analogues of graph decomposition theorems have been developed in~\cite{Local2sep,LocalKsep,canonicalGraphDec,Ponyhof}.
This perspective has also led to structural characterizations and canonical decompositions of locally chordal graphs~\cite{abrishamiGlobalStructureLocally2025,abrishamiLocallyChordalGraphs2025,jacobsCanonicalTreedecompositionsChordal2025}.
Local separators have applications in network analysis~\cite{Sarah} and distributed approximation~\cite{BGPW25}.

Our construction of the cyclic parts of the decomposition uses covering graphs to combine local separators into global separators.
This suggests a local analogue of our theorem that, under appropriate local connectivity assumptions, approximately represents local separations together with their crossing structure.

\subsection{Related work}\label{subsec:related-work}

The tangle-tree theorem of Robertson and Seymour~\cite{GMX} states that every graph admits a tree-decomposition distinguishing all its maximal tangles.
Canonical constructions and alternative proofs have been developed in~\cite{confing,CDHH13CanonicalAlg,CarmesinToTshort,entanglements}.
Refinements address the structure of the parts containing tangles~\cite{CDHH13CanonicalParts,AlbrechtsenOptimalToT}.
Erde~\cite{RefiningToT} showed how to further decompose the parts containing no tangle of the prescribed order, and Albrechtsen~\cite{AlbrechtsenRefiningToTinASS} extended this approach to abstract separation systems.
Related work studies $k$-blocks and tree-decompositions displaying them~\cite{Kblocks,CG14:isolatingblocks,AlbrechtsenRefiningTDCdisplayKblocks}.

The theory extends to profiles and abstract separation systems~\cite{ProfilesNew,StructuralASS,FiniteSplinters}, and to infinite settings~\cite{infinitesplinter,jacobs2023efficiently}.
Its algorithmic applications include computing with tangles~\cite{ComputingTangles} and clustering~\cite{diestel2017tangles,elbracht2020tangles}; see~\cite{GroheTangles3,TangleBook} for broader accounts.
Our emphasis in this paper is on an approximate representation of sufficiently broad separations together with their crossing structure.
\smallskip

Compact descriptions of separations also support data structures for pairwise connectivity queries.
For edge-connectivity, a Gomory--Hu tree~\cite{GH61} encodes all pairwise minimum-cut values in $O(n)$ space and can be preprocessed to answer queries in constant time.
For vertex-connectivity, Cohen, Di Battista, Kanevsky, and Tamassia~\cite{ReinventingWheel}, building on~\cite{KTDC91}, obtained an $O(k^3n)$-space representation for pairwise $(k+1)$-connectivity queries in $k$-connected graphs.
Pettie and Yin~\cite{StructureMinimumVertexCuts} achieved $O(kn)$ words of space and constant query time.
Pettie, Saranurak, and Yin~\cite{PSY26} proved an $\Omega(kn)$-bit lower bound for general graphs, which Nutov~\cite{Nutov26} extended to $k$-connected graphs.
Our theorem complements these results by describing how separations of small order fit together in a cactus structure, up to bounded edit-distance.
\smallskip

Other structural descriptions of graphs at low connectivity include decompositions of $3$-connected graphs~\cite{4tangles,Ugo,grohe2016quasi} and the characterization by Oporowski, Oxley, and Thomas~\cite{typical3con4con} of the typical subgraphs of $3$-connected and $4$-connected graphs.
A complementary question asks which substructures of a $k$-connected graph can be deleted or contracted while preserving $k$-connectivity.
For $3$-connectivity, this includes work of Ando, Enomoto, and Saito~\cite{ando1987contractible} and of Kriesell~\cite{KriesellNonEdges,KriesellContractibleSubgraph,KriesellAlmostAll,KriesellTriangleFree,kriesell2008number,KriesellVertexSuppression}.
The corresponding questions for $4$-connectivity are studied in~\cite{Number4con,generatingInternally4con,Martinov,MartionvTwo,Quasi4conStructure,chain4con,Splitting4con,Slater4con}.
\smallskip

The structure of $k$-separations in $k$-connected matroids has also been studied extensively.
Results for $2$-connectivity appear in~\cite{Structure2sepsMatroids,cunningham_edmonds_1980,Seymour2sepMatroid}.
Oxley, Semple, and Whittle~\cite{Structure3sepsMatroids} introduced flowers to study the structure of $3$-separations in $3$-connected matroids.
The case of $4$-connectivity was studied in~\cite{Structure4sepsMatroids,4connectedMatroids}.
Aikin and Oxley~\cite{Aikin_flowers} generalized flowers to arbitrary orders and classified their structure.
In a framework covering both graphs and matroids, Clark and Whittle~\cite{TanglesTreesFlowers} developed flower representations using an equivalence that can identify separations differing on arbitrarily many elements.
Our theorem gives a uniform bound on the edit-distance to representatives.
\smallskip

For directed graphs, Bowler, Gut, Hatzel, Kawarabayashi, Muzi and Reich~\cite{Flo} proved a Tutte-type decomposition theorem for $k=1$.
\smallskip

Our decomposition is closely related to \emph{unbreakable tree-decompositions}~\cite{CLPPS19,cygan2020randomized}, and \emph{lean tree-decompositions}~\cite{GMIV,T90}.
Roughly speaking, our cactus-decompositions are to unbreakable tree-decompositions what Tutte decompositions are to ear decompositions.
Indeed, deleting one edge from each cycle and adding its adhesion set to every bag along the resulting path converts our cactus-decompositions into unbreakable tree-decompositions of adhesion at most $k$.
We do not know a direct converse construction recovering the cyclic structure.

Unbreakable tree-decompositions support parameterized algorithms for \textsc{Minimum Bisection}~\cite{CLPPS19}, \textsc{Fair Bisection}~\cite{inamdar2023parameterized}, and model checking~\cite{LRSZ18,pilipczuk2022algorithms,schirrmacher2024model}.
Their construction has seen substantial improvements~\cite{cygan2020randomized,anand2025unbreakable}.
Korhonen~\cite{korhonen2024} obtained a parameterized linear-time algorithm for computing the stronger lean tree-decompositions.
This naturally suggests an analogous running-time bound for our cactus-decompositions.

\begin{conjecture}\label{conj99}
There exist computable functions $f,\eta$ such that, given $q,k \in \Nbb$ and an $n$-vertex $(q,k-1)$-unbreakable graph $G$, one can compute in time $f(q,k)\bigl(|V(G)|+|E(G)|\bigr)$ a cactus-decomposition displaying $(\leq k)$-separations such that every $(\leq k)$-separation of $G$ has edit-distance at most $\eta(q, k)$ to a displayed separation.
\end{conjecture}

\enlargethispage{2\baselineskip}

\centerline{\includegraphics[height=7\baselineskip]{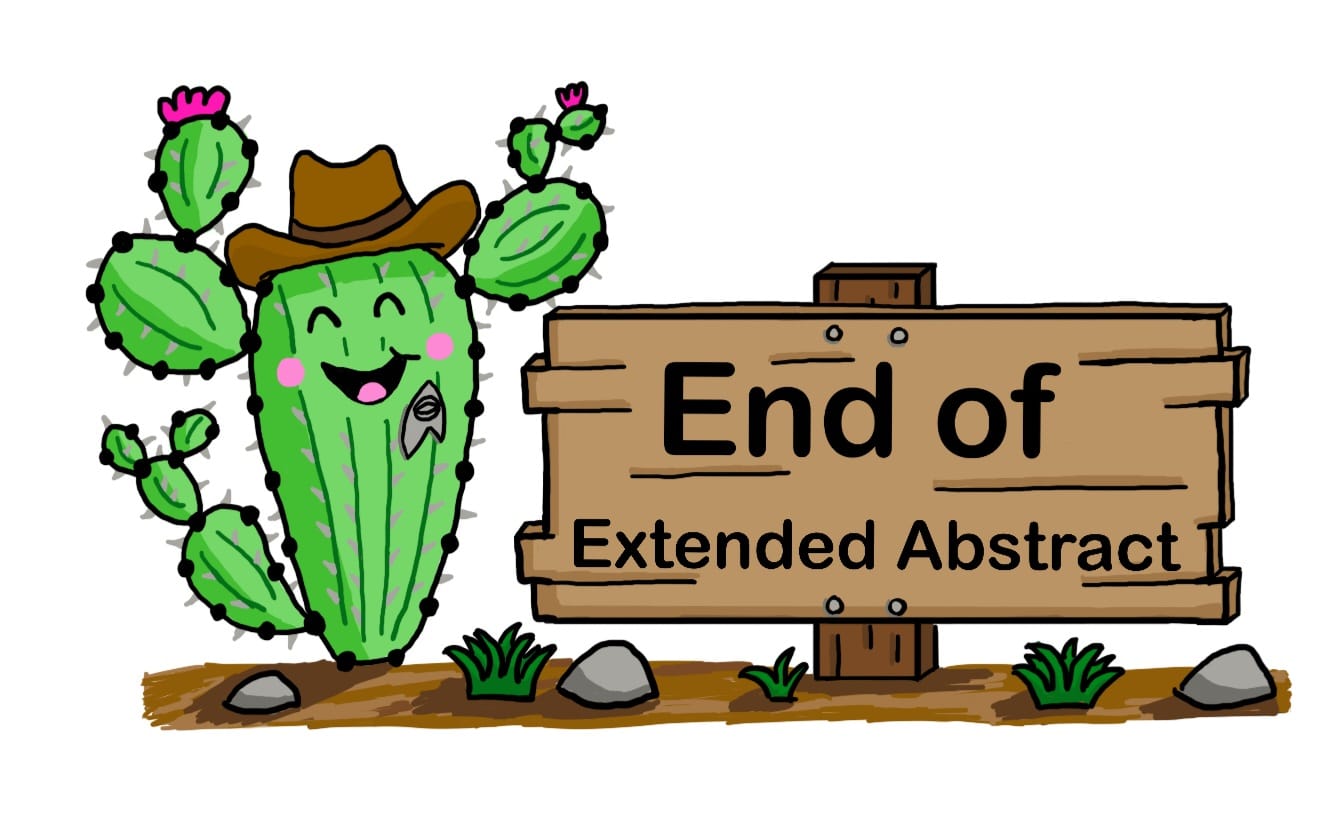}}

\newpage
\section{A detailed overview of the proof of \texorpdfstring{\cref{k-con-intro}}{Theorem 1.1}}\label{sec:overview}
Definitions and statements in this section are often simplified -- see the later sections for accurate information.

\subsection{Strategy}
\begin{figure}[ht]
    \centering
    \includegraphics[height=9\baselineskip]{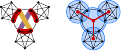}
    \caption{Left: The yellow and purple 2-separators cross, while the red 2-separators are totally-nested. Right: The Tutte-decomposition from a tree-decomposition perspective.}
    \label{fig:TutteFor2}
\end{figure}

Our proof is inspired by the pattern of Tutte's decomposition, which we briefly recall.
To decompose a $k$-connected graph for $k=0$ and $k=1$, it suffices to cut along all $k$-separators of the graph.
But already for $k=2$, this approach fails.
Indeed, different 2-separators can \emph{cross} in that they cut through each other, as happens in a cycle.
Then if we cut at one of the 2-separators, we lose the other.
Tutte realised that the solution is to cut the graph at all 2-separators that are \emph{totally-nested}, those that are not crossed by any 2-separator, and indeed no conflicts arise when cutting in this way; see \cref{fig:TutteFor2}.
The result is a tree-decomposition whose torsos are 3-connected, cycles or $K_2$'s.
It captures every proper $2$-separation~$s$ of the graph: either $s$ is induced by an edge of the decomposition tree, or $s$ stems from a 2-separation of a cyclic torso, or the separator of~$s$ is the bag of a node~$t$ with a $K_2$-torso and $s$ distributes the branches of the decomposition tree at~$t$ among its two sides.
Nodes with $K_2$-torsos become $P$-nodes in the cactus-formulation from the introduction (see \cref{fig:cactus} for an example). 
The extensions of Tutte's result to connectivity~$k\in\{3,4\}$~\cite{TriAiC,TetraDecomp} follow the same strategy: they cut the graph at all totally-nested tri-separations or tetra-separations, respectively, the modified $k$-separations from the introduction.
The resulting basic building blocks are quasi-$(k+1)$-connected graphs as well as graphs exhibiting a specific cyclic or star-like structure (see \cref{fig:flower99} for the cyclic case).
The latter are precisely the graphs in which every tri- or tetra-separation, respectively, is crossed by another one.
For small~$k$, such graphs can be captured with \emph{$k$-flowers}, a central notion in matroid theory.

We show that this correspondence extends to all~$k$, up to bounded errors.
For this, we let $G$ be a $(q,k-1)$-unbreakable graph for some $q,k\in\N$.
Our goal is to obtain a cactus-decomposition of~$G$ (\cref{thm:main-intro}).
The proof proceeds in two phases.
In the first phase, which consists of Steps~1--4 below and spans almost the entire paper, we construct a tree-decomposition of~$G$ along $k$-separations that are approximately totally-nested.
Its torsos are $(q',k)$-unbreakable for a constant~$q'$, or approximately have the structure of a $k$-flower.
These are our analogues of the 3-connected torsos on the one hand, and of the cyclic torsos and $P$-nodes on the other, in Tutte's decomposition.
This is our technical main result, \cref{mainTutte}, which we state at the end of Step~4.
In the second phase, which consists of Step~5, we refine this tree-decomposition into the desired cactus-decomposition: by slightly adjusting the flowers, and forging them into the tree-decomposition.

The brief overview in \cref{subsec:brief-overview} constructs the tree-decomposition first and the flowers afterwards.
Here we follow the order of the paper instead.
Steps~1 and~2 expand on the first two ideas listed there: Step~1 analyses how $k$-separations heavily cross with the help of three-dimensional corner diagrams, and Step~2 uses coverings to assemble heavily crossing $k$-separations into flowers.
Step~3 constructs the tree-decomposition, Step~4 expands on the third idea by analysing its torsos, and Step~5 refines it into the cactus-decomposition.

\subsection{Approximate terminology}
The \defn{breadth} of a separation $(A,B)$ is the minimum of $|A\sm B|$ and $|B\sm A|$.
We say that $(A,B)$ is \defn{$\beta$-broad} to mean that $(A,B)$ has breadth~$\ge \beta$.
In these terms, a graph is $(q,k)$-unbreakable if and only if every $(\le k)$-separation of it has breadth~$\le q$.
Two separations $(A_1,B_1)$ and $(A_2,B_2)$ are \defn{$\eps$-close} if all four sets $A_1\sm A_2$, $A_2\sm A_1$, $B_1\sm B_2$ and $B_2\sm B_1$ have size~$\le\eps$.
Closeness is our working version of the edit-distance from the introduction: $\eps$-close separations have edit-distance at most~$2\eps$ (see \cref{closeVSapprox}).
Closeness is quasi-transitive: if $r$ and~$s$ are $\eps$-close and $s$ and~$t$ are $\eps'$-close, then $r$ and~$t$ are $(\eps+\eps')$-close.
Let $S,S'$ be sets of separations of~$G$.
We say that $S$ \defn{$\eps$-approximates} $S'$ if every $\eps$-broad $s'\in S'$ is $\eps$-close to some separation in~$S$.

We write $X\,\defnMath{\se_\eps}\,Y$ for $X,Y\se V(G)$ to mean that all but $\le \eps$ vertices in $X$ are contained in~$Y$.
We define a relation \defn{$\leq_\eps$} on the set of separations of~$G$ by letting
\[
    (A,B) \,\defnMath{\leq_\eps}\, (C,D)\quad :\Longleftrightarrow\quad A\se_\eps C\text{ and }B\supseteq_\eps D.
\]
Two separations $(A,B)$ and $(C,D)$ of $G$ are \defn{$\eps$-nested} if, after possibly swapping the name $A$ with $B$ or $C$ with~$D$, we have $(A,B)\le_\eps (C,D)$.
For $\eps=0$, we drop the index: $(A,B)\le(C,D)$ means $A\se C$ and $B\supseteq D$, and $0$-nested separations are \defn{nested}; separations that are not nested \defn{cross}.
A set of separations is \defn{nested} if its elements are pairwise nested.

For the vast majority of the paper, we assume that $G$ is additionally \defn{$(\beta,k)$-narrow} for some $\beta\in\N$, meaning that for every $\beta$-broad $k$-separation $(A,B)$ of~$G$, every component of $G-(A\cap B)$ has neighbourhood equal to~$A\cap B$.
This technical property allows us to almost treat $G$ as if it were genuinely $k$-connected.
We remove this additional assumption in the finale of the paper without too much effort (\cref{sec:RemoveNarrow} and~\cref{sec:TutteNarrow}).
In particular, the key lemmas below assume that~$G$ is $(\beta,k)$-narrow, while \cref{mainTutte} holds for all $(q,k-1)$-unbreakable graphs~$G$.

\subsection{Step~1: Analysing how \texorpdfstring{$k$}{k}-separations heavily cross}\label{subsec:overview-step1}
\begin{figure}[ht]
    \centering
    \includegraphics[height=7\baselineskip]{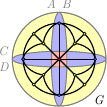}
    \caption{A crossing-diagram for two separations $(A,B)$ and $(C,D)$ of a graph~$G$. The four \defn{corners} are yellow, the four \defn{links} are blue, and the \defn{centre} is red.}
    \label{Overview:fig:2DcrossingDiagram}
\end{figure}

When two $k$-separations $(A,B)$ and $(C,D)$ of~$G$ cross, we can use a standard tool known as `crossing diagram' to analyse the way in which the two separations cross; see \cref{Overview:fig:2DcrossingDiagram}.
The crossing-diagram reveals, among other things, that $(C,D)$ induces a \defn{tri-partition} $(X,\{P,Q\})$ of $A\cap B$ into two links $P,Q$ -- namely $(A\cap B)\cap (C\sm D)$ and $(A\cap B)\cap (D\sm C)$ -- and the centre $X:=(A\cap B)\cap (C\cap D)$.

In general, we have little control over what tri-partition $(C,D)$ induces on~$A\cap B$.
Curiously, we show that this changes drastically if $(C,D)$ $\gamma$-crosses~$(A,B)$ with~$\gamma\approx q+k$, 
where \defn{$\gamma$-crossing} roughly means the same as crossing plus that all four corners have size~$\ge\gamma$.
So $\gamma$-crossing separations are essentially what we called heavily crossing in \cref{subsec:brief-overview}.
A set $S$ of separations of $G$ is \defn{uniformly $\gamma$-crossing} if every $s=(A,B)\in S$ admits a tri-partition $\cP(s)$ of $A\cap B$ such that every $(C,D)\in S$ that $\gamma$-crosses $s$ induces the tri-partition $\cP(s)$ on $A\cap B$.
Let $S_0$ denote the set of all $q$-broad $k$-separations of~$G$.
We show that heavy crossings within~$S_0$ are consistent: for every $s\in S_0$, all separations in~$S_0$ that $\gamma$-cross~$s$ induce the same centre and the same two links on the separator of~$s$.

\begin{keylemma}[Rough version of \cref{HugeSeparationsUniformlyCrossing}]\label{OverviewUniformlyCrossing}
    The set $S_0$ is uniformly $\gamma$-crossing.
\end{keylemma}

To prove \cref{OverviewUniformlyCrossing}, we analyse how a given $(A,B)\in S_0$ can be $\gamma$-crossed simultaneously by two other separations from~$S_0$.
For this, we introduce a new method -- a three-dimensional corner diagram -- and use the submodularity property of connectivity in a three-dimensional way.
See \cref{Overview:fig:3Dcube} for an illustrative example of how three separations can cross three-dimensionally.

\subsection{Step~2: Finding the flowers}\label{subsec:overview-step2}
\begin{figure}[ht]
    \centering
    \includegraphics[height=16\baselineskip]{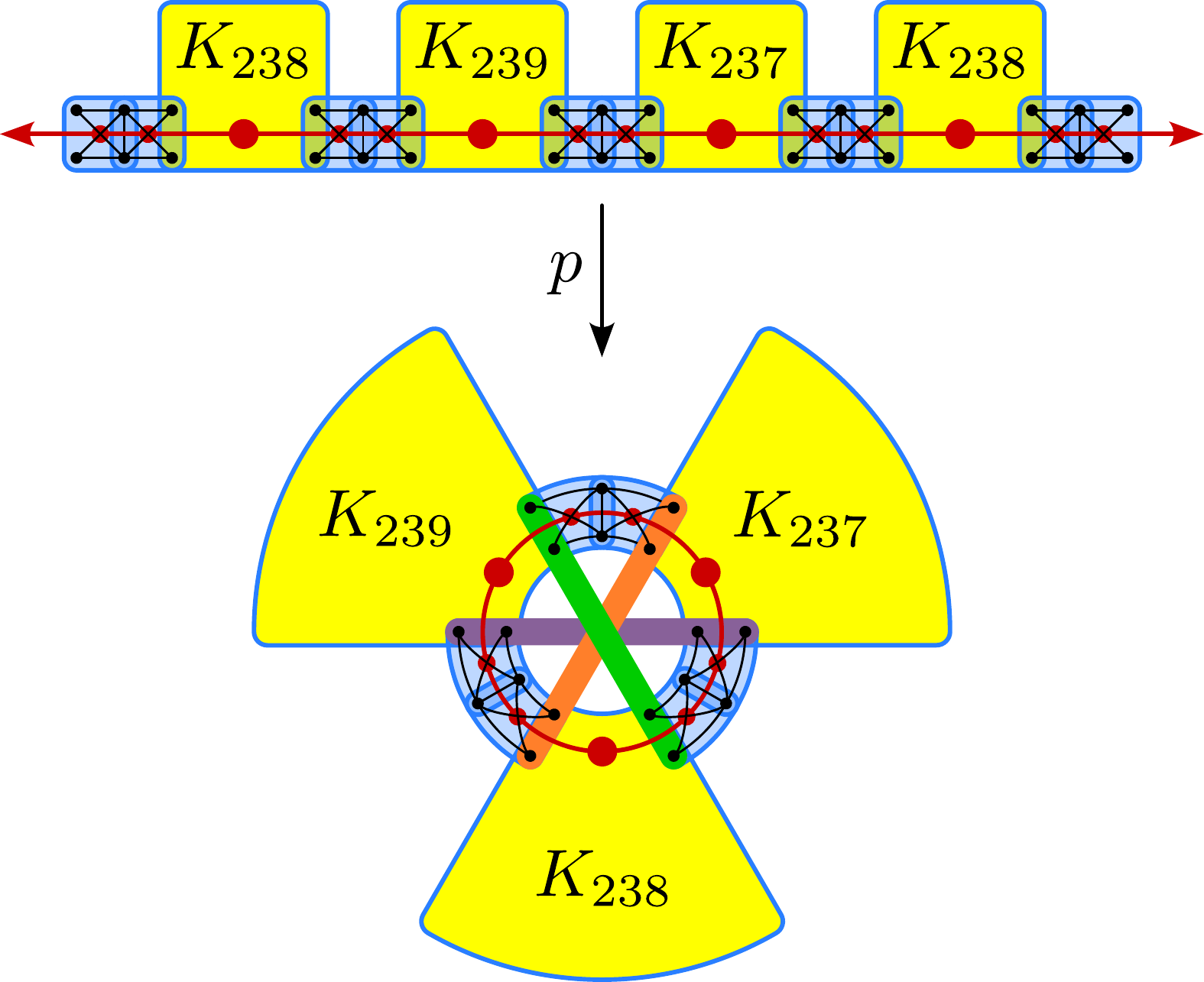}
    \caption{Bottom: The purple, orange and green separators of $G-X$ stem from the same asterisk~$S$ of~$G$. The indicated cycle-decomposition equals the daisy~$F(S)$. Top: The tree-decomposition $(D,\cW)$ of~$\hat G$.}
    \label{Overview:fig:Covering}
\end{figure}

Let $S$ be a nonempty set of separations of~$G$ that is \defn{symmetric} in that it contains for every $(A,B)\in S$ also its \defn{inverse}~$(A,B)^*:=(B,A)$.
We say that $S$ is \defn{$\gamma$-crossing-connected} if the reflexive and transitive closure of the union of the relations `$\gamma$-crossing' on~$S$ and $\{\,(s,s^*):s\in S\,\}$ has exactly one partition class.
Every symmetric set of separations admits a unique partition into classes that are maximally $\gamma$-crossing-connected, and these classes are again symmetric.
If the set is uniformly $\gamma$-crossing, then so are its classes, and we call them its \defn{$\gamma$-asterisks}.
By \cref{OverviewUniformlyCrossing}, this applies to~$S_0$ (see \cref{Overview:fig:Covering} for one of its $\gamma$-asterisks).
A~$\gamma$-asterisk is \defn{trivial} precisely if it consists of only one separation and its inverse.

We capture the nontrivial $\gamma$-asterisks with flowers.
Flowers were introduced by Oxley, Semple and Whittle in their work on decomposing $3$-connected matroids~\cite{Structure3sepsMatroids}.
One of the striking discoveries of that paper is that the left-over pieces of such decompositions, which they called `flowers', are in fact highly structured: by their dichotomy theorem, every flower is an `anemone' or a `daisy'.
Aikin and Oxley extended flowers from $3$-separations to $k$-separations for arbitrary~$k$~\cite{Aikin_flowers}. 
In graphs, flowers can be described in terms of cactus-decompositions, which we now define formally.
A \defn{cactus-decomposition} of~$G$ is a pair $\cC=(C,\cV)$ of a strong cactus~$C$, some of whose nodes are designated as \defn{$P$-nodes}, and a family $\cV=(V_t:t\in V(C))$ of sets of vertices of~$G$, its \defn{bags}, such that $G=\bigcup_{t\in V(C)}G[V_t]$ and, for every vertex~$v$ of~$G$, the nodes~$t$ with $v\in V_t$ induce a connected subgraph of~$C$.
Every \defn{cut}~$(U,W)$ of~$C$, a bipartition of~$V(C)$ in which either part may be empty, \defn{induces} the separation $(\bigcup_{t\in U}V_t,\bigcup_{t\in W}V_t)$ of~$G$.
A separation of~$G$ is \defn{displayed} by~$\cC$ if it is induced by a cut~$(U,W)$ of~$C$ such that $E_C(U,W)$ is empty, consists of a single edge, of two edges on a common cycle of~$C$, or of edges incident with a single $P$-node of~$C$.
Tree-decompositions are the cactus-decompositions whose cactus is a tree without $P$-nodes, so their displayed separations are those induced by single edges.
A \defn{$k$-daisy} of~$G$ is a cactus-decomposition of~$G$ whose cactus is a single cycle and all of whose displayed separations have order~$\le k$; see \cref{fig:flower99}.
A \defn{$k$-anemone} of~$G$ is a cactus-decomposition of~$G$ whose cactus is a star with a $P$-node as its central node, such that the bag of the central node, the \defn{central bag}, has size~$\le k$ and is included in all other bags, the \defn{leaf-bags}; so all displayed separations have order~$\le k$.
We call $k$-daisies and $k$-anemones \defn{$k$-flowers}.
(The formal definition of $k$-daisies in \cref{sec:Flowers} adds a technical condition.)
In practice, $k$-daisies capture the situation where many $k$-separators cut through each other as happens in a cycle, while $k$-anemones capture the situation where a single $k$-separator of~$G$ separates three or more components.
The \defn{centre} of a $k$-flower is the intersection of all its bags; for a $k$-anemone, this is its central bag.
For a node~$t$ of the cactus, the bag~$V_t$ is the \defn{petal} of the $k$-flower at~$t$, the displayed separation $(V_t,\bigcup_{u\neq t}V_u)$ that cuts it off is the \defn{petal-separation} at~$t$, and the \defn{petal-star} of the $k$-flower is the set of all its petal-separations.

\begin{keylemma}[Rough version of \cref{keylem:flower-approximation}]\label{OverviewFlowerTheorem}
    For every nontrivial $\gamma$-asterisk $S$ of~$S_0$ there is a $k$-flower $F(S)$ of~$G$ such that every separation in~$S$ is $\eps_F$-close to some displayed separation of~$F(S)$.
    Moreover, $F(S)$ is \defn{canonical}: every automorphism~$\varphi$ of~$G$ maps $F(S)$ onto $F(\varphi(S))$, bag by bag.
    Here, $\eps_F$ depends only on $q$ and~$k$, not on~$G$.
\end{keylemma}

To prove \cref{OverviewFlowerTheorem}, let $S$ be any nontrivial $\gamma$-asterisk of~$S_0$.
We choose $s=(A,B)\in S$ arbitrarily.
Let $(X,\{P,Q\}):=\cP(s)$, so every $(C,D)\in S$ that $\gamma$-crosses $s$ induces the tri-partition $(X,\{P,Q\})$ on $A\cap B$.
The centre~$X$ is the same for every choice of~$s\in S$ (\cref{coveringCentreAndLinksOfS}); we call it the \defn{common centre} of~$S$. 

If $X$ has size~$k$, then we are in an easy case: here we take $F(S)$ to be the \defn{$X$-anemone}, meaning that $F(S)$ is the $k$-anemone with $X$ as its central bag and a leaf-bag $V(K)\cup X$ for each component $K$ of $G-X$.
Standard arguments quickly verify that we are done in this case.

So we may assume that $|X| < k$.
Then we will find a suitable $k$-daisy to take as~$F(S)$.
To find~it, we build a covering $p\colon \hat G\to G-X$ (as known from general topology), which is depicted in \cref{Overview:fig:Covering}.
Here $G-X$ is the union of $G[A]-X$ and $G[B]-X$, which intersect precisely in the two links~$P$ and~$Q$, much like a circle is the union of two arcs that intersect precisely in their endpoints.
To unroll this circle, we arrange infinitely many copies of $G[A]-X$ and $G[B]-X$ alternately along a two-way infinite path and glue consecutive copies together along copies of~$P$ and~$Q$ in turn; this is~$\hat G$, and $p$ sends every vertex of a copy to the vertex it is a copy of (\cref{CoveringDefinition}).
It is not hard to show that the covering is independent of the choice of~$s$.
We recall from infinite graph theory that an \emph{end} is an equivalence class of one-way infinite paths, where two such paths are \emph{equivalent} if they are linked by infinitely many pairwise vertex-disjoint paths.
Just as $\R$ covers~$\mathbb{S}^1$ and $\R$ is two-way infinite,
also $\hat G$ covers $G-X$ and $\hat G$ is `two-way infinite': $\hat G$ has exactly two ends.
So, $G-X$ `looks like $\mathbb{S}^1$' from the perspective of~$\hat G$, and we will capture this circular structure of $G-X$ with the help of a $k$-daisy.

By standard arguments, $P$ and $Q$ have the same size~$\ell$, say.
Moreover, every separation of $\hat{G}$ that separates the two ends of $\hat{G}$ has order~$\ge\ell$.
We let $R$ be the (nonempty) set of all $\ell$-separations of $\hat{G}$ that separate the two ends of~$\hat{G}$.
To each separation $r\in R$ we assign its \defn{crossing-number} $\cross_R(r)$ which is the number of separations $r'\in R$ that cross~$r$.
We let $L$ be the set of all \defn{locally $R$-optimal} separations $r\in R$: those with $\cross_R(r')>\cross_R(r)$ for all $r'\in R$ that cross~$r$.
This ensures that the separations in $L$ are pairwise nested.
A short proof shows that $L$ is nonempty.
Hence the set $L$ is the recipe for a tree-decomposition $(D,\cW)$ of~$\hat{G}$.
Notably, the decomposition tree $D$ is a double-ray: a two-way infinite path.
We cannot directly project the tree-decomposition $(D,\cW)$ of~$\hat G$ along the covering map~$p$ to obtain a $2\ell$-daisy of~$G-X$.
Instead, we use a folding technique from~\cite{canonicalGraphDec} which exploits the high symmetry of the covering (as encoded in its group of deck transformations) to `fold' $(D,\cW)$ down to $G-X$ where it becomes a $2\ell$-daisy.
For this step, it is crucial that $(D,\cW)$ is canonical, as this ensures that the group of deck transformations admits a well-defined action on the tree-decomposition~$(D,\cW)$.
Finally, we add $X$ to all bags of the $2\ell$-daisy to obtain the desired $k$-daisy~$F(S)$; recall that $|X|+2\ell=k$.

Given an arbitrary $s'\in S$, we find an $\eps_F$-close displayed separation of $F(S)$ as follows.
First, we may assume without loss of generality that $s'=s$, where $s$ was used to construct the covering.
Then, we lift $s'$ to $\hat{G}$, where it splits into two $\ell$-separations~$s'_1,s'_2$.
Each $s'_i$ is contained in~$R$, but might not be in~$L$.
If some~$s'_i$ is not in~$L$, then it is crossed by a separation in~$R$ whose crossing number is not larger, and we replace~$s'_i$ by a corner-separation of the two that lies in~$R$ and has a lower crossing number than~$s'_i$.
We repeat this until $s'_i\in L$.
We show that all crossing numbers in~$R$ are bounded, so that there are only boundedly many replacements, and that any two crossing separations in~$R$ cross with two small opposite corners, so that every replacement changes~$s'_i$ by only a bounded number of vertices.
Hence the final $\ell$-separations $s'_i\in L$ are close to the original ones.
Combining the final $s'_1,s'_2$ and adding $X$ to both sides yields a displayed separation of~$F(S)$ that is $\eps_F$-close to~$s'$.

\subsection{Step~3: Finding the tree-decomposition}\label{subsec:overview-step3}
Let $S_1$ be the set of 
\begin{itemize}
    \item all petal-separations of the $k$-flowers found for the nontrivial $\gamma$-asterisks of~$S_0$ in Step~2, as well as their inverses, and
    \item all $k$-separations contained in trivial $\gamma$-asterisks of~$S_0$.
\end{itemize}
This follows the pattern of Tutte's decomposition, whose totally-nested 2-separations cut off the pieces attached to each cycle and leave the cycle as a torso (\cref{fig:TutteFor2}).
Similarly, many displayed separations of a flower cross each other, so we cannot cut along them; but every displayed separation is induced by a set of petals and hence nested with all petal-separations of the flower.
So cutting along the petal-separations separates the petals from each other without cutting through any displayed separation, and the flower then fits into the bag of the resulting node (Step~4).
It would ease our lives quite a lot if the separations in $S_1$ were pairwise nested, but unfortunately they usually are not.
They do, however, satisfy the following approximate version of Tutte's notion of totally-nested separations.
A proper $(\le k)$-separation of~$G$ is \defn{totally-$\eps$-nested} if it is $\eps$-nested with every proper $(\le k)$-separation of~$G$.
Indeed, we can use three-dimensional corner-diagrams to show that
\begin{keylemma}[Rough version of \cref{keylem:petal-separations-totally-nested}]
\label{OverviewEqSoneTotally}
The separations in $S_1$ are totally-$\gamma$-nested.
\end{keylemma}

\noindent (For \cref{OverviewEqSoneTotally} to hold, we have to slightly increase~$\gamma$ first, which is no problem.)
Since the separations in $S_1$ need not be pairwise nested, we look for the next best thing: a nested set $N$ of $k$-separations of~$G$ such that $N$ $\eps_N$-approximates~$S_1$, where $\eps_N$ is some constant independent of~$G$.

It is not clear to us if $N$ can be found within~$S_1$, but we manage to find a suitable~$N$ within a slightly larger set.
Let $M$ be the set of all totally-$\gamma$-nested $k$-separations of~$G$ (so $S_1\se M$).
Up to the bounded errors, these are the $k$-separations that are not heavily crossed, from which \cref{subsec:brief-overview} built the tree-decomposition.
Let $M'$ be the set of all totally-$\gamma'$-nested $k$-separations of~$G$ that are $\eps_N$-close to a separation in~$M$, where $\gamma'\approx k(k+\gamma)\approx\eps_N$.
Then let $N$ be a maximal nested subset of~$M'$, and let $(T,\cV)$ be the tree-decomposition that is uniquely determined by~$N$.
We will show more generally that
\begin{keylemma}[Rough version of \cref{keylem:max-nested-set}]\label{OverviewNapproximatesM}
    $N$ $\eps_N$-approximates~$M\supseteq S_1$.
    Moreover, every $\eps_N$-broad $k$-separation of~$G$ is $\eps_N$-close to some $k$-separation of~$G$ that is nested with~$N$.
\end{keylemma}
\begin{cor}[Rough version of \cref{canonical:max-nested-set}]
\label{OverviewCanonical}
    $N$ is $2\eps_N$-canonical.
\end{cor}
Here, a set~$S$ of separations of~$G$ is \defn{$\eps$-canonical} if for every automorphism~$\varphi$ of~$G$, the sets~$S$ and~$\varphi(S)$ $\eps$-approximate each other; for $\eps=0$, we call~$S$ \defn{canonical}.
Similarly, a tree-decomposition is (\defn{$\eps$-})\defn{canonical} if the set of separations induced by the edges of its decomposition tree is ($\eps$-)canonical.

For a node~$t$ of~$T$, each edge~$e$ of~$T$ at~$t$ induces the separation $(\bigcup_{u\in V(T_1)}V_u,\bigcup_{u\in V(T_2)}V_u)$ of~$G$, where $T_2$ is the component of $T-e$ that contains~$t$ and $T_1$ is the other component; these separations form the \defn{$t$-star} of~$(T,\cV)$.
To show \cref{OverviewNapproximatesM}, let any $(A,B)\in M$ be given.
If $(A,B)$ is $\gamma$-close to a separation in~$N$, we are done.
Otherwise, using that $(A,B)$ is totally-$\gamma$-nested, we find a unique $t$-star $\sigma$ of $(T,\cV)$ such that $(A,B)$ \defn{almost-$\gamma$-interlaces} $\sigma$, meaning that $s\le_\gamma (A,B)$ or $s\le_\gamma (B,A)$ for all $s\in \sigma$. 

If $(A,B)$ crosses at least about $3k/2$ of the separations in~$\sigma$, then we deduce that
$A\cap B$ is essentially a totally-nested separator of~$G$, and show with fairly standard arguments that $(A,B)$ is close to a separation in~$N$. 

Otherwise, $(A,B)$ crosses at most about $3k/2$ of the separations in~$\sigma$.
Then we `uncross' $(A,B)$ with those separations, which yields a $k$-separation $(A',B')$ that is $\eps_N$-close to the original $(A,B)$ but satisfies $s\le (A',B')$ or $s\le (B',A')$ (without the $\gamma$-error) for every $s\in\sigma$.
Essentially by maximality of~$N$, we have $(A',B')\in N$.
Hence $N$ $\eps_N$-approximates~$M\supseteq S_1$.

The `moreover' part of \cref{OverviewNapproximatesM} is proved by the same case analysis, applied to an arbitrary $\eps_N$-broad $k$-separation~$(A,B)$ of~$G$ instead of to a separation in~$M$.
Since every separation in~$N$ is totally-$\gamma'$-nested, $(A,B)$ is $\gamma'$-nested with every separation in~$N$.
What we can no longer conclude is that the resulting separation $(A',B')$ lies in~$N$: the maximality of~$N$ only applies to separations in~$M'$, and $(A',B')$ need not lie in~$M'$ as $(A,B)$ need not lie in~$M$.
However, $(A',B')$ will be nested with~$N$.

\subsection{Step~4: The structure of the torsos}\label{subsec:overview-step4}
Let $\sigma\se N$ be an arbitrary $t$-star of~$(T,\cV)$, and let $H$ denote the torso of~$V_t$.
We distinguish between two cases.

In the first case, every $k$-separation $(A,B)$ of~$G$ that `cuts through'~$H$ is sufficiently close to some separation in~$\sigma$, so that either $A\cap V(H)$ or $B\cap V(H)$ is small.
In this case, we use standard arguments to deduce that the torso $H$ is $(q',k)$-unbreakable for some constant $q'$ that depends only on $q$ and~$k$.

In the second case, there exists a $k$-separation $(A,B)$ of~$G$ such that $A\cap V(H)$ and $B\cap V(H)$ are large --  large enough that $(A,B)$ is not $\eps_N$-close to any separation in~$\sigma$, and hence not $\eps_N$-close to any separation in~$N$ either.
Then, since $N$ $\eps_N$-approximates~$M$ by \cref{OverviewNapproximatesM}, we know that $(A,B)$ must be $\gamma$-crossed by some other $k$-separation of~$G$, so $(A,B)$ belongs to a nontrivial $\gamma$-asterisk $S$ of~$S_0$.
By \cref{OverviewFlowerTheorem}, the $k$-flower $F(S)$ has a displayed separation~$d$ that is $\eps_F$-close to~$(A,B)$. 
If $d$ is close to some petal-separation of $F(S)$ (which is contained in~$S_1$), then $N$ approximates $d$ with some parameter slightly larger than~$\eps_N$ (by \cref{OverviewNapproximatesM} and quasi-transitivity of closeness), which allows us to push this outcome to the first case by tricking with the constants.
So we may assume that $d$ is sufficiently far from all petal-separations of~$F(S)$.
Then we employ the flower-structure of~$F(S)$ to construct another displayed separation $d'$ of $F(S)$ that $\delta$-crosses~$d$, where $\delta\in\N$ is much larger than all previous constants, and 
we define $F_t:=F(S)$.
Thus $F_t$ has \defn{spread} at least~$\delta$, meaning that two of its displayed separations $\delta$-cross.
As we may choose $\delta$ as large as we like, Step~5 will be able to exploit this.

\begin{figure}[ht]
    \centering
    \includegraphics[width=.9\linewidth]{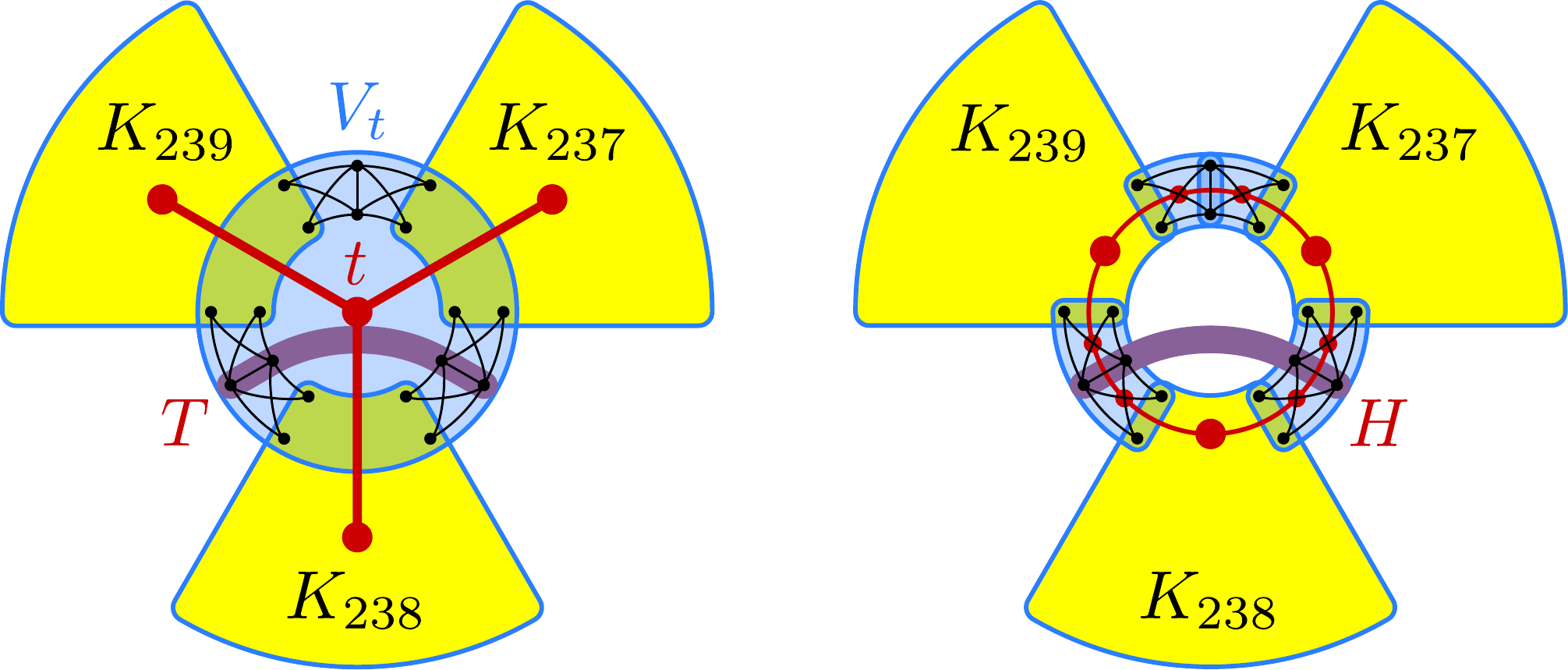}
    \caption{$G$ consists of the black vertices and edges together with the three yellow cliques $K_{237}$, $K_{238}$ and $K_{239}$. A tree-decomposition is indicated on the left: the decomposition tree~$T$ is red and the bags are circled in blue. The bag~$V_t$ of the central node~$t$ is refined by a daisy, whose decomposition cycle~$H$ is indicated on the right. The 4-separation with the purple 4-separator is a displayed separation of that daisy.}
    \label{fig:refines}
\end{figure}

It is not clear a priori that $F_t$ is well-defined, as $(A,B)$ was arbitrary.
However, by heavily exploiting the fact that $\delta$ is much larger than all previous constants, we can show that there is no competitor~$F(S')$ at the same node~$t$.
Finally, we show that $F_t$ fits into the bag~$V_t$ in the following sense.
A $k$-flower~$F$ \defn{$\eps$-refines} $(T,\cV)$ at~$t$ if the petal-star of~$F$ and the $t$-star of~$(T,\cV)$ $\eps$-approximate each other; see \cref{fig:refines}.
We show that $F_t$ $\eps$-refines $(T,\cV)$ at~$t$ for some $\eps\in\N$ that depends only on $q$ and~$k$.

Together, Steps~1--4 yield our technical main result, an approximate Tutte-type decomposition for all~$k$.

\begin{theorem}\label{mainTutte}
    For every $q,k\in\N$ there are $\gamma,q',\eps, \kappa\in\N$ such that, for every $(q,k-1)$-unbreakable graph~$G$, there exists a nested set $N$ of totally-$\gamma$-nested $(\leq k)$-separations of~$G$ with corresponding tree-decomposition $(T,\cV)$ of~$G$ such that:
    \begin{enumerate}
        \item\label{mainTutteTorsoStructure} every node $t$ of~$T$ either has a $(q',k)$-unbreakable torso, or comes with a $k$-flower $F_t$ of $G$ such that $F_t$ $\eps$-refines $(T,\cV)$ at~$t$;
        \item\label{mainTutteSeparationAnalysis} every $\eps$-broad $k$-separation of $G$ is $\eps$-close to either a separation in~$N$ or a displayed separation of some~$F_t$;
        \item\label{mainTutteCanonical} the family of all flowers~$F_t$ is $(q,k)$-canonical, and both the tree-decomposition $(T,\cV)$ and the set~$N$ are $\kappa$-canonical.
    \end{enumerate}
\end{theorem}

Informally, a family of flowers is $(q, k)$-canonical if it is canonical up to a bounded number of vertices in each flower, all of which lie in small components of $G$ minus the centre of the flower.
This weakening is the price for dropping the narrowness assumption at the end of the paper; for narrow graphs, the flowers are canonical (\cref{mainForNarrow}). 
See \cref{qkcanonical} for the formal definition.
See \cref{mainGeneral} for a technical strengthening. 

\subsection{Step~5: Refining the tree-decomposition into a cactus-decomposition}\label{subsec:overview-step5}
By \cref{mainTutte}~\ref{mainTutteSeparationAnalysis}, every sufficiently broad $k$-separation of~$G$ is close to a separation in~$N$ or to a displayed separation of one of the flowers~$F_t$.
The separations in~$N$ are displayed by~$(T,\cV)$, but those of the flowers need not be, and it is in order to display these as well that we pass from a tree to a cactus: we would like to substitute, for every node~$t$ that comes with a flower~$F_t$, the cactus of~$F_t$ for~$t$, which would turn $(T,\cV)$ into a cactus-decomposition. 
The obstacle is that $F_t$ is a flower of the whole graph~$G$, so a displayed separation~$d$ of~$F_t$ need not \defn{almost-interlace} the $t$-star~$\sigma$ of~$(T,\cV)$, that is, satisfy $s\le d$ or $s\le d^*$ for all $s\in\sigma$, and hence need not `live' inside the bag~$V_t$. 

Call a displayed separation of~$F_t$ \defn{active} if it is $\delta$-crossed by another displayed separation of~$F_t$, and \defn{passive} otherwise; here $\delta$ is the large constant from Step~4.
The passive ones cause no trouble. 
Indeed, those displayed by the empty cut are already displayed by~$(T,\cV)$.
Every other passive separation is close to a petal-separation of~$F_t$ or its inverse, or $\delta$-crossed by another displayed separation of~$F_t$ (\cref{FlowerClosePetalOrCrossed}), so every passive one is close to a petal-separation or its inverse; and these are approximated by~$N$ already, as $F_t$ $\eps$-refines $(T,\cV)$ at~$t$ by \cref{mainTutte}~\ref{mainTutteTorsoStructure}.

Since $F_t$ has spread at least~$\delta$, it has two active separations that $\delta$-cross each other, and we can show that all active separations lie in a single asterisk of the set of displayed separations of~$F_t$ (\cref{NewerUniqueAsterisk}).
This is where the spread bound from Step~4 comes in, and it makes the active separations crossing-connected.
Being $\delta$-crossed, every active~$d$ is also $\eps_N$-broad, so the `moreover' part of \cref{OverviewNapproximatesM} applies to it: there is a $k$-separation~$\tilde d$ of~$G$ that is $\eps_N$-close to~$d$ and nested with~$N$.
We fix such a~$\tilde d$ for every active~$d$ so that $d\mapsto\tilde d$ commutes with inversion, and we show that each~$\tilde d$ almost-interlaces~$\sigma$ and hence lives inside the bag~$V_t$.
As $\delta$ is much larger than~$\eps_N$, the separations~$\tilde d$ are crossing-connected as well.

To bundle the $k$-separations~$\tilde d$ into a flower inside~$V_t$, we let $\tilde G_t$ be the graph obtained from~$G$ by making $G[A]$ complete for every $(A,B)\in\sigma$.
Then $\tilde G_t$ again satisfies our assumptions on~$G$, each~$\tilde d$ is still a separation of~$\tilde G_t$ (this is what we needed the almost-interlacing for), and now \emph{every} separation of~$\tilde G_t$ almost-interlaces~$\sigma$.
So we may apply \cref{OverviewFlowerTheorem} a second time, to~$\tilde G_t$ in place of~$G$ and to the asterisk of~$\tilde G_t$ that contains the~$\tilde d$, and obtain a $k$-flower~$\tilde F_t$ that approximates all the~$\tilde d$ and all of whose displayed separations almost-interlace~$\sigma$.
Every~$A$ with $(A,B)\in\sigma$ is complete in~$\tilde G_t$, so it lies in a bag of~$\tilde F_t$, and we attach the branch of~$T$ that $(A,B)$ cuts off at the node of that bag.
Thus substituting~$\tilde F_t$ for~$t$ at every node~$t$ that comes with a flower, and intersecting the bags of~$\tilde F_t$ with~$V_t$, now turns $(T,\cV)$ into the desired cactus-decomposition.
Separations with a small side are close to a separation displayed by the empty cut.
Chasing an arbitrary sufficiently broad $k$-separation of~$G$ through the dichotomy above -- active separations via the corresponding~$\tilde d$ to a displayed separation of~$\tilde F_t$, passive ones back into~$N$ -- shows that the separations displayed by this cactus-decomposition approximate the set of \emph{all} $k$-separations of~$G$.
This proves \cref{thm:main-intro}.

\subsection{Organisation of the paper}

In \cref{sec:Flowers}, we introduce flowers for graphs and prove a dichotomy theorem for them.
In \cref{sec_3d_corner_new}, we introduce three-dimensional corner-diagrams, which we then use in \cref{sec:BroadSeparationsCrossUniformly} to show that broad $k$-separations are uniformly $\gamma$-crossing (\cref{OverviewUniformlyCrossing}).
In \cref{sec:covering}, we employ coverings from general topology to construct canonical flowers for asteriskal sets of separations, and thereby prove \cref{OverviewFlowerTheorem}.
In \cref{sec:petal-separations-totally-nested}, we show that every $k$-separation is $\gamma$-nested with the petal-separations of these flowers (\cref{OverviewEqSoneTotally}).
In \cref{sec:max-nested-set}, we construct the tree-decomposition $\cT$ for \cref{mainTutte} and prove \cref{OverviewNapproximatesM} and \cref{OverviewCanonical}.
In \cref{sec:TorsoStructure}, we analyse the structure of the torsos of~$\cT$, and show that each splitting star either has a torso of higher connectivity, in the sense of unbreakability, or is approximated by the petal-star of a flower.
In \cref{sec:main-theorem-for-narrow}, we prove \cref{mainTutte} under the additional assumption that the graph is $(\beta,k)$-narrow and establish \cref{thm:main-intro} in the narrow case.
In \cref{sec:RemoveNarrow}, we remove this assumption and prove \cref{thm:main-intro} in full generality.
Finally, in \cref{sec:TutteNarrow}, we remove this assumption to prove \cref{mainTutte} in the general case.
The appendix proves the unbreakable-minor lemma used to iterate our decomposition within minors of the original graph.

\section{Separations and Flowers}\label{sec:Flowers}

In this section we shall define flowers and prove some basic properties thereof.
Flowers were introduced by Oxley, Semple and Whittle when studying the structure of $3$-separations in $3$-connected matroids~\cite{Structure3sepsMatroids}. 
The case of $4$-connectivity was studied in~\cite{Structure4sepsMatroids,4connectedMatroids}.
Aikin and Oxley~\cite{Aikin_flowers} generalized flowers to arbitrary orders and classified their structure.
Clark and Whittle extended this concept to $k$-separations in vertically $k$-connected matroids~\cite{TanglesTreesFlowers}. 
When specialised to graphs this definition reads as follows (modulo some tiny details which are not relevant here). 
Given $k\in\N$ and a graph $G$, a $k$-flower (which we call a \defn{pseudo $k$-flower} in this paper) 
is a cyclic family $\cF = (P_i : i \in \Zbb_\ell)$ such that for every interval\footnote{It is common in the literature to only require this property for intervals of size at most two, and then deduce it for arbitrary intervals using submodularity.
Indeed, for this one only needs the mild conditions that all petals contain at least $k+1$ vertices and $G$ is $k$-connected.
Hence this is a minor technical detail, and the presentation given here is slightly easier for the reader to follow in the context of this paper.} $I \subseteq \Zbb_\ell$, we have that $(\bigcup_{i\in I} P_i, \bigcup_{j \notin I} P_j)$ is a $(\leq k)$-separation of~$G$.
Pseudo $k$-flowers were previously studied in the setting of $k$-connected graphs but in this paper we study them in the larger class of $(q,k-1)$-unbreakable graphs (for some $q\in \Nbb$).
Here we shall define flowers slightly differently than Whittle and Clark, and we shall show that these two definitions are equivalent for $k$-connected graphs in \cref{lem:flowers-defs-basically-agree} below, and in \cref{eg_add} we shall explain a technical difference between them.

\smallskip

In order to give our definition of flowers, we rely on the following standard notation.

\subsection{Separations}

Let~$G$ be a graph.
A \defn{separation} of~$G$ is a pair~$(A,B)$ with $G[A] \cup G[B] = G$ (so there is no edge in $G$ between $A\sm B$ and $B\sm A$).
The \defn{separator} of a separation~$s=(A,B)$ is the intersection $\defnMath{\sep(s)}:=A \cap B$.
We write $\defnMath{s^*}:=(B,A)$ for the \defn{inverse} of~$s$.
A set~$S$ of separations is \defn{symmetric} if $s^* \in S$ whenever $s \in S$.
The \defn{breadth} of a separation $(A,B)$ is the minimum of $|A\sm B|$ and $|B\sm A|$.
We say that $(A,B)$ is \defn{$b$-broad} to mean that $(A,B)$ has breadth~$\ge b$.
If $S$ is a set of separations, we denote by \defn{$S^*$} the set $\{s^* : s \in S\}$.
A separation~$(A_1,A_2)$ of a graph~$G$ is \defn{proper} if $A_i \sm A_{3-i} \neq \emptyset$ for $i \in \{1,2\}$.
Note that~$s$ is proper if and only if~$s$ is 1-broad.
A proper separation~$(A,B)$ of~$G$ is \defn{half-connected} if $G[A\sm B]$ or $G[B\sm A]$ is connected.
A separation~$(A_1,A_2)$ of~$G$ is \defn{tight} if for both $i \in \{1,2\}$ there is a component~$K_i$ of $G[A_i \sm A_{3-i}]$ with $N_G(K_i)=A_1 \cap A_2$.
Note that tight separations are proper.
A set~$S$ of separations of~$G$ is \defn{tight} if every separation in~$S$ is tight.
A separation $(A,B)$ of $G$ is \defn{narrow} if every component of $G - (A\cap B)$ has neighbourhood equal to~$A\cap B$.
Note that if~$(A,B)$ is narrow and proper, then it is tight.
Finally, for $b \in \mathbb{N}$, $G$ is \defn{$(b,k)$-narrow} if every $b$-broad $k$-separation is narrow.

Given two separations $s=(A,B)$ and $t=(C,D)$ of a graph~$G$, we define
\[
    s \,\defnMath{\wedge}\, t := (A \cap C, B \cup D)
    \quad\text{and}\quad
    s \,\defnMath{\vee}\, t := (A \cup C, B \cap D)~. 
\]
Note that $s \vee t = (s^* \wedge t^*)^*$.
Note that both $s \wedge t$ and $s \vee t$ again are separations of~$G$.
Moreover, the set of all separations of~$G$ together with~$\wedge$ and~$\vee$ forms a lattice.
Additionally, we define $\bigvee_{s \in \emptyset}s = (\emptyset, V(G))$ and $\bigwedge_{s \in \emptyset}s = (V(G), \emptyset)$.

\subsection{Stars and torsos}
Let~$G$ be a graph and~$S, \sigma$ be sets of separations of~$G$.
We say that~$\sigma$ is a \defn{star} if $s \leq t^*$ for all distinct $s,t \in \sigma$.
Let~$s$ be a separation of~$G$.
We say that~$s$ \defn{interlaces}~$\sigma$ if for every $t \in \sigma$, we have $t < s$ or $t < s^*$.
We say that~$s$ \defn{almost-interlaces}~$\sigma$ if for every $t \in \sigma$, we have $t \le s$ or $t \le s^*$.
We say that~$S$ (\defn{almost-})\defn{interlaces}~$\sigma$ if every separation $s \in S$ (almost-)interlaces~$\sigma$.
We say that~$S$ is \defn{consistent} if there do not exist $s,t \in S$ with $s^* \leq t$.
A star~$\sigma \subseteq S$ is \defn{splitting in~$S$} if it is consistent and no separation in~$S$ interlaces~$\sigma$.
The \defn{torso} of a star~$\sigma = \{(A_i,B_i)\}_{i \in I}$ is the graph obtained from $G[\bigcap_{i \in I} B_i]$ by making every separator $A_i \cap B_i$ into a clique.

\subsection{Graph-decompositions}
Let $G$ and $H$ be graphs, and let $\cV=(V_h:h\in V(H))$ be a family of sets~$V_h\se V(G)$.
Following~\cite{canonicalGraphDec}, the pair $(H,\cV)$ is called a \defn{graph-decomposition}, or \defn{$H$-decomposition}, of~$G$ if it satisfies the following two conditions:
\begin{enumerate}[label=(H\arabic*)]
    \item\label{itm:graph-decomposition-1} $G=\bigcup_{h\in H}G[V_h]$;
    \item\label{itm:graph-decomposition-2} for every $v\in V(G)$ the subgraph $H[\,\{h\in V(H):v\in V_h\}\,]$ is connected.
\end{enumerate}
The sets $V_h$ are the \defn{bags} of the $H$-decomposition.
For every edge $e = hh'$ of $H$, the set $V_h \cap V_{h'}$ is the \defn{adhesion-set at $e$}.
Whenever an $H$-decomposition is introduced as $(H,\cV)$, we tacitly assume that $\cV=(V_h:h\in V(H))$.
A \defn{cut} of~$H$ is a bipartition~$(U,W)$ of~$V(H)$, where either part is allowed to be empty.
Every cut $(U,W)$ of~$H$ \defn{induces} a separation $(A,B)$ of~$G$, where $A=\bigcup_{h\in U}V_h$ and $B=\bigcup_{h\in W}V_h$, see \cite[Lemma~3.4]{canonicalGraphDec}.
We also call~$(A,B)$ the separation \defn{induced by~$(U,W)$} or \defn{induced by~$U$}.
Moreover, the separator of $(A,B)$ equals $\bigcup_{uw\in F}(V_u\cap V_w)$ where $F=E_H(U,W)$.
We write $\defnMath{V\rest X}:=\bigcup_{h\in V(X)}V_h$ for all subgraphs~$X\se H$.

By abuse of notation, for every $\ell \geq 2$, we denote by~$\Zbb_\ell$ the graph that is the $\ell$-cycle with vertices in~$\Zbb_\ell$ arranged in the corresponding cyclic order.
For $\ell=2$, this consists of two vertices joined by two distinct parallel edges.

A \defn{cactus-decomposition} of a graph $G$ is a graph-decomposition $(H, \cV)$ of $G$ where $H$ is a strong cactus some of whose nodes are designated as \defn{$P$-nodes}.
Let~$\cG=(H,\cV)$ be a cactus-decomposition of a graph~$G$.
A separation of~$G$ is \defn{displayed} by~$\cG$ if it is induced by a cut~$(U,W)$ of~$H$ such that the set $F:=E_H(U,W)$ satisfies one of the following conditions:
\begin{enumerate}
    \item $|F|\leq1$;
    \item $|F|=2$ and both edges of~$F$ lie on a common cycle of~$H$;
    \item $|F|\geq 2$ and all edges of~$F$ are incident with a common $P$-node of~$H$.
\end{enumerate}
We denote the set of all displayed separations of~$\cG$ by~\defn{$D(\cG)$}.
Unless specified otherwise, we view tree-decompositions and cycle-decompositions as cactus-decompositions with no $P$-nodes.
Note that, if~$H=\Zbb_\ell$ then every displayed separation of~$\cG$ with both sides nonempty is induced by an interval of~$\Zbb_\ell$.
Moreover, if~$H$ is a star whose centre is a $P$-node, every cut~$(U,W)$ of~$H$ induces a displayed separation of~$\cG$.
We make the following definition for all decomposition graphs $H$ (but will only use it when $H$ is a cycle or a star).
For each node $t\in V(H)$, we call the separation of~$G$ induced by $(\{t\},V(H)-\{t\})$ the \defn{petal-separation} of $\cG$ at~$t$. 
The set $\defnMath{\pi(\cG)}$ of all petal-separations of~$\cG$ forms a star of separations, the \defn{petal-star} of~$\cG$; in particular, $\pi(\cG)$ is a nested set.
We refer to bag $V_t$ as the \defn{petal} at $t$.

\subsection{Flowers}
    Let $G$ be a graph and $k\in\N$.
        A \defn{$k$-flower} in~$G$ is a cycle-decomposition~$\cO$ of~$G$ such that every displayed separation has order~$\le k$. 
    A \defn{$k$-daisy} in~$G$ is a $k$-flower $\cD=(O,\cP)$ such that every narrow proper $k$-separation induced by some subset of~$V(O)$ is induced by the vertices of a path contained in~$O$.
    A \defn{weak $k$-anemone} in~$G$ is a $k$-flower such that every induced separation has order~$\leq k$.
    A \defn{(strong) $k$-anemone} in~$G$ is a $k$-flower $\cA=(O,\cP)$ such that the union of the adhesion-sets of~$\cA$ has size~$\leq k$.

\begin{eg}\label{eg_add}
Let $G$ be the $4$-wheel. We obtain a $3$-flower structure on $G$ by taking a cycle-decomposition with four bags that are the four triangles of this graph, with the appropriate cyclic ordering, see \cref{fig:w4flower}; denote it by $(\Zbb_4,(\Delta_i :  i\in \Zbb_4))$. Now $(\Delta_i\mid i\in \Zbb_4)$ is also a pseudo $3$-flower. 
If we delete the degree-four vertex from $\Delta_2$ and $\Delta_4$, we still have a pseudo $3$-flower, which is not a $3$-flower.

\begin{figure}[ht]
    \centering
    \includegraphics[height=6\baselineskip]{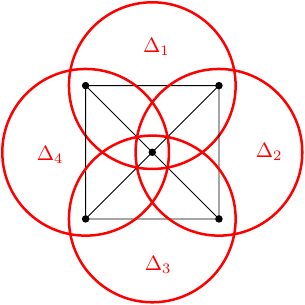}
    \caption{A flower on the $4$-wheel. Its bags $\Delta_1$ to $\Delta_4$ are highlighted in red.}
    \label{fig:w4flower}
\end{figure} 
\end{eg}

The following says that the definition of pseudo-flowers, which is the original definition of flower by Clark and Whittle for matroids specialised to graphs (modulo minor details), is basically the same as the definition of flowers we use in this paper, modulo technicalities arising from examples like the one of \cref{eg_add}.

\begin{obs}\label{lem:flowers-defs-basically-agree} 
    Given $k\in \Nbb$, for every graph $G$ the following holds:
    \begin{enumerate}
        \item\label{itm:flowers-defs-basically-agree-1} If~$(\Zbb_\ell,(P_o : o \in \Zbb_\ell))$ is a $k$-flower, then~$(P_o : o \in \Zbb_\ell)$ is a pseudo $k$-flower.
        \item\label{itm:flowers-defs-basically-agree-2} If~$G$ is $k$-connected and~$(P_o : o \in \Zbb_\ell)$ is a pseudo $k$-flower, then there is a $k$-flower $(\Zbb_\ell,(Q_o : o \in \Zbb_\ell))$ with $P_o\se Q_o$ for all $o\in \Zbb_\ell$.
    \end{enumerate}
\end{obs}

\begin{proof}
\ref{itm:flowers-defs-basically-agree-1} holds by definition of $k$-flower. 

\ref{itm:flowers-defs-basically-agree-2}. Assume that~$G$ is $k$-connected and $(P_o : o \in \Zbb_\ell)$ is a pseudo $k$-flower. We obtain $\Qcal=(Q_o : o \in \Zbb_\ell)$ from $(P_o : o \in \Zbb_\ell)$ by greedily adding vertices to the bags $P_o$ as long as the property of being a pseudo $k$-flower is preserved.
It suffices to show that $\cO=(\Zbb_\ell,\Qcal)$ is a graph-decomposition.
First we show~\ref{itm:graph-decomposition-2}.
Suppose for a contradiction that there is a vertex $v\in V(G)$ such that the subgraph~$H$ induced by $\{i \in \Zbb_\ell:v\in Q_i\}$ is not connected.
Since $H$ is not connected, the complement $\Zbb_\ell-H$ of the subgraph $H$ in the cycle $\Zbb_\ell$ has at least two nonempty components; pick two such components and denote their vertex-sets by $I$ and $J$.

\begin{claim}\label{containment}
 One of $\bigcup_{i\in I}Q_i$ and $\bigcup_{j\in J}Q_j$ contains at most $k-1$ vertices.
\end{claim}

\begin{cproof}
Denote the two intervals that make up the complement of $I\cup J$ in $\Zbb_\ell$ by $K_1$ and $K_2$.
By construction, for each of the four indices $x$ in the boundary of $K_1$ and $K_2$, the bag $Q_x$ contains the vertex $v$. 
Consider the separations $s=(A,B)$ and $t=(C,D)$ induced by $I\cup K_1$ and $I\cup K_2$, respectively.
Since $\Qcal$ is a pseudo $k$-flower, both these separations have order $\le k$.
By submodularity one of the separations $s \wedge t$ and $s \vee t$ has order $\le k$, say it is $s \wedge t = (A \cap C, B \cup D)$ (the other case is the same with the roles of `$I$' and `$J$' interchanged).
The separation $(\bigcup_{i\in I} Q_i, \bigcup_{m\notin I} Q_m)$ induced by~$I$ satisfies that its side $\bigcup_{i\in I} Q_i$ is contained in $A\cap C$, and its side $\bigcup_{m\notin I} Q_m$ is contained in $B\cup D$.
However, the containment $\bigcup_{i\in I} Q_i\se A\cap C$ is proper since the vertex $v$ is only in the latter set.
Hence the separation induced by $I$ has order $\le k-1$.
Since $G$ is $k$-connected, one of its sides contains at most $k-1$ vertices. 
Hence one of $\bigcup_{i\in I}Q_i$ or $\bigcup_{j\in J}Q_j$ contains at most $k-1$ vertices.
\end{cproof}

By \cref{containment} and by exchanging the roles of `$I$' and `$J$' if necessary, assume that $\bigcup_{i\in I}Q_i$ contains at most $k-1$ vertices.
We obtain $\Rcal$ from $\Qcal$ by adding the vertex $v$ to all $Q_i$ with $i\in I$.
To see that $\Rcal$ is a pseudo $k$-flower, let $M$ be an interval of $\Zbb_\ell$.
If $M$ and its complement are both not contained in $I$, then the separation induced by $M$ is the same for $\Rcal$ and $\Qcal$.
So assume that one of $M$ or its complement, say $M$, is a subset of $I$.
Then $\bigcup_{m\in M} R_m$ is contained in $\{v\}\cup\bigcup_{i\in I}Q_i$ and thus has size $\le k$.
Thus the separation induced by $M$ for $\Rcal$ has order $\le k$.
Hence $\Rcal$ is a pseudo $k$-flower.
This violates the maximality of $\Qcal$, as desired.

We now show \cref{itm:graph-decomposition-1}.
Pick any $i \in \mathbb{Z}_{\ell}$.
Since $(P_i, \bigcup_{j \neq i}P_j)$ is a separation of $G$, it follows that every vertex is contained in some bag~$P_o \se Q_o$.
Suppose that some edge~$uv$ lies in no bag~$Q_o$.
Let~$I_u, I_v$ be the set of indices of the bags~$Q_o$ containing~$u,v$, respectively.
By \cref{itm:graph-decomposition-2},~$I_u$ and~$I_v$ are disjoint nonempty connected subgraphs of the cycle~$\Z_\ell$.
The separation induced by~$I_u$ separates~$u$ from~$v$, contradicting that~$uv$ is an edge.

Hence $\cO=(\Zbb_\ell,\Qcal)$ is a graph-decomposition, and thus a $k$-flower.
\end{proof}

\begin{lem}\label{flower_is_anemone}
Let $\cO=(\Zbb_\ell,\cP)$ be a $k$-flower such that there is a partition  of $\Zbb_\ell$ into three nonempty intervals $I_1$, $I_2$ and $I_3$ such that all their corresponding induced separations have the same separator~$X$, and this separator has size exactly~$k$.
Then $\cO$ is a strong $k$-anemone.
\end{lem}

\begin{proof}
It suffices to show that for every $i\in \Zbb_\ell$, the separator of the separation $(P_i, \bigcup_{j\neq i} P_j)$ induced by $\{i\}$ is a subset of~$X$.
Suppose not for a contradiction.
Then there is $i\in \Zbb_\ell$ such that this separator contains a vertex $v\notin X$.
Since $\cO$ is a graph-decomposition, the vertex~$v$ is contained in one of the two bags adjacent to~$P_i$; these are the bags $P_{i-1}$ and $P_{i+1}$. 
By reversing the cyclic orientation of~$\Zbb_\ell$ if necessary, we assume that $v\in P_{i-1}$.
By relabelling if necessary, assume that $i\in I_1$, and that~$I_2$ is the interval just before~$I_1$, and~$I_3$ is the one just after.
Let~$M$ be the interval of $\Zbb_\ell$ obtained from $I_2$ by adding the subinterval of~$I_1$ from the initial element of~$I_1$ all the way to $i-1$ (if~$i$ is the initial element of $I_1$, this subinterval is empty). 
By construction, the vertex~$v$ is in the separator of the separation induced by~$M$.
Every vertex of~$X$ is in the separators of the separations induced by~$I_2$ and~$I_3$ by assumption.
Hence~$X$ is contained in $\bigcup_{j\in I_2} P_j$ and contained in $\bigcup_{j\in I_3} P_j$.
Thus~$X$ is contained in the separator of the separation induced by~$M$.
Thus the separator of the separation induced by~$M$ contains $X\cup\{v\}$, and thus has size $\geq k+1$.
Since $M$ is an interval, this contradicts the assumption that $\cO$ is a $k$-flower.
Hence for every $i\in \Zbb_\ell$ the separator of the separation $(P_i, \bigcup_{j\neq i} P_j)$ induced by $\{i\}$ is a subset of $X$, and so $\cO$ is a strong anemone, as desired.
\end{proof}

\begin{rem}
Since in $k$-connected graphs, all $k$-separations are narrow and there are no proper separations of order less than $k$, the notion of daisy we give here specialises to the well-known notion of daisies in $k$-connected graphs in the sense of Clark and Whittle~\cite{TanglesTreesFlowers}; below in \cref{eg99} we explain in detail the motivation behind our definition of `daisy'.
Thus the following theorem implies the well-known Flower-Dichotomy theorem~\cite{TanglesTreesFlowers} for graphs, showing that this central theorem lifts to our slightly more general setting.
\end{rem}

\begin{theorem}[Flower-Dichotomy Theorem]\label{thm:flower-dichotomy}
For every $k$, every $k$-flower is a strong $k$-anemone or a $k$-daisy. 
\end{theorem}

\begin{proof}
Let a $k$-flower $\cO=(\Zbb_\ell,(P_i\mid i\in \Zbb_\ell))$ be given, and assume that it is not a $k$-daisy.
For every $I \subseteq \Zbb_\ell$, we define the \emph{interval catalogue}~$\Mcal$ of~$I$ to be the set of maximal intervals that are contained in one of~$I$ or its complement $\Zbb_\ell\sm I$. Observe that the size of $\cM$ is twice the number of maximal intervals contained in $I$.
Since~$\cO$ is not a $k$-daisy, there is a narrow proper $k$-separation~$s$ induced by some~$I \subseteq \Zbb_\ell$ that is not induced by an interval.
We let~$I \subseteq \Zbb_\ell$ be such that it induces~$s$ and its interval catalogue~$\cM$ is of minimal size.
Denote the separator of~$s$ by~$X$.
Since~$I$ is a non-interval, the family~$\cM$ has size at least three (in fact at least four but we will not use that) and partitions~$\Zbb_\ell$.
We say that an interval $M\in \Mcal$ is \emph{bad} if for all $m\in M$ we have that the petal $P_m$ is contained in $X$.

\begin{claim}\label{catalogue-lemma}
If $M\in \cM$ is not bad, then the separation it induces has separator equal to $X$.
\end{claim}
\begin{cproof}
Let $M\in \cM$ be given and let~$t_M$ be the separation induced by~$M$.
Without loss of generality, assume that~$M$ is contained in~$I$.
Our first aim is to show $\sep(t_M) \subseteq X$. 
Let~$v \in \sep(t_M)$ be arbitrary. 
Then~$v$ is in $\bigcup_{i\in I} P_i$.
Moreover, there is some petal~$P_j$ with $j\notin M$ that contains~$v$.
If~$j \notin I$, then~$v$ is in~$X$ and we are done.
So assume~$j \in I$.
Since $\cO$ is a graph-decomposition, by the maximality of the interval $M$, there is some $j'\in\Zbb_\ell \sm I$ with $v\in P_{j'}$.
Hence $v$ is in $X$.
It follows that $\sep(t_M) \subseteq X$.

Since~$M$ is not bad, the set $\bigcup_{m\in M} P_m$ is not a subset of~$X$.
Since~$s$ is proper, its side $\bigcup_{j \notin I} P_j$ is not a subset of~$X$.
Since $M \subseteq I$, the set $\bigcup_{j\notin M} P_j$ is not a subset of~$X$.
Together with $\sep(t_M)\se X$ it follows that~$t_M$ is proper.
Since~$s$ is narrow, the separation~$t_M$ has the whole of~$X$ as its separator, as desired.
\end{cproof}

\begin{claim}\label{bad-claim}
If $M\in \cM$ is bad, then $I\Delta M$ induces $s$.
\end{claim}
\begin{cproof}    
Since~$s$ is proper, there are elements of~$\cM$, one contained in~$I$, the other in the complement of~$I$, that are not bad. By \cref{catalogue-lemma}, these two elements of $\cM$ induce separations that have the whole of $X$ in their separator. Thus, for every bad $M\in \Mcal$, the symmetric difference $I\Delta M$ induces the separation $s$ induced by $I$.
\end{cproof}

By \cref{bad-claim}, for every bad~$M \in \cM$, the index set~$I \Delta M$ induces~$s$ but has a strictly smaller interval catalogue.
Therefore, no $M \in \cM$ is bad.
Hence, since the separation induced by $I$ is narrow, by \cref{catalogue-lemma} the separator of every separation induced by some $M\in \Mcal$ is equal to~$X$. 

Thus we have constructed a partition of $\Zbb_\ell$ into at least three intervals such that for each of the separations induced by them, their separators are equal to $X$.
By merging adjacent intervals if necessary, we obtain such a family with exactly three elements. 
(Note that in a merge, the separator does not change as for each interval $J$ in our partition, $\left (\bigcup_{j\in J} P_j \right)\sm X$ has the whole of $X$ in its neighbourhood.) 

Hence we can apply \cref{flower_is_anemone} to deduce that the flower~$\cO$ is a strong anemone, as desired. 
\end{proof}

\begin{eg}\label{eg99}
In this example we explain the motivation behind our definition of \lq daisy\rq\ given above.
We approached this by attempting to have the most expressive definition of daisy that still allows for a dichotomy between daisies and anemones, and so that it specialises to the established definition of daisy in the $k$-connected case. 
The cycle-decomposition depicted in \cref{fig:narrow_daisy} is a $6$-flower that is not an anemone, and thus is a $6$-daisy by \allref{thm:flower-dichotomy}; in particular, the term `narrow' in the definition of daisy cannot be dropped.

\begin{figure}[ht]
    \centering
    \includegraphics[height=7\baselineskip]{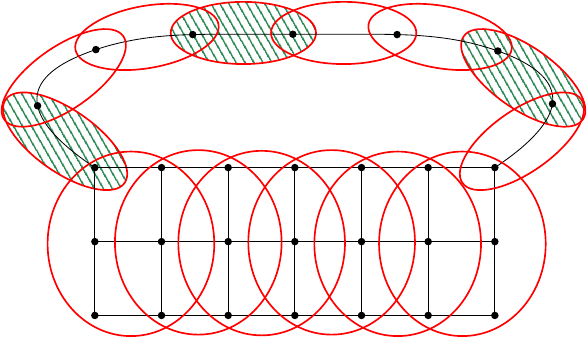}
    \caption{We depict a cycle-decomposition that is a $6$-flower.
    The separation that has the three bags shaded green on one side and all other bags on the other side has order exactly six but is not narrow.}
    \label{fig:narrow_daisy}
\end{figure} 
\end{eg}

\begin{eg}
The cycle-decomposition depicted in \cref{fig:not_ane} is a $3$-flower that is not an anemone, and thus is a daisy. The figure shows that the condition `order exactly $k$' in the definition of daisy cannot be relaxed to `order at most $k$' with the \allref{thm:flower-dichotomy} being true at the same time.

\begin{figure}[ht]
    \centering
    \includegraphics[height=8\baselineskip]{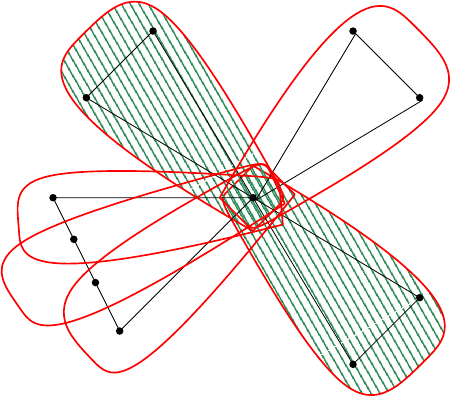}
    \caption{We depict a cycle-decomposition that is a $3$-flower. The separation that has the two bags shaded green on one side and all other bags on the other side has order one and is narrow.}
    \label{fig:not_ane}
\end{figure} 
\end{eg}

\begin{lem}\label{lem:flower-almost-interlaces-then-displayed-new}
    Let~$\cO=(\Zbb_\ell,(P_i\mid i\in \Zbb_\ell))$ be a $k$-flower of a graph~$G$ and let~$s$ be a narrow proper $(\leq k)$-separation that almost-interlaces the star of petal-separations of~$\cO$.
    Then~$s$ is equal to a separation induced by some $I\se \Zbb_\ell$.

    Moreover, if~$s$ has order exactly~$k$ and~$\cO$ is a daisy, then~$I$ can be chosen to be an interval. 
\end{lem}

\begin{proof}
    Let $s=:(C,D)$.
    Let~$I$ be the set of all indices $i\in \Zbb_\ell$ with $P_i\se C$.
    By construction $\bigcup_{i\in I} P_i\se C$.
    For every $j\in \Zbb_\ell \sm I$, we have that~$P_j$ is not a subset of~$C$.
    Since~$s$ almost-interlaces the star of petal-separations of~$\cO$, it follows that $P_j\se D$.
    Hence $\bigcup_{j\notin I} P_j\se D$.
    It follows that the separation induced by~$I$ is equal to~$s$ up to possibly reassigning vertices of the separator~$C\cap D$ to one of the sides.

    Since~$s$ is proper, the separation induced by~$I$ is proper too.
    Since~$s$ is narrow, the separator of the separation induced by~$I$ is equal to $\sep(s)$.
    Hence the separation induced by~$I$ is equal to~$s$.

    The `Moreover'-part is immediate from the definition of daisy.
    This completes the proof. 
\end{proof}

\subsection{Graphic-flowers}
Given a $k$-flower~$\cF=(\Z_\ell,\cP)$, we are interested in all narrow proper $k$-separations that are induced by~$\cF$, call them~$S_0$.
In the case of daisies, they are easy to spot: they are induced by some interval of~$\Z_\ell$.
However, in the case of anemones, they are induced by arbitrary bipartitions of~$\Z_\ell$ and are therefore no longer displayed by the cycle-decomposition~$\cF$.
Moreover, we could arrange the bags of~$\cP$ in an arbitrary order along the cycle.
Therefore, we have no canonical choice for constructing such $k$-anemones.
To circumvent both issues, we instead display all separations in~$S_0$ by a star-decomposition~$\cS$ in the following way.

\begin{definition}[Graphic-$k$-flower]
    A \defn{graphic-$k$-anemone} of a graph~$G$ is a star-decomposition~$\cS=(S,\cV)$ such that the central bag~$V_c$ has size~$\leq k$, and every leaf-bag contains~$V_c$.
    We define its centre~$c$ as a $P$-node, and declare that there are no other $P$-nodes.
    A \defn{graphic-$k$-daisy} of~$G$ is a $k$-daisy of~$G$.
    The term \defn{graphic-$k$-flower} summarizes graphic-$k$-anemones and graphic-$k$-daisies.
\end{definition}

For later use, we include the following definition:

\begin{definition}[Centre]
    The \defn{centre} of a graphic-$k$-flower is the intersection of all its bags.
\end{definition}

\begin{observation}
    \label{obs:daisy-as-graphic-daisy}
    Let~$\cO=(\Z_\ell,\cP)$ be a $k$-daisy of a graph~$G$.
    Then~$\cD=\cO$ is a graphic-$k$-daisy.
    Moreover, every narrow proper $k$-separation induced by~$\cO$ is a displayed separation of~$\cD$.
\end{observation}
\begin{proof}
    That~$\cD$ is a graphic-$k$-daisy is trivial.
    If~$s$ is a narrow proper $k$-separation induced by some $I \subseteq \Z_\ell$ in~$\cO$, then~$I$ can be chosen to be an interval (by definition of `$k$-daisy').
    Hence, $s$ is a displayed separation of~$\cD$.
\end{proof}

\begin{lem}\label{AnemoneAsDaisy}
    Let~$\cA=(S,\cV)$ be a graphic-$k$-anemone of a graph~$G$.
    Then there is a $k$-anemone~$\cO$ such that $D(\cO) \subseteq D(\cA)$ and $\pi(\cO)=\pi(\cA)$.
    Moreover, given two displayed separations $s,s' \in D(\cA)$ of~$\cA$, we can choose~$\cO$ such that $s,s' \in D(\cO)$.
\end{lem}
\begin{proof}
    We let~$O$ be a cycle on the nodes~$V(S)$.
    Then~$\cO:=(O,\cV)$ is a cycle-decomposition.
    It is easy to see that $\pi(\cO)=\pi(\cA)$.
    Let $s$ be a displayed separation of~$\cO$ induced by a cut~$(U, W)$ of $O$.
    Then,~$(U, W)$ induces~$s$ in~$\cA$, hence $s$ is displayed by $\cA$.

    For the `Moreover'-part, let $s,s' \in D(\cA)$  be displayed separations induced by cuts $(U,W)$ and~$(U',W')$ of~$S$, respectively.
    We arrange the nodes of $S$ cyclically in four consecutive blocks $U\cap U'$, $U\cap W'$, $W\cap W'$, and $W\cap U'$, omitting empty blocks and ordering the nodes within each block arbitrarily.
    Then, $s$ and $s'$ are displayed by $\cO$.
\end{proof}

In the following lemma, we prove that every separation displayed by a graphic-$k$-flower~$F$ with both sides nonempty is either close to a petal-separation or crossed significantly by a displayed separation of~$F$.

\begin{lem}\label{FlowerClosePetalOrCrossed}
    Let $F$ be a graphic-$k$-flower of a graph~$G$.
    Let $(A,B)$ be a displayed separation of~$F$ with both sides nonempty.
    Let $\eps\in\N$.
    Then we have at least one of the following:
    \begin{enumerate}
        \item\label{FlowerClosePetalOrCrossed:close} $(A,B)$ is $2\eps$-close to some petal-separation of~$F$ or its inverse; or
        \item\label{FlowerClosePetalOrCrossed:crossed} $F$ has a displayed separation that $(\eps-k)$-crisses~$(A,B)$.
    \end{enumerate}
\end{lem}
\begin{proof}
    By \cref{AnemoneAsDaisy},~$F$ has a cycle decomposition $(O,\cV)$ such that every 2-edge cut of~$O$ induces a separation displayed by~$F$ and there is a 2-edge-cut of~$O$ that induces~$(A,B)$. 
    Recall that $V\rest X:=\bigcup_{x\in V(X)}V_x$ for all subgraphs~$X\se O$.
    Let $P,Q\se O$ be paths with adjacent endvertices such that $A=V\rest P$ and $B=V\rest Q$.
    If $|A| \leq \eps$, then~$(A,B)$ is $\eps$-close to the petal-separation at any node of~$P$.
    Similarly, if $|B| \leq \eps$, then~$(A,B)$ is $\eps$-close to the inverse of the petal-separation at any node of~$Q$.
    Hence, we can assume $|A|,|B| > \eps$.
    
    Let $p$ be the first node on~$P$ such that $|V\rest Pp|>\eps$, which exists since $|A|=|V \rest P| > \eps$.
    If $|V\rest \mathring{p}P|\le\eps$, then $(A,B)$ is $2\eps$-close to the petal-separation $s_p$ of $F$ at~$p$, and so we have outcome~\cref{FlowerClosePetalOrCrossed:close}.
    Hence we have $|V\rest\mathring{p}P|>\eps$ in addition to $|V\rest Pp|>\eps$.
    Similarly, we find a node $q\in V(Q)$ such that $|V\rest Qq|>\eps$ and $|V\rest \mathring{q}Q|>\eps$.
    Let $C:=V\rest (Pp\cup Qq)$ and $D:=V \rest (\mathring{p}P\cup\mathring{q}Q)$.
    Then $(C,D)$ is a displayed separation of~$F$, and therefore has order at most~$k$.
    It remains to show that $(C,D)$ $(\eps-k)$-crisses $(A,B)$.
    By symmetry, it suffices to consider the corner for~$AC$.
    The $AC$-corner is the proper side of the displayed separation of $F$ corresponding to the path~$Pp$ of~$O$. 
    Since $F$ is a graphic-$k$-flower, this separation has order at most $k$.
    Therefore, the $AC$-corner has size $\geq |V \rest Pp| - k > \eps - k$.
    This gives outcome~\cref{FlowerClosePetalOrCrossed:crossed}.
\end{proof}

\section{The 3D Corner Lemma}\label{sec_3d_corner_new}

The aim of this section is to prove \cref{3dcorner}, which can be interpreted as a way of saying that the crossing pattern of separations is essentially 2-dimensional. In a sense it says that under certain assumptions the Venn-diagram of three crossing separations has lots of empty entries, which leads to a 2-dimensional representation of the Venn-diagram when empty regions of the Venn-diagram are omitted.
Our proof strategy is to first prove \cref{combinatorial_3D-corner-lemma}, a purely combinatorial lemma about the cube, and then apply it to the Venn-diagram of three crossing separations to deduce \cref{3dcorner}. 

\subsection{Some combinatorial results on the cube}

We remind the reader of the following standard notation.
If $G$ is a graph, a set $R \subseteq V(G)$ is a \defn{vertex cover} of $G$ if every edge of $G$ is incident with a vertex in $R$. 
The \defn{cube} is the CW complex with vertex set $\{0,1\}^3$ where for each $d\in \{0, 1, 2, 3\}$, its \defn{faces of dimension $d$} are the sets of vertices obtained by fixing $3-d$ coordinates and letting the remaining $d$ coordinates vary freely over $\{0,1\}$. Throughout this section, we denote it by \defn{$Q$}.
Two vertices of the cube are \defn{diagonal} if they differ in every coordinate; each vertex is diagonal to a unique vertex.
Given a partition of $V(Q)$ into two sets $A$ and $B$, the \defn{cut} $c=E(A,B)$ between these two vertex sets
consists of the edges with one endvertex in $A$ and the other in $B$.
$A$ and $B$ are the \defn{sides} of this cut.
A cut is \defn{atomic} if $A$ or $B$ is a singleton.
A cut $c=E(A,B)$ of the cube is called a \defn{coordinate cut} if there is a coordinate such that $A$ consists of all vertices that take the value zero in that coordinate and $B$ consists of all vertices taking the value one in that coordinate, or vice versa.

To state the next lemma, we further introduce the following notation.
We fix a function $\lambda$ from the 1D, 2D and 3D faces of the cube to non-negative integers.
The \defn{$\lambda$-value} of a cut $c$ is given by the sum $\lambda(c)=\sum \lambda(f)$ ranging over all faces $f$ which contain both endpoints of an edge in $c$.
For $k \in \mathbb{N}$, we say that $\lambda$ is a \defn{($k$-)flower} if the following holds:
\begin{enumerate}
    \item $\lambda$ assigns the value zero to all 2D faces; and 
    \item there is a pair of two diagonal vertices such that $\lambda$ assigns the value zero to all 1D faces incident with one of them; and 
    \item $\lambda$ assigns the same value $\ell\in \Nbb$ with $\ell \leq k/2$ to all remaining 1D faces; and
    \item $\lambda$ assigns the value $k-2\ell$ to the 3D face.
\end{enumerate}
We say that $\lambda$ is a \defn{$k$-corner-map} if every coordinate cut has $\lambda$-value $k$.

\begin{restatable}{lemma}{combthreeDcorner}\label{combinatorial_3D-corner-lemma}
Let $\lambda$ be a $k$-corner-map. Let $R$ be a vertex cover of the cube, and assume that every cut with at least one vertex of $R$ on each side has $\lambda$-value at least $k$. 
Then $\lambda$ is a flower. 
\end{restatable}

\begin{figure}[ht]
    \centering
    \includegraphics[height=4\baselineskip]{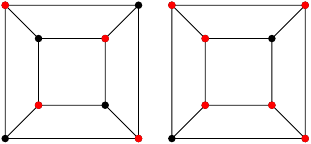}
    \caption{Left: a vertex cover of the cube by an independent set of size four. Right: a vertex cover of the cube by an anti-diagonal cycle.}
    \label{fig:cube}
\end{figure}

In order to prove this lemma, we first prove the following simple fact about the structure of vertex covers of the cube. 
A 6-cycle of the cube is \defn{anti-diagonal} if its complement consists of a pair of diagonal vertices, see \cref{fig:cube}.

\begin{lem}\label{cover-lemma}
A set $X \subseteq V(Q)$ is a vertex cover of the cube $Q$ if and only if $X$ contains an independent vertex set of size four, or induces an anti-diagonal 6-cycle.
\end{lem}

\begin{proof}
The reverse implication is trivial.
Conversely, assume that $X \subseteq V(Q)$ is a vertex cover of $Q$. 
Since we are done otherwise, assume that some vertex of $Q$ is not in $X$, call it $v$.
Since $X$ covers all edges of $Q$, all neighbours of $v$ must be in $X$. Denote by $w$ the unique vertex diagonal to $v$ in $Q$.
If $w\in X$ then the three neighbours of $v$ together with $w$ form an independent set of size four. If $w \notin X$ then all its three neighbours are in $X$, and so $X$ has size six and forms an anti-diagonal 6-cycle. 
\end{proof}

Motivated by \cref{cover-lemma}, we prove \cref{combinatorial_3D-corner-lemma} in the two special cases where $R$ is an independent set of size four and an anti-diagonal $6$-cycle. 
To prepare for this, we first prove the following. We abbreviate $\sum \lambda(f)$ where the sum ranges over all faces $f$ of dimension $d$ by \defn{$s(d)$}.

\begin{lem}\label{corner-count}
    Let $\lambda$ be a $k$-corner-map. Then,
    \[
    3s(3)+2s(2)+s(1)=3k.
    \]
\end{lem}

\begin{figure}[ht]
    \centering
\begin{tikzpicture}[scale=4, every node/.style={circle, draw, fill=white, inner sep=1pt}]
  \coordinate (A) at (0,0);    
  \coordinate (B) at (1,0);    
  \coordinate (C) at (1,1);    
  \coordinate (D) at (0,1);    

  \coordinate (E) at (0.3,0.3);  
  \coordinate (F) at (0.7,0.3);  
  \coordinate (G) at (0.7,0.7);  
  \coordinate (H) at (0.3,0.7);  

  \draw[color=red] (A) -- (B);
  \draw[color=green] (B) -- (C);
  \draw[color=red] (C) -- (D);
  \draw[color=green] (D) -- (A);

  \draw[color=red] (E) -- (F);
  \draw[color=green] (F) -- (G);
  \draw[color=red] (G)-- (H);
  \draw[color=green] (H) -- (E);

  \draw[color=blue] (A) -- (E);
  \draw[color=blue] (B) -- (F);
  \draw[color=blue] (C) -- (G);
  \draw[color=blue] (D) -- (H);

  \node at (A) [label={[fill=white, inner sep=1pt]center:001}] {};
  \node at (B) [label={[fill=white, inner sep=1pt]center:101}] {};
  \node at (C) [label={[fill=white, inner sep=1pt]center:100}] {};
  \node at (D) [label={[fill=white, inner sep=1pt]center:000}] {};
  \node at (E) [label={[fill=white, inner sep=1pt]center:011}] {};
  \node at (F) [label={[fill=white, inner sep=1pt]center:111}] {};
  \node at (G) [label={[fill=white, inner sep=1pt]center:110}] {};
  \node at (H) [label={[fill=white, inner sep=1pt]center:010}] {};
\end{tikzpicture}
    \caption{Illustration of the proof of \cref{corner-count}. The cuts that we consider are displayed in red, blue and green.
    The 3-dimensional face contains an edge from each cut, every 2-dimensional face contains an edge from two of the cuts, and every 1-dimensional face contains an edge from a single cut.}
    \label{fig:proof-corner-count}
\end{figure}
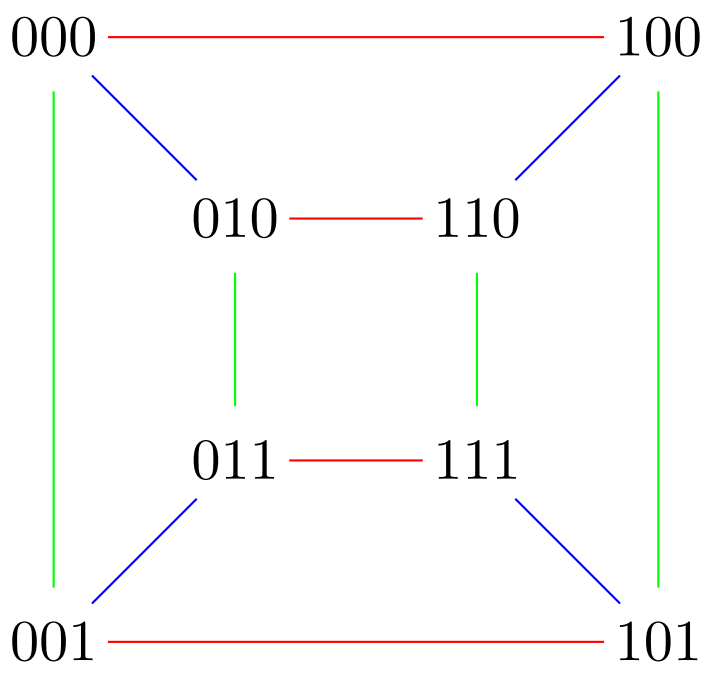

\begin{proof}
Let $c_1, c_2, c_3$ be the three coordinate cuts. 
By assumption, we have $\lambda(c_1) + \lambda(c_2) + \lambda(c_3) = 3k$.
The 3D face contributes to all three terms on the left hand side, every 2D face contributes to two terms and every 1D face contributes to one term (see \cref{fig:proof-corner-count}). Hence the left hand side is equal to $3s(3)+2s(2)+s(1)$, as desired. 
\end{proof}

\begin{lem}\label{combinatorial_3D-corner-lemma_c1}
Let $\lambda$ be a $k$-corner-map. Let $R$ be an independent set of size four of the cube, and assume that every atomic cut of a vertex in $R$ has $\lambda$-value at least $k$. 
Then $\lambda$ is a flower. 
\end{lem}

Given a vertex $r$, we abbreviate the atomic cut $c=E(\{r\},V-r)$ by \defn{$c(r)$}.

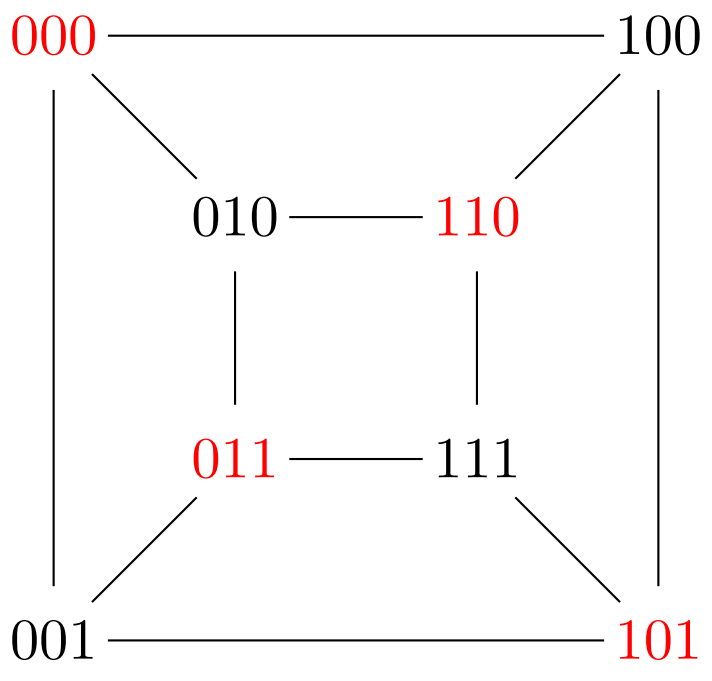
\begin{figure}[ht]
    \centering
\begin{tikzpicture}[scale=4, every node/.style={circle, draw, fill=white, inner sep=1pt}]
  \coordinate (A) at (0,0);    
  \coordinate (B) at (1,0);    
  \coordinate (C) at (1,1);    
  \coordinate (D) at (0,1);    

  \coordinate (E) at (0.3,0.3);  
  \coordinate (F) at (0.7,0.3);  
  \coordinate (G) at (0.7,0.7);  
  \coordinate (H) at (0.3,0.7);  

  \draw (A) -- (B);
  \draw (B) -- (C);
  \draw (C) -- (D);
  \draw (D) -- (A);

  \draw (E) -- (F);
  \draw (F) -- (G);
  \draw (G)-- (H);
  \draw (H) -- (E);

  \draw (A) -- (E);
  \draw (B) -- (F);
  \draw (C) -- (G);
  \draw (D) -- (H);

  \node at (A) [label={[fill=white, inner sep=1pt]center:001}] {};
  \node at (B) [label={[fill=white, text=red, inner sep=1pt]center:101}] {};
  \node at (C) [label={[fill=white, inner sep=1pt]center:100}] {};
  \node at (D) [label={[fill=white, text=red, inner sep=1pt]center:000}] {};
  \node at (E) [label={[fill=white, text=red, inner sep=1pt]center:011}] {};
  \node at (F) [label={[fill=white, inner sep=1pt]center:111}] {};
  \node at (G) [label={[fill=white, text=red, inner sep=1pt]center:110}] {};
  \node at (H) [label={[fill=white, inner sep=1pt]center:010}] {};
\end{tikzpicture}
    \caption{Illustration of the proof of \cref{combinatorial_3D-corner-lemma_c1}.
    The vertices in $R$ are displayed in red. We consider the atomic cuts at the vertices of $R$. The 3-dimensional face contains an edge from each of them, every 2-dimensional face contains an edge from exactly two of them, and every 1-dimensional face contains an edge from exactly one of them.}
    \label{fig:proof-combinatorial_3D-corner-lemma_c1}
\end{figure}

\begin{proof}[Proof of \cref{combinatorial_3D-corner-lemma_c1}]
By assumption, we have $\sum_{r\in R} \lambda(c(r))\geq 4k$.
The 3D face contributes to all four terms on the left hand side, every 2D face contributes to two terms and every 1D face contributes to one term (see \cref{fig:proof-combinatorial_3D-corner-lemma_c1}). Hence the left hand side is equal to $4s(3)+2s(2)+s(1)$, so \[
4s(3)+2s(2)+s(1)\geq 4k.
\]
Subtracting the equation from \cref{corner-count} gives $s(3)\geq k$. Since $s(3)$ contributes to the $\lambda$-value of every coordinate cut, and the $\lambda$-values of coordinate cuts are equal to $k$, by non-negativity of $\lambda$, we get that $s(3)=k$ and all other values contributing to the $\lambda$-value of some coordinate cut must be zero. Since every face contributes to some coordinate cut, $\lambda$ assigns the value $k$ to the 3D face and zero to all other faces, so $\lambda$ is a flower.
\end{proof}

\begin{lem}\label{combinatorial_3D-corner-lemma_c2}
Let $\lambda$ be a $k$-corner-map. Let $R$ be an anti-diagonal cycle of the cube, and assume that every cut with at least one vertex of $R$ on each side has $\lambda$-value at least~$k$.
Then $\lambda$ is a flower. 
\end{lem}

\begin{figure}[ht]
    \centering
\begin{tikzpicture}[scale=4, every node/.style={circle, draw, fill=white, inner sep=1pt}]
  \coordinate (A) at (0,0);    
  \coordinate (B) at (1,0);    
  \coordinate (C) at (1,1);    
  \coordinate (D) at (0,1);    

  \coordinate (E) at (0.3,0.3);  
  \coordinate (F) at (0.7,0.3);  
  \coordinate (G) at (0.7,0.7);  
  \coordinate (H) at (0.3,0.7);  

  \draw (A) -- (B);
  \draw (B) -- (C);
  \draw (C) -- (D);
  \draw (D) -- (A);

  \draw (E) -- (F);
  \draw (F) -- (G);
  \draw (G)-- (H);
  \draw (H) -- (E);

  \draw (A) -- (E);
  \draw (B) -- (F);
  \draw (C) -- (G);
  \draw (D) -- (H);

  \node at (A) [label={[fill=white, text=red, inner sep=1pt]center:001}] {};
  \node at (B) [label={[fill=white, text=red, inner sep=1pt]center:101}] {};
  \node at (C) [label={[fill=white, text=red, inner sep=1pt]center:100}] {};
  \node at (D) [label={[fill=white, inner sep=1pt]center:000}] {};
  \node at (E) [label={[fill=white, text=red, inner sep=1pt]center:011}] {};
  \node at (F) [label={[fill=white, inner sep=1pt]center:111}] {};
  \node at (G) [label={[fill=white, text=red, inner sep=1pt]center:110}] {};
  \node at (H) [label={[fill=white, text=red, inner sep=1pt]center:010}] {};
\end{tikzpicture}
    \caption{Illustration of the proof of \cref{combinatorial_3D-corner-lemma_c2}.
    The vertices in $R$ are displayed in red. We consider the atomic cuts at the vertices of $R$. The 3-dimensional face contains an edge from each of them, every 2-dimensional face contains an edge from exactly three of them. Every 1-dimensional face incident to two vertices of $R$ contains an edge from exactly two of them, and all other 1-dimensional faces contain an edge from exactly one of them.}
    \label{fig:proof-combinatorial_3D-corner-lemma_c2}
\end{figure}
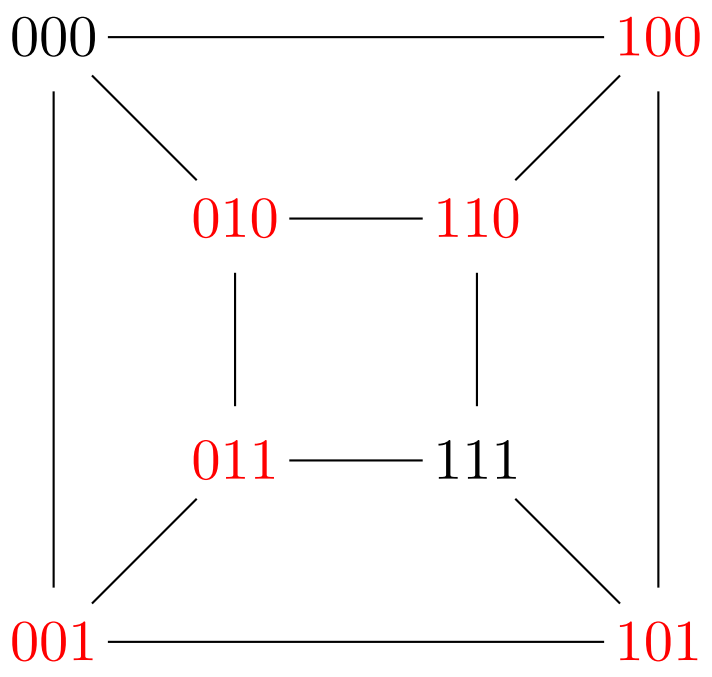

\begin{proof}
This proof follows the same strategy as the proof of \cref{combinatorial_3D-corner-lemma_c1}, with slightly more complicated computations.
By assumption, we have $\sum_{r\in R} \lambda(c(r))\geq 6k$.
The 3D face contributes to all six terms on the left hand side, every 2D face contributes to three terms, every 1D face incident to two vertices in $R$ contributes to two terms, and all other 1D faces contribute to exactly one term since $R$ is a vertex cover (see \cref{fig:proof-combinatorial_3D-corner-lemma_c2}).
Let $s^+(1) = \sum \lambda(f)$ where the sum ranges over all 1D faces with both endpoints in $R$, and $s^-(1) = \sum \lambda(f)$ where the sum ranges over all other 1D faces.
Hence, the left hand side is equal to $6s(3)+3s(2)+2s^+(1)+s^-(1)$, so \begin{equation}\label{eq:comb-3D-c2}
    6s(3)+3s(2)+2s^+(1)+s^-(1)\geq 6k.
\end{equation}

Subtracting twice the equation from \cref{corner-count}, and using that $s(1)=s^+(1)+s^-(1)$, we get 
\[
-s(2)-s^-(1)\geq 0
\]

Since $\lambda$ is non-negative, it follows that $s(2)=0$ and $s^-(1)=0$, so $\lambda$ satisfies the first two conditions in the definition of a flower. Furthermore, \cref{eq:comb-3D-c2} must be an equality. Hence, all the inequalities we used to derive \cref{eq:comb-3D-c2} must also be equalities. Thus, for every $r \in R$, we have $\lambda(c(r)) = k$.
Denote the $\lambda$-value of the 3D face by $x$, and let $k' \coloneqq k-x$. 
Now enumerate the edges (that is, 1D faces) that are on the anti-diagonal cycle by $e_0,\ldots,e_5$ according to a cyclic orientation of that cycle. 
Since $\lambda(c(r)) = k$ for every $r \in R$, we get $\lambda(e_i)+\lambda(e_{i+1})=k'$ for every $i \in \mathbb{Z}_6$. 
Combining these equations at consecutive vertices yields $\lambda(e_i)=\lambda(e_{i+2})$ for all $i\in \Zbb_6$. By renaming the edges if necessary, assume that $\lambda(e_0)\leq \lambda(e_1)$.
Let $u, v \in V(Q)$ be the endpoints of $e_1$, and consider the cut $c \coloneqq E(\{u, v\}, V(Q) \sm \{u, v\})$. Both sides of $c$ contain vertices of $R$, so $\lambda(c) \geq k$.
Since all 2D faces and all 1D faces incident to a vertex outside of $R$ have $\lambda$-value zero, we have $\lambda(c) = x + \lambda(e_0) + \lambda(e_2)$.
Therefore, $\lambda(e_0)+\lambda(e_2)\geq k'$ so $\lambda(e_0)\geq k'/2$. Thus $\lambda(e_0)=\lambda(e_1) = k'/2 \coloneqq \ell$. Hence $\lambda$ assigns the same value $\ell \leq k/2$ at all edges of the anti-diagonal cycle, and value $x = k-2\ell$ to the 3D face. This concludes the proof that $\lambda$ is a flower. 
\end{proof}

We now have all we need to prove \cref{combinatorial_3D-corner-lemma}, which we restate for convenience.

\combthreeDcorner*

\begin{proof}
By \cref{cover-lemma}, it suffices to prove it when $R$ is an independent set of size four or an anti-diagonal 6-cycle. This is done in \cref{combinatorial_3D-corner-lemma_c1,combinatorial_3D-corner-lemma_c2}, respectively.
\end{proof}

\begin{figure}[ht]
\centering
\begin{tikzpicture}[scale=4, every node/.style={circle, draw, fill=white, inner sep=1pt}]
  \coordinate (A) at (0,0);    
  \coordinate (B) at (1,0);    
  \coordinate (C) at (1,1);    
  \coordinate (D) at (0,1);    

  \coordinate (E) at (0.3,0.3);  
  \coordinate (F) at (0.7,0.3);  
  \coordinate (G) at (0.7,0.7);  
  \coordinate (H) at (0.3,0.7);  

  \draw (A) -- (B);
  \draw (C) -- (B);
  \draw[very thick, blue] (C) -- (D);
  \draw (A) -- (D);

  \draw[very thick, red] (E) -- (F);
  \draw (G) -- (F);
  \draw (G) -- (H);
  \draw (E) -- (H);

  \draw (A) -- (E);
  \draw (B) -- (F);
  \draw (C) -- (G);
  \draw (D) -- (H);

  \node at (A) [label={[fill=white, inner sep=1pt]center:001}] {};
  \node at (B) [label={[fill=white, inner sep=1pt]center:101}] {};
  \node at (C) [label={[fill=white, inner sep=1pt]center:100}] {};
  \node at (D) [label={[fill=white, inner sep=1pt]center:000}] {};
  \node at (E) [label={[fill=white, inner sep=1pt]center:011}] {};
  \node at (F) [label={[fill=white, inner sep=1pt]center:111}] {};
  \node at (G) [label={[fill=white, inner sep=1pt]center:110}] {};
  \node at (H) [label={[fill=white, inner sep=1pt]center:010}] {};

\fill[pattern=north east lines, pattern color=red] 
    (0.3,0.3) -- (.3,0.7) -- (.7,.7) -- (0.7,.3) -- cycle;

\fill[pattern=north east lines, pattern color=red] 
    (0.3,0.3) -- (.7,0.3) -- (1,0) -- (0,0) -- cycle;

\fill[pattern=north east lines, pattern color=blue] 
    (0.3,0.7) -- (.7,0.7) -- (1,1) -- (0,1) -- cycle;

\fill[even odd rule, pattern=north east lines, pattern color=blue] 
    (-.15,-.15) rectangle (1.15,1.15) 
    (1,1) rectangle (0,0);
  
\end{tikzpicture}
    \caption{Illustration of \cref{comb-crossing-with-2-nested-empty-links}.
    The vertices $u$ and $v$ are the endpoints of the red edge.
    The faces in red (plain or dashed) have $\lambda$-value $0$ by assumption.
    \cref{comb-crossing-with-2-nested-empty-links} then states that the faces in blue (plain or dashed) also have $\lambda$-value $0$.}
    \label{fig:statement-comb-crossing-with-2-nested-empty-links}
\end{figure}
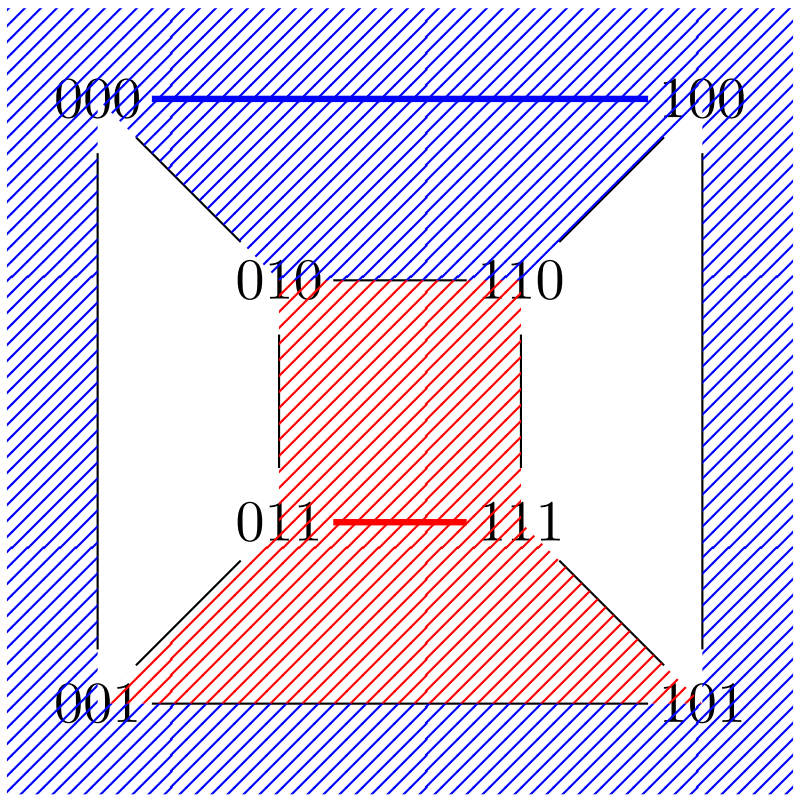

We now prove additional results in the same flavour, which we will later use to analyse how triples of separations can cross.
A cut $c = E(A, B)$ is a \defn{pair cut} if either $A$ or $B$ consists only of two adjacent vertices of $Q$. We then say that these two vertices are \defn{contained} in the pair cut. The statement of the following lemma is illustrated in \cref{fig:statement-comb-crossing-with-2-nested-empty-links}.

\begin{lem}\label{comb-crossing-with-2-nested-empty-links}
Let $\lambda$ be a $k$-corner-map.
Suppose that there exist adjacent vertices $u, v \in V(Q)$ such that all 1D and 2D faces which contain $u$ and $v$ have $\lambda$-value $0$, and all pair cuts which contain exactly one of $u, v$ have $\lambda$-value $k$.
Let $u', v' \in V(Q)$ be the vertices which are diagonal to $u, v$, respectively. Then, all 1D and 2D faces which contain $u'$ and $v'$ have $\lambda$-value $0$.
\end{lem}

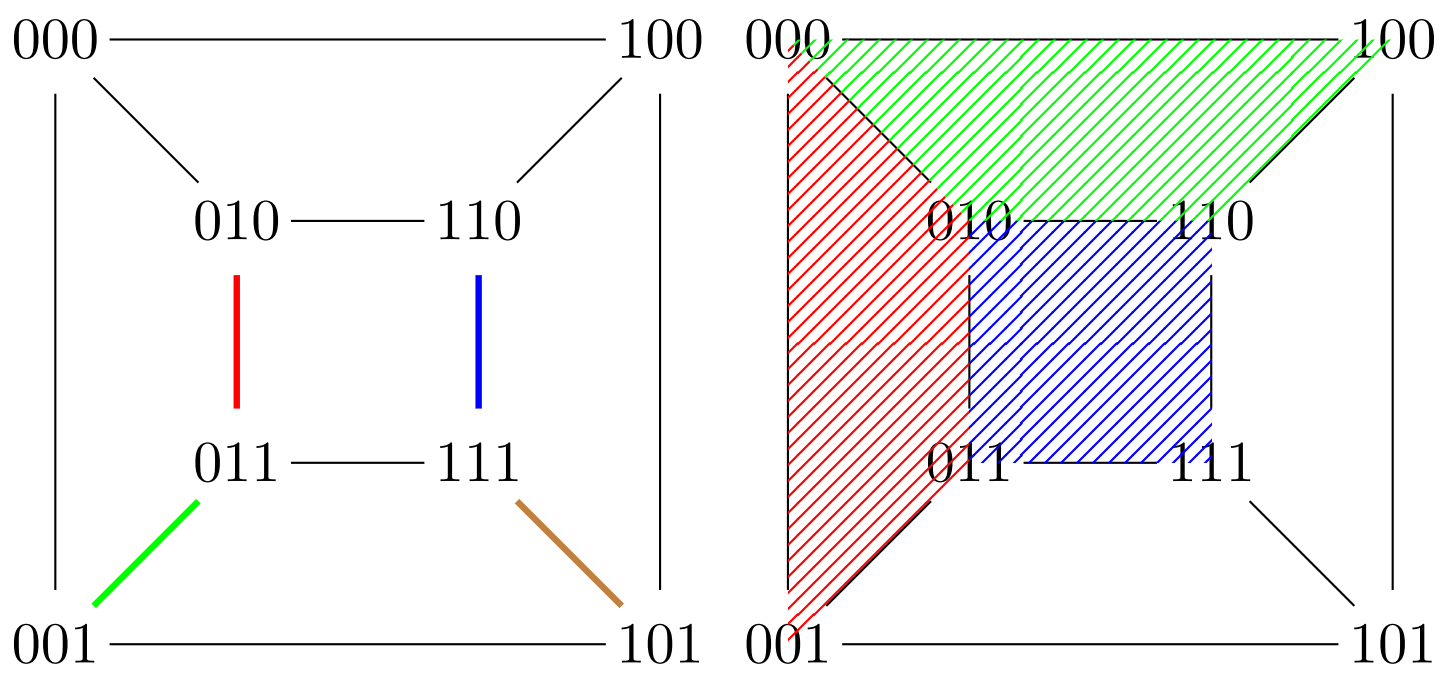
\begin{figure}[ht]
\centering
\begin{tikzpicture}[scale=4, every node/.style={circle, draw, fill=white, inner sep=1pt}]
  \coordinate (A) at (0,0);    
  \coordinate (B) at (1,0);    
  \coordinate (C) at (1,1);    
  \coordinate (D) at (0,1);    

  \coordinate (E) at (0.3,0.3);  
  \coordinate (F) at (0.7,0.3);  
  \coordinate (G) at (0.7,0.7);  
  \coordinate (H) at (0.3,0.7);  

  \draw (A) -- (B);
  \draw (C) -- (B);
  \draw (C) -- (D);
  \draw (A) -- (D);

  \draw (E) -- (F);
  \draw[very thick, blue] (G) -- (F);
  \draw (G) -- (H);
  \draw[very thick, red] (E) -- (H);

  \draw[very thick, green] (A) -- (E);
  \draw[very thick, brown] (B) -- (F);
  \draw (C) -- (G);
  \draw (D) -- (H);

  \node at (A) [label={[fill=white, inner sep=1pt]center:001}] {};
  \node at (B) [label={[fill=white, inner sep=1pt]center:101}] {};
  \node at (C) [label={[fill=white, inner sep=1pt]center:100}] {};
  \node at (D) [label={[fill=white, inner sep=1pt]center:000}] {};
  \node at (E) [label={[fill=white, inner sep=1pt]center:011}] {};
  \node at (F) [label={[fill=white, inner sep=1pt]center:111}] {};
  \node at (G) [label={[fill=white, inner sep=1pt]center:110}] {};
  \node at (H) [label={[fill=white, inner sep=1pt]center:010}] {};
\end{tikzpicture}
\begin{tikzpicture}[scale=4, every node/.style={circle, draw, fill=white, inner sep=1pt}]
  \coordinate (A) at (0,0);    
  \coordinate (B) at (1,0);    
  \coordinate (C) at (1,1);    
  \coordinate (D) at (0,1);    

  \coordinate (E) at (0.3,0.3);  
  \coordinate (F) at (0.7,0.3);  
  \coordinate (G) at (0.7,0.7);  
  \coordinate (H) at (0.3,0.7);  

  \draw (A) -- (B);
  \draw (C) -- (B);
  \draw (C) -- (D);
  \draw (A) -- (D);

  \draw (E) -- (F);
  \draw (G) -- (F);
  \draw (G) -- (H);
  \draw (E) -- (H);

  \draw (A) -- (E);
  \draw (B) -- (F);
  \draw (C) -- (G);
  \draw (D) -- (H);

  \node at (A) [label={[fill=white, inner sep=1pt]center:001}] {};
  \node at (B) [label={[fill=white, inner sep=1pt]center:101}] {};
  \node at (C) [label={[fill=white, inner sep=1pt]center:100}] {};
  \node at (D) [label={[fill=white, inner sep=1pt]center:000}] {};
  \node at (E) [label={[fill=white, inner sep=1pt]center:011}] {};
  \node at (F) [label={[fill=white, inner sep=1pt]center:111}] {};
  \node at (G) [label={[fill=white, inner sep=1pt]center:110}] {};
  \node at (H) [label={[fill=white, inner sep=1pt]center:010}] {};
  \fill[pattern=north east lines, pattern color=red] 
    (0,0) -- (.3,0.3) -- (.3,.7) -- (0,1) -- cycle;

  \fill[pattern=north east lines, pattern color=blue] 
    (0.3,0.3) -- (.3,0.7) -- (.7,.7) -- (0.7,.3) -- cycle;

  \fill[pattern=north east lines, pattern color=green] 
    (0,1) -- (.3,0.7) -- (.7,.7) -- (1,1) -- cycle;
  
\end{tikzpicture}
    \caption{Illustration of the proof of \cref{comb-crossing-with-2-nested-empty-links}.
    On the left, the colored edges correspond to the left side of the separations in $\Ccal_1$.
    On the right, the colored facets correspond to the left side of the separations in $\Ccal_2$.
    For each face $f$ of the cube, we compare the number of colored edges that intersect $f$ without containing $f$, with the number of colored facets that intersect $f$ without containing $f$, counting the red facet twice.
    This number is the same, except for the zones which are highlighted in \cref{fig:statement-comb-crossing-with-2-nested-empty-links}.}
    \label{fig:proof-comb-crossing-with-2-nested-empty-links}
\end{figure}

\begin{proof}
    Let $\mathcal{C}_1$ be the set of all pair cuts which contain exactly one of $u, v$, and let $\mathcal{C}_2$ be the multiset which contains the coordinate cut which contains the edge $uv$ with multiplicity two, and the other two coordinate cuts with multiplicity one.
    By assumption, we have $\sum_{c \in \mathcal{C}_1}\lambda(c) = 4k = \sum_{c \in \mathcal{C}_2} \lambda(c)$.
    We now analyse the contribution of each face to the left hand side and the right hand side, see \cref{fig:proof-comb-crossing-with-2-nested-empty-links}.
    The 3D face contributes to four terms on both sides.
    The 2D faces which contain the edge $uv$ have $\lambda$-value $0$.
    The 2D faces $f_1, f_2$ which contain the edge $u'v'$ contribute to two terms on the left hand side and to three terms on the right hand side.
    The other 2D faces contribute to two terms on both sides.
    The 1D face which contains $u$ and $v$ has $\lambda$-value $0$.
    The 1D face $f_3$ which contains $u'$ and $v'$ does not contribute to the left hand side, and contributes to two terms on the right hand side.
    The 1D faces which are incident to exactly one of $u, v, u', v'$ contribute to one term on both sides.
    The other 1D faces contribute to two terms on both sides.
    Thus, we get \[0 = \lambda(f_1) + \lambda(f_2) +2\lambda(f_3).\]

    Since $\lambda$ is non-negative, this concludes the proof.
\end{proof}

\begin{lemma} \label{lem:comb-3Dcrossing-lem-1}
    Let $\lambda$ be a $k$-corner-map.
    Let $u, v \in V(Q)$ be two vertices of $Q$ which are at distance $2$, and let $w$ be one of their common neighbours.
    Let $x$ be the third neighbour of $w$.
    Suppose that the atomic cuts $c(u)$ and $c(v)$ have $\lambda$-value at least $k$, and that the pair cuts which contain respectively $u$ and $w$, and $v$ and $w$ have $\lambda$-value at most $k$.
    Then, all one- and two-dimensional faces which contain both $w$ and $x$ have $\lambda$-value $0$.
\end{lemma}

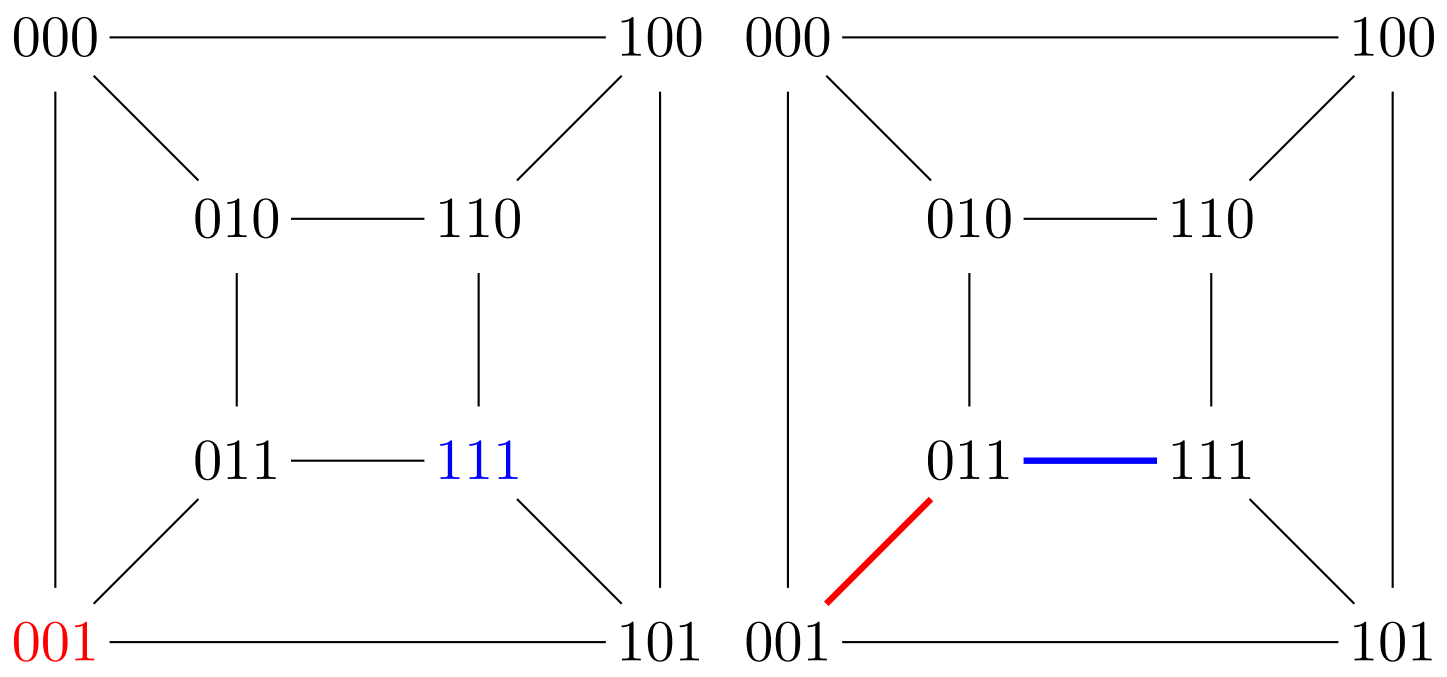
\begin{figure}[ht]
\centering
\begin{tikzpicture}[scale=4, every node/.style={circle, draw, fill=white, inner sep=1pt}]
  \coordinate (A) at (0,0);    
  \coordinate (B) at (1,0);    
  \coordinate (C) at (1,1);    
  \coordinate (D) at (0,1);    

  \coordinate (E) at (0.3,0.3);  
  \coordinate (F) at (0.7,0.3);  
  \coordinate (G) at (0.7,0.7);  
  \coordinate (H) at (0.3,0.7);  

  \draw (A) -- (B);
  \draw (C) -- (B);
  \draw (C) -- (D);
  \draw (A) -- (D);

  \draw (E) -- (F);
  \draw (G) -- (F);
  \draw (G) -- (H);
  \draw (E) -- (H);

  \draw (A) -- (E);
  \draw (B) -- (F);
  \draw (C) -- (G);
  \draw (D) -- (H);

  \node at (A) [label={[fill=white, text=red, inner sep=1pt]center:001}] {};
  \node at (B) [label={[fill=white, inner sep=1pt]center:101}] {};
  \node at (C) [label={[fill=white, inner sep=1pt]center:100}] {};
  \node at (D) [label={[fill=white, inner sep=1pt]center:000}] {};
  \node at (E) [label={[fill=white, inner sep=1pt]center:011}] {};
  \node at (F) [label={[fill=white, text=blue, inner sep=1pt]center:111}] {};
  \node at (G) [label={[fill=white, inner sep=1pt]center:110}] {};
  \node at (H) [label={[fill=white, inner sep=1pt]center:010}] {};
\end{tikzpicture}
\begin{tikzpicture}[scale=4, every node/.style={circle, draw, fill=white, inner sep=1pt}]
  \coordinate (A) at (0,0);    
  \coordinate (B) at (1,0);    
  \coordinate (C) at (1,1);    
  \coordinate (D) at (0,1);    

  \coordinate (E) at (0.3,0.3);  
  \coordinate (F) at (0.7,0.3);  
  \coordinate (G) at (0.7,0.7);  
  \coordinate (H) at (0.3,0.7);  

  \draw (A) -- (B);
  \draw (C) -- (B);
  \draw (C) -- (D);
  \draw (A) -- (D);

  \draw[very thick, color=blue] (E) -- (F);
  \draw (G) -- (F);
  \draw (G) -- (H);
  \draw (E) -- (H);

  \draw[very thick, color=red] (A) -- (E);
  \draw (B) -- (F);
  \draw (C) -- (G);
  \draw (D) -- (H);

  \node at (A) [label={[fill=white, inner sep=1pt]center:001}] {};
  \node at (B) [label={[fill=white, inner sep=1pt]center:101}] {};
  \node at (C) [label={[fill=white, inner sep=1pt]center:100}] {};
  \node at (D) [label={[fill=white, inner sep=1pt]center:000}] {};
  \node at (E) [label={[fill=white, inner sep=1pt]center:011}] {};
  \node at (F) [label={[fill=white, inner sep=1pt]center:111}] {};
  \node at (G) [label={[fill=white, inner sep=1pt]center:110}] {};
  \node at (H) [label={[fill=white, inner sep=1pt]center:010}] {};  
\end{tikzpicture}
    \caption{Illustration of the proof of \cref{lem:comb-3Dcrossing-lem-1}. 
    $u$ and $v$ are the colored vertices, $w$ is the common endpoint of the colored edges.
    On the left, the colored vertices correspond to the left side of the separations in $\Ccal_1$.
    On the right, the colored edges correspond to the right side of the separations in $\Ccal_2$.
    For each face $f$ of the cube, we compare the number of colored vertices that are contained in $f$, with the number of colored edges that intersect $f$ without containing $f$.}
    \label{fig:proof-comb-3Dcrossing-lem-1}
\end{figure}

\begin{proof}
    Let $\mathcal{C}_1 = \{c(u), c(v)\}$, and $\mathcal{C}_2$ contain the pair cuts which contain respectively $u$ and $w$, and $v$ and $w$.
    By assumption, we have $\sum_{c \in \mathcal{C}_1}\lambda(c) \geq 2k \geq \sum_{c \in \mathcal{C}_2} \lambda(c)$.
    We now analyse the contribution of each face to the left hand side and the right hand side, see \cref{fig:proof-comb-3Dcrossing-lem-1}.
    The 3D face contributes to two terms on both sides.
    The 2D face which contains $u$ and $v$ contributes to two terms on both sides.
    The other 2D faces $f_1, f_2$ which contain $x$ contribute to one term on the left hand side and to two terms on the right hand side.
    The 2D face which contains neither $u$ nor $v$ contributes to neither side.
    The last two 2D faces contribute to one term on both sides.
    The 1D face $f_3$ which contains $w$ and $x$ does not contribute to the left hand side, and contributes to two terms on the right hand side.
    The edges which are incident to one of $u, v$ contribute to one term on both sides.
    The other edges contribute to neither side.
    Thus, we get \[0 \geq \lambda(f_1) + \lambda(f_2) +2\lambda(f_3).\]

    Since $\lambda$ is non-negative, this concludes the proof.
\end{proof}

\subsection{Corner diagrams}

Later, we will see how to apply the previous results to deduce properties on how three separations can cross. Before this, we review some notions about corner diagrams, which depict how two separations cross.
Let us start by recalling some notions about corner diagrams.
Consider a graph $G$ and let $r = (A, B)$ and $s = (C, D)$ be two separations of $G$.
These separations induce a partition of $V(G)$, which we represent using a \emph{corner diagram}, see \cref{Overview:fig:2DcrossingDiagram}.
The corresponding partition of $V(G)$ has 9 parts, each of which corresponds to a zone of \cref{Overview:fig:2DcrossingDiagram}.

We now introduce the terminology that we will use to talk about these different parts. Since we will later work with corner diagrams for triples of separations, we need to introduce a new terminology to avoid ambiguity.
Each part is of the form $X \cap Y$ for $X \in \{A \setminus B, A \cap B, B \setminus A\}$ and $Y \in \{C \setminus D, C \cap D, D \setminus C\}$.
We then define the \defn{$AC$-corner} to be the set $(A \setminus B) \cap (C \setminus D)$, and define the \defn{$AD$-corner}, the \defn{$BC$-corner} and the \defn{$BD$-corner} analogously.
The \defn{extended $AC$-corner} is the set $A \cap C$, and the \defn{extended $AD$-corner}, the \defn{extended $BC$-corner} and the \defn{extended $BD$-corner} are defined analogously.
We define the \defn{$rC$-link} to be the set $(A \cap B) \cap (C \setminus D)$, and define the \defn{$rD$-link}, the \defn{$As$-link} and the \defn{$Bs$-link} analogously.
Finally, we define the \defn{$rs$-centre} as the set $(A \cap B) \cap (C \cap D)$.
We denote the size of the $AC$-corner by \defn{$|AC|$}, and define $|AD|, |BC|$ and $|BD|$ similarly.
We denote the size of the $rC$-link by \defn{$|rC|$}, and define $|rD|, |As|$ and $|Bs|$ similarly.
Finally, we denote the size of the $rs$-centre by \defn{$|rs|$}.
A \defn{corner separation} is any separation of the form $(A \cap C, B \cup D)$, up to exchanging $A$ and $B$, and/or $C$ and $D$.
Note that the corner separations are $r \wedge s$, $r \wedge s^*$, $r^* \wedge s$ and $r^* \wedge s^*$.
For simplicity, we call the separation $(A \cap C, B \cup D)$ the \defn{corner separation $AC$}, and similarly for the other corner separations.
The corners $AC$ and $BD$ are \defn{opposite}, and the corners $AD$ and $BC$ as well.
The following result is a fundamental property of the orders of opposite corner separations.

\begin{lemma}[Submodularity, {see e.g. \cite{bibel}}] \label{lem:submodularity}
    Let $(A, B)$ and $(C, D)$ be two separations of a graph.
    Then, \[|(A \cap C) \cap (B \cup D)| + |(A \cup C) \cap (B \cap D)| = |A \cap B| + |C \cap D|.\]
    In particular, if $(A, B)$ and $(C, D)$ are two $k$-separations, then the orders of the opposite corner separations $AC$ and $BD$ sum to $2k$.
\end{lemma}

\subsection{Approximate notions on separations}

Let $G$ be a graph and $\eps\in\N$.
We write $X\,\defnMath{\se_\eps}\,Y$ for $X,Y\se V(G)$ to mean that all but $\le \eps$ vertices in $X$ are contained in~$Y$; that is, $|X\sm Y|\le \eps$.
For $\eps=0$, this defaults to $X\se Y$.
We define $X\,\defnMath{\approx_\eps}\,Y$ to mean $X \se_\eps Y$ and $Y \se_\eps X$.
We define a relation \defn{$\leq_\eps$} on the set of separations of~$G$ by letting
\[
    (A,B) \,\defnMath{\leq_\eps}\, (C,D)\quad :\Longleftrightarrow\quad A\se_\eps C\text{ and }B\supseteq_\eps D.
\]
As usual, we write $(A,B) \,\defnMath{<_\eps}\, (C,D)$ if $(A,B) \leq_\eps (C,D)$ but $(C,D) \nleqslant_\eps (A,B)$.
We write 
\[
    (A,B) \,\defnMath{\approx_\eps}\, (C,D)\quad :\Longleftrightarrow\quad (A,B) \leq_\eps (C,D)\text{ and }(C,D) \leq_\eps (A,B).
\]
If $(A,B)\approx_\eps (C,D)$, then we say that $(A,B)$ and $(C,D)$ are \defn{$\eps$-close}.
Note that this is equivalent to $A \approx_\eps C$ and $B \approx_\eps D$.
Note that neither~$\leq_\eps$ nor~$\approx_\eps$ nor~$<_\eps$ needs to be transitive.
However, we have:

\begin{lem}\label{approxLeQuasiTransitive}
    Let $G$ be a graph.
    Let $r\le_\eps s \le_\delta t$ be separations of~$G$ where $\eps,\delta\in\N$.
    Then $r\le_{\eps + \delta} t$.\qed
\end{lem}

\begin{lem}\label{approxLeQuasiInfimum}
    Let $G$ be a graph.
    Let $r,s,t$ be separations of~$G$ such that
    $r\le_\eps s$ and $r \le_\delta t$, where $\eps,\delta\in\N$.
    Then $r\le_{\eps + \delta} s \wedge t$.\qed
\end{lem}

We say that $(A,B)$ and $(C,D)$ \defn{have edit-distance~$\le\eps$} if $|A\triangle C|\le\eps$ and $|B\triangle D|\le\eps$.

\begin{obs}\label{closeVSapprox}
    Let~$s$ and~$t$ be two separations in a graph~$G$, and let~$\eps\in\N$.
    \begin{enumerate}
        \item If $s$ and $t$ are $\eps$-close, then $s$ and $t$ have edit-distance~$\le 2\eps$.
        \item If $s$ and $t$ have edit-distance~$\le\eps$, then $s$ and $t$ are $\eps$-close.\qed
    \end{enumerate}
\end{obs}

Two separations $(A,B)$ and $(C,D)$ of $G$ are \defn{$\eps$-nested} if, after possibly swapping the name $A$ with $B$ or $C$ with~$D$, we have $(A,B)\le_\eps (C,D)$.
If $(A,B)$ and $(C,D)$ are not $\eps$-nested, then they \defn{$\eps$-cross}.
Note that for $\eps=0$, the notions of $\eps$-nested and $\eps$-cross coincide with the standard notions of nested and cross.
A set~$S$ of separations of~$G$ is \defn{$\eps$-nested} if the separations in~$S$ are pairwise $\eps$-nested.

\subsection{3D corner diagrams}

\begin{figure}[ht]
    \centering
    \begin{tikzpicture}[scale=.66, every node/.style={circle, draw, fill=white, inner sep=1pt, transform shape}]
  \coordinate (A) at (0,0);    
  \coordinate (B) at (6,0);    
  \coordinate (C) at (6,6);    
  \coordinate (D) at (0,6);    

  \coordinate (E) at (1.8,1.8);  
  \coordinate (F) at (4.2,1.8);  
  \coordinate (G) at (4.2,4.2);  
  \coordinate (H) at (1.8,4.2);  

  \draw (A) -- (B) -- (C) -- (D) -- cycle;

  \draw (E) -- (F) -- (G) -- (H) -- cycle;

  \draw (A) -- (E);
  \draw (B) -- (F);
  \draw (C) -- (G);
  \draw (D) -- (H);

  \node at (A) [label={[fill=white, inner sep=1pt]center:ADF}] {};
  \node at (B) [label={[fill=white, inner sep=1pt]center:BDF}] {};
  \node at (C) [label={[fill=white, inner sep=1pt]center:BDE}] {};
  \node at (D) [label={[fill=white, inner sep=1pt]center:ADE}] {};
  \node at (E) [label={[fill=white, inner sep=1pt]center:ACF}] {};
  \node at (F) [label={[fill=white, inner sep=1pt]center:BCF}] {};
  \node at (G) [label={[fill=white, inner sep=1pt]center:BCE}] {};
  \node at (H) [label={[fill=white, inner sep=1pt]center:ACE}] {};
  \node[draw=none,fill=none] at (3,.9) [label={[fill=white, inner sep=1pt]center:ABCDF}] {};
  \node[draw=none,fill=none] at (6,3) [label={[fill=white, inner sep=1pt]center:BDEF}] {};
\end{tikzpicture}
    \caption{The cube on vertex set $\{A, B\} \times \{C, D\} \times \{E, F\}$. The labels of some of the faces are indicated.}
    \label{fig:3D-corner}
\end{figure}
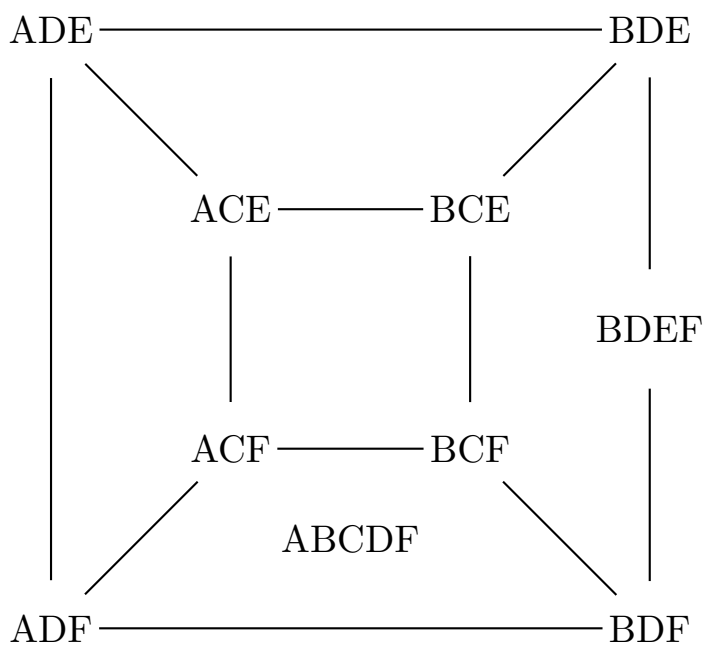

Let $G$ be a graph and let $(A,B)$, $(C,D)$ and $(E,F)$ be three $k$-separations of $G$.
Consider the cube on vertex set $\{A, B\} \times \{C, D\} \times \{E, F\}$.
Label each face $f$ of the cube by the set of letters in $\{A, B, C, D, E, F\}$ which appear in the names of the vertices in $f$ (for instance, the 1D face containing $ACE$ and $ACF$ has label $ACEF$), see \cref{fig:3D-corner}.
Let $\varphi$ be the function which maps every vertex $v \in V(G)$ to the face with label $\{Z \in \{A, B, C, D, E, F\} : v \in Z\}$.
For every 1D, 2D or 3D face $f$ of the cube, let $\lambda(f) = |\varphi^{-1}(f)|$.
We say that $\lambda$ is the \defn{visualisation} of $(A, B), (C, D)$ and $(E, F)$.
We think of $\lambda$ as a 3-dimensional generalization of corner diagrams.
We now prove two simple but key properties of $\lambda$. The first one simply states that $\lambda$ is a $k$-corner-map. The second one states that $\lambda$ correctly measures the order of the separations that can arise from the 3D corner diagram.

\begin{lemma}\label{lem:visualisation-is-corner-map}
    Let $(A, B), (C, D)$ and $(E, F)$ be three $k$-separations of a graph and let $\lambda$ be their visualisation. Then, $\lambda$ is a $k$-corner-map.
\end{lemma}

\begin{proof}
    By definition, $\lambda$ is a function from the 1D, 2D and 3D faces of the cube to the non-negative integers.
    By symmetry, it suffices to prove that the coordinate cut \[c \coloneqq E(\{ACE, ACF, ADE, ADF\}, \{BCE, BCF, BDE, BDF\})\] has $\lambda$-value $k$.
    It follows immediately from the definition of $\lambda$ and of $c$ that $\lambda(c) = |A \cap B| = k$.
\end{proof}

\begin{lemma}\label{lem:lambda-correct-order}
    Let $(A, B), (C, D)$ and $(E, F)$ be three separations of a graph and let $\lambda$ be their visualisation.
    Let $\mathcal{S} \subseteq \{A, B\} \times \{C, D\} \times \{E, F\}$, and let $A' \subseteq V(Q)$ be the set of vertices of $Q$ whose label is an element of $\mathcal{S}$ and let $c = E(A', V(Q) \sm A')$.
    Then, the separation \[s \coloneqq {\left(\bigcup_{(X, Y, Z) \in \mathcal{S}} X \cap Y \cap Z, \bigcup_{(X, Y, Z) \notin \mathcal{S}} X \cap Y \cap Z\right)}\] has order $\lambda(c)$.
\end{lemma}

\begin{proof}
    Let $v \in V(G)$. We have \begin{align*}
        \text{$v$ contributes to $\lambda(c)$ } \iff &\text{ the face $\varphi(v)$ of $Q$ intersects both $A'$ and $V(Q) \sm A'$} \\
        \iff &\text{ there exist $(X_1, Y_1, Z_1) \in \mathcal{S}$ such that $v \in X_1 \cap Y_1 \cap Z_1$, }\\ &\text{ and $(X_2, Y_2, Z_2) \notin \mathcal{S}$ such that $v \in X_2 \cap Y_2 \cap Z_2$} \\
        \iff &\text{ $v \in \sep(s)$.}
    \end{align*}
    Furthermore, if $v$ contributes to $\lambda(c)$ then $v$ contributes by exactly one. This concludes the proof.
\end{proof}

We now introduce the terminology that we will use to talk about the various zones in the Venn diagram of three separations.
Let $r \coloneqq (A, B)$, $s \coloneqq (C, D)$ and $t \coloneqq (E, F)$ be three separations.
We define the \defn{$ACE$-corner} to be the set $\varphi^{-1}(ACE) = (A \sm B) \cap (C \sm D) \cap (E \sm F)$, and define the 7 other corners analogously. The corners correspond to vertices of $Q$. 
We define the \defn{$rCE$(-edgy)-link} to be the set $\varphi^{-1}(rCE) := (A \cap B) \cap (C \sm D) \cap (E \sm F)$, and define the 11 other edgy links analogously. The edgy links correspond to edges of $Q$. For instance, the $rCE$-link corresponds to the edge $ABCE$.
We define the \defn{$rsE$(-facial)-link} to be the set $\varphi^{-1}(rsE) := (A \cap B) \cap (C \cap D) \cap (E \sm F)$, and define the other 5 facial links analogously. The facial links correspond to facets of $Q$. For instance, the $rsE$-link corresponds to the face labelled $ABCDE$.
We define the \defn{$rst$-centre} to be the set $\varphi^{-1}(rst) := (A \cap B) \cap (C \cap D) \cap (E \cap F)$. The $rst$-centre corresponds to the entire cube $Q$.

A \defn{3D corner separation} is any separation of the form $(A \cap C \cap E, B \cup D \cup F)$, up to exchanging $A$ and $B$, and/or $C$ and $D$, and/or $E$ and $F$.
For simplicity, we call the separation $(A \cap C \cap E, B \cup D \cup F)$ the \defn{3D corner separation $ACE$}, and similarly for the other 3D corner separations.

\subsection{Proof of the 3D Corner Lemma}

We now prove a result about how three well-behaved separations can cross in a $(q,k-1)$-unbreakable graph. To state this result, we first need some definitions.

Let $G$ be a graph.
Assume that two separations $s = (A,B)$ and $t = (C,D)$ of $G$ cross.
Let $P$ denote the $sC$-link, let $X$ denote the $st$-centre, and let $Q$ denote the $sD$-link.
We say that $(C,D)$ \defn{induces} $(X,\{P,Q\})$ on~$A\cap B$.
By a slight abuse of notation, we call~$(X,\{P,Q\})$ a \defn{tripartition} of~$A \cap B$.

Let $S$ be a set of separations of~$G$, and let $\gamma\in\N$.
We say that $S$ is \defn{uniformly $\gamma$-crossing} if for each $(A,B)\in S$ there is a tripartition $(X,\{P,Q\})$ of $A\cap B$ such that, whenever a $(C,D)\in S$ $\gamma$-crosses $(A,B)$, then $(C,D)$ induces the tripartition $(X,\{P,Q\})$ on $A\cap B$.
We abbreviate uniformly 0-crossing to \defn{uniformly crossing}.

If two $k$-separations cross, they cross \defn{equitably} if all four corner separators have size exactly $k$.
A set $S$ of $k$-separations is \defn{equitable} if whenever two elements of $S$ cross, they cross equitably.
Given $\varepsilon \in \mathbb{N}$, we say that two separations \defn{$\varepsilon$-criss} if their four corners all have size more than $\varepsilon$.
A set $S$ of separations is \defn{$\varepsilon$-criss-compatible} if whenever two separations in $S$ cross, they $\varepsilon$-criss.

The central statement of this section is the following.

\begin{restatable}[3D Corner Lemma]{keylemma}{cornerlemma}\label{3dcorner}
Let $G$ be a $(q,k-1)$-unbreakable graph and let $S$ be an equitable $(2q+k)$-criss-compatible triple of $k$-separations of $G$.
Then,~$S$ is uniformly crossing.
\end{restatable}

We can prove that the same result holds under the weaker assumption that no two separations in $S$ cross with two opposite corners of size at most $2q+k$, but this requires some extra work and we do not need this stronger result.
We note that in the statement of \cref{3dcorner}, we can drop the assumption of~$S$ being equitable (by \cref{equitable_calc}).

We start with some preliminary results.

\begin{lem}\label{to-cover}
Let $G$ be a $(q,k-1)$-unbreakable graph and let $S$ be an equitable $(2q+k)$-criss-compatible triple of $q$-broad $k$-separations of $G$.
Let $\lambda$ be the visualisation map of $S$ and let $R$ be the set of vertices $XYZ$ of $Q$ such that $|X \cap Y \cap Z| \geq q+k$.
Suppose that the separations in $S$ pairwise cross.
Then, $R$ is a vertex cover of $Q$.
\end{lem}

\begin{proof}
    By contradiction, suppose that some edge of $Q$ is not covered by $R$.
    Without loss of generality, suppose that this edge is the edge $(ACE, ACF)$, see \cref{fig:3D-corner}.
    Therefore, $|A \cap C| \leq |A \cap C \cap E| + |A \cap C \cap F| \leq 2(q+k)$.
    Since $S$ is equitable and $(A, B)$ and $(C, D)$ cross, the 2D corner separation $(A \cap C, B \cup D)$ has order exactly $k$, so the $AC$-corner has size at most $2q+k$, a contradiction.
\end{proof}

\begin{lem}\label{lem:cut-R-either-side-lambda}
Let $G$ be a $(q,k-1)$-unbreakable graph and let $S$ be a triple of $k$-separations of $G$.
Let $\lambda$ be the visualisation map of $S$ and let $R$ be the set of vertices $XYZ$ of $Q$ such that $|X \cap Y \cap Z| \geq q+k$.
Then, every cut of $Q$ with one vertex of $R$ on each side has $\lambda$-value at least $k$.
\end{lem}

\begin{proof}
    Let $c$ be a cut of $Q$ with one vertex of $R$ on each side.
    Then, $c$ lifts to a separation $s$ of $G$ as in \cref{lem:lambda-correct-order}.
    By definition of $R$ and since there is a vertex of $R$ on each side of $c$, both sides of $s$ have size at least $q+k$.
    Since $G$ is $(q,k-1)$-unbreakable, $s$ has order at least $k$.
    By \cref{lem:lambda-correct-order}, $\lambda(c)$ is the order of this separation, so $\lambda(c) \geq k$.
\end{proof}

\begin{lem}\label{lem:flower-unif-crossing}
    Let $G$ be a graph and let $S$ be a triple of $k$-separations of $G$.
    Let $\lambda$ be the visualisation map of $S$.
    If $\lambda$ is a flower then $S$ is uniformly crossing.
\end{lem}

\begin{proof}
    Write $S = \{(A, B), (C, D), (E, F)\}$ and let $r = (A, B)$, $s = (C, D)$ and $t = (E, F)$.
    Since $\lambda$ is a flower, without loss of generality, we can assume that all the 1D faces incident to $ACE$ or $BDF$ have $\lambda$-value zero, that all other 1D faces have the same $\lambda$-value $\ell \leq k/2$, that all the 2D faces have $\lambda$-value zero and that the 3D face has $\lambda$-value $k-2\ell$, see \cref{fig:proof-flower-like}.
    By definition of $\lambda$, this means that the $rsE$-, $rsF$-, $rCt$-, $rDt$-, $Ast$- and $Bst$-links are empty, as well as the $AsE$-, $ACt$-, $rCE$-, $rDF$-, $BsF$- and $BDt$-links.

    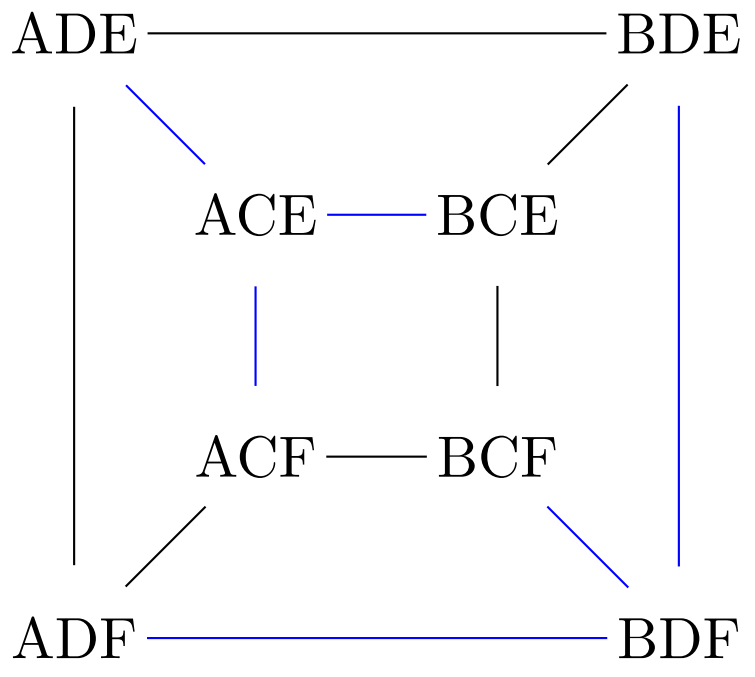
\begin{figure}[ht]
    \centering
\begin{tikzpicture}[scale=4, every node/.style={circle, draw, fill=white, inner sep=1pt}]
  \coordinate (A) at (0,0);    
  \coordinate (B) at (1,0);    
  \coordinate (C) at (1,1);    
  \coordinate (D) at (0,1);    

  \coordinate (E) at (0.3,0.3);  
  \coordinate (F) at (0.7,0.3);  
  \coordinate (G) at (0.7,0.7);  
  \coordinate (H) at (0.3,0.7);  

  \draw[blue] (A) -- (B);
  \draw[blue] (B) -- (C);
  \draw (C) -- (D);
  \draw (D) -- (A);

  \draw (E) -- (F);
  \draw (F) -- (G);
  \draw[blue] (G) -- (H);
  \draw[blue] (H) -- (E);

  \draw (A) -- (E);
  \draw[blue] (B) -- (F);
  \draw (C) -- (G);
  \draw[blue] (D) -- (H);

  \node at (A) [label={[fill=white, inner sep=1pt]center:ADF}] {};
  \node at (B) [label={[fill=white, inner sep=1pt]center:BDF}] {};
  \node at (C) [label={[fill=white, inner sep=1pt]center:BDE}] {};
  \node at (D) [label={[fill=white, inner sep=1pt]center:ADE}] {};
  \node at (E) [label={[fill=white, inner sep=1pt]center:ACF}] {};
  \node at (F) [label={[fill=white, inner sep=1pt]center:BCF}] {};
  \node at (G) [label={[fill=white, inner sep=1pt]center:BCE}] {};
  \node at (H) [label={[fill=white, inner sep=1pt]center:ACE}] {};
\end{tikzpicture}
    \caption{Illustration of the proof of \cref{lem:flower-unif-crossing}. 
    Every 2D face has $\lambda$-value $0$, as well as all blue 1D faces.
    All pairwise centres are then equal to the $rst$-centre.
    Furthermore, all pairwise links are equal to some non-blue edge.
    }
    \label{fig:proof-flower-like}
\end{figure}

    Whenever two distinct separations in $S$ cross, their centre is the union of the $rst$-centre and two facial links. 
    However, we saw that all facial links are empty.
    Therefore, all three pairwise centres are equal to the $rst$-centre.

    Similarly, for the crossing separations $r = (A, B)$ and $s =(C, D)$, the $rC$-link is the union of the $rCE$-link, the $rCF$-link and the $rCt$-link, which is simply the $rCF$-link, and the $rD$-link is the union of the $rDE$-link, the $rDF$-link and the $rDt$-link, which is simply the $rDE$-link.
    For the crossing separations $r = (A, B)$ and $t = (E, F)$, the $rF$-link is the union of the $rCF$-link, the $rDF$-link and the $rsF$-link, which is simply the $rCF$-link, and the $rE$-link is the union of the $rCE$-link, the $rDE$-link and the $rsE$-link, which is simply the $rDE$-link.
    This proves that the separations $s$ and $t$ uniformly cross the separation $r$. Similar arguments conclude the proof that $S$ is uniformly crossing.
\end{proof}

\begin{lem}\label{crossing-with-2-nested}
    Let $G$ be a graph and let $S = \{r, s, t\}$ be an equitable set of $k$-separations of $G$. 
    Assume that $r, s$ cross, $r, t$ cross and $s \leq t$.
    Then, $S$ is uniformly crossing.
\end{lem}

\begin{proof}
    Since the separations $s$ and $t$ are nested, it suffices to prove that the crossing separations $r, s$, and $r, t$ induce the same tri-partition of the separator of $r$ into one centre and two links.

    Let $r \eqqcolon (A, B)$, $s \eqqcolon (C, D)$ and $t \eqqcolon (E, F)$.
    Let $\lambda$ be the visualisation of $S$. By \cref{lem:visualisation-is-corner-map}, $\lambda$ is a $k$-corner-map.
    The vertices $ACF$ and $BCF$ are adjacent in $Q$ (see \cref{fig:3D-corner}).
    Since $C \subseteq E$ and $F \subseteq D$, we have that the $rCt$-link, the $rsF$-link and the $rCF$-link are empty.
    Thus, by definition of $\lambda$, the 1D and 2D faces which contain $ACF$ and $BCF$ have $\lambda$-value $0$.

    Since $S$ is equitable, and $(A, B)$ and $(C, D)$ cross, and $(A, B)$ and $(E, F)$ cross, we have that the $AC$-, $AF$-, $BC$- and $BF$-corner separations have order exactly $k$.
    Thus, by \cref{lem:lambda-correct-order}, all pair cuts which contain exactly one of $ACF, BCF$ have $\lambda$-value $k$.
    Therefore, by \cref{comb-crossing-with-2-nested-empty-links}, all 1D and 2D faces which contain $ADE$ and $BDE$ have $\lambda$-value $0$.
    By definition of $\lambda$, this means that the $rDt$-link, the $rsE$-link and the $rDE$-link are also empty.
    See \cref{fig:statement-comb-crossing-with-2-nested-empty-links} for an illustration of the setting.

    Since the $rCt$-link, the $rDt$-link, the $rsE$-link and the $rsF$-link are empty, the two pairs of crossing separations $(A, B)$ and $(C, D)$, and $(A, B)$ and $(E, F)$ have the same centre, namely the $rst$-centre.

    Furthermore, since the $rCt$-link, the $rsE$-link, the $rCF$-link and the $rDE$-link are empty, it follows that $rCE$ is both the $rC$-link of the crossing separations $(A, B)$ and $(C, D)$, and the $rE$-link of the crossing separations $(A, B)$ and $(E, F)$.
    Similarly, $rDF$ is both the $rD$-link of the crossing separations $(A, B)$ and $(C, D)$, and the $rF$-link of the crossing separations $(A, B)$ and $(E, F)$.
    This concludes the proof.
\end{proof}

Combining these results, we can now prove \cref{3dcorner}, which we restate for convenience.

\cornerlemma*

\begin{proof}[Proof of \cref{3dcorner}]
Write $S = \{(A, B), (C, D), (E, F)\}$.
Let $\lambda$ be the visualisation map of $S$.
By \cref{lem:visualisation-is-corner-map}, $\lambda$ is a $k$-corner-map.
If no separation in $S$ crosses the other two, then we are done. 
Hence assume that one of them, say $(A,B)$, crosses the other two.
Since $S$ is $(2q+k)$-criss-compatible, it follows that all separations in~$S$ have breadth~$\geq 4q+2k \geq q$.

\casen{Case 1: $(C,D)$ and $(E,F)$ cross.} 
Then the three separations $(A,B)$, $(C,D)$ and $(E,F)$ are pairwise crossing.
Let $R$ be the set of vertices $XYZ$ of $Q$ such that $|X \cap Y \cap Z| \geq q+k$.
Then, $R$ is a vertex cover of $Q$ by \cref{to-cover} and every cut of $Q$ with one vertex of $R$ on each side has $\lambda$-value at least $k$ by \cref{lem:cut-R-either-side-lambda}.
Thus, $\lambda$ is a flower by \cref{combinatorial_3D-corner-lemma}.
Therefore, $S$ is uniformly crossing by \cref{lem:flower-unif-crossing}.

\casen{Case 2: $(C,D)$ and $(E,F)$ are nested.} Without loss of generality, suppose that $C \subseteq E$ and $F \subseteq D$. Then, by \cref{crossing-with-2-nested}, $S$ is uniformly crossing.
\end{proof}

\section{Broad separations are uniformly crossing}
\label{sec:BroadSeparationsCrossUniformly}

We now use the 3D Corner Lemma to show that sufficiently strong crossings impose consistent centres and links. This consistency will allow us to organise crossing separations into asteriskal sets of separations in the next section.

The main goal of this section is to prove the following statement.

\begin{keylemma}\label{HugeSeparationsUniformlyCrossing}
    Let~$G$ be a $(q,k-1)$-unbreakable graph, and define $\gamma_{\ref{HugeSeparationsUniformlyCrossing}}(q,k):=3q+3k$.
    Then the set of all $q$-broad $k$-separations of~$G$ is uniformly $\gamma_{\ref{HugeSeparationsUniformlyCrossing}}$-crossing.
\end{keylemma} 

We will also prove some additional statements which will be useful in later sections.

\begin{convention}\label{keepBoundsSimple}
    Often we will define bounds of the form $\eps_{X}(q,k)$ where $q,k\in\N$ determine the degree of $(q,k-1)$-unbreakability and $X$ is the number of a lemma or whatever that introduces the bound.
    When we use these bounds for some $(q,k-1)$-unbreakable graph, we sometimes drop `$(q,k)$' and just write $\eps_X$ to reduce notational clutter, provided that $q$ and $k$ are clear from context.
\end{convention}

\subsection{Uniformly crossing sets of separations}

Recall that a separation $(A,B)$ of a graph~$G$ is \defn{$\beta$-broad} if $|A \setminus B| \geq \beta$ and $|B \setminus A| \geq \beta$. 
Moreover, recall that two separations $s$ and $t$ \defn{$\eps$-criss} if all four corners have size~$>\eps$.

\begin{lem}\label{equitable_calc}
Let~$G$ be a $(q,k-1)$-unbreakable graph, and let~$s,s'$ be two $q$-crissing $k$-separations of~$G$.
Then~$s$ and~$s'$ cross equitably.
\end{lem}

\begin{proof}
Assume~$s$ and~$s'$ do not cross equitably.
Then, one of their corner separations has order less than~$k$.
By $(q,k-1)$-unbreakability of~$G$, the corresponding corner contains at most $q$ vertices. So,~$s$ and~$s'$ do not $q$-criss.
\end{proof}

\begin{observation}\label{crissCross}
    Let $s=(A,B)$ and $t=(C,D)$ be two $(\le k)$-separations in a graph~$G$, and let~$\eps\in\N$.
    \begin{enumerate}
        \item\label{crissCross:smallCornerImpliesNested} If the $AD$-corner has size~$\leq \eps-\ell$ and both the $At$-link and $sD$-link have size~$\leq\ell$, then $s \leq_\eps t$.
        \item\label{crissCross:crossImpliesCriss} If $s$ and $t$ $\eps$-cross, then they $(\eps-k)$-criss.
        \item\label{crissCross:crossImpliesCriss2} If $s$ and $t$ $\eps$-cross and all links have size~$\le\ell$, then $s$ and $t$ $(\eps-\ell)$-criss.
        \item\label{crissCross:crissImpliesCross} If $s$ and $t$ $\eps$-criss with~$\eps>0$, then they also $\eps$-cross.
        \item\label{crissCross:crossLargeCorner} If $s$ and $t$ $\eps$-cross then $|A \cap C|, |A \cap D|, |B \cap C|, |B \cap D| > \eps$.
    \end{enumerate}
\end{observation}
\begin{proof}
    \ref{crissCross:smallCornerImpliesNested}.
    We have $|A \sm C| \leq |A \sm (B \cup C)| + |(A \cap B) \sm C| \leq (\eps-\ell)+\ell=\eps$.
    Similarly, $|D \sm B| \leq \eps$.
    It follows that $A \subseteq_\eps C$ and $D \subseteq_\eps B$.
    Hence $(A,B) \leq_\eps (C,D)$.
    
    \ref{crissCross:crossImpliesCriss2}.
    If~$s=:(A,B)$ and~$t=:(C,D)$ do not $(\eps-\ell)$-criss, then there is a corner of size~$\leq \eps-\ell$, say the $AD$-corner.
    By \cref{crissCross:smallCornerImpliesNested}, $s \leq_\eps t$, and they are $\eps$-nested.

    \ref{crissCross:crossImpliesCriss} follows from~\ref{crissCross:crossImpliesCriss2} with $\ell:=k$.

    \ref{crissCross:crissImpliesCross}.
    We show the contrapositive.
    Let~$s=:(A,B)$ and~$t=:(C,D)$ be $\eps$-nested, say $(A,B) \leq_\eps (C,D)$.
    Then $|(A \cap D) \sm (B \cup C)| \leq |A \sm C| \leq \eps$.
    That is,~$(A,B)$ and~$(C,D)$ do not $\eps$-criss.

    \ref{crissCross:crossLargeCorner}.
    If~$|A \cap C| \leq \eps$, say, then $|A \sm D|\leq \eps$ and $|C \sm B| \leq \eps$.
    Then $(A,B) \leq_\eps (D,C)$, in particular they do not $\eps$-cross.
\end{proof}

\begin{lem}
\label{lem:vee-wedge-is-close}
    Let~$s,t$ be two separations of a graph~$G$ with $s \leq_\gamma t$.
    Then $s \wedge t \approx_\gamma s$ and $t \approx_\gamma s \vee t$.
\end{lem}
\begin{proof}
    Write $s=:(A,B)$ and $t=:(C,D)$.
    Then $s \leq_\gamma t$ implies $|A \sm C| \leq \gamma$ and $|D \sm B| \leq \gamma$.
    It follows that $|A \sm (A \cap C)| \leq \gamma$ and $|(B \cup D) \sm B| \leq \gamma$.
    Therefore, $A \subseteq_\gamma A \cap C \se A$ and $B \supseteq_\gamma B \cup D \supseteq B$.
    Hence, $s \wedge t=(A \cap C,B \cup D) \approx_\gamma (A,B)=s$.

    Similarly, $s \leq_\gamma t$ implies $t^* \leq_\gamma s^*$.
    By the first part, we have $(s \vee t)^* = s^* \wedge t^* \approx_\gamma t^*$.
    This implies $t \approx_\gamma s \vee t$.
\end{proof}

\begin{lem}\label{LargeCrossedImpliesBroad}
    If two $(\le k)$-separations of a graph~$G$ $\eps$-cross, then both are $2(\eps-k)$-broad.
\end{lem}
\begin{proof}
    Let $(A,B)$ and $(C,D)$ be $(\le k)$-separations of~$G$ that $\eps$-cross.
    Then every corner has size~$\geq \eps-k$ by \cref{crissCross}, so $|A \sm B|, |B \sm A| \geq 2 (\eps-k)$.
\end{proof}

\begin{lem}\label{lem:meet-join-k-separation}
    Let $\gamma,\ell\in\N$.
    Let~$G$ be a graph and let~$s,t$ be two $(\ell + \gamma)$-broad $k$-separations of~$G$ with $s \leq_\gamma t$.
    Then both $s \wedge t$ and $s \vee t$ are $\ell$-broad.
    Moreover, if~$G$ is $(q,k-1)$-unbreakable and $\ell>q$ then both $s \wedge t$ and $s \vee t$ are $\ell$-broad $k$-separations of~$G$.
\end{lem}
\begin{proof}
    Let~$s=:(A,B)$, $t=:(C,D)$, $s \wedge t=:(X_1,Y_1)=:r_1$ and $s \vee t=:(X_2,Y_2)=:r_2$.
    Since $s \leq_\gamma t$, we have $|D \sm B| \leq \gamma$.
    So,
    \begin{align*}
        |Y_1 \sm X_1| &= |(B \cup D) \sm (A \cap C)| \geq |B \sm A| \geq \ell+\gamma \geq \ell~,~\text{and}\\
        |X_1 \sm Y_1| &= |(A \cap C) \sm (B \cup D)| = |A\sm B| - |(A \cap D) \sm B| \geq |A \sm B| - |D \sm B| \geq (\ell+\gamma)-\gamma=\ell~.
    \end{align*}
    Similarly, $|X_2 \sm Y_2| \geq \ell$ and $|Y_2 \sm X_2| \geq \ell$.
    So, both~$r_1$ and~$r_2$ are $\ell$-broad.

    Moreover, assume that~$G$ is $(q,k-1)$-unbreakable.
    Then, from $\ell > q$ it follows that both~$r_1$ and~$r_2$ have order~$\geq k$.
    From submodularity, it follows that the order of both~$r_1$ and~$r_2$ is exactly~$k$.
\end{proof}

\begin{lem}\label{CloseVsCross}
    Let $G$ be a graph.
    Let $s,s',t$ be separations of~$G$.
    Let $\gamma,\eps\in\N$ with $\gamma\ge\eps$.
    Assume that $s$ and $s'$ $\gamma$-cross while $s$ and $t$ are $\eps$-close.
    Then $t$ and $s'$ $(\gamma-\eps)$-cross.
\end{lem}
\begin{proof}
    Since $s$ and $t$ are $\eps$-close, we may assume without loss of generality that $s\le_\eps t\le_\eps s$.
    Assume for a contradiction that $t$ and $s'$ fail to $(\gamma-\eps)$-cross.
    Then $t\le_{\gamma-\eps} s'$, say.
    With $s\le_\eps t$ and \cref{approxLeQuasiTransitive}, we find that $s\le_\gamma s'$, contradicting that $s$ and~$s'$ $\gamma$-cross.
\end{proof}

\begin{lem}\label{LowBreadthClose}
    Let $s$ and $t$ be separations of a graph~$G$.
    Assume that $s$ has order~$\le k$, breadth~$\beta$, and that $s$ is $\eps$-close to~$t$.
    Then $t$ has breadth~$\beta'$ with~$\beta-\eps\le \beta'\le \beta+\eps+k$.
\end{lem}
\begin{proof}
    Let $s=:(A,B)$ and $t=:(C,D)$.
    Without loss of generality, $|A\sm B|=\beta$ and $|B \sm A| \geq \beta$.
    First, we show $\beta'\le \beta+\eps+k$.
    Since~$s$ is $\eps$-close to~$t$, we have $s \approx_\eps t$.
    Hence $|C\sm A|\le\eps$.
    From
    \[
        C\sm D\se C\se (C\sm A)\cup A\se (C\sm A)\cup (A\cap B)\cup (A\sm B)
    \]
    we obtain $|C\sm D|\le \eps+k+\beta$.

    Now, we show $\beta'\ge \beta-\eps$.
    From $|A \sm B| \geq \beta$ and $|D \sm B| \leq \eps$, we obtain
    \[
        |C \sm D| \geq |(A \sm B) \cap (C \sm D)| = |A \sm B| - |D \sm B| \geq \beta - \eps~.
    \]
    Similarly, from $|B \sm A| \geq \beta$ and $|C \sm A| \leq \eps$ we obtain $|D \sm C| \geq \beta - \eps$.
    Hence,~$t$ is $(\beta-\eps)$-broad.
\end{proof}

A set~$S$ of separations of a graph~$G$ is \defn{$\eps$-triple uniformly crossing} if for every triple of elements of $S$, either the triple is uniformly crossing, or two separations do not $\eps$-cross.
The next result essentially says that to prove that a set of separations is uniformly $\gamma$-crossing, it suffices to prove that it is $\varepsilon$-triple uniformly crossing, for some $\varepsilon < \gamma$.

\begin{lem}\label{lem:tripleUniform-implies-uniform-2}
    Let~$G$ be a $(q,k-1)$-unbreakable graph, let $\eps \in \Nbb$ and 
    \[\gamma \geq \gamma_{\ref{lem:tripleUniform-implies-uniform-2}}(\eps, q, k) \coloneqq \max\{\,\eps+k+q\, , \, 2k \,\}~.\]
    Then, any $\eps$-triple uniformly crossing set of $k$-separations of $G$ is uniformly $\gamma$-crossing.
\end{lem}
\begin{proof}
    Let $S$ be an $\eps$-triple uniformly crossing set of $k$-separations of $G$, and consider an arbitrary separation $r \in S$. 
    Let $s, t\in S$ be two separations which $\gamma$-cross~$r$.
    Our aim is to show that~$r$ is crossed by both of them with the same tri-partition of the separator $\sep(r)$ into centre and links. 
    Since $S$ is $\eps$-triple uniform, either the triple $\{r, s, t\}$ is uniformly crossing, or two of its separations do not $\eps$-cross.
    Since we are done in the first case, assume that two of its separations do not $\eps$-cross.
    Since $\gamma \geq \eps$ and both $s$ and $t$ $\gamma$-cross $r$, we can assume that $s$ and~$t$ are $\eps$-nested, say $s \leq_\eps t$.

    By \cref{lem:vee-wedge-is-close}, we have $s \wedge t \approx_\eps s$.
    Since $r$ $\gamma$-crosses $s$ and~$t$, all three separations are $2(\gamma - k)$-broad by \cref{LargeCrossedImpliesBroad}, with $2(\gamma-k) \geq \gamma$ (since $\gamma \geq 2k$). 
    By \cref{lem:meet-join-k-separation} and since $\gamma-\eps > q$, the separation~$s \wedge t$ is a $(\gamma-\eps)$-broad separation of order exactly~$k$.

    By \cref{CloseVsCross}, $r$ and~$s \wedge t$ $(\gamma-\eps)$-cross.
    By \cref{crissCross}, $r$ and~$s \wedge t$ $(\gamma-\eps-k)$-criss.
    By \cref{equitable_calc} and since $\gamma-\eps-k\geq q$, the separations~$r$ and $s \wedge t$ cross equitably.

    Since $s$ and $t$ both $(\gamma-k)$-criss $r$ (by \cref{crissCross}) with $\gamma-k \geq q$, by \cref{equitable_calc}, the sets $\{r, s, s \wedge t\}$ and $\{r, t, s \wedge t\}$ are equitable.
    Since $s$ and $s \wedge t$ are nested, the triple 
    $\{r, s, s \wedge t\}$ is uniformly crossing by \cref{crossing-with-2-nested}. 
    Similarly, $\{r, t, s \wedge t\}$ is uniformly crossing.
    
    Consider the tri-partition of the separator $\sep(r)$ obtained when crossing $r$ with $s \wedge t$. 
    Since $\{r, s, s \wedge t\}$ is uniformly crossing, this is the same as when crossing $r$ with $s$, and since $\{r, t, s \wedge t\}$ is uniformly crossing, this is the same as when crossing $r$ with $t$. This completes the proof. 
\end{proof}

\begin{lem}\label{lem:largeSeparation-is-tripleUniform}
Let $G$ be a $(q,k-1)$-unbreakable graph, and define $\varepsilon = \varepsilon(q, k) \coloneqq 2q+2k$.
Then, the set of $q$-broad $k$-separations of $G$ is $\varepsilon$-triple uniformly crossing. 
\end{lem}

\begin{proof}
    Let $S$ be a triple of $q$-broad $k$-separations of $G$.
    If two separations in $S$ are $\eps$-nested then we are done.
    Otherwise, every pair in~$S$ $\eps$-crosses, hence $(\eps-k)$-crisses by \cref{crissCross} (with $\eps-k=2q+k$).
    Thus, $S$ is equitable by \cref{equitable_calc}.
    Then, $S$ is uniformly crossing by \cref{3dcorner}.
\end{proof}

\begin{proof}[Proof of \cref{HugeSeparationsUniformlyCrossing}]
    This follows immediately from \cref{lem:tripleUniform-implies-uniform-2} and \cref{lem:largeSeparation-is-tripleUniform}.
\end{proof}

\subsection{Empty links in 3D corner diagrams} \label{subsec:more-results}

We now show some additional results that will be used in \cref{sec:petal-separations-totally-nested}.

\begin{corollary} \label{cor:3D-crossing1}
    Let $G$ be a $(q,k-1)$-unbreakable graph.
    Let $r=(A, B)$, $s=(C, D)$ and $t=(E, F)$ be three $k$-separations of~$G$ such that the 3D corners $ADF$, $BCF$ and $BDE$ have size more than~$q$.
    Then, the $Ast$-link, the $rCt$-link, the $rsE$-link, the $ACt$-link, the $AsE$-link and the $rCE$-link are empty.
\end{corollary}

\begin{proof}
    Since $G$ is $(q,k-1)$-unbreakable, the 3D corner separations $ADF$, $BCF$ and $BDE$ have order at least $k$.
    Similarly, the 2D corner separations $AD$, $AF$, $DF$, $BC$, $BF$, $CF$, $BD$, $BE$ and $DE$ have order at least $k$.
    By submodularity, the 2D corner separations $BC$, $BE$, $CE$, $AD$, $AE$, $DE$, $AC$, $AF$ and $CF$ have order at most $k$.
    Let $\lambda$ be the visualisation map of $r, s$ and $t$.
    By \cref{lem:visualisation-is-corner-map}, $\lambda$ is a $k$-corner-map.
    By \cref{lem:lambda-correct-order}, the $\lambda$-values correspond to the order of the separations of $G$.
    See \cref{fig:3D-corner} for an illustration of the setting.
    Applying \cref{lem:comb-3Dcrossing-lem-1} with $u = ADF$, $v = BCF$ and $w = ACF$ (which we can do since the $AF$ and $CF$ corner separations have order at most $k$), we get that the $Ast$-link, the $rCt$-link and the $ACt$-link are empty.
    Applying \cref{lem:comb-3Dcrossing-lem-1} with $u = ADF$, $v = BDE$ and $w = ADE$ (which we can do since the $AD$ and $DE$ corner separations have order at most $k$), we get that the $Ast$-link, the $rsE$-link and the $AsE$-link are empty.
    Finally, applying \cref{lem:comb-3Dcrossing-lem-1} with $u = BCF$, $v = BDE$ and $w = BCE$ (which we can do since the $BC$ and $BE$ corner separations have order at most $k$), we get that the $rCt$-link, the $rsE$-link and the $rCE$-link are empty.
\end{proof}

\begin{corollary}\label{cor:3D-crossing3}
    Let $r=(A, B)$, $s=(C, D)$ and $t=(E, F)$ be three $k$-separations of a graph, such that the 3D corner separations $ADF$, $ACF$, $ACE$ and $BCE$ have order at least $k$, and the 2D corner separations $AC$, $AF$ and $CE$ have order at most $k$.
    Then, the $Ast$-link, the $rCt$-link, the $rsE$-link, the $rsF$-link, the $AsE$-link and the $rCF$-link are empty.
\end{corollary}

\begin{proof}
    Let $\lambda$ be the visualisation map of $r, s$ and $t$.
    By \cref{lem:visualisation-is-corner-map}, $\lambda$ is a $k$-corner-map.
    By \cref{lem:lambda-correct-order}, the $\lambda$-values correspond to the order of the separations of $G$.
    See \cref{fig:3D-corner} for an illustration of the setting.
    Applying \cref{lem:comb-3Dcrossing-lem-1} with $u = ACF$, $v = BCE$ and $w = ACE$ (which we can do since the $AC$ and $CE$ corner separations have order at most $k$), we get that the $Ast$-link, the $rsE$-link and the $AsE$-link are empty.
    Similarly, applying \cref{lem:comb-3Dcrossing-lem-1} with $u = ADF$, $v = ACE$ and $w = ACF$ (which we can do since the $AF$ and $AC$ corner separations have order at most $k$), we get that the $rCt$-link, the $rsF$-link and the $rCF$-link are empty.
\end{proof}

\section{From asterisks to flowers via coverings}
\label{sec:covering}

The purpose of this section is to prove \cref{keylem:flower-approximation}.
For this we need the following notation.

\begin{figure}[ht]
    \centering
    \includegraphics[height=9\baselineskip]{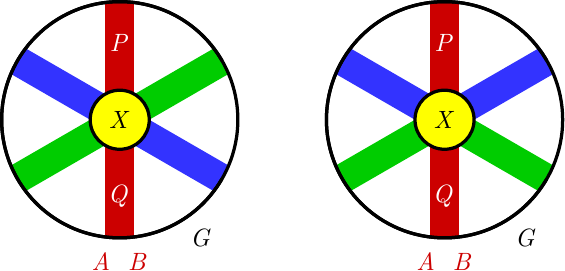}
    \caption{Each white chunk contains $\ge\gamma$ vertices. The red, green and blue separations all include $X$ in their separators. Together, they form a $\gamma$-asteriskal set of separations.}
    \label{fig:Asterisk}
\end{figure}

We say that a set~$S$ of separations is \defn{$\gamma$-crossing-connected} if $S\neq\emptyset$ and the reflexive and transitive closure of the union of the relations `$\gamma$-crossing' on~$S$ and $\{\,(s,s^*):s\in S\,\}$ has exactly one partition class.
Equivalently,~$S$ is $\gamma$-crossing-connected if $S\neq\emptyset$ and for any two separations $s,t \in S$ with $s \neq t^*$ there is a sequence $s=s_0,s_1,\ldots,s_\ell=t$ of separations in~$S$ of some length~$\ell \in\mathbb{N}$ such that~$s_i$ $\gamma$-crosses~$s_{i-1}$ for every $i \in [\ell]$.
We call $S$ \defn{$\gamma$-asteriskal} if $S$ is 
\begin{itemize}
    \item uniformly $\gamma$-crossing,
    \item $\gamma$-crossing-connected, and
    \item symmetric;
\end{itemize}
see \cref{fig:Asterisk}.
A~$\gamma$-asteriskal set~$S$ is \defn{trivial} if it consists of a single separation and its inverse.
Assume now that $S$ is uniformly $\gamma$-crossing and symmetric, but not necessarily $\gamma$-crossing-connected.
Then~$S$ admits a unique partition into maximal $\gamma$-asteriskal subsets of~$S$.
We refer to these maximal subsets as the \defn{$\gamma$-asterisks of~$S$}.

Let $(H,\cV)$ be a graph-decomposition.
For an edge $e=tt'\in E(H)$ we write $\defnMath{V_e}:=V_t\cap V_{t'}$.

\begin{definition}[Compatible]
\label{def:F-compatible}
    Let $\mathcal{S} = (S, \cV)$ be a graphic-$k$-anemone and let $Z \coloneqq V_c$, where $c$ is the centre of $S$.
    A separation $s = (A, B)$ is \defn{$\gamma$-compatible with~$\cS$} if, whenever $s$ $\gamma$-crosses a petal-separation $t$ of $\cS$ then $s$ and $t$ cross with centre $Z$ and four empty links.
    
    Let~$\cD=(O,(V_p)_{p \in O})$ be a graphic-$k$-daisy of a graph~$G$ and let $Z:= \bigcap_{p \in O} V_p$.
    A separation $s=(A,B)$ is \defn{$\gamma$-compatible with~$\cD$}
    if, whenever $s$ $\gamma$-crosses a petal-separation~$t$ of~$\cD$, then all of the following assertions hold:
    \begin{enumerate}
        \item $Z = V_e \cap V_f$ and $Z$ is the $st$-centre, where~$e,f$ are the edges incident with the node of~$O$ corresponding to~$t$; and
        \item $V_e \sm Z$ and $V_f \sm Z$ are the $At$-link and the $Bt$-link; and
        \item all links of the two separations~$s,t$ are of the same size.
    \end{enumerate}
\end{definition}

Intuitively, we require that whenever $s$ $\gamma$-crosses a petal-separation $t$ of the flower, this crossing happens in the way we expect it to, namely with the centre the ``centre of the flower'', the $t$-links being $V_e \sm V_f$ and $V_f \sm V_e$, and with all links of the same size.

\smallskip

Let $\mathcal{D} = (O, \mathcal{V})$ be a graphic-$k$-daisy.
Let $e \neq f$ be two edges of $O$, and let $X, Y$ be the two connected components of $O \sm \{e, f\}$.
Let $s$ be the separation induced by $(X, Y)$.
Then, $s$ is a displayed separation of $\cD$.
We say that $s$ is \defn{displayed by $\{e, f\}$}. Note that the separation $s^*$ is also displayed by $\{e, f\}$.

\begin{definition}[$q$-lean]
Let $\mathcal{F} = (H, \cV)$ be a graphic-$k$-flower of a graph $G$.
For every $e \in E(H)$, let $V_e$ be the adhesion-set at $e$.
Let $Z = \bigcap_{h \in V(H)}V_h$.
If $\mathcal{F}$ is a graphic-$k$-daisy, we say that $\mathcal{F}$ is \defn{$q$-lean} if for every $(>q)$-broad separation of~$G$ displayed by a set of edges $\{e, f\}$, the sets $V_e \sm Z$ and $V_f \sm Z$ are disjoint.
If $\mathcal{F}$ is a graphic-$k$-anemone, we say that~$\mathcal{F}$ is \defn{$q$-lean} if $|Z| = k$.
\end{definition}

\begin{definition}[Clover-property]
    Let $F=(O,\cV)$ be a graphic-$k$-daisy of a graph~$G$.
    We say that $F$ has the \defn{clover-property} if the centre $X := \bigcap_{o \in V(O)} V_o$ has size~$<k$ and~$F$ has four adhesion-sets $V_{e_0},V_{e_1},V_{e_2},V_{e_3}$ with $V_{e_i} \cap V_{e_j} = X$ for all $i \neq j$.
\end{definition}

\begin{definition}[Canonical]
Let~$S$ be a set of separations of a graph~$G$.
We say that~$S$ is \defn{canonical} if for every automorphism~$\varphi$ of~$G$ we have $\varphi(S)=S$.
Similarly, we call a tree-decomposition~$\cT=(T,\cV)$ \defn{canonical} if the set of separations that are induced by an edge of~$T$ is canonical.
\end{definition}

\begin{definition}[$\psi$-isomorphic]
\label{defn:psi-isomorphic-graph-decomposition}
Let $\cH = (H, \cV)$ be a graph-decomposition of a graph $G$ and $\cH' = (H', \cV')$ be a graph-decomposition of a graph $G'$ and let $\psi\colon V(G) \to V(G')$ be an isomorphism between $G$ and~$G'$. 
We say that $\cH$ and $\cH'$ are \defn{$\psi$-isomorphic} if there exists an isomorphism $\phi \colon V(H) \to V(H')$ such that for every $h \in V(H)$, we have $\psi(V_h) = V'_{\phi(h)}$.
\end{definition}

\begin{restatable}[Flower Approximation Lemma]{keylemma}{flowerapprox}\label{keylem:flower-approximation}
For all $q, k \in \N$, there exist constants\\
$\beta_{\ref{keylem:flower-approximation}}(q,k)\ge q$, $\gamma_{\ref{keylem:flower-approximation}}(q,k)$ and $\eps_{\ref{keylem:flower-approximation}}(q,k)$ such that the following holds.
Let $G$ be a $(q,k-1)$-unbreakable graph and $\gamma>\gamma_{\ref{keylem:flower-approximation}}$.
Let $S$ be a nontrivial $\gamma$-asteriskal set of separations of~$G$ that consists of $\beta_{\ref{keylem:flower-approximation}}$-broad narrow $k$-separations of~$G$.
Then there is a $q$-lean graphic-$k$-flower $F:=F_{\ref{keylem:flower-approximation}}(S)$ of $G$ such that:
\begin{enumerate}
    \item\label{item:keylem:flower-approximation-1} every separation in~$S$ 
    has edit-distance $\le\eps_{\ref{keylem:flower-approximation}}$ to some displayed separation of~$F$; 
    \item\label{item:keylem:flower-approximation-2} for every $\eps \geq \eps_{\ref{keylem:flower-approximation}}$, every $k$-separation of~$G$ that is $\eps$-compatible with~$F$ is $\eps$-nested with all the petal-separations of~$F$;  
    \item\label{item:keylem:flower-approximation-3} if~$F$ is a graphic-$k$-daisy, then $F$ has the clover-property; and
    \item\label{flower-canonical}
    The graphic-$k$-flower $F_{\ref{keylem:flower-approximation}}(S)$ is canonical: whenever $\psi$ is an automorphism of~$G$, the graphic-$k$-flowers $F_{\ref{keylem:flower-approximation}}(S)$ and $F_{\ref{keylem:flower-approximation}}(\psi(S))$ are $\psi$-isomorphic.
\end{enumerate}
\end{restatable}

 The constants will be (for $q, k \in \N$): 
    \begin{align*}
        \beta_{\ref{keylem:flower-approximation}}(q,k)&{}:=\beta_{\ref{lem:coveringDisjointPaths}}(q,k)\\
        \gamma_{\ref{keylem:flower-approximation}}(q,k)&{}:=\max\,\{\,\gamma_{\ref{coveringCycleLong}}(q,k),\gamma_{\ref{coveringLinkSize}}(q,k)\,\}\\
        \eps_{\ref{keylem:flower-approximation}}(q,k)&{}:=\max\,\{\,\eps_{\ref{coveringDaisy}}(q,\ell): 0\le \ell\le k/2\,\}
    \end{align*}

We start by establishing some basic properties of asteriskal sets of separations in \cref{subsec:prop:asterisks}.
In particular, if~$S$ is a $\gamma$-asteriskal set of separations, then for each $s\in S$ there is a tripartition $(X_s,\{P_s,Q_s\})$ of $\sep(s)$ such that every $s'\in S$ that $\gamma$-crosses $s$ induces $(X_s,\{P_s,Q_s\})$ on~$\sep(s)$. Furthermore, there exists a set~$X$ and an integer~$\ell$ such that for every $s \in S$ we have $X_s = X$ and $|P_s| = \ell = |Q_s|$.

When~$X$ is the separator of all separations in $S$, the graphic-$X$-anemone (\cref{defn:X-anemone}) satisfies the outcome of \cref{keylem:flower-approximation}. This easy case is handled in \cref{subsec:anemone}.

The heart of the proof lies in the case where the separations in~$S$ do not satisfy $X = \sep(s)$. Then, we build a graphic-$k$-daisy that satisfies the outcome of \cref{keylem:flower-approximation}.
The construction roughly goes as follows.
\begin{itemize}[leftmargin=1em]
    \item First, we consider a covering $(\hat G, p)$ of $G-X$, where~$\hat G$ is an infinite graph and $p\colon V(\hat G)\to~V(G-X)$ is a local isomorphism.
    To construct this covering, we consider an arbitrary separation $s = (A, B) \in S$, and obtain~$\hat G$ by iteratively gluing copies of $G[A-X]$ and $G[B-X]$ on~$P_s$ and~$Q_s$ (see \cref{fig:CoveringDefinition}).
    This construction is detailed in \cref{subsec:covering}, where we also demonstrate some of its key properties. In particular, the covering does not depend on the separation $s \in S$ we initially picked. Furthermore, the graph~$\hat G$ has exactly two ends~$\omega_1$ and~$\omega_2$, and every separation of~$\hat G$ that distinguishes them has order at least~$\ell$.
    \item Then, denote by~$R$ the set of all $\ell$-separations $(U, W)$ of~$\hat G$ such that~$\omega_1$ lives in~$U$ and~$\omega_2$ lives in~$W$.
    Say that a separation $s \in R$ is \emph{locally $R$-optimal} if for every separation $r \in R$ that crosses~$s$, we have $\cross_R(s) < \cross_R(r)$, where~$\cross_R(\cdot)$ denotes the number of separations in~$R$ that cross a given separation.
    Let~$L$ be the set of all locally $R$-optimal separations.
    Note that the separations in~$L$ are pairwise nested.
    In \cref{subsec:L-approximates-R}, we prove that every separation in~$R$ is crossed by a bounded number of separations in~$R$.
    From this, we deduce that every separation $r \in R$ is at bounded edit-distance from a separation $s(r) \in L$.
    \item Finally, since the separations in~$L$ are pairwise nested, they induce a unique (canonical) tree-decomposition of~$\hat G$. This tree-decomposition of~$\hat G$ projects (via~$p$) to a (canonical) cycle-decomposition of $G-X$. 
    In \cref{subsec:daisy}, we show that after adding~$X$ to all bags, this cycle-decomposition is a graphic-$k$-daisy of~$G$ that satisfies the outcome of \cref{keylem:flower-approximation}.
\end{itemize}

\subsection{Centres and links of asteriskal sets of separations} \label{subsec:prop:asterisks}

In this section, we prove some basic properties of asteriskal sets of separations.

\begin{lem}\label{coveringCentreAndLinksOfS}
    Let $G$ be a graph and $\gamma\in\N$.
    Let~$S$ be a nontrivial $\gamma$-asteriskal set of separations of~$G$.
    Then:
    \begin{enumerate}
        \item\label{coveringCentreAndLinksOfS:trip} For every $s\in S$ there is a unique tripartition $(X_s,\{P_s,Q_s\})$ of $\sep(s)$ such that every $s'\in S$ that $\gamma$-crosses $s$ induces $(X_s,\{P_s,Q_s\})$ on~$\sep(s)$.
        \item\label{coveringCentreAndLinksOfS:centre} There is a unique $X\se V(G)$ such that $X=X_s$ for all $s\in S$.
    \end{enumerate}
\end{lem}
\begin{proof}
    \ref{coveringCentreAndLinksOfS:trip}.
    By $|S| \geq 2$ and since~$S$ is $\gamma$-crossing-connected, for every $s \in S$ there is some $s' \in S$ that $\gamma$-crosses~$s$.
    As $S$ is uniformly $\gamma$-crossing, the tripartition $(X_s,\{P_s,Q_s\})$ is unique.

    \ref{coveringCentreAndLinksOfS:centre} follows from \ref{coveringCentreAndLinksOfS:trip} and from~$S$ being $\gamma$-crossing-connected.
\end{proof}

We refer to $X$ as in \cref{coveringCentreAndLinksOfS}~\cref{coveringCentreAndLinksOfS:centre} as the \defn{common centre} of~$S$.

\begin{lem}\label{coveringLinkSize}
    Let $G$ be a $(q,k-1)$-unbreakable graph.
    Let $\gamma\in\N$ such that $\gamma\ge\gamma_{\ref{coveringLinkSize}}(q,k):=q+k$.
    Let~$S$ be a nontrivial $\gamma$-asteriskal set of separations of~$G$ that consists of $k$-separations.
    For every $s\in S$ let $(X_s,\{P_s,Q_s\})$ be as in \cref{coveringCentreAndLinksOfS}~\cref{coveringCentreAndLinksOfS:trip}.
    Then there is a unique $\ell\in\N$ such that $|P_s|=|Q_s|=\ell$ for all $s\in S$.
\end{lem}
\begin{proof}
    Let~$X$ be the common centre of~$S$ and define $\ell := (k-|X|)/2$.
    We show that $|P_s|=|Q_s|=\ell$ for all $s \in S$.
    Assume otherwise, that there is a separation~$s \in S$ with $|Q_s| < \ell < |P_s|$, say.
    Let~$s'$ be a separation in~$S$ that $\gamma$-crosses~$s$.
    Say, $|Q_{s'}| \leq \ell \leq |P_{s'}|$.
    Let~$K$ be the corner that is adjacent to~$Q_s$ and~$Q_{s'}$, and let~$(C,D)$ be the $K$-corner-separation with $K=C \sm D$, say.
    Then~$(C,D)$ has order $|Q_s| + |Q_{s'}| + |X| < 2\ell+|X| \leq k$.
    Since~$s$ and~$s'$ $\gamma$-cross, they $(\gamma-k)$-criss by \cref{crissCross}.
    Therefore, $|C \sm D| = |K| > \gamma-k \geq q$ and $|D \sm C| > 3(\gamma -k) \geq q$ (using the assumption that $\gamma\ge q+k$).
    Hence, $(C,D)$ witnesses that~$G$ is not $(q,k-1)$-unbreakable, a contradiction.
\end{proof}

We refer to~$\ell$ as in \cref{coveringLinkSize} as the \defn{common link size} of the $\gamma$-asteriskal set of separations~$S$.
We remark that this definition requires $\gamma\ge\gamma_{\ref{coveringLinkSize}}(q,k)$.

\subsection{Constructing the graphic-anemone} \label{subsec:anemone}

If the common link size is zero, all separators are equal to the common centre $X$. In this section, we show that the graphic-$X$-anemone displays every separation in $S$.

\begin{definition}[Graphic-$X$-anemone]
\label{defn:X-anemone}
    Let~$X$ be a set of vertices of a graph~$G$ and let~$\cC$ be the set of components of $G-X$.
    The \defn{graphic-$X$-anemone} of~$G$ is the graphic-$|X|$-anemone $\cA=(S,\cV)$, where~$S$ is a star with~$\cC$ as set of leaves, $V_c=X$ is the central bag, and $V_C=V(C)\cup X$ for every component~$C\in\cC$.
\end{definition}

\begin{lem}\label{coveringAnemone}
    Let $G$ be a $(q,k-1)$-unbreakable graph and let~$\gamma\in\N$.
    Let~$S$ be a nontrivial $\gamma$-asteriskal set of separations of~$G$ that consists of $k$-separations.
    Let~$X$ denote the common centre of~$S$ and assume that $|X|=k$.
    Let $\cA$ be the graphic-$X$-anemone.
    Then, $\cA$ is a $q$-lean graphic-$k$-anemone that displays every separation in~$S$, and for every $\eps\in \Nbb$, every separation which is $\eps$-compatible with $\cA$ is $\eps$-nested with all the petal-separations of~$\cA$.
\end{lem}
\begin{proof}
    $\cA$ is a $q$-lean graphic-$k$-anemone since $|X| = k$.
    Let $s\in S$.
    Then $\sep(s)=X$ as $X$ is the common centre of~$S$, and has size $k = |\sep(s)|$.
    Hence $s$ is displayed by~$\cA$.

    Let $\eps \in \Nbb$ and let $s$ be a separation which is $\eps$-compatible with $\cA$.
    By contradiction, suppose that $s$ $\eps$-crosses some petal-separation~$t$ of $\cA$.
    Since $s$ is $\eps$-compatible with $\mathcal{A}$, we have that $s$ and $t$ cross with centre $X$ and four empty links.
    Thus, $X = \sep(s)$ so~$s$ is displayed by~$\cA$, so cannot cross $t$.
\end{proof}

\begin{rem}\label{coveringAnemonecanonical}
    The graphic-$k$-anemone $\cA=\cA(S)$ in \cref{coveringAnemone} is canonical: 
    whenever $\psi$ is an automorphism of~$G$,
    the graphic-$k$-anemones
    $\cA(S)$ and $\cA(\psi(S))$ are $\psi$-isomorphic.
\end{rem}

\subsection{The covering construction}\label{subsec:covering}

If the common link size is $\ell > 0$, we build a covering $(\hat G, p)$ of~$G-X$.
In this section, we describe the construction of this covering (see \cref{CoveringDefinition}) and analyse some of its properties.
We show that~$\hat G$ has exactly two ends and that every separation of~$\hat G$ that distinguishes them has order at least~$\ell$ (see \cref{lem:inducedCoveringEnds}).
We also show that the covering $(\hat G, p)$ is canonical given $(G, S)$ (see \cref{lem:covering-equivalent}).

A \defn{ray} is a one-way infinite path.
A \defn{tail} of a ray means a subray.
An \defn{end} of a graph~$G$ is an equivalence class of rays in~$G$, where two rays in~$G$ are considered \defn{equivalent} if they are linked by infinitely many pairwise disjoint paths.
We say that a ray of an infinite graph~$G$ \defn{lives} in a set~$A \subseteq V(G)$ if it has a tail contained in~$A$.
An end~$\omega$ of an infinite graph~$G$ \defn{lives} in a set~$A \subseteq V(G)$ if every ray contained in~$\omega$ lives in $A$.
Note that, given a separation~$(A,B)$ of~$G$ of finite order, any end~$\omega$ of~$G$ lives either in~$A$ or in~$B$, but not in both.
A separation~$s$ of~$G$ \defn{distinguishes} two ends~$\omega_1$ and~$\omega_2$ if they live in distinct sides of~$s$.
We say that~$s$ distinguishes~$\omega_1$ and~$\omega_2$ \defn{efficiently} if~$s$ is of smallest order among all separations that distinguish~$\omega_1$ and~$\omega_2$.

\begin{figure}[ht]
    \centering
    \includegraphics[height=14\baselineskip]{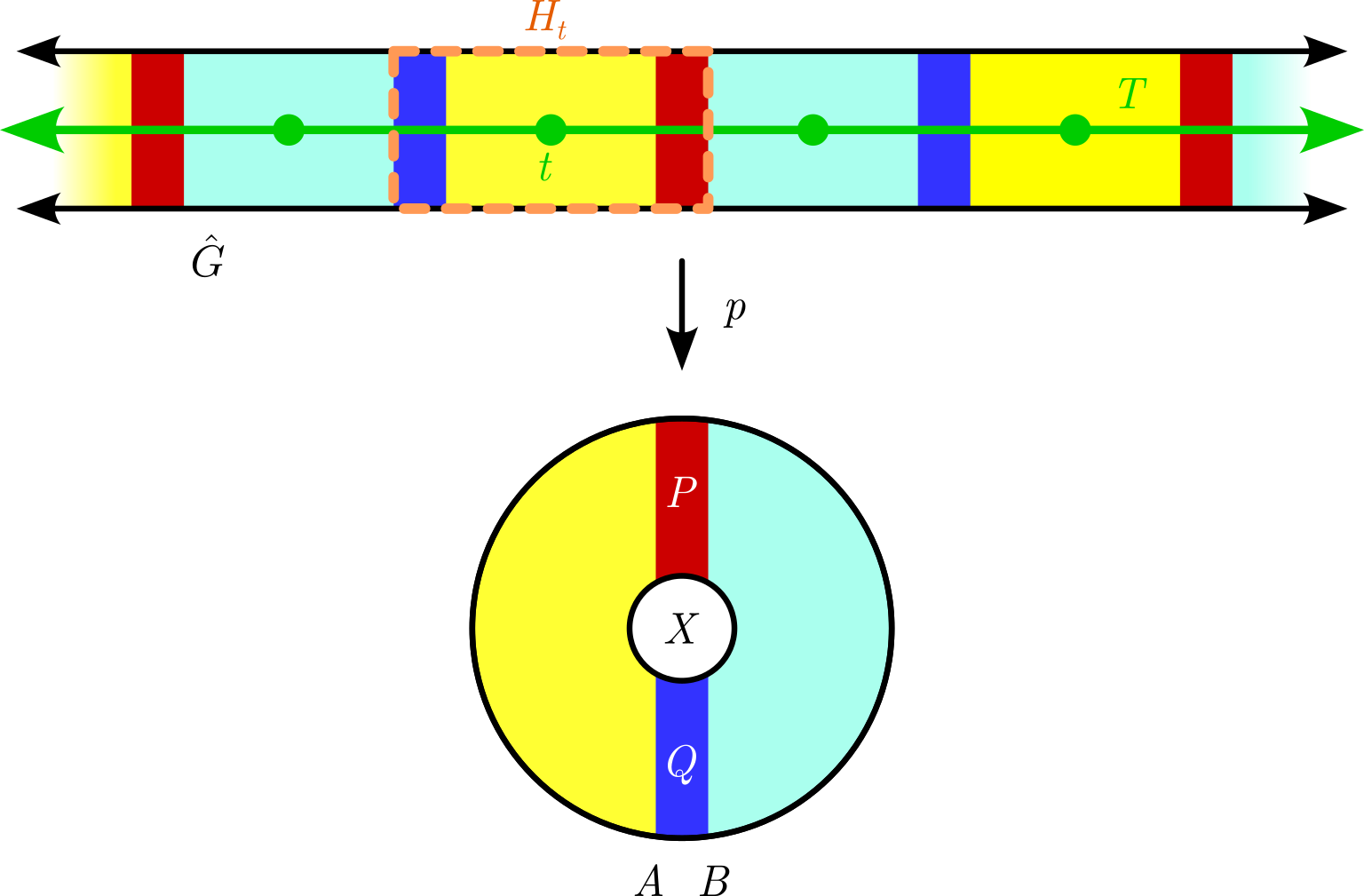}
    \caption{The covering $p\colon \hat{G}\to G-X$ that is constructed in \cref{CoveringDefinition}. On the one hand, $H_A$ consists of yellow, red and blue. On the other hand, $H_B$ consists of turquoise, red and blue.}
    \label{fig:CoveringDefinition}
\end{figure}

\begin{construction}[Covering]\label{CoveringDefinition}
This definition is supported by \cref{fig:CoveringDefinition}.
Let $(A,B)$ be a narrow proper $k$-separation of~$G$.
Let $(X,\{P,Q\})$ be a tripartition of $A\cap B$ such that $|X|<k$, $|P|=|Q|$ and there is no edge in between~$P$ and~$Q$.
Let $H_A:=G[A]-X$ and $H_B:=G[B]-X$.
Let $T$ be a double-ray $\ldots t(-1)\, t(0)\, t(1)\ldots$, viewed as a bipartite graph with classes~$T_A$ and~$T_B$.
For each $t\in T_A$ let~$H_t$ be a copy of~$H_A$, and for each $t\in T_B$ let~$H_t$ be a copy of~$H_B$.
Let~$P_t$ and~$Q_t$ denote the copies of~$P$ and~$Q$ in~$H_t$, respectively.
We choose the graphs~$H_t$ so that they are pairwise disjoint except that they intersect as follows for every edge $t(i)\, t(i+1)$ of~$T$:
\begin{itemize}
    \item $P_{t(i)}=P_{t(i+1)}=:L(i)$ if $t(i)\in T_A$ and $t(i+1)\in T_B$, and
    \item $Q_{t(i)}=Q_{t(i+1)}=:L(i)$ if $t(i)\in T_B$ and $t(i+1)\in T_A$.
\end{itemize}
Hence $V(H_{t(i)})\cap V(H_{t(i+1)})=L(i)$ and $V(H_{t(i)})\cap V(H_{t(j)})=\emptyset$ whenever $|i-j|>1$.
Let $\hat G :=\bigcup_{t\in V(T)}H_t$.
Note that~$\hat G$ is connected since $(A,B)$ is narrow.
Let $V_t:=V(H_t)$.
Then $(T,\cV)$ with $\cV=(V_t:t\in T)$ is a tree-decomposition of~$\hat G$ with pairwise disjoint adhesion-sets~$L(i)$.
Moreover, the map $p\colon \hat G\to G-X$ that sends every vertex in a graph~$H_t$ to the vertex of $G-X$ it stems from is a covering map (as there is no edge in between~$P$ and~$Q$).
We call~$p$ the \defn{covering induced by}~$(A,B)$ and $(X,\{P, Q\})$, as~$p$ is unique (up to equivalence of coverings).
By a slight abuse of notation, we also call~$\hat G$ the \defn{induced covering}.
We call $(T,\cV)$ the \defn{standard decomposition} of~$\hat G$.
\end{construction}

The following lemma will help us prove that every separation that distinguishes the two ends of~$\hat G$ has order at least $\ell$.

\begin{lem}\label{lem:coveringDisjointPaths}
    Let~$G$ be a $(q,k-1)$-unbreakable graph.
    Let $\beta:=\beta_{\ref{lem:coveringDisjointPaths}}(q,k):=2q+k$.
    Let~$(A,B)$ be a $\beta$-broad $k$-separation of~$G$ and let $(X,\{P,Q\})$ be a tripartition of~$A \cap B$ with $|P|=|Q|=:\ell$.
    Then there are~$\ell$ disjoint $P$--$Q$ paths in~$G[A]-X$.
\end{lem}
\begin{proof}
    Assume otherwise.
    Recall that $H_A=G[A]-X$.
    Then by Menger's Theorem, there is a $(<\ell)$-separation $(C,D)$ of~$H_A$ with $P \subseteq C$ and $Q \subseteq D$.
    Then $(C_1,D_1):=(C\cup X,D \cup B)$ and $(C_2,D_2):=(D\cup X,C\cup B)$ are $(<k)$-separations of~$G$.
    Since~$G$ is $(q,k-1)$-unbreakable and $|D_i \sm C_i| \geq |B \sm A| > q$, we must have $|C_i \sm D_i| \leq q$ for $i \in \{1,2\}$.
    Then, however, 
    \[|A \sm B| \le |C_1 \sm D_1| + |C_2 \sm D_2| + |C \cap D| < 2q+\ell\le \beta~,\]
    a contradiction.
\end{proof}

\begin{lem}\label{lem:inducedCoveringEnds}
    Let $G$ be a $(q,k-1)$-unbreakable graph.
    Let~$(A,B)$ be a $\beta_{\ref{lem:coveringDisjointPaths}}(q,k)$-broad narrow $k$-separation of~$G$ and let~$(X,\{P,Q\})$ be a tripartition of~$A \cap B$ with $|X|<k$ and $|P|=|Q|=:\ell$ and no edge in between~$P$ and~$Q$.
    Let~$\hat G$ be the covering induced by~$(A,B)$ and~$(X,\{P,Q\})$.
    Then:
    \begin{enumerate}
        \item\label{itm:inducedCoveringEnds-1} $\hat G$ has exactly two ends~$\omega_1$ and~$\omega_2$;
        \item\label{itm:inducedCoveringEnds-2} Every separation of~$\hat G$ distinguishing both ends has order~$\geq \ell$.
    \end{enumerate}
\end{lem}
\begin{proof}
    \ref{itm:inducedCoveringEnds-1}.
    Since~$(A,B)$ is narrow, and in particular tight, there is a component of $G[A\sm B]$ with neighbourhood equal to $A\cap B$.
    Hence every two vertices in~$P \cup Q$ are joined by a path in~$H_A:=G[A]-X$, and similarly in~$H_B:=G[B]-X$.
    Therefore, if two rays~$R_1, R_2$ in~$\hat G$ follow~$T$ in the same direction, there are infinitely many vertex-disjoint $R_1$--$R_2$ paths in~$\hat G$, so~$R_1$ and~$R_2$ lie in a common end.
    As every ray in~$\hat G$ follows~$T$ in some direction,~$\hat G$ has at most two ends.

    Conversely, recall that~$(A,B)$ is tight and~$|P|=|Q| \geq 1$.
    Hence, there are two rays~$R_1$ and~$R_2$ in~$\hat G$ that follow~$T$ in opposite directions.
    Then any adhesion-set~$L(i)$ of~$(T, \cV)$ separates~$R_1$ and~$R_2$.
    Hence,~$\hat G$ has at least two ends.

    \ref{itm:inducedCoveringEnds-2}.
    By \cref{lem:coveringDisjointPaths}, there are~$\ell$ vertex-disjoint $P$--$Q$ paths contained in~$H_A:=G[A]-X$.
    Similarly, there are~$\ell$ vertex-disjoint $P$--$Q$ paths in~$H_B:=G[B]-X$.
    The copies of these paths in~$\hat G$ combine to~$\ell$ vertex-disjoint double-rays between the two ends of~$\hat G$.
    Hence every separator that distinguishes the two ends has order at least~$\ell$.
\end{proof}

Recall that two coverings $p_1\colon \hat G_1\to G$ and $p_2\colon \hat G_2\to G$ are \defn{equivalent} if there is a graph-isomorphism $f\colon \hat G_1\to \hat G_2$ such that $p_1=p_2\circ f$.

\begin{lem}\label{lem:covering-equivalent}
    Let~$G$ be a $(q,k-1)$-unbreakable graph and $\gamma\in\N$.
    Let~$S$ be a nontrivial $\gamma$-asteriskal set of separations of~$G$ that consists of narrow $k$-separations, where $\gamma\ge\gamma_{\ref{coveringLinkSize}}(q,k)$.
    For each~$s \in S$, let~$(X,\{P_s,Q_s\})$ be as in \cref{coveringCentreAndLinksOfS}, where~$X$ is the common centre of~$S$.
    Assume that $|X|<k$.
    Then:
    \begin{enumerate}
        \item\label{lem:covering-equivalent:exist} Every $s\in S$ and $(X,\{P_s,Q_s\})$ induce a covering $p_s\colon \hat G_s\to G-X$.
        \item\label{lem:covering-equivalent:eq} All coverings $p_s$ for $s\in S$ are pairwise equivalent.
    \end{enumerate}
\end{lem}
\begin{proof}
    \cref{lem:covering-equivalent:exist}.\
    Let $s\in S$.
    Since $S$ is nontrivial, there is $s'\in S$ that $\gamma$-crosses~$s$.
    As~$S$ is uniformly $\gamma$-crossing, $s'$ induces $(X,\{P_s,Q_s\})$ on the separator of~$s$.
    Since~$s'$ is a separation, it follows that~$G$ contains no edge in between~$P_s$ and~$Q_s$.
    Therefore, the premises of \cref{CoveringDefinition} are met, so~$s$ and $(X,\{P_s,Q_s\})$ induce a covering of~$G-X$.

    \cref{lem:covering-equivalent:eq}.\
    Let $s\neq s'$ be separations in~$S$.
    Since~$S$ is $\gamma$-crossing-connected, it suffices to prove it under the assumption that~$s$ and~$s'$ $\gamma$-cross.
    Write $(A,B):=s$ and $(C,D):=s'$.
    Since~$S$ is uniformly $\gamma$-crossing, we have $P_{s'}\se A\sm B$ and $Q_{s'}\se B\sm A$, say.
    Adding $P_{s'}$ and $Q_{s'}$ to \cref{fig:CoveringDefinition}, it is straightforward to see that the covering induced by~$s$ and~$(X,\{P_s,Q_s\})$ is equivalent to the covering induced by~$s'$ and~$(X,\{P_{s'},Q_{s'}\})$.
\end{proof}

Let~$S$ and $X$ be as in the assumptions of \cref{lem:covering-equivalent}.
By \cref{lem:covering-equivalent}, each $s\in S$ induces a covering of~$G-X$, and all these coverings are equivalent.
Hence, we may speak of \defn{the covering induced by the $\gamma$-asteriskal set of separations~$S$}.

\subsection{Approximating by locally-optimal separations}\label{subsec:L-approximates-R}

Given a set~$S$ of separations of a graph~$G$ and a separation $s \in S$, we denote by~$\defnMath{\cross_S(s)}$ the \defn{$S$-crossing-number} of~$s$, that is, the number of separations in~$S$ that cross~$s$.
Bounds for the $S$-crossing-number will usually use the letter~$\chi$.
A separation $s\in S$ is \defn{locally $S$-optimal} if $\cross_S(s)<\cross_S(r)$ for every $r\in S$ that crosses~$s$.

Throughout the next two sections, we will work in the following setting.

\begin{setting}\label{set:covering}
    Let~$G$ be a $(q,k-1)$-unbreakable graph, and let $\gamma\in\N$ such that $\gamma\ge\gamma_{\ref{coveringLinkSize}}(q,k)$.
    Let~$S$ be a nontrivial $\gamma$-asteriskal set of separations of~$G$ that consists of $\beta_{\ref{lem:coveringDisjointPaths}}(q,k)$-broad narrow $k$-separations.
    Let~$X$ denote the common centre of~$S$ and assume that $|X|<k$.
    Let $\ell$ denote the common link size of~$S$.
    Let~$p\colon \hat G \to G-X$ denote the covering of~$G-X$ induced by~$S$.
    By \cref{lem:inducedCoveringEnds},~$\hat G$ has exactly two ends~$\omega_1$ and~$\omega_2$.
    Let~$R$ be the set of all $\ell$-separations $(U,W)$ of~$\hat G$ such that $\omega_1$ lives in~$U$ and $\omega_2$ lives in~$W$.
    Let $L$ be the set of all locally $R$-optimal separations in~$R$, which is nested by definition.
\end{setting}

In this section, we show that every separation $r \in R$ is at bounded edit-distance from some separation $s(r) \in L$ (see \cref{coveringLemmaFive}). 
To do so, we first show that every separation $r \in R$ is only crossed by a bounded number of separations $r' \in R$ (see \cref{mu-bound-cov}, which immediately follows from \cref{coveringLemmaTwo,conn-cov}). 
Then, we argue that whenever two separations in~$R$ cross, they cross with two small opposite corners (see \cref{lem:corner-in-cover}). 
To prove \cref{coveringLemmaFive}, we then argue that starting from a separation $r \in R$, it suffices to ``wiggle'' it a bounded number of times by a small amount to obtain a separation in~$L$.

We will use the following result of Elbracht, Kneip and Teegen. Recall that a (potentially infinite) graph is \defn{locally finite} if every vertex has finite degree.
\begin{lem}\label{HalinForTightSeparations}
    \cite[Proof of Proposition~6.2]{infinitesplinter}
    Let $G$ be a locally finite connected graph.
    Let $(A,B)$ be a tight finite-order separation of~$G$, and let $k\in\N$.
    Then only finitely many tight separations of~$G$ of order~$\le k$ cross~$(A,B)$.
\end{lem}

We now show that~$R$ is closed under taking meets and joins, and that these operations can be used to create separations with smaller $R$-crossing-number.

\begin{lem}[{\cite[Proof of Corollary~3.3]{entanglements}}]
\label{lem:crossing-number-gets-smaller}
    Let~$r$ be a separation of a graph~$G$, let~$S$ be a set of separations of~$G$ and let~$s \in S$.
    Assume that both $\cross_S(r)$ and~$\cross_S(s)$ are finite and that~$r$ and~$s$ cross.
    Then \[
        \cross_S(r \wedge s) + \cross_S(r \vee s) < \cross_S(r) + \cross_S(s)~.
    \]
\end{lem}
\begin{proof}
    If a separation $t \in S$ satisfies $t \leq r$ or $t \leq s$, then we have $t \leq r \vee s$ since $r,s \leq r \vee s$.
    Similarly, $r \leq t$ or $s \leq t$ implies $r\wedge s \leq t$ since $r \wedge s \leq r,s$.
    Hence, if a separation~$t \in S$ is nested with~$r$ or~$s$, then it is nested with~$r \wedge s$ or $r \vee s$.
    Moreover, if~$t \in S$ is nested with both~$r$ and~$s$, then it is nested with both $r \wedge s$ and $r \vee s$ (since $r$ and $s$ cross).
    It follows that $\cross_S(r \wedge s) + \cross_S(r \vee s) \leq \cross_S(r) + \cross_S(s)$.

    Since $s \in S$ crosses~$r$, the separation~$s$ is counted once on the right side of the inequality.
    Since~$s$ is nested with~$r \wedge s$ and~$r \vee s$, it is not counted on the left side of the inequality.
    So, we have $\cross_S(r \wedge s) + \cross_S(r \vee s) < \cross_S(r) + \cross_S(s)$.
\end{proof}

\begin{lem}\label{coveringCornersAreEfficient}
    Assume \cref{set:covering}.
    Let~$s_1,s_2$ be two separations in~$R$.
    Let $c_1:=s_1 \wedge s_2$ and $c_2:=s_1\vee s_2$.
    Then both $c_1$ and $c_2$ are contained in~$R$.
    Moreover, if~$s_1$ and~$s_2$ cross and $\cross_R(s_1)\le\cross_R(s_2)$, then there is $i\in [2]$ such that $\cross_R(c_i)<\cross_R(s_2)$.
\end{lem}
\begin{proof}
    If~$s_1$ and~$s_2$ do not cross then we have $s_1 \leq s_2$ or $s_2 \leq s_1$ since~$\omega_1$ lives on the left side of each~$s_i$ and~$\omega_2$ lives on the right side of each~$s_i$.
    Therefore, $\{c_1,c_2\}=\{s_1,s_2\}\se R$ and we are done.
    
    Hence, we may assume that $s_1$ and $s_2$ cross. 
    Both $c_1$ and $c_2$ distinguish $\omega_1$ and~$\omega_2$, hence both have order~$\ge\ell$ by \cref{lem:inducedCoveringEnds}.
    Since $c_1$ and $c_2$ are opposite corner-separations, it follows from submodularity that both~$c_i$ have order equal to~$\ell$.
    Moreover,~$\omega_1$ lives on the left side of each~$c_i$, and~$\omega_2$ lives on the right side of each~$c_i$.
    Hence both $c_i$ are contained in~$R$.

    Assume now that $\cross_R(s_1)\le\cross_R(s_2)$.
    Note that every separation $s \in R$ is tight since it efficiently distinguishes $\omega_1$ and $\omega_2$.
    Thus, both $\cross_R(s_i)$ are finite by \cref{HalinForTightSeparations}.
    By \cref{lem:crossing-number-gets-smaller}, we have
    \[
        \cross_R(c_1)+\cross_R(c_2)<\cross_R(s_1)+\cross_R(s_2)\le 2\cross_R(s_2).
    \]
    Hence there is $i\in [2]$ such that $\cross_R(c_i)<\cross_R(s_2)$.
\end{proof}

The next result is our first step towards proving that $\cross_R(\cdot)$ is bounded.

\begin{lem}\label{coveringLemmaTwo}
    Assume \cref{set:covering}.
    Let $\beta\in\N$ and let $\chi_{\ref{coveringLemmaTwo}}(\beta,\ell):=3^{2\beta+3\ell}$.
    Let $s=(A,B)\in R$, and let~$R_s$ be the set of all $(C,D)\in R$ that cross $(A,B)$.
    Then at least one of the following holds:
    \begin{enumerate}
        \item\label{itm:coveringLemmaTwo-1} $|R_s|< \chi_{\ref{coveringLemmaTwo}}(\beta,\ell)$, or
        \item\label{coveringLemmaTwoOutcomeUW} there is a $(<2\ell)$-separation $(U,W)$ of~$\hat{G}$ such that $|U\sm W|\ge\beta$ and both $\omega_i$ live in the side~$W$.
    \end{enumerate}
\end{lem}

We first give some intuition about the proof.
To prove that~$R_s$ has bounded size, it suffices to prove that for every tripartition $A \cap B = P \sqcup X \sqcup Q$, the set $R_s(P, X, Q)$ of all separations $t = (C, D) \in R_s$ that cross $s=(A, B)$ with~$P$ as the $sD$-link,~$X$ as the $st$-centre, and~$Q$ as the $sC$-link, has bounded size. 
Fix such a tripartition $A \cap B = P \sqcup X \sqcup Q$, and consider the set $R_s(P, X, Q)$. 
Consider the supremum and the infimum of all separations in $R_s(P, X, Q)$, and look at the region~$F$ between them (this is the union of the small sides of the green and yellow separations in \cref{fig:CoveringIteratedInfSup}).
Observe that each separation in $R_s(P, X, Q)$ is entirely defined by its restriction to~$F$. Thus, if~$|F|$ is small then $|R_s(P, X, Q)|$ is small. If~$|F|$ is large, then one of these two separations has a small side which is large, and thus satisfies the second condition.

\begin{proof}[Proof of \cref{coveringLemmaTwo}]
    The set~$R_s$ is finite by \cref{HalinForTightSeparations}.
    Let~$T$ denote the set of all tripartitions $(P,X,Q)$ of $A\cap B$.
    For each $(P,X,Q)\in T$, let $R_s(P,X,Q)$ denote the set of all $t = (C,D)\in R_s$ that cross~$(A,B)$ with~$P$ as the $sD$-link,~$X$ as the $st$-centre, and~$Q$ as the $sC$-link.
    Now choose $(P,X,Q)\in T$ with $R_s(P,X,Q)$ of maximum size.
    If $R_s(P,X,Q)$ is empty, then $|R_s|=0$ and we trivially have~\cref{itm:coveringLemmaTwo-1}.
    Hence we may assume that $R_s(P,X,Q)$ is not empty.

    Let $s_1=(X_1,Y_1)$ be the infimum of~$R_s(P,X,Q)$ and let~$s_2=(X_2,Y_2)$ be the supremum of~$R_s(P,X,Q)$.
    By \cref{coveringCornersAreEfficient},~$s_1$ and~$s_2$ are both in~$R$.
    Moreover, for every~$s \in R_s(P,X,Q)$ we have $s_1 \leq s \leq s_2$ by definition.
    
    \begin{claim}\label{coveringLemmaTwoInfSupLinks}
        Let $i\in [2]$ and consider the two separations $s_i=(X_i,Y_i)$ and $s=(A,B)$.
        Then the $sY_i$-link~$L(Y_i)$ equals~$P$, the $ss_i$-centre equals~$X$, and the $sX_i$-link~$L(X_i)$ equals~$Q$.
    \end{claim}
    \begin{cproof}
        Say $i=1$.
        Since $(P,X,Q)$ is a tripartition of $A\cap B$, it suffices to show the three inclusions 
        \[
            P\se L(Y_1),\quad X\se X_1\cap Y_1,\quad \text{and}\quad  Q\se L(X_1).
        \]
        We have $Q\se X_1\sm Y_1$ as $s_1$ is the infimum of $R_s(P,X,Q)$, and $Q\se A\cap B$ by definition.
        Hence $Q\se L(X_1)$.
        Similarly, $P\se L(Y_1)$.
        Finally,~$X$ is included in all separators of the separations in $R_s(P,X,Q)$, so~$X$ must also be included in the separator of the infimum of all these separations.
        Hence $X\se X_1\cap Y_1$.
    \end{cproof}

    The following definitions are supported by \cref{fig:CoveringIteratedInfSup}.
    Let
    \begin{align*}
        t_1:=(U_1,W_1) &:= s_1^* \wedge s_2 \wedge s = (Y_1,X_1) \wedge (X_2,Y_2) \wedge (A,B)\quad\text{and}\\
        t_2:=(U_2,W_2) &:= s_1^* \wedge s_2 \wedge s^* = (Y_1,X_1) \wedge (X_2,Y_2) \wedge (B,A)~.
    \end{align*}

    \begin{figure}[ht]
        \centering
        \includegraphics[height=6\baselineskip]{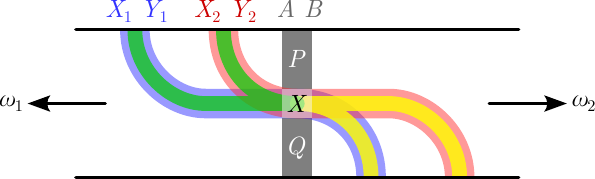}
        \caption{The situation in the proof of \cref{coveringLemmaTwo}. Here $U_1\cap W_1$ is green and $U_2\cap W_2$ is yellow.}
        \label{fig:CoveringIteratedInfSup}
    \end{figure}

    \begin{claim}\label{claim:coveringLemmaTwo}
        For both $i=1,2$ we have that $(U_i,W_i)$ is a $(<2\ell)$-separation such that both ends~$\omega_1$ and~$\omega_2$ live in~$W_i$.
    \end{claim}
    \begin{cproof}
        Recall that $\omega_1$ lives in~$X_1$ and~$\omega_2$ lives in~$Y_2$.
        Hence, both~$\omega_1$ and~$\omega_2$ live in~$X_1 \cup Y_2 \subseteq W_i$.
        Assume for a contradiction that $|U_1 \cap W_1| \geq 2\ell$, say.
        Let $Z:=(X_1\cap Y_1)\cup (X_2\cap Y_2)$.
        By the definition of infimum and supremum, we have
        \[
            U_1\cap W_1\se A\cap Z \quad \text{and} \quad U_2\cap W_2\se B\cap Z~,
        \]
        see \cref{fig:CoveringIteratedInfSup}.
        Hence
        \[
            |U_1\cap W_1|+|U_2\cap W_2|\le |A\cap Z|+|B\cap Z|\overset{\text{\cref{coveringLemmaTwoInfSupLinks}}}{=}|Z|+|X|\le 2\ell~.
        \]
        Now $|U_1\cap W_1|\ge 2\ell$ implies that $U_2\cap W_2$ is empty.
        Since both ends of~$\hat G$ live on the side~$W_2$ and since $\hat{G}$ is connected, the side~$U_2$ must be empty.
        By the definition of infimum we have
        \begin{equation}\label{eq:U2empty}
            Y_1\cap X_2\cap B=U_2=\emptyset~.
        \end{equation}
        Let $s'=(I,J)\in R_s(P,X,Q)$.
        Then $(X_1,Y_1)\le (I,J)\le (X_2,Y_2)$ by the definition of $(X_i,Y_i)$.
        Hence, 
        \[
            I\cap J\se Y_1\cap X_2\overset{\cref{eq:U2empty}}{\se} A\sm B.
        \]
        We have $(A, B), (I, J) \in R$, so $\omega_1$ lives in the $AI$-corner and $\omega_2$ lives in the $BJ$-corner.
        Since $I \cap J \subseteq A \setminus B$, the $ss'$-centre and the $Bs'$-link are empty.
        The $BJ$-corner-separation separates $\omega_1$ and $\omega_2$, so it has order $\geq \ell$ by \cref{lem:inducedCoveringEnds}. Therefore, the $sJ$-link has size $\ell$, so the $sI$-link is empty.
        Since $\hat{G}$ is connected, this implies that the $BI$-corner is also empty, which contradicts that $s$ and $s'$ cross.
    \end{cproof}

    By \cref{claim:coveringLemmaTwo}, we may assume that both $(U_i,W_i)$ (for $i=1,2$) are $(<2\ell)$-separations such that both ends $\omega_1,\omega_2$ live in~$W_i$.
    If for some $i\in [2]$ we have $|U_i\sm W_i|\ge\beta$, then we have \cref{coveringLemmaTwoOutcomeUW} and are done.
    So we may assume that $|U_i\sm W_i|<\beta$ for both~$i\in [2]$.

    Let $F:=Y_1\cap X_2$.
    Then
    \[
    F\se (X_1\cap Y_1)\cup (X_2\cap Y_2)\cup (U_1\sm W_1)\cup (U_2\sm W_2)
    \]
    which implies $|F|<2\beta+2\ell$.
    Every $(C,D)\in R_s(P,X,Q)$ satisfies
    \[
    (X_1,Y_1)\le (C,D)\le (X_2,Y_2).
    \]
    Hence $(C,D)$ is uniquely determined by $(C\cap F,D\cap F)$.
    There are at most $3^{|F|}$ separations of~$\hat G[F]$.
    Thus
    \[
    |R_s(P,X,Q)|<3^{2\beta+2\ell}.
    \]
    We also have $|T|\le 3^\ell$.
    Therefore,
    \[
    |R_s|\le |R_s(P,X,Q)|\cdot |T|<3^{2\beta+3\ell}=\chi_{\ref{coveringLemmaTwo}}(\beta,\ell)
    \]
    which completes the proof.
\end{proof}

The next result is our second step towards proving that $\cross_R(\cdot)$ is bounded.

\begin{lem}\label{conn-cov}
    Assume \cref{set:covering}.
    Let~$(U,W)$ be a $(<2\ell)$-separation of~$\hat{G}$ such that both~$\omega_1$ and~$\omega_2$ live in~$W$.
    Then $|U \setminus W| < \beta_{\ref{conn-cov}}(q,\ell):= 4\ell (q + 2\ell)$.
\end{lem}

We first give a quick proof sketch. 
Let $(A, B) \in S$ be any separation, and let $(X, \{P, Q\})$ be its corresponding tripartition.
We consider the standard decomposition of~$\hat G$ induced by $(A, B)$ and $(X, \{P, Q\})$ and analyse how~$U$ interacts with each bag~$V_t$ (recall that~$\hat G[V_t]$ is a copy either of $G[A-X]$ or of $G[B-X]$, with a copy~$P_t$ of~$P$ and a copy~$Q_t$ of~$Q$).
Let~$V_t$ be a bag such that $U \cap V_t \neq \emptyset$.
Consider the separation $(U_t, W_t) = (U \cap V_t, (W \cap V_t) \cup P_t \cup Q_t)$ (this is the yellow separation in \cref{fig:conn-cov:UtWt}).
We first argue that this separation has order $<2\ell$.
This separation (of $G[A-X]$, say) gives a $(<k)$-separation of~$G$ by adding~$X$ to the separator and $B \setminus A$ to the side of~$W_t$.
This then implies that $U_t = U \cap V_t$ is small.
Then, we argue that if $U \cap V_t \neq \emptyset$ then $U \cap W \cap V_t \neq \emptyset$.
Since non-consecutive bags of the tree-decomposition are disjoint, there can only be at most~$2k$ bags~$V_t$ such that $U \cap V_t \neq \emptyset$.
This immediately implies that~$U$ has bounded size.

\begin{proof}
    Let~$(A,B)\in S$ be arbitrary and let $(X,\{P,Q\})$ be the tripartition that all separations in~$S$ that $\gamma$-cross $(A, B)$ induce on~$A\cap B$.
    Consider the tree-decomposition $(T,\cV)$ of~$\hat G$ induced by~$(A,B)$ and~$(X,\{P,Q\})$.
    Without loss of generality, assume $T=\Zbb$.
    For $t \in V(T)$, we define $C_t:= V_t$ and $D_t:=\bigcup_{s \in V(T-t)} V_s$.
    Observe that~$(C_t,D_t)$ is a separation of~$\hat G$ of order~$|P\cup Q|=2\ell$ for every $t \in V(T)$.
    We define $(U_t,W_t):=(U,W) \wedge (C_t,D_t)$, see \cref{fig:conn-cov:UtWt}.

    \begin{figure}[ht]
        \centering
        \includegraphics[height=6\baselineskip]{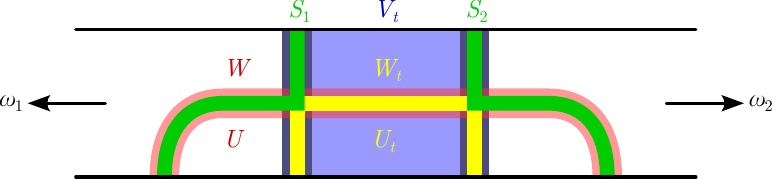}
        \caption{The interaction of $(U,W)$ with $V_t$ and the location of $(U_t,W_t)$}
        \label{fig:conn-cov:UtWt}
    \end{figure}

    Moreover, we let~$E_t:=\bigcup_{t' \leq t} V_{t'}$ and $F_t:=\bigcup_{t' > t} V_{t'}$.
    Observe that for every~$t \in T$,~$(E_t,F_t)$ is an $\ell$-separation such that~$\omega_1$ lives in~$E_t$ and~$\omega_2$ lives in~$F_t$.
    In particular,~$(E_t,F_t) \in R$.

    \begin{claim}
    \label{claim:conn-cov-1}
        $(U_t,W_t)$ has order~$<2\ell$ for every~$t \in V(T)$.
    \end{claim}
    \begin{cproof}        
        Consider the separation~$s_1:=(E_{t-1},F_{t-1}) \wedge (W,U)$ and let~$S_1$ be its separator, see \cref{fig:conn-cov:UtWt}.
        Note that~$s_1$ distinguishes~$\omega_1$ from~$\omega_2$, so $|S_1| \geq \ell$ by \cref{lem:inducedCoveringEnds}.
        Similarly, the separation $s_2:=(E_t,F_t) \vee (U,W)$ with separator~$S_2$ (see again \cref{fig:conn-cov:UtWt}) distinguishes~$\omega_1$ and~$\omega_2$, so $|S_2| \geq \ell$ by \cref{lem:inducedCoveringEnds}.
        
        Consider $(X_t,Y_t) := (U,W) \vee (C_t,D_t)$.
        The separator of~$(X_t,Y_t)$ is the disjoint union of $S_1$ and~$S_2$, see \cref{fig:conn-cov:UtWt}. Note that the union is indeed disjoint since $S_1 \subseteq E_{t-1}$ and $S_2 \subseteq F_t$, and $E_{t-1} \cap F_t = \emptyset$. 
        Hence, $|X_t \cap Y_t| = |S_1| + |S_2| \geq 2\ell$.
        By submodularity, we get
        \[
            |U_t \cap W_t| = |U \cap W| + |C_t \cap D_t| - |X_t \cap Y_t| \leq (2\ell-1) + 2\ell - 2\ell < 2\ell
        \]
        which completes the proof.
    \end{cproof}

    \begin{claim}\label{claim:conn-cov-2}
        For every $t\in V(T)$ we either have $U_t = \emptyset$ or $U \cap W \cap V_t \neq \emptyset$.
    \end{claim}
    \begin{cproof} 
        We assume that $U \cap W \cap V_t = \emptyset$ and show that $U_t = \emptyset$.
        By construction of~$\hat{G}$, we either have $p(V_t) \subseteq A$ or $p(V_t) \subseteq B$, say the former.
        Since~$(A,B)$ is tight, there is a component~$K$ of~$G[A \setminus B]$ with $N_G(K) = A \cap B$.
        Let~$K'$ be the induced subgraph of $\hat{G}[V_t]$ such that $p(K') = K$.
        Similarly, let $P',Q' \subseteq V_t$ be such that $p(P')=P$ and $p(Q')=Q$.
        See \cref{fig:CoveringUWempty} for a depiction of the situation.
        Since $U_t \cap W_t = (U \cap V_t) \cap (W \cup D_t)$ and $U \cap W \cap V_t = \emptyset$, it follows that $U_t \cap W_t = U \cap V_t \cap D_t \subseteq P' \cup Q'$.
        Since $N_{G}(K) = A \cap B$, we have $N_{\hat{G}}(K') = P' \cup Q'$.
        Since $U \cap W \cap V_t = \emptyset$, we either have $V(K') \subseteq U \setminus W$ or $V(K') \subseteq W \setminus U$.
        
        \casen{Case $V(K') \subseteq U \setminus W$.}
        Since $N_{\hat{G}}(K')=P' \cup Q'$, we have $P' \cup Q' \subseteq U$ so $P' \cup Q' \subseteq U_t \cap W_t$.
        With $|P'|=|P| = \ell = |Q| = |Q'|$, we get $|U_t \cap W_t| \geq 2\ell$, contradicting \cref{claim:conn-cov-1}.

        \casen{Case $V(K') \subseteq W \setminus U$.}
        Since $N_{\hat{G}}(K') = P' \cup Q'$, we have $P' \cup Q' \subseteq W$.
        Since we also have $P' \cup Q' \subseteq V_t$ and since $U \cap W \cap V_t = \emptyset$, it follows that $P' \cup Q' \subseteq W \setminus U$.
        Hence, $U_t \cap W_t = \emptyset$.
        Since~$\hat{G}[V_t]$ is connected, $U_t \setminus W_t$ is empty.
        Hence, $U_t = \emptyset$.
    \end{cproof}

    \begin{figure}[ht]
        \centering
        \includegraphics[height=6\baselineskip]{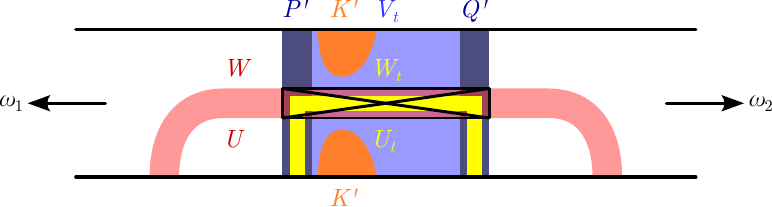}
        \caption{The two possible locations of~$K'$. There are no vertices in the crossed out box as it equals $U\cap W\cap V_t$.}
        \label{fig:CoveringUWempty}
    \end{figure}

    We define $Y:=\{t \in T \mid U_t \neq \emptyset\}$, and we let $Y' \subseteq Y$ contain every other~$t \in Y$ (along the double-ray~$T$).
    In particular, $|Y| \leq 2|Y'|$, and for every $t \neq t' \in Y'$ we have $V_t \cap V_{t'} = \emptyset$.
    Since~$G$ is $(q,k-1)$-unbreakable, every $(<k)$-separation in~$G$ has a proper side with~$\le q$ vertices.
    Note that for every~$t \in V(T)$, the covering map~$p$ is injective on~$U_t$.
    Moreover, there is a (unique) separation $(U_t'\cup X,W_t'\cup X)$ of~$G$ such that $p(U_t)=U_t'$ and $p(U_t \cap W_t)=U_t' \cap W_t'$.
    Note that $(U_t'\cup X,W_t'\cup X)$ has order~$<2\ell+|X| \leq k$.
    By symmetry, suppose that $p(V_t)\subseteq A$.
    Then $B\setminus A\subseteq W_t'\setminus U_t'$ and
    $|B\setminus A|>q$, so $(q,k-1)$-unbreakability yields
    $|U_t'\setminus W_t'|\le q$.
    Hence $|U_t|=|U_t'|<q+2\ell$.
    Moreover, by \cref{claim:conn-cov-2} there is an injection from~$Y'$ to~$U \cap W$, e.g.\ by mapping each $t \in Y'$ to a vertex in~$U \cap W \cap V_t$.
    In particular, $|Y'| \leq |U \cap W| < 2\ell$. 
    Putting all inequalities together, we get
    \[
        |U \setminus W| \leq \left| \bigcup_{t \in Y} U_t \right|
        \leq (q+2\ell) \cdot |Y| \leq (q+2\ell) \cdot 2|Y'| \leq (q+2\ell) \cdot 2|U \cap W| < 4\ell (q+2\ell)=\beta_{\ref{conn-cov}}(q,\ell)~.
    \]
\end{proof}

We can now prove that $\cross_R(\cdot)$ is bounded.

\begin{lem}\label{mu-bound-cov}
    Assume \cref{set:covering}.
    Let $\chi_{\ref{mu-bound-cov}}(q,\ell) := \chi_{\ref{coveringLemmaTwo}}(\beta_{\ref{conn-cov}}(q,\ell),\ell)$.
    Then $\cross_R(s) \leq \chi_{\ref{mu-bound-cov}}(q,\ell)$ for all $s \in R$.
\end{lem}
\begin{proof}
    By \cref{coveringLemmaTwo}, applied with $\beta := \beta_{\ref{conn-cov}}(q,\ell)$, we have at least one of the following two cases.

    \casen{Case \cref{coveringLemmaTwo}~\ref{itm:coveringLemmaTwo-1}}: $\cross_R(s) \leq \chi_{\ref{coveringLemmaTwo}}(\beta,\ell)$.
    In this case we are done.

    \casen{Case \cref{coveringLemmaTwo}~\ref{coveringLemmaTwoOutcomeUW}}: there is a $(<2\ell)$-separation $(U,W)$ of~$\hat{G}$ with $|U \setminus W| \geq \beta$ and both~$\omega_i$ live in the side~$W$.
    In particular, we have~$|U \setminus W| \geq \beta_{\ref{conn-cov}}(q,\ell)$, contradicting \cref{conn-cov}.
\end{proof}

We now prove that any two crossing separations in $R$ cross with two small opposite corners.

\begin{lemma}
\label{lem:corner-in-cover}
    Assume \cref{set:covering}. 
    Let~$(A_1,A_2), (B_1,B_2) \in R$ be two crossing separations.
    Then, the $A_1B_2$-corner and the $A_2B_1$-corner both have size less than $\beta_{\ref{conn-cov}}(q,\ell)$.
\end{lemma}

\begin{figure}[ht]
    \centering
    \includegraphics[height=6\baselineskip]{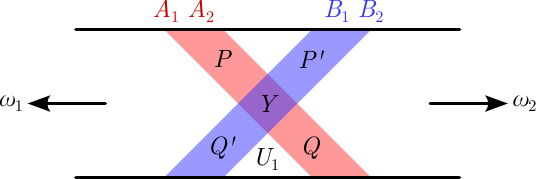}
    \caption{The situation in the proof of \cref{lem:corner-in-cover}.}
    \label{fig:CoveringInfimumSupremum}
\end{figure}

\begin{proof}
    Let~$s:=(A_1,A_2)$ and $t:=(B_1,B_2)$.
    Let~$(U_1,W_1):=(A_1,A_2) \wedge (B_2,B_1)$ and~$(U_2,W_2):=(A_2,A_1)\wedge (B_1,B_2)$.
    This proof is supported by \cref{fig:CoveringInfimumSupremum}.

    \begin{claim}
        For $i \in \{1,2\}$,~$(U_i,W_i)$ is a $(<2\ell)$-separation of~$\hat{G}$ such that both~$\omega_1$ and~$\omega_2$ live in~$W_i$.
    \end{claim}

    \begin{cproof}
        By symmetry, it suffices to prove it for $i=1$.
        First, recall that $\omega_1$ lives in~$B_1$ and~$\omega_2$ lives in~$A_2$.
        Hence, both~$\omega_i$ live in $B_1 \cup A_2 = W_1$.
        
        Now, consider the crossing separations~$(A_1,A_2)$ and~$(B_1,B_2)$.
        Let~$P$ be the $sB_1$-link, let~$Q$ be the $sB_2$-link, let~$P'$ be the $A_2t$-link, let~$Q'$ be the $A_1t$-link and let~$Y$ be the $st$-centre.
        Recall that $|P|+|Q|+|Y|=\ell=|P'|+|Q'|+|Y|$.
        The order of~$(U_1,W_1)$ is 
        \[
            |U_1 \cap W_1| = |Q|+|Q'|+|Y| = 2\ell - |P| - |P'| - |Y|~.
        \]
        Hence, if~$P$ or~$P'$ or~$Y$ is nonempty, then~$(U_1,W_1)$ has order $<2\ell$.
        Otherwise, $P=P'=Y=\emptyset$.
        Since~$\hat G$ is connected, the $A_2 B_1$-corner is empty.
        But then $(A_1,A_2)$ and $(B_1,B_2)$ are nested, a contradiction.
    \end{cproof}
    
    Thus, for every $i \in \{1, 2\}$, we have $|U_i \setminus W_i| < \beta_{\ref{conn-cov}}(q,\ell)$ by \cref{conn-cov}, which concludes the proof.
\end{proof}

We are now ready to prove the main result of this section, that every separation in $R$ is at bounded edit-distance from a separation in $L$.

\begin{lem}\label{coveringLemmaFive}
    Assume \cref{set:covering}.
    Let $\beta:=\beta_{\ref{conn-cov}}(q,\ell)$ and
    $\eps_{\ref{coveringLemmaFive}}(q,\ell):=(\beta+\ell)\cdot \chi_{\ref{mu-bound-cov}}(q,\ell)$.
    Then for every $r\in R$ there is $s(r)\in L$ such that $r$ and $s(r)$
    have edit-distance $\le (\beta+\ell)\cdot\cross_R(r)\le \eps_{\ref{coveringLemmaFive}}(q,\ell)$.
\end{lem}
\begin{proof}
    We proceed by induction on $\cross_R(r)$.

    If $\cross_R(r)=0$, then $r\in L$.
    And if $r\in L$, then we are done (with $s(r):=r$).
    Hence we may assume now that $\cross_R(r)>0$ and that some $t\in R$ witnesses $r\notin L$, i.e.\ that $t\in R$ crosses $r$ with $\cross_R(t)\le\cross_R(r)$.

    By \cref{coveringCornersAreEfficient}, there are two opposite corner-separations $c,d$ of $r,t$ that are both contained in~$R$.
    Moreover, the lemma yields $\cross_R(c)<\cross_R(r)$, say.
    By the induction hypothesis, there is $s\in L$ such that $s$ and $c$ have edit-distance $\le (\beta+\ell)\cdot\cross_R(c)$.

    \begin{claim}\label{coveringLemmaFiveClose}
        The separations $r$ and $c$ have edit-distance $<\beta+\ell$.
    \end{claim}
    \begin{cproof}
        Write $r=:(A,B)$, $t=:(C,D)$ and $c=:(A\cap C,B\cup D)$, say.
        Consider how $r$ and $t$ cross.
        By \cref{lem:corner-in-cover}, the $AD$-corner has size $<\beta_{\ref{conn-cov}}(q,\ell) = \beta$.
        The vertex-set $A \Delta (A \cap C) = A \setminus C$ is equal to the union of the $AD$-corner and the $rD$-link.
        Thus, $|A \Delta (A \cap C)| < \beta + \ell$.
        Similarly, $B \Delta (B \cup D) = D \setminus B$ is equal to the union of the $AD$-corner and the $At$-link.
        Thus, $|B \Delta (B \cup D)| < \beta + \ell$.
        Hence, $r$ and~$c$ have edit-distance $<\beta + \ell$.
    \end{cproof}
    
    The edit-distance between $r$ and $s$ is at most the edit-distance between $r$ and $c$ plus the edit-distance between $c$ and $s$, so by \cref{coveringLemmaFiveClose} and the induction hypothesis this is strictly less than
    \[
        (\beta+\ell)+(\beta+\ell)\cdot \cross_R(c)\le (\beta+\ell)\cdot\cross_R(r)\le (\beta+\ell)\cdot \chi_{\ref{mu-bound-cov}}(q,\ell)= \eps_{\ref{coveringLemmaFive}}(q,\ell)~,
    \]
    where we used \cref{mu-bound-cov} implicitly.
\end{proof}

\subsection{Constructing the graphic-daisy}\label{subsec:daisy}
In this section, we show that the separations in~$L$ induce a path-decomposition of~$\hat G$, which projects to a cycle-decomposition of $G-X$.
Then, we argue that adding~$X$ to all bags of this cycle-decomposition yields a graphic-$k$-daisy $\cD$ of~$G$, such that every separation $s \in S$ is at bounded edit-distance from a separation displayed by~$\cD$.

Recall that a \defn{deck transformation} of the covering $p\colon \hat G \to G$ from \cref{set:covering} is an automorphism~$\varphi$ of~$\hat G$ that commutes with~$p$ in that $p=p\circ\varphi$.
The deck transformations form a subgroup of the automorphism group of~$\hat G$, called the \defn{deck transformation group}.

\begin{construction}
\label{TDCfromNestedSetConstruct}
    Let $N$ be a nested symmetric set of proper separations of a graph~$G$.
    We define a graph $T$ by taking the splitting stars of~$N$ as the nodes of~$T$, and joining two nodes $\sigma_1\neq\sigma_2$ by an edge in~$T$ if there is $(A_1,A_2)\in N$ such that $(A_{3-i},A_i)\in\sigma_i$ for both~$i=1,2$.
    For each node $\sigma$ we put $V_\sigma:=\bigcap\,\{\,B:(A,B)\in \sigma\,\}$.
    Write $\defnMath{\cT(N)}:=(T,\cV)$.
\end{construction}

The following lemma was known to \cite{infinitesplinter}, who provide all the proof ingredients. In the proof, we rely on \cite{canonicalGraphDec} for convenience.

\begin{lem}[{\cite{infinitesplinter}}]
\label{TDCfromNesteSetVerify}
    Let $G$ be a locally finite connected graph and $k\in\N$.
    Let $N$ be a nested symmetric set of proper $(\le k)$-separations of~$G$.
    If $G$ is infinite, we additionally assume that all separations in $N$ are tight.
    \begin{enumerate}
        \item $\cT(N)$ is a tree-decomposition of~$G$ of adhesion at most $k$.
        \item There is a natural order-preserving bijection between the (directed) edges of~$T$ and $(N,\le)$.
        \item The torso of a splitting star $\sigma$ of~$N$ equals the torso of~$V_\sigma$ when viewed as a bag of the tree-decomposition~$\cT(N)$.
    \end{enumerate}
\end{lem}
\begin{proof}
    We combine the following results from \cite{canonicalGraphDec}:
    the discussion at the beginning of §7.1, Lemma 7.2 and Lemma~7.5.
\end{proof}

\begin{lem}\label{covering:LisomZ}
    Assume \cref{set:covering}.
    The set~$L$ with the usual ordering on separations is order-isomorphic to~$\Z$.
\end{lem}
\begin{proof}
    The definition of locally-$R$-optimal immediately yields that $L$ is nested.
    Since each $(A,B)\in L$ has $\omega_1$ living in $A$ and $\omega_2$ living in~$B$,
    we find that $(L,\le )$ is a linear order.
    Note that all separations in $L \cup L^*$ are tight since they efficiently distinguish~$\omega_1$ and~$\omega_2$.
    By \cref{TDCfromNesteSetVerify}, $L\cup L^*$ uniquely determines a tree-decomposition $(H,\cW)$ of~$\hat G$.
    Moreover, $(L\cup L^*,\le)$ is order-isomorphic to the set of directed edges of~$H$.
    As $(L,\le)$ is a linear order, $H$ must be a (possibly infinite) path.
    
    We claim that $H$ is a double-ray, and assume for a contradiction that it is not.
    Then $(L,\le )$ has a maximal element $(U,W)$, say.
    Let $v\in U\cap W$ be arbitrary.
    Consider the standard decomposition $(T,\cV)$ of~$\hat G$.
    Since the adhesion-sets of $(T,\cV)$ are pairwise disjoint, and since $\omega_1$ and $\omega_2$ correspond to the two ends of the double-ray $T$,
    a distance argument shows that there is an edge $t_1 t_2$ of $T$ whose induced separation $(A_1,A_2)\in R$ satisfies $U\cap W\se A_1\sm A_2$.
    Let $v'\in A_2\sm A_1$ such that $p(v')=p(v)$.
    Note that $v'\in W\sm U$ since $U \cap W$ does not separate $v'$ from~$\omega_2$.
    Let $\varphi$ be a deck transformation of~$p$ such that $\varphi(v)=v'$.
    Since $\varphi$ is a deck transformation, it fixes $\omega_i$ for both~$i$ by the construction of~$p$.
    Hence $(U',W'):=(\varphi(U),\varphi(W))\in L$.
    Since $(L,\le )$ is a linear order, and since $v'\in U'\cap W'$ satisfies $v'\in W\sm U$, we must have $(U,W)<(U',W')$.
    This contradicts the maximality of $(U,W)$ in~$(L,\le )$.
\end{proof}

\begin{lem}\label{coveringCycleLong}
    Assume \cref{set:covering} with 
    \[
        \gamma > \gamma_{\ref{coveringCycleLong}}(q,k) := \max\,\{\,2\eps_{\ref{coveringLemmaFive}}(q,k')+2k'+k+\beta_{\ref{conn-cov}}(q,k'):0\le k'\le k/2\,\}~.
    \]
    The set~$L$ with the usual ordering on separations is order-isomorphic to~$\Z$ by \cref{covering:LisomZ}, which yields an enumeration $L=\{t_i:i\in\Z\}$.
    Then all of the following hold.
    \begin{enumerate}
        \item\label{itm:coveringCycleLong-1} For every nontrivial deck transformation $\varphi$ of~$p$, we have that $\varphi(t_i)=t_j$ implies $|i-j|\ge 4$.
        \item\label{itm:coveringCycleLong-2} There are $s_0,s_1,s_2,s_3\in L$ such that the following holds.
        Let $\Gamma\cdot s_j=\{\cdots<s_j^{-1}<s_j^0<s_j^1<\cdots\}$ with~$s_j^0=s_j$, where~$\Gamma$ is the deck transformation group of~$p$.
        Then we have $s_0^i < s_1^i < s_2^i < s_3^i < s_0^{i+1}$ for every~$i \in\Zbb$, and the separators of the~$s_j^i$ are pairwise disjoint.
    \end{enumerate}
\end{lem}
\begin{proof}
    By the assumptions of \cref{set:covering},
    the covering $p$ arises from a separation $x\in S$ that is $\gamma$-crossed by another separation~$y\in S$.
    The four links of $x,y$ determine a two-way infinite sequence of separations $r_k=(X_k,Y_k)\in R$ such that
    \begin{enumerate}
        \item\label{FourSequenceIncreases} $\cdots < r_{-1}<r_0<r_1<\cdots$ and 
        \item\label{FourSequenceFourOrbits} for every nontrivial deck transformation $\varphi$ of~$p$, we have that $\varphi(r_k)=r_\ell$ implies $|k-\ell|\ge 4$.
    \end{enumerate}
    By \cref{coveringLemmaFive}, for every $r_k$ there is an $s_k=(U_k,W_k)\in L$ such that $r_k$ and $s_k$ have edit-distance~$\le\eps_{\ref{coveringLemmaFive}}(q,\ell)$.
    We choose the $s_k$ so that, 
    for every deck transformation $\varphi$ of~$p$, we have that $\varphi(r_i)=r_j$ implies $\varphi(s_i)=s_j$.

    \begin{claim}\label{coveringCycleLongIncreasingSequence}
        We have
        \[
            \cdots < s_{-1} < s_0 < s_1 < \cdots 
        \]
        That is to say that $r_k\mapsto s_k$ is an order-preserving injection $\{r_k:k\in\Z\}\to L$.
        Moreover, the separators of two consecutive~$s_k$ are disjoint.
    \end{claim}

    \begin{cproof}
        For the first part, assume for a contradiction that there are $i < j$ with $s_i\not < s_j$.
        Since $(L, \leq)$ is a linear order by \cref{covering:LisomZ}, we have $s_j\le s_i$.
        In particular, $(U_j,W_j)\le (U_i,W_i)$.
        The following calculation is supported by \cref{fig:CoveringCycleLong}.
        We obtain
        \[
            Y_i\cap X_j\se (U_i\sm X_i)\cup (W_j\sm Y_j)\cup (X_i\cap Y_i)\cup (X_j\cap Y_j).
        \]
        Recall $|X_i \cap Y_i|=\ell=|X_j \cap Y_j|$.
        Moreover,~$r_k$ and $s_k$ having edit-distance $\le\eps_{\ref{coveringLemmaFive}}(q,\ell)$ implies $|U_i \setminus X_i| \leq \eps_{\ref{coveringLemmaFive}}(q,\ell)$ and $|W_j \setminus Y_j| \leq \eps_{\ref{coveringLemmaFive}}(q,\ell)$.
        Hence $|Y_i\cap X_j|\le 2\eps_{\ref{coveringLemmaFive}}(q,\ell)+2\ell$.
        But then one of the four corners of $x,y$ has size $\le 2\eps_{\ref{coveringLemmaFive}}(q,\ell)+2\ell\le\gamma_{\ref{coveringCycleLong}}(q,k)-k<\gamma-k$, contradicting that $x,y$ $(\gamma-k)$-criss by~\cref{crissCross}.

        For the `moreover'-part, assume for a contradiction that the separators of two consecutive~$s_k$ intersect, say~$s_0$ and~$s_1$.
        Recall that $s_0 < s_1$.
        Consider the separation $s':=s_0\vee s_1^*$, which is a separation of order $<|\sep(s_0)|+|\sep(s_1)|=2\ell$.
        Observe that both ends~$\omega_1$ and~$\omega_2$ live in the left side of~$s'$.
        By \cref{conn-cov} applied to~$s'$, we have $|(W_0 \sm U_0) \cap (U_1 \sm W_1)| < \beta_{\ref{conn-cov}}(q,\ell)$.
        Therefore $|W_0 \cap U_1| < \beta_{\ref{conn-cov}}(q,\ell) + 2\ell$.
        Since~$s_i$ and~$r_i$ have edit-distance $\leq \eps_{\ref{coveringLemmaFive}}(q,\ell)$, they are $\eps_{\ref{coveringLemmaFive}}(q,\ell)$-close by \cref{closeVSapprox}. 
        So, we have $|Y_0 \cap X_1| < 2\eps_{\ref{coveringLemmaFive}}(q,\ell) + \beta_{\ref{conn-cov}}(q,\ell) + 2\ell < \gamma-k$.
        But then one of the four corners of~$x,y$ has size~$< \gamma-k$, contradicting that $x,y$ $(\gamma-k)$-criss by~\cref{crissCross}.
    \end{cproof}

    \begin{figure}[ht]
        \centering
        \includegraphics[height=6\baselineskip]{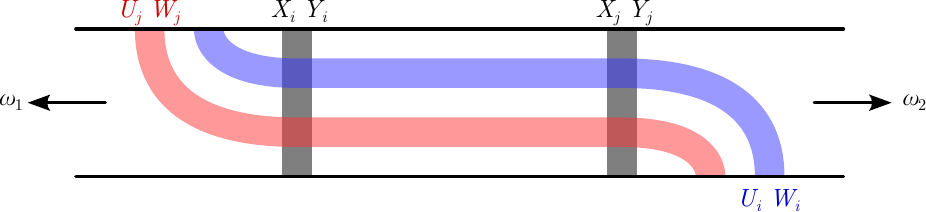}
        \caption{The separations $(X_i,Y_i)\le (X_j,Y_j)$ and $(U_j,W_j)\le (U_i,W_i)$}
        \label{fig:CoveringCycleLong}
    \end{figure}

    \begin{claim}\label{coveringCycleLongDistanceFour}
        For every nontrivial deck transformation $\varphi$ of~$p$, we have that $\varphi(s_i)=s_j$ implies $|i-j|\ge 4$.
    \end{claim}
    \begin{cproof}
        Assume that $\varphi(s_i)=s_j$.
        Since $\varphi$ is a deck transformation, $\varphi$ fixes the two ends $\omega_1$ and $\omega_2$ of~$\hat G$.
        In particular, $\varphi$ acts on $\hat G$ by `shifting' but without `reflection' with respect to the two ends~$\omega_1$ and~$\omega_2$.
        Hence there is a unique index $k\in\Z$ such that $r_k=\varphi(r_i)$.
        Then $s_k=\varphi(s_i)$ by the choice of~$s_k$.
        Combining this with the assumption that $\varphi(s_i)=s_j$ yields $s_k=s_j$.
        Then $k=j$ by \cref{coveringCycleLongIncreasingSequence}.
        Hence $\varphi(r_i)=r_k=r_j$.
        Therefore, $|i-j|\ge 4$ by~\cref{FourSequenceFourOrbits}.
    \end{cproof}

    Finally, assume that some nontrivial deck transformation $\varphi$ of~$p$ satisfies $\varphi(t_i)=t_j$ for some $i,j\in\Z$.
    By \cref{coveringCycleLongIncreasingSequence} there are unique indices $I,J\in\Z$ such that
    \[
        s_I\le t_i < s_{I+1} \quad \text{and} \quad s_J\le t_j < s_{J+1}.
    \]
    Since $L$ is canonical, $\varphi(t_i)=t_j$ implies $\varphi(s_I)=s_J$.
    Then $|I-J|\ge 4$ by \cref{coveringCycleLongDistanceFour}.
    Therefore, $|i-j|\ge 4$ as desired.
    This gives~\ref{itm:coveringCycleLong-1}.
    \ref{itm:coveringCycleLong-2} follows from~\cref{coveringCycleLongIncreasingSequence} .
\end{proof}

\begin{corollary}\label{coveringIncreasingSequence}
    Assume \cref{set:covering} with $\gamma>\gamma_{\ref{coveringCycleLong}}(q,k)$.
    Let $r_1<r_2<r_3\in R$ be induced by edges of the standard tree-decomposition of~$\hat G$.
    Assume that there is a deck transformation $\varphi$ of~$p$ such that $\varphi(r_1)=r_3$.
    Then there are $s_1,s_2,s_3\in L$ such that:
    \begin{enumerate}
        \item $r_i$ and $s_i$ have edit-distance~$\le\eps_{\ref{coveringLemmaFive}}(q,\ell)$ for every $i\in [3]$,
        \item $s_1<s_2<s_3$, and
        \item $\varphi(s_1)=s_3$.
    \end{enumerate}
\end{corollary}
\begin{proof}
    We use \cref{coveringLemmaFive} to find $s_1$ and~$s_2$ such that $r_i$ and $s_i$ have edit-distance~$\le\eps_{\ref{coveringLemmaFive}}$ for $i=1,2$.
    We let $s_3:=\varphi(s_1)$.
    Since $\varphi$ is an automorphism, we find that $r_3$ and $s_3$ also have edit-distance~$\le\eps_{\ref{coveringLemmaFive}}$.
    It remains to show that $s_1 < s_2 < s_3$.
    For this, we closely follow the proof of \cref{coveringCycleLongIncreasingSequence}.
\end{proof}

\begin{lem}\label{coveringInjective}
    Assume \cref{set:covering} with $\gamma > \gamma_{\ref{coveringCycleLong}}(q,k)$.
    Let $(T,\cV):=\cT(L\cup L^*)$ be the tree-decomposition of~$\hat G$ that $L\cup L^*$ defines by \cref{TDCfromNesteSetVerify}.
    Note that $T$ is a double-ray by \cref{covering:LisomZ}.
    Let $t\in V(T)$ with neighbouring edges $e,f$.
    Assume that $p(V_e)$ and $p(V_f)$ are disjoint.
    Then $p\rest V_t$ is injective.
\end{lem}
\begin{proof}
    Assume for a contradiction that there are two distinct vertices $u,v\in V_t$ such that $p(u)=p(v)=:w$.
    Let $\Gamma$ denote the group of deck transformations of~$p$.
    Let $\varphi\in\Gamma$ such that $\varphi(u)=v$.
    Since $(T,\cV)$ is canonical, there is $t'\in V(T)$ such that $\varphi(V_t)=V_{t'}$.
    By \cref{coveringCycleLong}, we have $t'\neq t$.
    Since $\Gamma$ acts transitively on the fibres of~$p$, and since $(T,\cV)$ is canonical, we have that as soon as one of $u,v$ lies in some adhesion-set of~$(T,\cV)$, then so does the other.
    Hence, it suffices to distinguish the following cases:

    \casen{Case~1: $u\in V_e$ and $v\in V_f$ (or vice versa).}
    This case is not possible, as $p(V_e)$ and $p(V_f)$ are disjoint by assumption.

    \casen{Case~2: $u,v\in V_f$, say.}
    Let $e',f'$ be the two edges of $T$ that are incident to~$t'$.
    Since $(T,\cV)$ is canonical, the deck transformation $\varphi$ acts on~$T$, and we name $e',f'$ so that $\varphi$ moves $e$ to $e'$ and $f$ to $f'$.
    By the construction of~$p$, all deck transformations fix $\omega_i$ for both~$i$, so we have $e<f\le e'<f'$ in the order of the double-ray~$T$, say.
    Now, on the one hand, we have $v\in V_f$ by assumption, and on the other hand, $u\in V_f$ implies $\varphi(u)\in V_{f'}$ and, in particular, $v\in V_{f'}$.
    But then $v\in V_{e'}$ as $f\le e'<f'$.
    Hence $u\in V_e$, in addition to $u\in V_f$.
    This contradicts our assumption that $p(V_e)$ and $p(V_f)$ are disjoint.

    \casen{Case~3: $u,v\in V_t\sm (V_e\cup V_f)$.}
    This case is impossible.
    Indeed, since $v=\varphi(u)$ we immediately find that $t'=t$ is the only option, contradicting that we have $t'\neq t$.
\end{proof}

\begin{lem}\label{separationsRandomInequality}
    Let $A,B,C,D$ be four sets.
    Then
    \[
        (A\cap B)\triangle (C\cap D)\se (A\triangle C)\cup (B\triangle D).
    \]
\end{lem}
\begin{proof}
    Let $v\in (A\cap B)\triangle (C\cap D)$.
    By symmetry, we may assume without loss of generality that $v\in (A\cap B)\sm (C\cap D)$.
    Then $v\notin C$ or $v\notin D$.
    If $v\notin C$, then $v\in A\triangle C$ and we are done.
    If $v\notin D$, then $v\in B\triangle D$ and we are done.
\end{proof}

\begin{lem}\label{coveringDaisy}
    Assume \cref{set:covering} with $\gamma> \gamma_{\ref{coveringCycleLong}}(q,k)$.
    Let $\eps_{\ref{coveringDaisy}}(q,\ell):=2\eps_{\ref{coveringLemmaFive}}(q,\ell)$.
    Then~$G$ has a $q$-lean graphic-$k$-daisy~$\cD=(O,\cP')$ with the following properties.
    \begin{enumerate}
        \item\label{coveringDaisy:adhesionSets} $\mathcal{D}$ has centre $X$ and has the clover property.
        \item\label{coveringDaisy:ItemDisplay} Every separation in~$S$ has edit-distance $\le \eps_{\ref{coveringDaisy}}$ to some displayed separation of~$\cD$.        \item\label{coveringDaisy:ItemCompatible} For every $\varepsilon \geq \eps_{\ref{coveringDaisy}}$, every $k$-separation of~$G$ that is $\eps$-compatible with~$\cD$ is $\eps$-nested with all the petal-separations of~$\cD$.
    \end{enumerate}
\end{lem}

\begin{proof}
    The set $L$ uniquely determines a tree-decomposition $(T,\cV)$ of~$\hat G$ such that $L \cup L^*$ equals the set of separations induced by the oriented edges of $T$, see \cref{TDCfromNesteSetVerify}.
    The tree~$T$ is a double-ray whose ends correspond to the ends~$\omega_1$ and~$\omega_2$ of~$\hat{G}$, respectively, by \cref{covering:LisomZ}.
    We orient the edges of $T$ from $\omega_1$ to~$\omega_2$.
    As $(T,\cV)$ is canonical, it defines via the covering map $p$ a graph-decomposition $\cO=(O,\cP)$ of~$G-X$ with family $\cP=(P_o:o\in V(O))$ of bags as in \cite[Construction~3.8 and Lemma~3.9]{canonicalGraphDec}.

    \begin{claim}\label{coveringDaisy:PreDaisy}
        $\cO$ is a graphic-$2\ell$-daisy of~$G-X$, and~$\cO$ has four pairwise disjoint adhesion-sets.
    \end{claim}
    \begin{cproof}
        First, we show that~$\cO$ is a $2\ell$-flower of~$G-X$.
        Since $T$ is a double-ray, $O$ is a cycle of length~$\ge 4$ by \cref{coveringCycleLong}.
        Let $(C,D)$ be a 2-edge-cut of~$O$ with cross-edges $c_i d_i$ so that $c_i\in C$ and $d_i\in D$ for $i=1,2$.
        Let $(C',D')$ be the displayed separation of~$\cO$ that is induced by~$(C,D)$.
        Then 
        \[
            C'\cap D'=\bigcup_{i=1,2}(P_{c_i}\cap P_{d_i})
        \]
        by \cite[Lemma~3.4]{canonicalGraphDec}.
        Each intersection $P_{c_i}\cap P_{d_i}$ equals the projection (with regard to~$p$) of some adhesion-set of $(T,\cV)$.
        Since every adhesion-set of $(T,\cV)$ has size~$\ell$, the separator $C'\cap D'$ has size~$\le 2\ell$.
        Hence,~$\cO$ is a $2\ell$-flower.

        Now we show that~$\cO$ is a graphic-$2\ell$-daisy.
        By \cref{thm:flower-dichotomy} and since every $2\ell$-daisy is a graphic-$2\ell$-daisy, it suffices to show that~$\cO$ is not a $2\ell$-anemone.
        From \cref{coveringCycleLong}~\ref{itm:coveringCycleLong-2} we get that $\cO$ has four pairwise disjoint adhesion-sets~$P_{e_0},\ldots,P_{e_3}$.
        Since~$\cO$ has four pairwise disjoint adhesion-sets, the union of the adhesion-sets of~$\cO$ has size~$\geq 4\ell > 2\ell$.
        So~$\cO$ is not a $2\ell$-anemone, and hence~$\cO$ is a graphic-$2\ell$-daisy, which completes the proof.
    \end{cproof}
    
    For each $P_o\in\cP$ let $P'_o:=P_o\cup X$, and let $\cP':=(P'_o:o\in V(O))$.
    Put $\cD:=(O,\cP')$.
    Then~$\cD$ is a graphic-$k$-daisy of~$G$ as $\cO$ is a graphic-$2\ell$-daisy of~$G-X$ by \cref{coveringDaisy:PreDaisy}.

    \begin{claim}\label{coveringDaisy:qlean}
        $\cD$ is $q$-lean.
    \end{claim}
    \begin{cproof}
        For an edge~$e$ of~$O$, we denote by~$P_e$ the adhesion-set of~$e$ in~$\cO$, and by~$P'_e$ the adhesion-set of~$e$ in~$\cD$.
        Let~$Z:=\bigcap_{o \in V(O)}P'_o$.
        Let~$s=(A,B)$ be a $(>q)$-broad separation of~$G$ displayed by a set of edges~$\{e,f\}$ of~$O$ in~$\cD$.
        We need to show that the sets $P'_e \sm Z$ and $P'_f \sm Z$ are disjoint.

        Since~$\cO$ has four pairwise disjoint adhesion-sets by \cref{coveringDaisy:PreDaisy}, we have $\bigcap_{o \in V(O)}P_o = \emptyset$ and therefore $Z=\bigcap_{o \in V(O)}P'_o = X$.
        Let~$k'$ be the order of~$s$, and note that $k' \geq k$ since~$G$ is $(q,k-1)$-unbreakable.
        Since $X \subseteq \sep(s)$, we have that $s':=(A \sm X, B \sm X)$ is a $(k'-|X|)$-separation of~$G-X$ that is displayed by~$\{e,f\}$ in~$\cO$.
        In particular, the order of~$s'$ is $k'-|X| \geq k-|X| = 2\ell$.
        Since every adhesion-set of~$\cO$ has size exactly~$\ell$, the adhesion-sets for~$e$ and~$f$ in~$\cO$ are disjoint.
        Hence, $P_e=P'_e \sm Z$ and $P_f=P'_f \sm Z$ are disjoint.
    \end{cproof}

    \begin{claim}\label{coveringDaisy:Display}
        Every $(A,B)\in S$ has edit-distance $\le \eps_{\ref{coveringDaisy}}(q,\ell)$ to some displayed separation of~$\cD$.
    \end{claim}
    \begin{cproof}
    Let $(A,B)\in S$ be arbitrary.
    By \cref{lem:covering-equivalent}, we may assume without loss of generality that~$p$ is induced by~$(A,B)$ and a unique tripartition~$\cQ$ of~$A \cap B$.
    Let $(H,\cW)$ denote the standard tree-decomposition of~$\hat G$ (induced by $(A,B)$ and~$\cQ$).
    Recall that~$H$ is a double-ray whose two ends reflect the two ends~$\omega_1$ and~$\omega_2$ of~$\hat G$.
    We consider the edges of~$H$ as oriented from~$\omega_1$ to~$\omega_2$.
    Now let $e,e',e''$ be three consecutive edges on the oriented double-ray~$H$.
    Let~$h$ be the node of~$H$ shared by~$e$ and~$e'$.
    Let~$h'$ be the node of~$H$ shared by~$e'$ and~$e''$.

    Let $(X_e,Y_e)\in R$ denote the separation induced by~$e$ and $(H,\cW)$, and define $(X_{e'},Y_{e'})$ and $(X_{e''},Y_{e''})$ similarly.
    Then $e<e'<e''$ implies
    \[
        (X_e,Y_e)<(X_{e'},Y_{e'})<(X_{e''},Y_{e''}).
    \]
    Moreover, there is a deck transformation~$\varphi$ of~$p$ such that
    \[
        (\varphi(X_e),\varphi(Y_e))=(X_{e''},Y_{e''})
    \]
    by the construction of~$p$.
    So by \cref{coveringIncreasingSequence}, there are three separations
    \[
        (U_f,W_f)<(U_{f'},W_{f'})<(U_{f''},W_{f''})
    \]
    in $L$ induced by oriented edges $f<f'<f''\in E(T)$ such that $(X_e,Y_e)$ and $(U_f,W_f)$ have edit-distance $\le \eps_{\ref{coveringLemmaFive}}(q,\ell)$ and so on, and such that
    \[
        (\varphi(U_f),\varphi(W_f))=(U_{f''},W_{f''}).
    \]
    
    Without loss of generality, we have $A=p(W_h)\sqcup X$ and $B=p(W_{h'})\sqcup X$.
    Let $\hat C:=\mathring{f}T\mathring{f}' \subseteq V(T)$ and $\hat D:=\mathring{f}' T \mathring{f}'' \subseteq V(T)$.
    Let~$\Gamma$ denote the deck transformation group of~$p$, which acts on $T$ as $(T,\cV)$ is canonical.
    By construction, $V(O)$ consists of the orbits of the $\Gamma$-action on~$V(T)$.
    Let
    \[
        C:=\{\,\Gamma\cdot t \ :t \ \in\hat C\,\} \subseteq V(O)\qquad\text{and}\qquad D:=\{\,\Gamma\cdot t \ : t \ \in\hat D\,\} \subseteq V(O)~.
    \]
    Then $(C,D)$ is a 2-edge-cut of~$O$ with the two cross-edges $\Gamma\cdot f=\Gamma\cdot f''$ and $\Gamma\cdot f'$.
    Let $(C',D')$ denote the displayed $(\le k)$-separation of $\cD$ that corresponds to~$(C,D)$.
    We have
    \begin{equation}\label{eq:Bside:coveringDaisy}
        B\sm X=p(W_{h'})=p(Y_{e'}\cap X_{e''})
    \end{equation}
    and
    \begin{equation}\label{eq:Ddash:coveringDaisy}
        D'\sm X=\bigcup_{o\in D}P_o=p\Big(\bigcup_{d\in\hat D}V_d\Big)=p(W_{f'}\cap U_{f''}).
    \end{equation}
    Combining $X\se B\cap D'$ with \cref{eq:Bside:coveringDaisy,eq:Ddash:coveringDaisy}, we find that
    \begin{align*}
        |B\triangle D'|&{}=|(B\sm X)\triangle (D'\sm X)|\\
        &{}=|p(Y_{e'}\cap X_{e''})\triangle p(W_{f'}\cap U_{f''})|\\
        &{}\le |(Y_{e'}\cap X_{e''})\triangle (W_{f'}\cap U_{f''})|\\
        &{}\overset{(\ast)}{\le} |Y_{e'}\triangle W_{f'}|+|X_{e''}\triangle U_{f''}|\\
        &{}\le 2\eps_{\ref{coveringLemmaFive}}(q,\ell) = \eps_{\ref{coveringDaisy}}(q,\ell)
    \end{align*}
    where we use \cref{separationsRandomInequality} for~$(\ast)$.
    A similar argument shows that $|A\triangle C'|\le \eps_{\ref{coveringDaisy}}(q,\ell)$.
    Hence $(A,B)$ and $(C',D')$ have edit-distance $\le \eps_{\ref{coveringDaisy}}(q,\ell)$.
    \end{cproof}

    \begin{claim}\label{coveringDaisy:compatible}
        For every $\varepsilon \geq \eps_{\ref{coveringLemmaFive}}$, every $k$-separation $(A,B)$ of~$G$ that is $\eps$-compatible with~$\cD$ is $\eps$-nested with all the petal-separations of~$\cD$.
    \end{claim}
    \begin{cproof}
        Let $\eps \geq \eps_{\ref{coveringLemmaFive}}$ and let $r:=(A, B)$ be a $k$-separation of $G$ that is $\eps$-compatible with~$\cD$.
        Assume for a contradiction that $r$ $\eps$-crosses a petal-separation $s:=(U,W)$ of~$\cD$ where $U=P'_o$ for some node $o\in O$.
        Let $o^-$ and $o^+$ be the two neighbours of $o$ on~$O$.
        Abbreviate $e:=oo^-$ and $f:=oo^+$.
        Recall that $P'_e$ and $P'_f$ denote the adhesion-sets of $\cD$ at $e$ and~$f$, respectively.
        Then, by the definition of $\eps$-compatible with~$\cD$, we have all of the following:
        \begin{enumerate}
            \item $X = P'_e \cap P'_f$ and $X$ is the centre of~$(A,B)$ and~$(U,W)$;
            \item $P'_e \sm X$ and $P'_f \sm X$ are the $As$-link and the $Bs$-link, respectively (after possibly swapping the names of $A$ and~$B$); and
            \item all links of the two separations~$r,s$ are of the same size.
            Since $X$ is the centre of $(A,B)$ and $(U,W)$, all links have size~$\ell$.
        \end{enumerate}
        Let $A':=A\cap P_o$ and $B':=B\cap P_o$.
        Then $(A',B')$ is an $\ell$-separation of~$G[P_o]$.
        Since $X=P'_e\cap P'_f$, we have that $P_e=P'_e\sm X$ is the $As$-link and $P_f=P'_f\sm X$ is the $Bs$-link.
        Therefore, we have $P_e\se A'\sm B'$ and $P_f\se B'\sm A'$ on the one hand, and we have that $P_e$ and $P_f$ are disjoint on the other hand.
        
        Let $t\in V(T)$ such that $o=\Gamma\cdot t$, which exists by the construction of $\cO$ (which gave rise to $\cD$ by taking the union of every bag with~$X$) from~$(T,\cV)$ by folding via~$p$.
        Let $t^-,t^+$ be the two neighbours of~$t$ in~$T$ with $o^\pm=\Gamma\cdot t^\pm$.
        Then let $\hat e:=\{t^-,t\}$ and $\hat f:=\{t,t^+\}$.
        We assume without loss of generality that the edge $\hat e$ is oriented from $t^-$ to~$t$ in the fixed orientation of the edges of the double ray~$T$ (from $\omega_1$ to $\omega_2$), and similarly that $\hat f$ is oriented from $t$ to~$t^+$.
        Note that $P_o=p(V_t)$.
        Let $\hat A:=p^{-1}(A')\cap V_t$ and $\hat B:=p^{-1}(B')\cap V_t$. 
        Then $(\hat A,\hat B)$ is an $\ell$-separation of~$\hat G[V_t]$ with $V_{\hat e}\se \hat A\sm \hat B$ and $V_{\hat f}\se \hat B\sm\hat A$, as $p\rest V_t$ is injective by \cref{coveringInjective}.
        So there is a unique $\ell$-separation $\bar r=(\bar A,\bar B)\in R$ with $\bar A\cap V_t=\hat A$ and $\bar B\cap V_t=\hat B$.
        
        By \cref{coveringLemmaFive}, there is $s(\bar r)=(C,D)\in L$ such that $\bar r$ and $s(\bar r)$ have edit-distance~$\le \eps_{\ref{coveringLemmaFive}}(q,\ell)$.
        Since $s(\bar r)\in L$, we have that $s(\bar r)$ is induced in $(T,\cV)$ by an (oriented) edge~$\hat g$ of~$T$.
        Note that $\hat g$ is oriented from $\omega_1$ to $\omega_2$, since $\bar r\in R$ has $\omega_1$ living on its left side and $\omega_2$ on its right side, and since $\bar r$ and $s(\bar r)$ have finite edit-distance.

        We claim that we may assume without loss of generality that $\hat g=\hat e$, so that $s(\bar r)\le (\bar A,\bar B)$.
        Assume first that $\hat g\notin\{\hat e,\hat f\}$.
        Then $\hat e$ lies in between $\hat g$ and $\hat f$ on~$T$, say, so $\hat g < \hat e < \hat f$ in the fixed orientation of the edges of~$T$.
        Hence $s(\bar r)\le s(\hat e)\le \bar r$, where $s(\hat e)\in L$ is the separation induced in $(T,\cV)$ by the oriented edge~$\hat e$.
        Hence $\bar r$ and $s(\hat e)$ have edit-distance~$\le \eps_{\ref{coveringLemmaFive}}(q,\ell)$, so we may replace $\hat g$ with $\hat e$ in this case.
        The other cases work analogously or by renaming.
        Hence we may assume now that $\hat g=\hat e$.
        
        Since $s(\bar r)$ and $\bar r$ have edit-distance~$\le \eps_{\ref{coveringLemmaFive}}$, 
        and since $\hat A\sm V_{\hat e}\se C\triangle\bar A$,
        we have $|\hat A\sm V_{\hat e}|\le \eps_{\ref{coveringLemmaFive}}$.
        Since $P_o=p(V_t)$ and $A'=p(\hat A)$, we have $|A'\sm P_e|\le \eps_{\ref{coveringLemmaFive}}$.
        Therefore, the $AU$-corner (for $r$ and $(U,W)$) has size~$\le \eps_{\ref{coveringLemmaFive}}-\ell$ (where we deduct $\ell$ for the $rU$-link).
        But $r$ and $s=(U,W)$ $\eps$-cross by assumption with all links of size~$\ell$, so they $(\eps-\ell)$-criss by \cref{crissCross}, with $\eps \geq \eps_{\ref{coveringLemmaFive}}$, a contradiction.
    \end{cproof}

    \cref{coveringDaisy:PreDaisy} gives~\cref{coveringDaisy:adhesionSets}.
    \cref{coveringDaisy:Display} gives~\cref{coveringDaisy:ItemDisplay}.
    \cref{coveringDaisy:compatible} implies \cref{coveringDaisy:ItemCompatible} since $\eps_{\ref{coveringDaisy}}\ge \eps_{\ref{coveringLemmaFive}}$.
\end{proof}

\begin{rem}\label{coveringDaisyCanonical}
    The graphic-$k$-daisy $\cD=\cD(S)$ constructed in \cref{coveringDaisy} is canonical:
    whenever $\psi$ is an automorphism of~$G$,
    the graphic-$k$-daisies
    $\cD(S)$ and $\cD(\psi(S))$ are $\psi$-isomorphic.
\end{rem}

We can now conclude the proof of \cref{keylem:flower-approximation}, which we restate for convenience.

\flowerapprox*

\begin{proof}[Proof of \cref{keylem:flower-approximation}.]
    The $\gamma$-asteriskal set of separations~$S$ has a common centre~$X$ by \cref{coveringCentreAndLinksOfS}, and it has a common link size~$\ell$ by \cref{coveringLinkSize}.

    If $|X|=k$, then we apply \cref{coveringAnemone} and \cref{coveringAnemonecanonical} to get a $q$-lean graphic-$k$-anemone~$\cA$ that displays every separation in~$S$ and such that every separation that is $\eps \geq \eps_{\ref{keylem:flower-approximation}}$-compatible with $\cA$ is $\eps$-nested with all the petal-separations of~$\cA$. 
        
    If $|X|<k$, then we may assume \cref{set:covering} and are done by \cref{coveringDaisy} and \cref{coveringDaisyCanonical}.
\end{proof}

\section{Petal-separations are compatible}
\label{sec:petal-separations-totally-nested}

The aim of this section is to prove \cref{keylem:petal-separations-totally-nested}.
In \cref{sec:main-theorem-for-narrow}, we will apply this lemma to deduce that the petal-separations of the graphic-$k$-flowers constructed in \cref{keylem:flower-approximation} are totally-$\gamma'$-nested.

Let~$S,S'$ be two sets of separations of a graph~$G$.
We say that~$S$ \defn{$\eps$-approximates}~$S'$ if every $\eps$-broad separation of~$S'$ is $\eps$-close to a separation in~$S$.
A graphic-$k$-flower~$F$ of a graph~$G$ \defn{$\eps$-approximates} a set~$S$ of separations of~$G$ if the set of displayed separations of~$F$ $\eps$-approximates~$S$.

\begin{restatable}{keylemma}{petalseptotnested}\label{keylem:petal-separations-totally-nested}
    Let~$G$ be a $(q,k-1)$-unbreakable graph, let $\eps \in \N$ and 
    \[
        \gamma\geq\gamma_{\ref{keylem:petal-separations-totally-nested}}(q,k,\eps):=2q+3k+2\eps+1
        \quad\text{and}\quad
        \gamma'\geq\gamma'_{\ref{keylem:petal-separations-totally-nested}}(q,k):=3q+3k+1~.
    \]
    Let~$A$ be a nontrivial $\gamma$-asteriskal set of $k$-separations of~$G$.
    Let~$F$ be a $q$-lean graphic-$k$-flower of~$G$ that $\eps$-approximates~$A$.
    Then every $k$-separation of~$G$ is $\gamma'$-compatible with~$F$.
\end{restatable}

We now give a high-level sketch of the proof.
We focus here on the case where $F$ is a graphic-daisy, which is the most difficult.
Let $r$ be a $k$-separation of $G$, and let $p$ be a petal-separation of~$F$ that $\gamma'$-crosses $r$.
Let $e, f$ be the edges of $O$ such that $p$ is displayed by $\{e, f\}$.
Since $p$ and $r$ $\gamma'$-cross, we have that~$p$ is $2(q+k+1)$-broad (\cref{LargeCrossedImpliesBroad}).
Since~$A$ is nontrivial and $\gamma$-asteriskal and since~$F$ $\varepsilon$-approximates~$A$, there are displayed separations $s_1, s_2$ of~$F$ which are $(2q+2k+1)$-crossing such that $p = s_1 \wedge s_2$ holds (see \cref{lem:petal-separation-as-meet}).
We then study the 3D-corner diagram of the separations $s_1, s_2$ and~$r$ to deduce that the $pr$-centre is the $s_1s_2$-centre, and that the $p$-links are links of the crossing separations $s_1$ and $s_2$ (see \cref{lem:crosses-corner-then-induces-tripartition}).
Finally, using that $F$ is $q$-lean, we get that the $s_1s_2$-centre is $Z$, and that the two aforementioned links are $V_e \sm V_f$ and $V_f \sm V_e$ (see \cref{lem:lean-daisy-correct-links}).

\begin{lem}\label{lem:close-separations-still-cross}
    Let~$\gamma \geq 2\eps$ and let~$s$ and~$t$ be two $\gamma$-crossing separations of a graph~$G$.
    Let~$s'$ and~$t'$ be two separations of~$G$ that are $\eps$-close to~$s$ and~$t$, respectively.
    Then~$s'$ and~$t'$ $(\gamma-2\eps)$-cross.
\end{lem}
\begin{proof}
    Let $\gamma'=\gamma-2\eps$, and suppose by contradiction that $s'$ and $t'$ are $\gamma'$-nested. Without loss of generality, assume that $s' \leq_{\gamma'} t'$.
    Then $s \approx_\eps s' \leq_{\gamma'} t' \approx_\eps t$, so we have $s \leq_{\gamma'+2\eps} t$ by \cref{approxLeQuasiTransitive}, which is a contradiction.
\end{proof}

\begin{lem}\label{lem:petal-separation-as-meet}
    Let $G$ be a $(q,k-1)$-unbreakable graph, let $\beta > q$ and $\gamma \geq q+k$.
    Let~$F$ be a graphic-$k$-flower of~$G$ and let~$s$ and~$t$ be two displayed separations of~$F$ that $\gamma$-cross.
    Let~$p$ be a $\beta$-broad petal-separation of~$F$.
    Then there are two $(\gamma-k)$-broad displayed separations~$s'$ and~$t'$ of~$F$ that $\min\{\beta-1,\gamma-k\}$-cross, and such that $p=s' \wedge t'$.
\end{lem}

\begin{proof}
    By \cref{AnemoneAsDaisy}, there is a cycle-decomposition $\cO=(O,\cV)$ of~$G$ together with two 2-edge-cuts of~$O$ that induce~$s$ and~$t$ in~$G$, respectively.
    Without loss of generality, we can assume that $O=\Zbb_\ell$, where $\ell=|O|$, and that \[
        p=\left(V_0, \bigcup_{i \neq 0} V_i\right) \text{ and } 
        s=(A,B)=\left(\bigcup_{i=i_2}^{i_1-1} V_i, \bigcup_{i=i_1}^{i_2-1}V_i\right)\text{ and }
        t=(C,D)=\left(\bigcup_{i=j_2}^{j_1-1} V_i, \bigcup_{i=j_1}^{j_2-1}V_i\right)
    \] with $0 < i_1 < j_1 < i_2 < j_2 \leq \ell$.
    Then we define
    \[
        s':=(A',B'):=\left(\bigcup_{i=i_2}^{\ell} V_i, \bigcup_{i=1}^{i_2-1}V_i\right)\text{ and }
        t':=(C',D'):=\left(\bigcup_{i=0}^{j_1-1} V_i, \bigcup_{i=j_1}^{\ell-1}V_i\right)~.
    \]
    Since~$p$ is $\beta$-broad, the $A'C'$-corner has size~$\geq \beta$.
    The $BC$-corner has size~$> \gamma-k$ by \cref{crissCross}~\ref{crissCross:crossImpliesCriss} and is contained in the $B'C'$-corner.
    Similarly, the $AD$-corner has size~$> \gamma-k$ by \cref{crissCross}~\ref{crissCross:crossImpliesCriss} and is contained in the $A'D'$-corner.
    Moreover, the $BD$-corner is equal to the $B'D'$-corner, which has size~$> \gamma-k$ by \cref{crissCross}~\ref{crissCross:crossImpliesCriss}.
    Thus, all four corners of $s'$ and $t'$ have size $> \min\{\beta-1,\gamma-k\}$, so $s'$ and $t'$ $\min\{\beta-1,\gamma-k\}$-cross by \cref{crissCross}~\ref{crissCross:crissImpliesCross}.
    In particular, all four corners of $s'$ and $t'$ have size $> q$.
    Since $G$ is $(q,k-1)$-unbreakable, the corresponding corner separations all have order at least $k$. Since $s'$ and $t'$ both have order at most $k$, submodularity implies that all four corner separations have order exactly $k$.
    Moreover, $p$ has order $k$ since it is $\beta$-broad with $\beta > q$.
    The separations $p$ and $s' \wedge t'$ have the same second side, and the first side of $p$ is contained in that of $s' \wedge t'$.
    This containment cannot be strict, otherwise $p$ and $s' \wedge t'$ would not have the same order.
    Thus, $p=s' \wedge t'$.

    Finally, since the $B'C'$-corner and the $A'D'$-corner have size~$> \gamma-k$, both~$s'$ and~$t'$ are $(\gamma-k)$-broad.
\end{proof}

\begin{lem}\label{lem:crosses-corner-then-induces-tripartition-1}
    Let~$G$ be a $(q,k-1)$-unbreakable graph.
    Let~$s_1=(A_1,B_1)$, $s_2=(A_2,B_2)$ be two $(q+k+1)$-crossing $k$-separations of~$G$.
    Denote by~$L_i$ the $A_is_{3-i}$-link, let~$X$ denote the $s_1s_2$-centre and let $p=s_1 \wedge s_2$.
    Let $r=(C,D)$ be a $k$-separation of~$G$ such that the $A_1B_2D$-corner, the $A_1A_2D$-corner, the $A_1A_2C$-corner and the $B_1A_2C$-corner each have size~$> q$.
    Then~$r$ induces $(X,\{L_1,L_2\})$ on~$\sep(p)$. 
\end{lem}

\begin{proof}
    First, we show that the prerequisites for applying \cref{cor:3D-crossing3} are fulfilled for the separations~$s_1$,~$s_2$ and~$r$.
    Note that the 3D corner-separations $A_1B_2D$, $A_1A_2D$, $A_1A_2C$ and $B_1A_2C$ have order $\geq k$ since the corresponding corners each have size~$> q$ and~$G$ is $(q,k-1)$-unbreakable.
    Now consider the two separations~$s_1$ and~$r$.
    Then the $A_1D$-corner and the $B_1C$-corner both have size $> q$.
    Hence, since~$G$ is $(q,k-1)$-unbreakable and by submodularity, the corresponding corner-separations have order exactly~$k$.
    Similarly, the $A_2C$-corner-separation and the $B_2D$-corner-separation of the separations~$s_2$ and~$r$ both have order exactly~$k$.
    Now, consider the two separations~$s_1$ and~$s_2$.
    Since they are $(q+k+1)$-crossing, each (2D) corner has size~$>q$ by \cref{crissCross}~\ref{crissCross:crossImpliesCriss}.
    By submodularity and since~$G$ is $(q,k-1)$-unbreakable, each 2D corner-separation of~$s_1$ and~$s_2$ has order exactly~$k$.
    In particular, the 2D corner-separation for~$A_1A_2$ has order exactly~$k$.
    
    So we can apply \cref{cor:3D-crossing3} to~$s_1$,~$s_2$ and~$r$.
    Then, we get that the $A_1s_2r$-link, the $s_1A_2r$-link, the $s_1s_2C$-link, the $s_1s_2D$-link, the $A_1s_2C$-link and the $s_1A_2D$-link are empty.
    Then, $L_1 = (A_2 \cap B_2) \setminus B_1$ is the union of the $A_1s_2r$-link, $A_1s_2C$-link and $A_1s_2D$-link, so~$L_1$ is equal to the $A_1s_2D$-link and thus $L_1 \subseteq D \setminus C$.
    Similarly, $L_2 = (A_1 \cap B_1) \setminus B_2$ is the union of the $s_1A_2r$-link, $s_1A_2C$-link and $s_1A_2D$-link, so~$L_2$ is equal to the~$s_1A_2C$-link and thus $L_2 \subseteq C \setminus D$.
    Finally, $X = A_1 \cap B_1 \cap A_2 \cap B_2$ is the union of the $s_1s_2r$-centre, $s_1s_2C$-link and $s_1s_2D$-link, so~$X$ is equal to the $s_1s_2r$-centre and thus $X \subseteq C \cap D$.
    Therefore,~$r$ induces $(X,\{L_1,L_2\})$ on~$\sep(p)$.
\end{proof}

\begin{lem}\label{lem:crosses-corner-then-induces-tripartition-2}
    Let~$G$ be a $(q,k-1)$-unbreakable graph.
    Let~$s_1=(A_1,B_1)$, $s_2=(A_2,B_2)$ be two $(q+k+1)$-crossing $k$-separations of~$G$. 
    Let $r=(C,D)$ be a $k$-separation of~$G$ such that the $A_1A_2D$-corner, the $A_1B_2C$-corner, the $B_1A_2C$-corner, and the $B_1B_2D$-corner each have size~$> q$.
    Then $\sep(s_1)=\sep(s_2)=\sep(r)$.
\end{lem}

\begin{proof}
    By \cref{cor:3D-crossing1} applied to the separations~$s_1$,~$s_2$ and~$r^*$, we get that the $A_1s_2r$-link, the $s_1A_2r$-link, the $s_1s_2D$-link, the $A_1A_2r$-link, the $A_1s_2D$-link and the $s_1A_2D$-link are empty.
    So, the neighbourhood of the $A_1A_2D$-corner is contained in the $s_1s_2r$-centre.
    Since the $A_1A_2D$-corner has size~$>q$ and~$G$ is~$(q,k-1)$-unbreakable, the~$A_1A_2D$-corner-separation has order~$\ge k$.
    Since~$s_1$,~$s_2$ and~$r$ each have order~$k$, it follows that $\sep(s_1)=\sep(s_2)=\sep(r)$.
\end{proof}

\begin{lem}\label{lem:crosses-corner-then-induces-tripartition}
    Let~$G$ be a $(q,k-1)$-unbreakable graph and let $\delta:=\delta(q,k):=2q+2k+1$ and $\delta':=\delta'(q,k):=3q+3k+1$.
    Let~$s_1=(A_1,B_1)$ and~$s_2=(A_2,B_2)$ be two $\delta$-crossing $k$-separations of~$G$.
    Let~$X$ be the $s_1s_2$-centre and let~$L_i$ be the $A_is_{3-i}$-link.
    Assume that $p:=s_1 \wedge s_2$ has order~$k$ and is $\delta'$-crossed by a $k$-separation $r=(C,D)$ of~$G$.
    Then~$r$ induces~$(X,\{L_1,L_2\})$ on~$\sep(p)$.
\end{lem}

\begin{proof}
    Since~$p$ and~$r$ $(q+k+1)$-cross, the $A_1A_2C$-corner and the $A_1A_2D$-corner both have size~$>q$ by \cref{crissCross}~\ref{crissCross:crossImpliesCriss}.
    Since~$s_1$ and~$s_2$ $\delta$-cross, the 2D corners for~$A_1B_2$ and for~$B_1A_2$ both have size~$\geq \delta-k$.
    Since the $A_1B_2r$-link has size~$\leq |\sep(r)| \leq k$, one of the corners for~$A_1B_2C$ or~$A_1B_2D$ has size $\geq (\delta-2k)/2 > q$.
    Similarly, one of the corners for~$B_1A_2C$ or~$B_1A_2D$ has size $> q$.
    Hence, we can distinguish the following cases.
    
    \casen{Case~1: the $A_1B_2D$-corner and the $B_1A_2C$-corner have size $>q$.} Then the claim follows from \cref{lem:crosses-corner-then-induces-tripartition-1}.
    
    \casen{Case~2: the $A_1B_2C$-corner and the $B_1A_2D$-corner have size $>q$.} This case is symmetric to the previous case.

    \casen{Case~3: the $A_1B_2C$-corner and the $B_1A_2C$-corner have size $>q$, and the $A_1B_2D$-corner and the $B_1A_2D$-corner have size $\leq q$.} 
    First, we check that the prerequisites for applying \cref{lem:crosses-corner-then-induces-tripartition-2} are fulfilled.
    It remains to show that the $B_1B_2D$-corner has size~$>q$.
    Since~$r$ $\delta'$-crosses~$p$, each corresponding 2D corner has size $\geq \delta'-k$ by \cref{crissCross}~\ref{crissCross:crossImpliesCriss}.
    Consider the two-dimensional $(B_1 \cup B_2)D$-corner, which is equal to the union of the three-dimensional $A_1B_2D$-corner, the $B_1A_2D$-corner, the $B_1B_2D$-corner, the $B_1s_2D$-link and the $s_1B_2D$-link.
    Both the $B_1s_2D$-link and the $s_1B_2D$-link have size~$\leq k$, and both the $A_1B_2D$-corner and the $B_1A_2D$-corner have size $\le q$.
    So, the $B_1B_2D$-corner has size $\ge \delta' - k - 2k - 2q > q$.
    Hence, we can apply \cref{lem:crosses-corner-then-induces-tripartition-2} and get that $\sep(s_1)=\sep(s_2)=\sep(r)=X$.
    In particular, $X=\sep(p)$ and $L_1=\emptyset=L_2$ and~$r$ induces~$(X,\{L_1,L_2\})$ on~$\sep(p)$.

    \casen{Case~4: the $A_1B_2D$-corner and the $B_1A_2D$-corner have size $>q$, and the $A_1B_2C$-corner and the $B_1A_2C$-corner have size $\leq q$.} This case is symmetric to the previous case.
\end{proof}

\begin{lem}\label{lem:lean-anemone-correct-links}
    Let $G$ be a $(q,k-1)$-unbreakable graph and let $F$ be a $q$-lean graphic-$k$-anemone of $G$, with centre $Z$.
    Let $p$ be a $(>q)$-broad petal-separation of $F$ and let $s_1 = (A_1, B_1)$ and $s_2 = (A_2, B_2)$ be displayed separations of $F$ that $q$-criss such that $p = s_1 \wedge s_2$.
    Then, the $s_1s_2$-centre is $Z$.
\end{lem}

\begin{proof}
    Since $F$ is $q$-lean, we have $|Z| = k$.
    Since $s_1$ and $s_2$ $q$-criss, they are both $(>q)$-broad. 
    Since $G$ is $(q,k-1)$-unbreakable, $s_1$ and $s_2$ are $k$-separations.
    Since they are displayed by $F$, and since $Z$ is the centre of $F$, we have $\sep(s_1) = \sep(s_2) = Z$.
    Thus, the $s_1s_2$-centre is $Z$.
\end{proof}

\begin{lem}\label{lem:lean-daisy-correct-links}
    Let $G$ be a $(q,k-1)$-unbreakable graph and let $F \eqqcolon (O, \cV)$ be a $q$-lean graphic-$k$-daisy of $G$.
    Let $Z \coloneqq \bigcap_{o \in V(O)}V_o$ be the centre of $F$.
    Let $p$ be a $(>q)$-broad petal-separation of $F$ and let $s_1 = (A_1, B_1)$ and $s_2 = (A_2, B_2)$ be displayed separations of $F$ that $q$-criss such that $p = s_1 \wedge s_2$.
    Denote by~$L_i$ the $A_is_{3-i}$-link.
    For every $e \in E(O)$, let $V_e$ be the adhesion-set at $e$.
    Let $e,f \in E(O)$ be the edges of $O$ such that $p$ is displayed by $\{e, f\}$.
    Then, the $s_1s_2$-centre is $V_e \cap V_f = Z$ and $\{L_1, L_2\} = \{V_e \sm Z, V_f \sm Z\}$.
\end{lem}

\begin{proof}
    Since $s_1$ and $s_2$ are displayed separations of $F$ which cross, and $p = s_1 \wedge s_2$, there exist edges $e', f' \in E(O)$ such that the edges $e, f, e', f'$ appear in this order along $O$, $s_1$ is displayed by $\{e, e'\}$, and $s_2$ is displayed by $\{f, f'\}$.
    Since $s_1$ and $s_2$ $q$-criss, all their corner separations are $(>q)$-broad.
    These separations are displayed respectively by $\{e, f\}$, $\{f, e'\}$, $\{e', f'\}$, and $\{f', e\}$.
    Since $s_1$ and $s_2$ $q$-criss, they are $(>q)$-broad.
    Since $F$ is $q$-lean, we then have that the sets $V_{e} \sm Z, V_{f} \sm Z, V_{e'} \sm Z$ and $V_{f'} \sm Z$ are pairwise disjoint.
    Thus, the $s_1s_2$-centre is $Z$, and $L_1 = V_{f} \sm Z$ and $L_2 = V_{e} \sm Z$.
\end{proof}

We now prove \cref{keylem:petal-separations-totally-nested}, which we restate for convenience.

\petalseptotnested*

\begin{proof}
    Let $r$ be a $k$-separation of $G$. We show that $r$ is $\gamma'$-compatible with $F$.
    Let $p$ be a petal-separation of $F$ which $\gamma'$-crosses $r$.
    Note that~$p$ is $2(\gamma'-k)$-broad by \cref{LargeCrossedImpliesBroad}.
    Since $2(\gamma'-k) > q$ and $G$ is $(q,k-1)$-unbreakable, $p$ has order exactly $k$.
    Since~$A$ is nontrivial and $\gamma$-asteriskal, there are two separations~$r_1$ and~$r_2$ in~$A$ that $\gamma$-cross.
    Note that both~$r_i$ are $2(\gamma-k)$-broad by \cref{LargeCrossedImpliesBroad}, in particular they are $\eps$-broad since $2(\gamma-k) \geq \eps$.
    Since~$F$ $\eps$-approximates~$A$, there are two displayed separations~$r_1'$ and~$r_2'$ of~$F$ that are $\eps$-close to~$r_1$ and~$r_2$, respectively.
    Then, $r_1'$ and~$r_2'$ $(\gamma-2\eps)$-cross by \cref{lem:close-separations-still-cross}.
    By \cref{lem:petal-separation-as-meet}, there are $\delta$-crossing displayed separations~$s_1$ and~$s_2$ of~$F$ such that $p=s_1 \wedge s_2$, where 
    \[
        \delta = \min\{2\gamma'-2k-1, \gamma-2\eps-k\} \geq 2(q+k)+1~.
    \]
    Write $s_i=:(A_i,B_i)$ and let~$L_i$ be the $A_is_{3-i}$-link of the two separations~$s_1$ and~$s_2$ and let~$X$ denote the $s_1s_2$-centre.
    Recall that~$p$ is $\gamma'$-crossed by~$r$, where $\gamma' \geq 3q+3k+1$.
    By \cref{lem:crosses-corner-then-induces-tripartition},~$r$ induces the tripartition~$(X,\{L_1,L_2\})$ on~$p$. 
    Since $s_1$ and $s_2$ are $\delta$-crossing, they are $q$-crissing.

    If $F$ is a graphic-$k$-anemone, let $Z$ denote its centre. Then, by \cref{lem:lean-anemone-correct-links}, we get $X = Z$. 
    Since $F$ is $q$-lean, we have $|Z| = k$.
    Since $p$ and $r$ are $k$-separations, they cross with centre $Z$ and four empty links.
    Thus, $r$ is $\gamma'$-compatible with $F$.

    If $F$ is a graphic-$k$-daisy, write $F= (O, \cV)$ and let $Z = \bigcap_{o \in V(O)}V_o$.
    Let $e, f$ be the edges of $O$ such that $p$ is displayed by $\{e, f\}$.    
    By \cref{lem:lean-daisy-correct-links}, we have $X = V_e \cap V_f = Z$, and $\{L_1, L_2\} = \{V_e \sm Z, V_f \sm Z\}$. This proves the first two conditions of \cref{def:F-compatible}.

    Since $p$ and $r$ $\gamma'$-cross, they $(\gamma'-k)$-criss, with $\gamma'-k > q$.
    Thus, their four corner separations are $(>q)$-broad.
    Since $G$ is $(q,k-1)$-unbreakable, these four separations have order $\ge k$.
    By submodularity, since $p$ and $r$ are $k$-separations, these four corner separations have order exactly~$k$.
    This implies that all the links of the two separations $p, r$ are of the same size.
    Therefore, $r$ is $\gamma'$-compatible with $F$.
\end{proof}

 \section{Approximation by a maximal nested set}
\label{sec:max-nested-set}

Following the strategy outlined in the overview, we construct a tree-decomposition from an approximately nested family of separations. We allow bounded perturbations before choosing a maximal nested family, and show that this family still approximates the original separations.

Recall that a graph~$G$ is $(b,k)$-narrow if every $b$-broad $k$-separation is narrow.

Let~$S$ be a set of separations of a graph~$G$.
A separation~$s$ is \defn{$S$-totally-$\eps$-nested} if $s \in S$ and $s$ is $\eps$-nested with every separation in~$S$.
A separation of~$G$ is \defn{totally-$\eps$-nested} if it is $S$-totally-$\eps$-nested where~$S$ is the set of all proper $(\leq k)$-separations of~$G$.
A separation of~$G$ is \defn{narrow-$\eps$-nested} if it is $S$-totally-$\eps$-nested where~$S$ is the set of all narrow proper $(\leq k)$-separations of~$G$.
Here, $k$ will be clear from context.
We abbreviate $S$-totally-$0$-nested to \defn{$S$-totally-nested}, and totally-$0$-nested to \defn{totally-nested}, and narrow-$0$-nested to \defn{narrow-nested}.

In this section, we prove the following result.
\begin{restatable}{keylemma}{maxnestedset}\label{keylem:max-nested-set}
    Let~$G$ be a $(q,k-1)$-unbreakable $(b,k)$-narrow graph and let $\gamma \in \N$.
    Let 
    \begin{align*}
        \delta&:=\delta_{\ref{keylem:max-nested-set}}(k,b,\gamma) := 3\delta_{\ref{lem:max-nested-set}}(k,b,\gamma) + 2\gamma+2k \quad \text{and}\\
        \gamma'&:=\gamma'_{\ref{keylem:max-nested-set}}(k,b,\gamma):= \delta+\gamma \quad\text{and}\\
        \beta&:=\beta_{\ref{keylem:max-nested-set}}(q,k,b):= \beta_{\ref{lem:max-nested-set}}(q,k,b) \quad\text{and}\\
        \eps&:=\eps_{\ref{keylem:max-nested-set}}(q,k,b,\gamma):= \max\{\, \eps_{\ref{lem:max-nested-set}}(q,k,b,\gamma)\, , \, \delta+\beta\,\} \quad\text{and}\\
        \delta'&:=\delta'_{\ref{keylem:max-nested-set}}(k,b,\gamma) := \delta_{\ref{lem:max-nested-set}}(k,b,\gamma') \quad\text{and}\\
        \eps'&:=\eps'_{\ref{keylem:max-nested-set}}(q,k,b,\gamma):=\eps_{\ref{lem:max-nested-set}}(q,k,b,\gamma')~.
    \end{align*}
    Let~$S$ be the set of all narrow-$\gamma$-nested $k$-separations of~$G$.
    Let~$S'$ be the set of all half-connected narrow-$\gamma'$-nested $\beta$-broad $k$-separations of~$G$ that are $\delta$-close to some separation in~$S$.
    Let~$N \subseteq S'$ be a maximal nested set of separations in~$S'$.
    
    Then~$N$ $\eps$-approximates~$S$.
    Moreover, for every $\eps'$-broad $k$-separation of~$G$ there is a $\delta'$-close $k$-separation of~$G$ that is nested with~$N$.
\end{restatable}

Let~$S$ be a set of separations of a graph~$G$.
We say that~$S$ is \defn{$\eps$-canonical} if for every automorphism~$\varphi$ of~$G$, the sets~$S$ and~$\varphi(S)$ $\eps$-approximate each other.
So if~$S$ is $\eps$-canonical for $\eps=0$, then~$S$ is canonical.
Similarly, we call a tree-decomposition~$\cT=(T,\cV)$ \defn{$\eps$-canonical} if the set of separations that are induced by an edge of~$T$ is $\eps$-canonical.

\begin{corollary}\label{canonical:max-nested-set}
    Assume the premises of \cref{keylem:max-nested-set}.
    Then the set~$N$ is $2\eps$-canonical.
\end{corollary}
\begin{proof}
    Let~$\varphi$ be an automorphism of~$G$ and let~$N':=\varphi(N)$.
    By \cref{keylem:max-nested-set},~$N$ $\eps$-approximates~$S$.
    Similarly,~$N'=\varphi(N)$ $\eps$-approximates~$\varphi(S)=S$.
    Since every separation in~$S'$ is $\delta$-close to some separation in~$S$, the set~$S$ $\delta$-approximates~$S'$.
    Recall that~$N\subseteq S'$.
    Similarly, we have~$N'=\varphi(N) \subseteq \varphi(S')=S'$.
    Hence,~$S$ $\delta$-approximates~$N$ and~$N'$.

    To summarize,~$N$ $\eps$-approximates~$S$, which in turn $\delta$-approximates~$N'$.
    So, $N$ $(\eps+\delta)$-approximates~$N'$ by \cref{approxLeQuasiTransitive}.
    Similarly,~$N'$ $(\eps+\delta)$-approximates~$N$.
    It follows that~$N$ and~$N'$ $(\eps+\delta)$-approximate each other, hence~$N$ is $(\eps+\delta)$-canonical, so $2\eps$-canonical.
\end{proof}

A proper separation~$(A,B)$ of~$G$ is \defn{half-connected} if $G[A\sm B]$ or $G[B \sm A]$ is connected.

\medskip

We now give a high-level overview of the proof. 
Let~$s$ be a narrow-$\gamma$-nested separation of~$G$.
By definition, $s$ is $\gamma$-nested with every separation in~$S'$. 
Thus, either~$s$ is $\gamma$-close to some separation of~$N$ or there is a unique splitting star~$\sigma$ of~$N$ such that~$s$ $\gamma$-interlaces $\sigma$ (\cref{defn:eps-interlaces} and \cref{lem:gamma-close-or-gamma-interlaces}).
In the first case, we are done, so suppose that~$s$ $\gamma$-interlaces~$\sigma$.
In that case, either~$s$ crosses few separations of~$\sigma$, or each side of~$s$ has many connected components (\cref{lem:few-crossings-or-many-components}).
If~$s$ crosses few separations of $\sigma$ then we can ``wiggle'' (see \cref{lem:wiggle-lemma}) $s$ to a $k$-separation~$s' \in S'$ that almost-interlaces~$\sigma$ and is $\delta$-close to~$s$.
Then, $s'$ is nested with every separation in $N$, so $s' \in N$ by maximality of $N$.
If each side of~$s$ has many connected components, then every half-connected separation~$s'$ with separator $\sep(s)$ is not crossed by a narrow proper $(\leq k)$-separation of~$G$ (\cref{lem:many-components-implies-tot-nested}), so it essentially belongs to~$N$.
Then, since~$s$ is narrow-$\gamma$-nested, it is $\delta$-close to some~$s' \in N$ (\cref{FlowerClosePetalOrCrossed}), and we are done. 

\subsection{Splitting stars and approximate interlacing}

\begin{lemma}
    \label{lem:assymmetric-nested-implies-not-broad}
    Let~$\gamma \in \Nbb$, let~$G$ be a graph, and let~$s$ be a separation of~$G$.
    If $s \leq_{\gamma} s^*$ then~$s$ has breadth~$\leq \gamma$.
\end{lemma}
\begin{proof}
    Let~$s=:(A,B)$.
    Then $s \leq_{\gamma} s^*$ implies $A \subseteq_{\gamma} B$.
    Hence, $|A\sm B| \leq \gamma$.
\end{proof}

Given a set~$N$ of separations of~$G$, an \defn{orientation} of~$N$ is a subset $M \subseteq N$ such that for every $s \in N$ exactly one of~$s$ and~$s^*$ is in~$M$.
The following result is folklore; we include a proof for completeness.

\begin{lem}
\label{lem:consistent-then-splitting-star}
    Let~$\vec N$ be a nested consistent finite set of proper separations of a graph~$G$.
    Then the set of all $\leq$-maximal separations in~$\vec N$ is the unique splitting star in~$\vec N$.
\end{lem}

\begin{proof}
    Let $\sigma$ be the set of all $\leq$-maximal separations in~$\vec N$.
    First, we argue that $\sigma$ is a splitting star in $\vec N$.
    Let $s \neq t \in \sigma$.
    Since $\vec N$ is nested, we have $s \leq t$ or $s \leq t^*$ or $t \leq s$ or $t^* \leq s$.
    Since the elements of $\sigma$ are $\leq$-maximal in $\vec N$, we have either $s \leq t^*$ or $t^* \leq s$.
    Since $\vec N$ is consistent, we must have $s \leq t^*$, so $\sigma$ is a star.
    Now we show that~$\sigma$ is splitting in~$\vec N$.
    Assume the contrary.
    Then there exists $s \in \vec N$ that interlaces $\sigma$.
    Thus, $s \notin \sigma$, so~$s$ is not $\leq$-maximal in~$\vec N$.
    Hence there exists $t \in \sigma$ such that $s < t$.
    Since $s$ interlaces $\sigma$, we have either $t < s$ or $t < s^*$.
    The first outcome contradicts the $\leq$-maximality of~$t$, so we have $t < s^*$.
    Thus, $s < t < s^*$, which contradicts that~$s$ is proper.

    Now, let $\tau$ be a splitting star in $\vec N$.
    We first show that $\sigma \subseteq \tau$.
    By contradiction, suppose that there exists $s \in \sigma \sm \tau$.
    Since $\tau$ is splitting, $s$ does not interlace $\tau$ so there exists $t \in \tau$ such that $t \not< s$ and $t \not< s^*$.
    Since $\vec N$ is nested, we have either $s \leq t$ or $s^* \leq t$. 
    The first outcome would contradict the $\leq$-maximality of $s$ and the second that $\vec N$ is consistent.
    Thus, $\sigma \subseteq \tau$.
    Now, we show that $\tau \se \sigma$.
    By contradiction, suppose that there exists $t \in \tau \sm \sigma$.
    Since $t \notin \sigma$, there exists $s \in \sigma$ such that $t \leq s$.
    Since $\tau$ is a star and $s \in \sigma \subseteq \tau$, we have $s \leq t^*$.
    Thus, $t \leq s \leq t^*$, which contradicts that $t$ is proper.
    Therefore, $\tau = \sigma$, which concludes the proof.
\end{proof}

It is folklore that if a proper separation~$s$ is nested with every separation from a finite symmetric nested set~$N$ of proper separations, either $s \in N$ or~$s$ interlaces a unique splitting star of~$N$.
\cref{lem:gamma-close-or-gamma-interlaces} is a generalization of this statement to the approximate setting.
For this, we need the following definition.

\begin{definition}[Approximate version of interlacing]
\label{defn:eps-interlaces}
    Let~$\sigma$ be a star of separations of a graph~$G$, and let~$s$ be a separation of~$G$.
    We say that~$s$ \defn{$\eps$-interlaces}~$\sigma$ if for every $t \in \sigma$, we have $t <_\eps s$ or $t <_\eps s^*$.
    We say that~$s$ \defn{almost-$\eps$-interlaces}~$\sigma$ if for every $t \in \sigma$, we have $t \le_\eps s$ or $t \le_\eps s^*$.
    Now let~$S$ be a set of separations of~$G$.
    We say that~$S$ (\defn{almost-})\defn{$\eps$-interlaces}~$\sigma$ if every separation in~$S$ (almost-)$\eps$-interlaces~$\sigma$.
\end{definition}

\begin{lemma}[{Generalization of~\cite[Corollary V.9]{TriAiC}}]
\label{lem:gamma-close-or-gamma-interlaces}
    Let~$G$ be a graph and let~$N$ be a symmetric nested finite set of proper separations of~$G$.
    Let~$s$ be a $(2\gamma+1)$-broad separation of~$G$ that is $\gamma$-nested with every separation in~$N$.
    Then either
    \begin{enumerate}
        \item\label{itm:gamma-close-or-gamma-interlaces-1} $s$ is $\gamma$-close to some $n \in N$, or
        \item\label{itm:gamma-close-or-gamma-interlaces-2} $O(s):=\{t \in N \mid t \leq_\gamma s \text{ or } t \leq_\gamma s^*\}$ is a consistent orientation of~$N$, and the set of $\leq$-maximal elements in~$O(s)$ is the unique splitting star~$\sigma$ of~$N$ such that~$s$ $\gamma$-interlaces~$\sigma$.
    \end{enumerate}
\end{lemma}

\begin{proof}
    Assume that \ref{itm:gamma-close-or-gamma-interlaces-1} fails, i.e.~$s$ is not $\gamma$-close to any $n \in N$.

    \begin{claim}\label{claim:gamma-close-or-gamma-interlaces-1}
        For all $n \in N$, we have $n \in O(s)$ if and only if $n^* \notin O(s)$.
    \end{claim}
    \begin{cproof}
        Note that since~$s$ is $\gamma$-nested with~$N$, we either have $n \in O(s)$ or $n^* \in O(s)$ for every $n \in N$.
        On the other hand, assume both $n \in O(s)$ and~$n^* \in O(s)$ for some $n \in N$.
        Then either $n \leq_\gamma s$ or $n \leq_\gamma s^*$, and either $n^* \leq_\gamma s$ or $n^* \leq_\gamma s^*$.
        Since~$s$ is not $\gamma$-close to~$n$ and not $\gamma$-close to~$n^*$, we either have $s^* \leq_\gamma n \leq_\gamma s$ or $s \leq_\gamma n \leq_\gamma s^*$.
        Without loss of generality, we can assume the latter.
        Then $s \leq_{2\gamma} s^*$, which implies 
        that~$s$ has breadth~$\leq 2\gamma$ by \cref{lem:assymmetric-nested-implies-not-broad}, a contradiction.
    \end{cproof}

    \begin{claim}\label{claim:gamma-close-or-gamma-interlaces-2}
        If $n_1 \in N$ and~$n_2 \in O(s)$ are such that $n_1 \leq n_2$, then $n_1 \in O(s)$.
    \end{claim}
    \begin{cproof}
        Indeed, $n_2 \in O(s)$ implies $n_2 \leq_{\gamma} s$ or $n_2 \leq_{\gamma} s^*$.
        Since $n_1 \leq n_2$, we also have $n_1 \leq_{\gamma} s$ or $n_1 \leq_{\gamma} s^*$.
        So, $n_1 \in O(s)$.
    \end{cproof}
    
    \begin{claim}\label{claim:gamma-close-or-gamma-interlaces-3}
        Let~$\sigma$ be a splitting star of~$N$.
        Then $\sigma \subseteq O(s)$ if and only if~$s$ $\gamma$-interlaces~$\sigma$.
    \end{claim}
    \begin{cproof}
    \[
    \sigma \subseteq O(s) \;\Longleftrightarrow\; \forall n \in \sigma,\; n \leq_\gamma s \text{ or } n \leq_\gamma s^* \;\Longleftrightarrow\; s \text{ $\gamma$-interlaces }\sigma~.
    \]
    \end{cproof}

    From \cref{claim:gamma-close-or-gamma-interlaces-1,claim:gamma-close-or-gamma-interlaces-2}, it follows that~$O(s)$ is consistent.
    By \cref{lem:consistent-then-splitting-star}, there is a unique splitting star~$\sigma$ in~$O(s)$.
    By \cref{claim:gamma-close-or-gamma-interlaces-1,claim:gamma-close-or-gamma-interlaces-2}, it follows that $\sigma$ is the only splitting star of $N$ contained in $O(s)$.
    From \cref{claim:gamma-close-or-gamma-interlaces-3}, we finally get that~$\sigma$ is the unique splitting star of~$N$ that is $\gamma$-interlaced by~$s$.
\end{proof}

We now establish two results about how narrow separations can cross other separations.

\begin{lemma}\label{lem:link-empty-or-component-in-corners}
    Let~$s=(A,B)$ and~$t=(C,D)$ be crossing $k$-separations in a graph~$G$. 
    Assume that~$s$ is narrow and proper.
    Then either
    \begin{enumerate}
        \item\label{itm:link-empty-or-component-in-corners-1} the $At$-link is nonempty, or
        \item\label{itm:link-empty-or-component-in-corners-2} $\sep(s)=\sep(t)$ and $G-\sep(t)$ has a component in the $AC$-corner and a component in the $AD$-corner.
    \end{enumerate}
\end{lemma}
\begin{proof}
    Assume that~\ref{itm:link-empty-or-component-in-corners-1} fails, i.e.\ the $At$-link is empty.
    First, we claim that both the $AC$-corner and the $AD$-corner are nonempty.
    Suppose for contradiction that the $AC$-corner is empty, say.
    Then the $AD$-corner is nonempty.
    Since~$s$ is narrow, the $sC$-link is empty.
    So~$s$ and~$t$ are nested, a contradiction.

    Hence, both the $AC$-corner and the $AD$-corner are nonempty.
    Since~$(A,B)$ is narrow, both the $sC$-link and the $sD$-link are empty, so the $st$-centre has size~$k$.
    Hence, $\sep(s)=\sep(t)$ and we have outcome~\ref{itm:link-empty-or-component-in-corners-2}.
\end{proof}

\begin{lem}\label{lem:many-components-implies-tot-nested}
    Let~$G$ be a graph and let~$s$ be a narrow half-connected proper separation of~$G$.
    Assume that $G-\sep(s)$ has $\geq k+1$ components.
    Then~$s$ is not crossed by any narrow proper $(\leq k)$-separation of~$G$.
\end{lem}
\begin{proof}
    Let~$s=(A,B)$ and assume that there is a narrow $(\leq k)$-separation~$t=(C,D)$ that crosses~$s$.
    \begin{claim}
         The $sC$-link or the $sD$-link is empty.
    \end{claim}
    \begin{cproof}
        Assume for contradiction that there is a vertex~$v$ in the $sC$-link and a vertex~$w$ in the $sD$-link.
        Since $s$ is narrow and $G-\sep(s)$ has $\geq k+1$ components, there are~$k+1$ independent $v$--$w$ paths in~$G$.
        Each of those paths intersects the separator of~$t$, contradicting $|\sep(t)| \leq k$.
    \end{cproof}

    Assume without loss of generality that the $sC$-link is empty.
    Consider any component~$K$ of ~$C \sm D$.
    Say,~$K$ is contained in the $AC$-corner.
    Since~$t$ is narrow, we have $N(K)=\sep(t)$ so the $Bt$-link is empty.
    Since~$s$ and~$t$ cross, the $BC$-corner must be nonempty.
    Since~$t$ is narrow, the $At$-link is empty, so the separator of~$t$ is contained in the $st$-centre; i.e.\ $\sep(t) \subseteq \sep(s)$.

    Since~$s$ is half-connected, either $G[A \sm B]$ or $G[B \sm A]$ is connected, say the former.
    Then the $AD$-corner is empty.
    Since~$s$ and~$t$ cross, there is a vertex~$u$ in the $sD$-link.
    Then, however,~$u$ does not have a neighbour in~$A \sm B$, contradicting that~$s$ is tight.
\end{proof}

We now demonstrate some properties of the operations $\wedge$ and $\vee$, and how they interact with $\leq_{\gamma}$. The proofs are elementary.

\begin{lem}\label{lem:corners-when-crossing-two-nested-sepns}
    Let $s=(A,B)$, $t=(C,D)$ and $s'=(E,F)$ be separations of a graph~$G$, and $\gamma \in \Nbb$.
    Assume $s \leq t^*$.
    Let $(X_1,Y_1):=s \vee s'$ and $(X_2,Y_2):=s^* \wedge s'$.
    Then all of the following hold:
    \begin{enumerate}
        \item\label{itm:corners-when-crossing-two-nested-sepns-1} the $CF$-corner is equal to the $CY_1$-corner and equal to the $CY_2$-corner;
        \item\label{itm:corners-when-crossing-two-nested-sepns-2} the $CE$-corner is equal to the $CX_1$-corner and equal to the $CX_2$-corner;
        \item\label{itm:corners-when-crossing-two-nested-sepns-3} If~$t \leq_\gamma s'$, then both $t \leq_\gamma s \vee s'$ and $t \leq_\gamma s^* \wedge s'$.
        \item\label{itm:corners-when-crossing-two-nested-sepns-4} If~$t \leq_\gamma (s')^*$, then both $t \leq_\gamma (s \vee s')^*$ and $t \leq_\gamma (s^* \wedge s')^*$.
    \end{enumerate}
\end{lem}

\begin{proof}
    \ref{itm:corners-when-crossing-two-nested-sepns-1} and~\ref{itm:corners-when-crossing-two-nested-sepns-2}.
    By symmetry, it suffices to show~\ref{itm:corners-when-crossing-two-nested-sepns-1}.
    Since $s \leq t^*$, we have $A \subseteq D$ and $C \subseteq B$.
    Note that the $CF$-corner is $(C \cap F) \setminus (D \cup E)$.
    We have $X_1 = A \cup E$ and $Y_1 = B \cap F$. 
    Thus, the $CY_1$-corner is 
    \begin{align*}
        (C \cap Y_1) \setminus (D \cup X_1) &= (C \cap B \cap F) \setminus (D \cup A \cup E) &&\\
            &= (C \cap F) \setminus (D \cup E) &&\text{(since $C \subseteq B$ and $A \subseteq D$).}
    \end{align*}

    We also have $X_2 = B \cap E$ and $Y_2 = A \cup F$. Thus, the $CY_2$-corner is
    \begin{align*}
        (C \cap Y_2) \setminus (D \cup X_2) &= (C \cap (A \cup F)) \setminus (D \cup (B \cap E)) && \\
            &= (C \cap F) \setminus (D \cup (B \cap E)) && \text{(since $A \subseteq D$).} \\
            &= (C \cap F) \setminus (D \cup E) && \text{(since $C \subseteq B$).}
    \end{align*}

    \ref{itm:corners-when-crossing-two-nested-sepns-3} and~\ref{itm:corners-when-crossing-two-nested-sepns-4}.
    By symmetry, it suffices to show~\ref{itm:corners-when-crossing-two-nested-sepns-3}.
    We have $C \subseteq_\gamma E \subseteq A \cup E$ and $D \supseteq_\gamma F \supseteq B \cap F$.
    So, $t=(C,D) \leq_\gamma (A \cup E, B \cap F) = s \vee s'$.
    Similarly, with $C \subseteq B$ and $C \subseteq_\gamma E$ and $D \supseteq_\gamma F$ and $D \supseteq A$, we get $t = (C,D) \leq_\gamma (B \cap E, A \cup F) = s^* \wedge s'$.
\end{proof}

\begin{lem}\label{lem:shifting-preserves-interlacing}
    Let~$\sigma$ be a star of separations of a graph~$G$,
    let $p\in\sigma$ and let $\gamma\geq 0$ be an integer.
    If a separation~$s$ almost-$\gamma$-interlaces~$\sigma\sm\{p\}$,
    then both $s\vee p$ and $s\wedge p^*$ almost-$\gamma$-interlace~$\sigma$.
\end{lem}
\begin{proof}
    Let $t \in \sigma\sm\{p\}$.
    Since~$\sigma$ is a star, we have $p \leq t^*$.
    We have $t \leq_\gamma s$ or $t \leq_\gamma s^*$ since~$s$ almost-$\gamma$-interlaces~$\sigma \sm\{p\}$.
    Then \cref{lem:corners-when-crossing-two-nested-sepns}~\ref{itm:corners-when-crossing-two-nested-sepns-3} and~\ref{itm:corners-when-crossing-two-nested-sepns-4} imply $t \leq_\gamma s \vee p$ or $t \leq_\gamma (s\vee p)^*$.
    Similarly, we get $t \leq_\gamma s \wedge p^*$ or $t \leq_\gamma (s\wedge p^*)^*$.
    Moreover, we obviously have $p \leq s \vee p$ and $p \leq (s\wedge p^*)^*$.
    This proves the statement.
\end{proof}

The following observation is a direct consequence of \cref{approxLeQuasiTransitive}.

\begin{obs}\label{obs:totally-nested-close-separation}
    Let~$s$ and~$t$ be $(\leq k)$-separations of a graph~$G$.
    \begin{enumerate}
        \item\label{itm:totally-nested-close-separation-1} If~$s$ is totally-$\gamma$-nested and~$t$ is proper and $s \approx_\eps t$, then~$t$ is totally-$(\gamma+\eps)$-nested.
        \item\label{itm:totally-nested-close-separation-2} If~$s$ is narrow-$\gamma$-nested and~$t$ is narrow and proper and $s \approx_\eps t$, then~$t$ is narrow-$(\gamma+\eps)$-nested.
    \end{enumerate}
\end{obs}
\begin{proof}
    \ref{itm:totally-nested-close-separation-1}.
    Note that~$t$ is a proper $(\leq k)$-separation.
    By \cref{approxLeQuasiTransitive},~$t$ is $(\gamma+\eps)$-nested with every proper $(\leq k)$-separation of~$G$.
    Hence,~$t$ is totally-$(\gamma+\eps)$-nested.

    \ref{itm:totally-nested-close-separation-2}.
    Note that~$t$ is a narrow proper $(\leq k)$-separation.
    By \cref{approxLeQuasiTransitive},~$t$ is $(\gamma+\eps)$-nested with every narrow proper $(\leq k)$-separation of~$G$.
    Hence,~$t$ is narrow-$(\gamma+\eps)$-nested.
\end{proof}

\begin{lem}
\label{lem:almost-interlaces-then-nested-with-N}
    Let~$\eps \in \N$, let~$G$ be a graph, let~$N$ be a nested set of separations of~$G$ and let~$\sigma \subseteq N$ be a splitting star of~$N$.
    Let~$s$ be a separation of~$G$.
    If~$s$ almost-$\eps$-interlaces~$\sigma$, then~$s$ is $\eps$-nested with every separation in~$N$.
    Moreover, if~$s$ $\eps$-interlaces~$\sigma$, then for every~$n \in N$ we have $s <_\eps n$ or $s^* <_\eps n$ or $s <_\eps n^*$ or $s^* <_\eps n^*$.
\end{lem}
\begin{proof}
    Let~$n \in N$ be arbitrary.
    Since~$\sigma$ is a splitting star and $N$ is nested, there is some $n' \in \sigma$ such that $n \leq n'$ or $n^* \leq n'$, say the former.
    Since~$s$ almost-$\eps$-interlaces~$\sigma$, we have $n' \leq_\eps s$ or $n' \leq_\eps s^*$, say the former.
    Then, $n \leq n' \leq_\eps s$ implies $n \leq_\eps s$ by \cref{approxLeQuasiTransitive}.
    Hence,~$s$ is $\eps$-nested with every separation in~$N$.

    For the `moreover'-part, it suffices to replace `$\leq_\eps$' by `$<_\eps$' in the above argument.
\end{proof}

\subsection{Proof of \texorpdfstring{\cref{keylem:max-nested-set}}{Key Lemma 8.1}}

We start by analysing how a separation can interact with a star of narrow separations.

\begin{lemma}\label{lem:few-crossings-or-many-components}
    Let~$\ell \in \N$.
    Let~$S$ be the set of narrow proper $k$-separations of a graph~$G$.
    Let $\sigma \subseteq S$ be a star and let $s \in S$.
    Then either
    \begin{enumerate}
        \item\label{itm:few-crossings-or-many-components-1} $s$ crosses $< k + \ell$ many separations in~$\sigma$, or
        \item\label{itm:few-crossings-or-many-components-2} each proper side of~$s$ has $\geq \ell$ components.
    \end{enumerate}
\end{lemma}
\begin{proof}
    Let $\sigma =: \{(A_i,B_i) \mid i \in I\}$ and let~$s=:(C,D)$.
    Assume that~\ref{itm:few-crossings-or-many-components-1} fails, i.e.~$s$ crosses~$\geq k + \ell$ many separations in~$\sigma$.
    Let~$\sigma_1$ be the set of all $(A_i,B_i) \in \sigma$ such that the $sA_i$-link is nonempty.
    Since $A_j \sm B_j$ is disjoint from $A_{j'} \sm B_{j'}$ whenever $j \neq j'$, the star~$\sigma_1$ contains $\leq k$ separations.
    Let~$\sigma_2$ be the set of all $(A_i,B_i) \in \sigma \sm \sigma_1$ that cross~$s$.
    So, $|\sigma_2| \geq \ell$.
    By \cref{lem:link-empty-or-component-in-corners}, every separation~$(A_i,B_i) \in \sigma_2$ has separator $A_i \cap B_i = C \cap D$ and $G-(C \cap D)$ has a component in the $A_i C$-corner and a component in the $A_i D$-corner.
    Since the $(A_i \sm B_i)$'s are pairwise disjoint, each proper side of~$s$ has~$\geq \ell$ components.
\end{proof}

The next two lemmas turn almost-$\gamma$-interlacing into almost-interlacing by ``wiggling'' a separation to accommodate the members of the star one at a time. These adjustments preserve its order while keeping its distance from the original separation bounded.

\begin{lem}\label{lem:shift-to-interlace}
    Let $G$ be a $(q,k-1)$-unbreakable graph and let $\gamma \in \N$.
    Let $r$ be a $(q+\gamma+1)$-broad $k$-separation and let $p$ be a $(2q+k+1)$-broad $k$-separation.
    If $p\le_\gamma r$ or $p\le_\gamma r^*$, then there is a $k$-separation $r'\in\{r\vee p,\ r\wedge p^*\}$ that is $\gamma$-close to~$r$.
    In particular, $p\le r'$ or $p\le (r')^*$.
\end{lem}

\begin{proof}
    Suppose first that $\gamma<q+k$.
    Since $2q+k+1\geq q+\gamma+1$, both $p$ and $r$
    are $(q+\gamma+1)$-broad.
    After possibly reversing~$r$, assume that $p\leq_\gamma r$.
    By \cref{lem:meet-join-k-separation}, $r\vee p$ has order~$k$,
    and by \cref{lem:vee-wedge-is-close}, it is
    $\gamma$-close to~$r$.
    Undoing any reversal of~$r$ gives the required separation
    in $\{r\vee p,\ r\wedge p^*\}$.

    Hence we may assume that $\gamma\geq q+k$.
    Write $r\eqqcolon(A,B)$ and $p\eqqcolon(C,D)$.
    By symmetry, we may orient $r$ so that $p \leq_{\gamma} r$.
    We first show that we can further choose this orientation so that the $AC$-corner has size $>q$.

    Suppose that the $AC$-corner has size at most $q$.
    Then, $|A \sm D| = |A \cap (C \sm D)| \leq |A \cap B| + |(A \sm B) \cap (C \sm D)| \leq k+q \leq \gamma$, and similarly $|C \sm B| \leq \gamma$ so $p \leq_{\gamma}r^*$.
    Since $p$ is $(2q+k+1)$-broad and $r$ has order $k$, we also have that the $BC$-corner has size $>q$.
    Thus, replacing $r$ by $r^*$ gives the desired orientation.

    Now, let $r' \coloneqq r \vee p$.
    Since $p\le_\gamma r$,~$r'$ is $\gamma$-close to~$r$ by \cref{lem:vee-wedge-is-close}.
    Since $r$ is $(q+\gamma+1)$-broad, \cref{LowBreadthClose} implies that $r'$ has breadth $> q$.
    The opposite corner-separation $r \wedge p$
    is also $(> q)$-broad: its first proper side is the $AC$-corner and its second proper side contains $B \sm A$.

    Since $G$ is $(q,k-1)$-unbreakable, it follows that both corner-separations have order at least $k$.
    By submodularity, they both have order exactly $k$.
\end{proof}

\begin{lem}\label{lem:wiggle-lemma}
    Let~$G$ be a $(q,k-1)$-unbreakable graph,
    $\gamma,\ell\in\Nbb$ with $\ell>q$.
    Let~$\sigma$ be a finite star of $(2q+k+1)$-broad
    $k$-separations of~$G$.
    Let $\beta:=\ell+\gamma\cdot|\sigma|$.
    Let~$s$ be a $\beta$-broad $k$-separation of~$G$
    that almost-$\gamma$-interlaces~$\sigma$.
    Then there is an $\ell$-broad $k$-separation~$s'$
    that almost-interlaces~$\sigma$ and is
    $(\gamma\cdot|\sigma|)$-close to~$s$.

    Moreover, if $\tau\supseteq\sigma$ is a star
    and~$s$ almost-interlaces $\tau\sm\sigma$,
    then~$s'$ can be chosen to almost-interlace~$\tau$.
\end{lem}
\begin{proof}
    We prove the assertion involving~$\tau$ by induction on~$|\sigma|$.
    The first assertion then follows by taking $\tau=\sigma$.
    If $\sigma=\emptyset$, then we take $s':=s$.
    Since $\beta=\ell$,~$s'$ is $\ell$-broad and almost-interlaces~$\tau=\tau\sm\sigma$ and is 0-close to~$s$.

    Now suppose that $\sigma\neq\emptyset$.
    We choose $p\in\sigma$ and define $\sigma'=\sigma\sm\{p\}$.
    Since $\beta\geq\ell+\gamma>q+\gamma$,
    \cref{lem:shift-to-interlace} gives a $k$-separation $r\in\{s\vee p,\ s\wedge p^*\}$ that is $\gamma$-close to~$s$.
    By \cref{LowBreadthClose},~$r$ is $(\beta-\gamma)$-broad, where $\beta-\gamma=\ell+\gamma\cdot|\sigma'|$.
    Moreover, $p\leq r$ or $p\leq r^*$.
    By \cref{lem:shifting-preserves-interlacing}, applied to~$\sigma$ (with parameter~$\gamma$) and to $\tau\sm\sigma'$ (with parameter~$0$),~$r$ almost-$\gamma$-interlaces~$\sigma'$ and almost-interlaces~$\tau\sm\sigma'$.

    The induction hypothesis applied to~$r$, $\sigma'$ and~$\tau$ now yields an $\ell$-broad $k$-separation~$s'$
    that almost-interlaces~$\tau$ and satisfies
    \[
        s'\approx_{\gamma|\sigma'|}r\approx_\gamma s~.
    \]
    By \cref{approxLeQuasiTransitive}, $s'\approx_{\gamma|\sigma|}s$, as required.
\end{proof}

Using the previous results, we now prove that if a broad $k$-separation~$s$ almost-$\gamma$-interlaces a star~$\sigma$ of broad separations, either~$s$ is close to a $k$-separation that almost-interlaces~$\sigma$, or the broad petal-separations in the graphic-$\sep(s)$-anemone are narrow-nested.

\begin{lem}\label{AlphaInterlacesCloseToInterlaces}
    Let~$G$ be a $(q,k-1)$-unbreakable $(b,k)$-narrow 
    graph, $\gamma\in\N$ and
    \begin{align*}
        \beta&:=\beta_{\ref{AlphaInterlacesCloseToInterlaces}}(q,k,b,\gamma) := \max\{\, q+1+\lfloor 3k/2\rfloor\cdot\gamma\,,\,b\,\}~\text{ and }\\
        \beta'&:=\beta'_{\ref{AlphaInterlacesCloseToInterlaces}}(q,k,b) := \max\{\,2q+k+1\,,\,b\,\}~\text{ and }\\ \eps&:=\eps_{\ref{AlphaInterlacesCloseToInterlaces}}(k,\gamma):=\lfloor 3k/2\rfloor\cdot \gamma~.
    \end{align*}
    Let~$\sigma$ be a star of $\beta'$-broad $k$-separations of~$G$. 
    Let~$s$ be a $\beta$-broad $k$-separation that almost-$\gamma$-interlaces~$\sigma$.
    Then we have one of the following two outcomes:
    \begin{enumerate}
        \item\label{AlphaInterlacesCloseToInterlaces-1} $s$ is $\eps$-close to some $k$-separation~$s'$ that almost-interlaces~$\sigma$.
        \item\label{AlphaInterlacesCloseToInterlaces-2} $s$ is a displayed separation of the graphic-$X$-anemone~$\cA$ with $X:=\sep(s)$, and the narrow proper petal-separations of~$\cA$ are narrow-nested.
        Moreover, $s$ is not half-connected.
    \end{enumerate}
\end{lem}

\begin{proof}
    The graph~$G$ is $(b, k)$-narrow and every separation in $\sigma$ is $b$-broad, so~$\sigma$ is a star of narrow $k$-separations of $G$. Similarly, $s$ is a narrow $k$-separation of $G$.
    Let $\ell := \lceil(k+1)/2\rceil$.
    We distinguish the following two cases. 
    
    \casen{Case~\ref{itm:few-crossings-or-many-components-1}: $s$ crosses $<k+\ell$ many separations in~$\sigma$.}
    If $s$ is $\gamma$-close to some $p\in\sigma\cup\sigma^*$,
    then $p$ almost-interlaces~$\sigma$, so taking $s':=p$ gives outcome~\ref{AlphaInterlacesCloseToInterlaces-1}, since $\gamma\leq\eps$. Hence assume that no such~$p$ exists.
    
    Let~$\sigma'$ consist of the separations in~$\sigma$ that cross~$s$.
    Then $|\sigma'|\leq k+\ell-1\leq \lfloor 3k/2\rfloor$.
    We show that $s$ almost-interlaces~$\sigma\sm\sigma'$.
    Otherwise, some $p\in\sigma\sm\sigma'$ satisfies neither $p\leq s$ nor $p\leq s^*$.
    Thus, $s \leq p$ or $s^* \leq p$.
    Without loss of generality, assume the former.
    Since $p\leq_\gamma s$ or $p\leq_\gamma s^*$, either $s\approx_\gamma p$, which contradicts our assumption, or $s\leq_\gamma s^*$, which contradicts that~$s$ is $\beta$-broad by \cref{lem:assymmetric-nested-implies-not-broad} and $\beta>\gamma$.

    Every member of~$\sigma'$ is $(2q+k+1)$-broad, and~$s$ is $\beta$-broad with $\beta\geq q+1+\gamma|\sigma'|$.
    So we can apply \cref{lem:wiggle-lemma} to~$\sigma'$ (with $\ell=q+1$ and $\tau=\sigma$).
    Its `Moreover'-part yields a $k$-separation~$s'$ that almost-interlaces~$\sigma$ and is $(\gamma|\sigma'|)$-close to~$s$.
    Since $\gamma|\sigma'|\leq \lfloor 3k/2\rfloor\cdot \gamma = \eps$, this gives outcome~\ref{AlphaInterlacesCloseToInterlaces-1}.

    \casen{Case~\ref{itm:few-crossings-or-many-components-2}: not Case~\ref{itm:few-crossings-or-many-components-1}.} By \cref{lem:few-crossings-or-many-components},
    each proper side of~$s$ has~$\geq \ell$ components.
    Let~$X:=\sep(s)$, and let~$\Acal$ be the graphic-$X$-anemone.
    Since $\sep(s)=X$, the separation~$s$ is a displayed separation of~$\Acal$, which is the first part of \ref{AlphaInterlacesCloseToInterlaces-2}. 
    Let $C_1,C_2,\ldots,C_m$ be the connected components of $G-X$.
    With $\ell = \lceil(k+1)/2\rceil$, we have $m \geq 2\ell \geq k+1$.
    Then, $s$ is not half-connected, otherwise \cref{lem:many-components-implies-tot-nested} would imply that $s$ is nested with every separation in $\sigma$, contradicting that we are not in Case~\ref{itm:few-crossings-or-many-components-1}.
    We define $s_i:=(V(C_i) \cup X,V(G) \sm V(C_i))$, which is a half-connected $k$-separation of~$G$ with separator $\sep(s_i)=X=\sep(s)$.
    Every narrow proper petal-separation of $\cA$ is some $s_i$.
    By \cref{lem:many-components-implies-tot-nested}, no narrow proper~$s_i$ is crossed by a narrow proper $(\leq k)$-separation of~$G$.
    Hence, every narrow proper~$s_i$ is narrow-nested, which yields the second part of outcome \ref{AlphaInterlacesCloseToInterlaces-2} and thus completes the proof.
\end{proof}

\begin{lem}\label{lem:two-opposite-corners}
    Let~$s$ and~$t$ be two $\beta$-broad $(\leq k)$-separations of a $(q,k-1)$-unbreakable graph~$G$, where $\beta > 2q+k$.
    Then there are two opposite corner-separations of order~$k$ and breadth~$\geq (\beta-k)/2$.
\end{lem}
\begin{proof}
    Let~$s=:(A,B)$ and~$t=:(C,D)$.
    First, we show that there are two opposite corners each of size~$\geq (\beta-k)/2$.
    Assume the contrary.
    Then either the $AC$-corner or the $BD$-corner has size $<(\beta-k)/2$, say the former.
    Similarly, either the $AD$-corner or the $BC$-corner has size $<(\beta-k)/2$, say the former.
    Then $|A \sm B| < (\beta-k)/2 + (\beta-k)/2 + |\sep(t)| \leq \beta$, contradicting that~$s$ is $\beta$-broad.

    Hence, there are two opposite corners each of size~$\geq (\beta-k)/2$.
    So, their corresponding corner-separations have breadth $\geq (\beta-k)/2 > q$.
    Since~$G$ is $(q,k-1)$-unbreakable, these separations have order~$\geq k$.
    By submodularity, they have order exactly~$k$.
\end{proof}

We shall use the following variant of {\cite[Lemma 4.4]{CG14:isolatingblocks}}; we include a proof here for completeness.

\begin{lem}\label{lem:narrow-half-connected-and-link-empty}
    Whenever two narrow proper separations cross with one empty link, they have the same separator.
    In particular, none of them can be half-connected.
\end{lem}
\begin{proof}
    Let~$s=(A,B)$ and $t=(C,D)$ be two narrow proper separations that cross with one empty link, say the $sC$-link.
    First, we show that there are two adjacent empty links.
    Suppose for a contradiction that no two adjacent links are empty.
    Then, the $At$-link and the $Bt$-link are nonempty, say they contain vertices $v$ and $w$, respectively.
    Since~$t$ is proper, there exists a connected component $K$ of $G - \sep(t)$ such that $V(K) \subseteq C \sm D$. Since~$t$ is narrow, $N(K) = \sep(t)$ so there exists a path between $v$ and $w$ all of whose internal vertices are in $C \sm D$.
    Since $v \in A \sm B$ and $w \in B \sm A$, this path intersects $\sep(s)$ at some vertex $u \in C \sm D$. This contradicts that the $sC$-link is empty.
    
    Hence there are two adjacent empty links, say the $At$-link and the $sC$-link. 
    Since~$s$ and~$t$ cross, the $AC$-corner is nonempty.
    Let $K$ be a component of the graph induced by the $AC$-corner.
    Now $K$ is a component of $G - \sep(s)$ and of $G - \sep(t)$.
    Since $s$ and $t$ are narrow, we have $\sep(s) = N(K) = \sep(t)$.
    
    To see the \lq In particular\rq-part, note that half-connected separations are nested with all separations with the same separator. 
\end{proof}

\begin{lem}
\label{lem:separation-with-three-components-crosses}
    Let $q,k,\beta\in \N$ with $\beta \geq \beta_{\ref{lem:separation-with-three-components-crosses}}(q,k) := 4q+k+1$, and let $G$ be a $(q,k-1)$-unbreakable graph with a $\beta$-broad narrow $(\leq k)$-separation~$s$ such that $G-\sep(s)$ has at least three components.
    Then every $\beta$-broad narrow $(\leq k)$-separation that crosses~$s$ has the same separator as~$s$.
    
    In particular, if~$s$ is half-connected then~$s$ is $S$-totally-nested, where~$S$ is the set of $\beta$-broad narrow $(\leq k)$-separations of~$G$.
\end{lem}

\begin{proof}
    Let $s=:(A,B)$ and let $t=:(C,D)$ be a $\beta$-broad narrow $(\leq k)$-separation that crosses~$s$.
    By symmetry assume that $G[A \sm B]$ has at least two components.
    By \cref{lem:two-opposite-corners} (with $\beta > 2q+k$), there are two opposite corner-separations of order~$k$ and breadth $\geq \lceil (\beta-k)/2 \rceil \geq 2q+1$.
    Without loss of generality, these are the $AC$-corner-separation and the $BD$-corner-separation.

    First we show that the corner-diagram of $s,t$ has an empty link.
    Suppose for a contradiction that all four links are nonempty.
    Thus since~$s$ is narrow, every component of $G[A \sm B]$ intersects the $At$-link.
    Let~$K$ be a component of $G[A \sm B]$ whose intersection with the $AC$-corner has minimal size.
    We define $s':=(A',B'):=(A \sm K, B \cup K)$, which is a separation of~$G$ with $\sep(s')=\sep(s)$.
    Note that since~$K$ intersects the $At$-link, the $A'C$-corner-separation has order~$<k$.
    Since $G[A \sm B]$ has at least two components and since we chose~$K$ minimal, the $A'C$-corner still contains~$\geq \lceil (2q + 1) / 2 \rceil > q$ vertices. This is a contradiction to the assumption that $G$ is $(q,k-1)$-unbreakable.

    Hence the corner diagram of $s$ and $t$ has an empty link, so $s$ and $t$ have the same separator by \cref{lem:narrow-half-connected-and-link-empty}, as desired. 
    The \lq In-particular\rq-part follows from the \lq In particular\rq-part of \cref{lem:narrow-half-connected-and-link-empty}.
\end{proof}

\begin{lemma}\label{lem:max-nested-set}
    Let~$G$ be a $(q,k-1)$-unbreakable $(b,k)$-narrow graph and let $\gamma \in \N$.
    Let 
    \begin{align*}
        \delta&:=\delta_{\ref{lem:max-nested-set}}(k,b,\gamma) := \max\{\lfloor 3k/2\rfloor\cdot \gamma\, , \, 2(\gamma+k+b)\} \quad \text{and}\\
        \beta&:=\beta_{\ref{lem:max-nested-set}}(q,k,b):= \max\{\, \beta'_{\ref{AlphaInterlacesCloseToInterlaces}}(q,k,b)\, , \, \beta_{\ref{lem:separation-with-three-components-crosses}}(q,k)\,\} \quad\text{and}\\
        \eps&:=\eps_{\ref{lem:max-nested-set}}(q,k,b,\gamma):= \delta + \beta~.
    \end{align*}
    Let~$N$ be a nested set of half-connected $\beta$-broad $k$-separations of~$G$.
    Let~$s$ be an $\eps$-broad $k$-separation of~$G$ that is $\gamma$-nested with every separation in~$N$.
    Then there is a $k$-separation~$s'$ of~$G$ that is nested with~$N$ and $\delta$-close to~$s$.
\end{lemma}
\begin{proof}
    Without loss of generality, we can assume that~$N$ is symmetric.

    Recall that~$s$ is $\gamma$-nested with every separation in~$N$.
    By \cref{lem:gamma-close-or-gamma-interlaces} (which we can apply since $s$ has breadth~$\geq \eps =\delta+\beta \geq 2\gamma+1$), either~$s$ is~$\gamma$-close to a separation in~$N$ or~$s$ $\gamma$-interlaces some splitting star~$\sigma$ of~$N$. 
    We are done in the former case (since $\delta \geq \gamma$), so we can assume the latter case.
    Observe that~$s$ has breadth $\eps \geq \beta_{\ref{AlphaInterlacesCloseToInterlaces}}(q,k,b,\gamma)$ and $\delta \geq \eps_{\ref{AlphaInterlacesCloseToInterlaces}}(k, \gamma)$.

    Note that every separation in~$\sigma$ is $\beta$-broad, hence
    $\beta'_{\ref{AlphaInterlacesCloseToInterlaces}}(q,k,b)$-broad.
    So, we can apply \cref{AlphaInterlacesCloseToInterlaces} to the separation~$s$ and the star~$\sigma$, and have one of the following two outcomes.

    \casen{Case~\ref{AlphaInterlacesCloseToInterlaces-1}: $s$ is $\delta$-close to some $k$-separation~$s'$ that almost-interlaces~$\sigma$.}
    Since~$s'$ almost-interlaces~$\sigma$, it is nested with every separation in~$N$ by \cref{lem:almost-interlaces-then-nested-with-N} and we are done.

    \casen{Case~\ref{AlphaInterlacesCloseToInterlaces-2}: $s$ is a displayed separation of the graphic-$X$-anemone~$\cA$ with $X:=\sep(s)$, and the narrow proper petal-separations of~$\cA$ are narrow-nested.}
    Let~$p_1,p_2,\ldots,p_m$ be the petal-separations of~$\cA$.
    By \cref{FlowerClosePetalOrCrossed} applied to~$\cA$ and~$s$ with $\eps':=\gamma+k+b$, we have one of the following two cases.

    \casen{Case~\ref{AlphaInterlacesCloseToInterlaces-2}~\ref{FlowerClosePetalOrCrossed:close}: $s$ is $2\eps'$-close to some~$p_i$ or its inverse.}
    Since $\delta \geq 2(\gamma+k+b) = 2\eps'$, this separation is $\delta$-close to~$s$. 
    Moreover, the separation~$p_i$ has breadth $\geq \eps - 2\eps' \geq \eps - \delta \geq \beta$ by \cref{LowBreadthClose}.
    Since~$p_i$ is $\beta$-broad, it is $b$-broad, hence narrow since $G$ is $(b, k)$-narrow.
    Thus, $p_i$ is narrow-nested by assumption.
    Furthermore, all separations in $N$ are $b$-broad, and $G$ is $(b, k)$-narrow, so all separations in~$N$ are narrow proper $(\leq k)$-separations of $G$.
    Therefore, $p_i$ is nested with~$N$, and thus $p_i^*$ is also nested with~$N$.
    We put~$s':=p_i$ or~$s':=p_i^*$ and are done.

    \casen{Case~\ref{AlphaInterlacesCloseToInterlaces-2}~\ref{FlowerClosePetalOrCrossed:crossed}: $\cA$ has a displayed separation that $(\eps'-k)$-crisses~$s$.}
    So,~$G-\sep(s)$ has at least three components.
    Moreover,~$s$ is $\beta_{\ref{lem:separation-with-three-components-crosses}}(q,k)$-broad and narrow (since it is $b$-broad and~$G$ is $(b,k)$-narrow).
    Similarly, every separation in~$N$ is $\beta_{\ref{lem:separation-with-three-components-crosses}}(q,k)$-broad and narrow.
    By \cref{lem:separation-with-three-components-crosses} and since every separation in~$N$ is half-connected,~$s$ is nested with~$N$.
    We put $s':=s$ and are done.
\end{proof}

We can now prove \cref{keylem:max-nested-set}, which we restate for convenience.
\maxnestedset*

\begin{proof}[Proof of \cref{keylem:max-nested-set}]
    First we show that~$N$ $\eps$-approximates~$S$.
    Let~$s \in S$ be arbitrary.
    If $s$ is not $\eps$-broad there is nothing to show, so we can assume that $s$ is $\eps$-broad.
    The graph~$G$ is $(b, k)$-narrow and all separations in $N$ are $b$-broad $k$-separations, so they are narrow proper $(\leq k)$-separations of~$G$.
    Since~$s$ is narrow-$\gamma$-nested, it is $\gamma$-nested with every separation in~$N$.
    By \cref{lem:max-nested-set}, there is a $k$-separation~$s'$ of~$G$ that is nested with~$N$ and $\delta_{\ref{lem:max-nested-set}}(k, b, \gamma)$-close to~$s$.
    By \cref{LowBreadthClose}, $s'$ has breadth at least $\eps - \delta_{\ref{lem:max-nested-set}}(k, b, \gamma) \geq \eps - \delta \geq \beta \geq b$, so $s'$ is narrow.
    By \cref{obs:totally-nested-close-separation},~$s'$ is narrow-$(\delta_{\ref{lem:max-nested-set}}(k, b, \gamma)+\gamma)$-nested.
    This implies that no $k$-separation of~$G$ $(\delta_{\ref{lem:max-nested-set}}(k, b, \gamma)+\gamma)$-crosses~$s'$.
    Indeed, if a $k$-separation $t$ of $G$ $(\delta_{\ref{lem:max-nested-set}}(k, b, \gamma)+\gamma)$-crosses~$s'$ then, using that $\delta_{\ref{lem:max-nested-set}}(k, b, \gamma) \geq b + k$, we have that $t$ is $b$-broad by \cref{LargeCrossedImpliesBroad}.
    But then $t$ is narrow since $G$ is $(b, k)$-narrow, which contradicts that $s'$ is narrow-$(\delta_{\ref{lem:max-nested-set}}(k, b, \gamma)+\gamma)$-nested.

    Let $\cA$ be the graphic-$\sep(s')$-anemone.
    By \cref{FlowerClosePetalOrCrossed}, there is a separation $p \in \pi(\cA) \cup \pi(\cA)^*$ that is $2(\delta_{\ref{lem:max-nested-set}}(k, b, \gamma)+\gamma+k)$-close to $s'$.
    If $s'$ is half-connected, we choose $p=s'$.
    Note that~$p$ is half-connected and $\delta$-close to $s \in S$.
    Since~$s$ is $\eps$-broad and $\eps\geq \delta+\beta$, we have that~$p$ is $\beta$-broad by \cref{LowBreadthClose}.
    Since~$\beta > q$ and~$G$ is $(q,k-1)$-unbreakable, the separation~$p$ has order~$k$.
    Since $\beta \geq b$ and~$G$ is $(b,k)$-narrow, the separation~$p$ is narrow.
    By \cref{obs:totally-nested-close-separation} and with $\gamma'= \delta+\gamma$,~$p$ is narrow-$\gamma'$-nested.
    Hence,~$p \in S'$.
    If $s'$ is half-connected, then $p=s'$; so~$p$ is nested with~$N$.
    Otherwise, $G-\sep(s') = G-\sep(p)$ has at least three components, so by \cref{lem:separation-with-three-components-crosses} and since $\beta \geq \beta_{\ref{lem:separation-with-three-components-crosses}}(q,k)$, the separation~$p$ is nested with~$N$.
    Since~$N$ is maximal, we have $p \in N$.
    This completes the proof that~$N$ $\eps$-approximates~$S$.

    We now show the `Moreover'-part.
    Let~$s$ be an $\eps'$-broad $k$-separation of~$G$.
    So~$s$ is $b$-broad and hence narrow since~$G$ is $(b,k)$-narrow.
    Since every separation in~$N$ is narrow-$\gamma'$-nested,~$s$ is $\gamma'$-nested with every separation in~$N$.
    By \cref{lem:max-nested-set} there is a $k$-separation~$s'$ of~$G$ that is nested with~$N$ and $\delta'$-close to~$s$.
\end{proof}

\section{The structure of the torsos}\label{sec:TorsoStructure}

In this section, we locate the flowers within the tree-decomposition. We show that asterisks of sufficiently large spread determine distinct splitting stars and that each such star and its corresponding petal-star approximate each other. We then prove that the remaining torsos have higher connectivity.

This enables us to introduce and prove a technical strengthening (\cref{keylemma:FlowerFinding} and \cref{keylemma:FlowerCriterion}) of the main result (\cref{mainTutte}), under two additional assumptions:
On the one hand, we assume that the graph~$G$ we are to decompose is $(b,k)$-narrow.
On the other hand, we assume that virtually all key lemmas from the previous sections have been combined, with suitably chosen constants, to create a setting (\ref{set:FlowerFinding}) in which we will work for the vast majority of this section.
Assuming \cref{set:FlowerFinding} helps us keep the paper fairly modular.
Later, in \cref{sec:main-theorem-for-narrow} we will show how to establish \cref{set:FlowerFinding} using the key lemmas from the previous sections.

\subsection{Preparations}

We start by establishing an equivalence between totally-nested separations and totally-nested separators in~\cref{lem:equivalence-totally-nested-set}.
Moreover, we strengthen \cref{AlphaInterlacesCloseToInterlaces} in the special case that the given star~$\sigma$ of separations is a petal-star of a graphic-flower (see \cref{NewerInterlacingPetalStar}).

Let~$G$ be a graph and let~$X \subseteq V(G)$.
Let~$S$ be a set of separations of~$G$.
We say that~$X$ is \defn{$S$-totally-nested} if no separation in~$S$ separates two vertices of~$X$. 

\begin{lem}\label{lem:equivalence-totally-nested-set}
    Let~$G$ be a graph, and let $X \subseteq V(G)$.
    Let~$S$ be the set of all narrow proper $k$-separations of~$G$, and let~$S_X \subseteq S$ be the set of half-connected separations in~$S$ with separator~$X$.  
    Assume $S_X \neq \emptyset$.
    Then the following assertions are equivalent:
    \begin{enumerate}
        \item\label{itm:equivalence-totally-nested-set-1} $X$ is $S$-totally-nested;
        \item\label{itm:equivalence-totally-nested-set-2} every separation in~$S_X$ is $S$-totally-nested.
        \item\label{itm:equivalence-totally-nested-set-3} some separation in~$S_X$ is $S$-totally-nested.
    \end{enumerate}
\end{lem}
\begin{proof}
    \ref{itm:equivalence-totally-nested-set-1}~$\Rightarrow$~\ref{itm:equivalence-totally-nested-set-2}. 
    Let~$X$ be $S$-totally-nested.
    Assume for a contradiction that there is a half-connected narrow proper separation~$s=(A,B)$ with separator $X=A \cap B$ that is crossed by some separation~$t=(C,D) \in S$.
    Assume without loss of generality that $G[A \sm B]$ is connected.

    First, we show that the $At$-link is nonempty.
    Suppose that it is empty.
    Since~$G[A \sm B]$ is connected, the $AC$-corner or the $AD$-corner is empty, say the former.
    Then the $AD$-corner is nonempty and there is a component of $G-X$ contained in the $AD$-corner.
    Since~$s$ is narrow, the $sC$-link is empty.
    By the fact that the $At$-link, the $sC$-link and the $AC$-corner are empty, it follows that~$s$ and~$t$ are nested, a contradiction.
    Hence the $At$-link is nonempty.
    
    Then, since~$s$ and~$t$ have the same order, either the $sC$-link or the $sD$-link is nonempty, say the $sC$-link.
    Since~$X$ is $S$-totally-nested, the $sD$-link is empty.
    Since~$s$ and~$t$ cross, either the $Bt$-link is nonempty or the $BD$-corner is nonempty.
    If the $Bt$-link is nonempty, every component of $G[D \sm C]$ has neighbourhood of size~$<k$, contradicting that~$t$ is narrow.
    If the $Bt$-link is empty, there is a component of $G[D \sm C]$ contained in the $BD$-corner. This component has a neighbourhood of size~$<k$, contradicting that~$t$ is narrow.

    \ref{itm:equivalence-totally-nested-set-2}~$\Rightarrow$~\ref{itm:equivalence-totally-nested-set-3}. This holds since $S_X \neq \emptyset$.

    \ref{itm:equivalence-totally-nested-set-3}~$\Rightarrow$~\ref{itm:equivalence-totally-nested-set-1}.
    Suppose that~$X$ is not $S$-totally-nested.
    Let~$t=(C,D) \in S$ be a separation that separates two vertices of~$X$.
    Let~$s$ be any narrow proper half-connected separation with separator~$X$.
    Since the $sC$-link and the $sD$-link are both nonempty, the separations~$s$ and~$t$ cross.
    Hence, no separation in~$S_X$ is $S$-totally-nested.
\end{proof}

\begin{lem}
    \label{lem:make-a-petal-star-broad}
    Let~$G$ be a graph and let~$\beta \in \Nbb$.
    Let~$F$ be a graphic-$k$-flower in~$G$ that has a $\beta$-broad displayed separation.
    Then there is a $k$-flower~$F'$ such that the following holds:
    \begin{enumerate}
        \item\label{itm:make-a-petal-star-broad-0} every separation induced by $F'$ is induced by $F$;
        \item\label{itm:make-a-petal-star-broad-1} every petal-separation of~$F'$ is $\beta$-broad;
        \item\label{itm:make-a-petal-star-broad-2} for every petal-separation~$t' \eqqcolon (C, D)$ of~$F'$ such that $|C \sm D| \geq \beta+2\eps$, there is a $\beta$-broad petal-separation~$t$ of~$F$ which is $\eps$-close to~$t'$, where $\eps:=\eps_{\ref{lem:make-a-petal-star-broad}}(\beta,k):=\beta+k$; and
        \item\label{itm:make-a-petal-star-broad-3} if a separation~$s$ of~$G$ almost-$\gamma$-interlaces the petal-star~$\pi$ of~$F$ for some $\gamma \in \Nbb$, then~$s$ almost-$\gamma'$-interlaces the petal-star of~$F'$ for $\gamma':=\gamma'_{\ref{lem:make-a-petal-star-broad}}(\beta,k,\gamma):=\max\{2\eps+\beta+k, \eps+\gamma\}$.
    \end{enumerate}
\end{lem}

\begin{proof}
    Let $t$ be a $\beta$-broad displayed separation of $F$.
    By \cref{AnemoneAsDaisy}, there is a cycle-decomposition $\cO=:(O,\cV)$ of $G$ which displays $t$, such that every displayed separation of $\cO$ is a displayed separation of $F$, and every petal-separation of $\cO$ is a petal-separation of $F$.
    Let~$B \subseteq V(O)$ be the set of nodes of~$O$ whose corresponding petal-separation is $\beta$-broad.
    Our strategy is to `merge' nodes of~$O$ that are not in~$B$ with adjacent nodes in order to obtain the desired $k$-flower~$F'$, as follows.
    
    We call a path $Q \subseteq O$ \emph{good} if the separation induced by~$(Q,O-Q)$ is $\beta$-broad, and~$Q$ does not have a proper subpath~$Q'$ such that the separation induced by $(Q',O-Q')$ is $\beta$-broad.
    Note that every good path either consists of a single node of~$B$ or is contained in $O-B$.
    
    Let~$\cQ$ be a maximal set of pairwise disjoint good paths.
    Since~$\cO$ displays the $\beta$-broad separation~$t$, and since~$t^*$ is a $\beta$-broad displayed separation of~$\cO$ too, we have $|\cQ| \ge 2$.
    Then, we can write $O=P_0Q_0P_1Q_1\ldots P_{m-1}Q_{m-1}$ with $m := |\cQ| \geq 2$, such that $V(P_i) \subseteq V(O-\bigcup \mathcal{Q})$ where $P_i$ is a (possibly empty) path and $\{\,Q_i:i\in \Z_m\,\}=\cQ$. 
    We define $F':=(\Zbb_m, \cV')$ with $\cV':=(V'_o : o \in \Zbb_m)$ where $V'_o:=V \rest (P_o \cup Q_o)$ (that is, for each~$o$ we `merge' the bags and nodes of~$\cO$ that correspond to the path $P_o \cup Q_o$).
    Observe that~$F'$ is a $k$-flower.
    We show that~$F'$ satisfies~\ref{itm:make-a-petal-star-broad-0}--\ref{itm:make-a-petal-star-broad-3}.

    \ref{itm:make-a-petal-star-broad-0}. By construction, every separation induced by $F'$ is induced by $\cO$, hence by $F$. 

    \ref{itm:make-a-petal-star-broad-1}. This follows from $m \geq 2$.

    \ref{itm:make-a-petal-star-broad-2}. Let~$t' \eqqcolon (C, D)$ be a petal-separation of~$F'$ such that $|C \sm D| \geq \beta+2\eps$, corresponding to some node $o \in \Zbb_m$, so $C = V'_o$.
    First, we show that~$Q_o$ consists of a single node of~$B$.
    Suppose not for a contradiction, so $Q_o$ is a path with at least two vertices.
    Then, there are nonempty paths~$Q_{o,1},Q_{o,2}\se Q_o$ such that $V(Q_o)=V(Q_{o,1}) \sqcup V(Q_{o,2})$.
    Since $Q_o$ is good, the separation induced by~$(Q_{o,i},O-Q_{o,i})$ has breadth~$<\beta$ (on the left side, for $i \in \{1,2\}$).
    Furthermore,~$F$ is a graphic-$k$-flower and every displayed separation of $\cO$ is a displayed separation of $F$, so the adhesion-set of the edge of~$Q_o \subseteq O$ joining~$Q_{o,1}$ and~$Q_{o,2}$ has size~$\leq k$.
    Thus, the separation induced by $(Q_o,O-Q_o)$ has breadth $< 2\beta+k$ (on the left side).
    Moreover, the separation induced by $(P_o,O-P_o)$ has breadth~$<\beta$ (on the left side), since otherwise we could find some good subpath of~$P_o$, contradicting the maximality of~$\cQ$.
    Then it similarly follows that the separation induced by $(P_oQ_o,O-P_oQ_o)$ has breadth (on the left side)~$<(2\beta+k)+\beta+k=3\beta+2k = \beta + 2\eps$ (where we use $\eps=\eps_{\ref{lem:make-a-petal-star-broad}}(\beta,k)=\beta+k$), a contradiction with $|C \sm D| \geq \beta + 2\eps$.

    Thus,~$Q_o$ consists of a single node~$b \in B$.
    Recall that the separation induced by $(P_o,O-P_o)$ has breadth~$< \beta$ (on the left side) and the adhesion-set of the edge of~$O$ joining~$P_o$ and~$Q_o$ has size~$\leq k$, as before.
    Then it follows that the petal-separation~$t'$ is $\eps$-close to the petal-separation of~$\cO$ (hence of $F$) corresponding to the node~$b$.

    \ref{itm:make-a-petal-star-broad-3}.
    Let~$s$ be a separation of~$G$ that almost-$\gamma$-interlaces the petal-star $\pi$ of~$F$ for some $\gamma \in \Nbb$.
    Let~$t' \eqqcolon (A', B')$ be an arbitrary petal-separation of~$F'$.
    We show that we have $t' \leq_{\gamma'} s$ or $t' \leq_{\gamma'} s^*$.
    If $|A' \sm B'| < \gamma' - k$, then this is trivially satisfied.
    Hence, suppose that~$|A' \sm B'| \geq \gamma'-k \geq \beta+2\eps$.
    
    By \ref{itm:make-a-petal-star-broad-2}, there is a petal-separation~$t$ of~$F$ which is $\eps$-close to~$t'$.
    Since~$s$ almost-$\gamma$-interlaces~$\pi$, we have $t' \approx_\eps t \leq_\gamma s$ or $t' \approx_\eps t \leq_\gamma s^*$.
    By \cref{approxLeQuasiTransitive}, we have $t' \leq_{\gamma+\eps} s$ or $t' \leq_{\gamma+\eps} s^*$.
    So,~$s$ almost-$\gamma'$-interlaces the petal-star of~$F'$.
\end{proof}

\begin{lem}\label{lem:broad-almost-interlaces-then-broad-displayed}
    Let $G$ be a graph, let $F$ be a graphic-$k$-flower of $G$ and let~$\pi$ be the petal-star of~$F$.
    Let $\beta', \gamma \in \N$, and let $\beta \geq \beta_{\ref{lem:broad-almost-interlaces-then-broad-displayed}}(k, \beta', \gamma) \coloneqq \max\{\beta' + \gamma, 2(\beta'+k)\}$.
    Suppose that there exists a~$\beta$-broad $k$-separation $s$ that almost-$\gamma$-interlaces $\pi$.
    Then, $F$ has a $\beta'$-broad displayed separation.
\end{lem}

\begin{proof}
    Write $F \eqqcolon (H, \cV)$ and $s \eqqcolon (A, B)$.
    We distinguish the following two cases.
    
    \casen{Case 1:} There exists a petal~$V_h$ with $h \in V(H)$ such that $|V_h|\geq\beta'+k$.
    Let $p \eqqcolon (C, D) \in \pi$ be the petal-separation at $h$, with $C = V_h$.
    Since $s$ almost-$\gamma$-interlaces $\pi$, we may assume that $(C, D) \leq_{\gamma} (A, B)$, so $|B \sm D| \leq \gamma$.
    Thus, $|D| \geq |B| - |B \sm D| \geq \beta + k - \gamma \geq \beta' + k$.
    Since $p$ is displayed by $F$, $p$ has order at most $k$, so $p$ is $\beta'$-broad.

    \casen{Case 2:} Otherwise, every petal $V_h$ with $h \in V(H)$ satisfies $|V_h| < \beta' + k$.
    Since $s$ is $\beta$-broad, we have $|V(G)| \geq 2\beta \geq 4(\beta'+k)$.
    Let $d \eqqcolon (C, D)$ be a displayed separation of $F$ with $|C| \geq \beta'+k$ and $|C|$ minimum under these conditions.
    By minimality of $|C|$ together with $|V_h|<\beta'+k$ for every $h \in V(H)$, we have $|C| < 2(\beta' + k)$.
    Thus, $|D| \geq 4(\beta'+k) - 2(\beta' + k) \geq \beta'+k$.
    Since $d$ is displayed by $F$, $d$ has order at most $k$, hence $d$ is $\beta'$-broad.
\end{proof}

\begin{lem}\label{NewerInterlacingPetalStar}
    Let~$G$ be a $(q,k-1)$-unbreakable $(b, k)$-narrow graph and $\gamma \in \N$. 
    Let 
    \begin{align*}
        \gamma' & \coloneqq \gamma'_{\ref{lem:make-a-petal-star-broad}}(\beta'_{\ref{AlphaInterlacesCloseToInterlaces}}(q,k,b),k, \gamma)~,\\
        \eps & \coloneqq \eps_{\ref{NewerInterlacingPetalStar}}(q,k,b,\gamma) \coloneqq  \eps_{\ref{AlphaInterlacesCloseToInterlaces}}(k, \gamma')~,\\
        \beta & \geq \beta_{\ref{NewerInterlacingPetalStar}}(q, k, b, \gamma) \coloneqq \max\{\,\beta_{\ref{AlphaInterlacesCloseToInterlaces}}(q, k, b, \gamma'),\enspace \beta_{\ref{lem:broad-almost-interlaces-then-broad-displayed}}(k, \beta'_{\ref{AlphaInterlacesCloseToInterlaces}}(q,k,b), \gamma),\enspace \eps+q+b+1\,\}~.
    \end{align*}
    Let~$F$ be a graphic-$k$-flower of $G$
    and let~$\pi$ be the petal-star of~$F$.
    Let~$s$ be a $\beta$-broad $k$-separation that almost-$\gamma$-interlaces~$\pi$.
    Then we have one of the following two outcomes:
    \begin{enumerate}
        \item\label{NewerInterlacingPetalStar-1} $s$ is $\eps$-close to some displayed separation of~$F$;
        \item\label{NewerInterlacingPetalStar-2} $X:=\sep(s)$ is $S$-totally-nested with $S$ the set of all narrow proper $k$-separations of~$G$, and~$s$ is a displayed separation of the graphic-$X$-anemone. Moreover, $s$ is not half-connected.
    \end{enumerate}
\end{lem}
\begin{proof}
    By \cref{lem:broad-almost-interlaces-then-broad-displayed}, $F$ has a $\beta'_{\ref{AlphaInterlacesCloseToInterlaces}}(q,k,b)$-broad displayed separation.
    By \cref{lem:make-a-petal-star-broad}, there is a $k$-flower $F'$ which satisfies the outcome of \cref{lem:make-a-petal-star-broad}.
    Let $\sigma$ be the star of petal-separations of~$F'$.
    By \cref{lem:make-a-petal-star-broad}~\ref{itm:make-a-petal-star-broad-1}, $\sigma$ is a star of $\beta'_{\ref{AlphaInterlacesCloseToInterlaces}}(q,k,b)$-broad separations of $G$.
    Since $G$ is $(q,k-1)$-unbreakable, all the separations in $\sigma$ are $k$-separations.
    By \cref{lem:make-a-petal-star-broad}~\ref{itm:make-a-petal-star-broad-3}, $s$ almost-$\gamma'$-interlaces $\sigma$.
    Since $\beta \geq \beta_{\ref{AlphaInterlacesCloseToInterlaces}}(q,k,b, \gamma')$, we can apply \cref{AlphaInterlacesCloseToInterlaces} to~$s$ and~$\sigma$ and get one of the following two outcomes:
    
    \casen{Case~\ref{AlphaInterlacesCloseToInterlaces-1}: $s$ is $\eps$-close to some $k$-separation~$s'$ that almost-interlaces~$\sigma$.}
    Since $s$ is $\beta$-broad with $\beta \geq \eps + q + b + 1$, we have that $s'$ is a proper separation of $G$ of breadth at least $q+b+1$. Since $G$ is $(b, k)$-narrow and $(q,k-1)$-unbreakable, we get that $s'$ is a narrow proper $k$-separation of $G$ that almost-interlaces~$\sigma$.
    By \cref{lem:flower-almost-interlaces-then-displayed-new}, $s'$ is induced by~$F'$, and hence by~$F$ by \cref{lem:make-a-petal-star-broad}~\ref{itm:make-a-petal-star-broad-0}.
    Since $s'$ is a narrow proper $k$-separation, it is displayed by $F$: this follows from the definition of a daisy in the daisy case, while every induced separation of an anemone is displayed.
    So, we have outcome~\ref{NewerInterlacingPetalStar-1}.

    \casen{Case~\ref{AlphaInterlacesCloseToInterlaces-2}: $s$ is a displayed separation of the graphic-$X$-anemone~$\cA$ with $X:=\sep(s)$, and the narrow proper petal-separations of~$\cA$ are narrow-nested. Moreover, $s$ is not half-connected.}
    In particular, the narrow proper petal-separations of~$\cA$ are $S$-totally-nested.
    Let $S_X \subseteq S$ be the set of half-connected separations in $S$ with separator $X$. 
    Since $s$ is $\beta$-broad with $\beta \geq q + b + 1$, and since $G$ is $(b, k)$-narrow and $(q,k-1)$-unbreakable, we get that $s$ is a narrow proper $k$-separation of $G$, with separator $X$. Thus, there exists a separation in $S$ with separator $X$, so $S_X \neq \emptyset$.
    Furthermore, every separation in~$S_X$ is a narrow proper petal-separation of $\cA$ or its inverse.
    Hence, every separation in $S_X$ is $S$-totally-nested. 
    Thus,~$X$ is $S$-totally-nested by \cref{lem:equivalence-totally-nested-set}.
    Then we have outcome~\ref{NewerInterlacingPetalStar-2}.
\end{proof}

\subsection{Setting}

We now introduce the setting in which we will work throughout this section, and under which assumptions we are essentially able to prove \cref{mainTutte} (see \cref{keylemma:FlowerFinding} and \cref{keylemma:FlowerCriterion}).
Later, in the proof of \cref{mainForNarrow}, we will show how the key lemmas from the previous sections combine to obtain the conditions required for this setting:
\cref{HugeSeparationsUniformlyCrossing} will show that~$S_0$ is uniformly $\gamma_\cA$-crossing, \cref{keylem:flower-approximation} constructs the flowers~$F_A$, \cref{keylem:petal-separations-totally-nested} essentially shows that the petal-separations of~$F_A$ are narrow-$\gamma_N$-nested, and \cref{keylem:max-nested-set} constructs the nested set~$N$ of separations.

\begin{setting}\label{set:FlowerFinding}
    Let $q,k,b,\beta,\gamma_N,\eps_N,\gamma_\cA,\eps_F\in\N$.
    Let $G$ be a $(q,k-1)$-unbreakable graph that is $(b,k)$-narrow.
    Let $S_0$ be the set of all $\beta$-broad $k$-separations of~$G$.
    Assume that $S_0$ is uniformly $\gamma_\cA$-crossing.
    Let $\cA\triv$ and $\cA$ consist of all trivial and nontrivial $\gamma_\cA$-asterisks of~$S_0$, respectively.
    Let $(F_A : A\in\cA )$ be a family of graphic-$k$-flowers.
    Let $S_1:=\bigcup\,\{\,\pi(F_A)\cup (\pi(F_A))^*:A\in\cA\,\}\cup\bigcup\cA\triv$.
    Let $N$ be a nested symmetric set of $k$-separations of~$G$.
    Assume that we have the following properties:
    \begin{defenum}
    \item\label{set:FlowerFinding:Inequalities} $\beta\ge \max\{2q+k,q+1\}$ and $\gamma_\cA\ge\beta\ge b$ and $\eps_N\ge\gamma_N\ge\beta+k$.
    \item\label{set:FlowerFinding:XAnemone} If the common centre $X$ of $A\in\cA$ has size~$k$, then $F_A$ is the graphic-$X$-anemone.
    \item\label{set:FlowerFinding:CloseFlower} For every $A\in\cA$ and $s\in A$ there is $d\in D(F_A)$ such that $s$ and $d$ are $\eps_F$-close.
    \item\label{set:FlowerFinding:S1totallyNested} Every $k$-separation of $G$ is $\gamma_N$-nested with all the petal-separations of all the graphic-$k$-flowers~$F_A$ for $A\in\cA$. 
    \item\label{set:FlowerFinding:NtotallyNested} The separations in $N$ are narrow-$\gamma_N$-nested.
    \item\label{set:FlowerFinding:NcontainsTotallyNested} Every narrow proper $k$-separation of~$G$ that is nested with every narrow proper $k$-separation of~$G$ is contained in~$N$. 
    \item\label{set:FlowerFinding:Napproximates} $N$ $\eps_N$-approximates~$S_1$.
\end{defenum}
\end{setting}

\begin{obs}\label{FlowerFindingBroadSeparationsNice}
    Assume \cref{set:FlowerFinding}.
    Then every $\beta$-broad $(\le k)$-separation $s$ of~$G$ is narrow and has order exactly~$k$.
\end{obs}
\begin{proof}
    Since $G$ is $(q,k-1)$-unbreakable and $\beta> q$ by \ref{set:FlowerFinding:Inequalities}, the order of $s$ must be~$\ge k$.
    Since $G$ is $(b,k)$-narrow and $\beta\ge b$ by \ref{set:FlowerFinding:Inequalities}, we have that $s$ is narrow.
\end{proof}

\subsection{High-spread asterisks interlace unique splitting stars}

The main objective in this subsection is to prove \cref{HighSpreadUniqueSplittingStar}, which roughly says that every asterisk~$A$ of large enough \emph{spread} (that is, the largest integer~$\gamma$ such that two separations in~$A$ $\gamma$-cross) $\eps$-interlaces a unique splitting star of~$N$.
One main step is to prove \cref{NewerUniqueAsterisk}, which essentially says that although a lot of $\delta'$-asterisks of $D(F)$ can `live' in a given flower $F$, there is a unique $\delta$-asterisk of $D(F)$ (with $\delta \approx \delta'/2$) that contains all those $\delta'$-asterisks.

\begin{lem}\label{lem:gamma-crossed-uniform}
    Let $G$ be a $(q,k-1)$-unbreakable graph and $\gamma\ge \gamma_{\ref{HugeSeparationsUniformlyCrossing}}(q,k)$.
    Let $S$ be a set of $(\leq k)$-separations of~$G$.
    Then $S' \coloneqq \{\,s\in S:s\text{ is $\gamma$-crossed by some }s'\in S\,\}$ is uniformly $\gamma$-crossing.
\end{lem}

\begin{proof}
    If $s \in S'$ then $s$ is $\gamma$-crossed by some $s' \in S$, so $s$ is $2(\gamma-k)$-broad by \cref{LargeCrossedImpliesBroad}.
    Since $\gamma \geq \gamma_{\ref{HugeSeparationsUniformlyCrossing}}(q,k) = 3q+3k$, we have $2(\gamma-k) > q$, so $s$ is $(>q)$-broad.
    Since $G$ is $(q,k-1)$-unbreakable, it follows that $s$ is a $k$-separation.
    Thus, $S'$ is uniformly $\gamma$-crossing by \cref{HugeSeparationsUniformlyCrossing}.
\end{proof}

\begin{definition}[$\cA(S,\gamma)$]\label{defn:cA-gamma-asterisks}
    Let $G$ be a $(q,k-1)$-unbreakable graph and $\gamma\ge \gamma_{\ref{HugeSeparationsUniformlyCrossing}}(q,k)$.
    Let $S$ be a symmetric set of $(\leq k)$-separations of~$G$.
    We denote by $\defnMath{\cA(S,\gamma)}$ the set of all $\gamma$-asterisks of~$\{\,s\in S:s\text{ is $\gamma$-crossed by some }s'\in S\,\}$, which is well-defined by \cref{lem:gamma-crossed-uniform}.
\end{definition}

\begin{lemma}\label{NewerUniqueAsterisk}
    Let $\delta,\delta'\in\N$ such that $\delta \geq \gamma_{\ref{HugeSeparationsUniformlyCrossing}}(q, k)$ and $\delta'\geq \delta'_{\ref{NewerUniqueAsterisk}}(k,\delta):= 2\delta +2k$.
    Let $G$ be a $(q,k-1)$-unbreakable graph, and let $F$ be a graphic-$k$-flower of~$G$ such that $\cA(D(F),\delta')\neq\emptyset$.
    Then there is a unique $A\in \cA(D(F),\delta)$ such that $\bigcup\cA(D(F),\delta')\se A$.
\end{lemma}

\begin{proof}
    It suffices to show that for every $A\neq A'\in \cA(D(F),\delta')$, there is a displayed separation of $F$ that $\delta$-crosses at least one separation in each of $A$ and~$A'$.
    Let $A\neq A'\in \cA(D(F),\delta')$.
    Let $r_1, s_1 \in A$ $\delta'$-cross and $r_2, s_2 \in A'$ $\delta'$-cross.
    Write $r_i \eqqcolon (A_i, B_i)$ and $s_i \eqqcolon (C_i, D_i)$ for $i=1,2$.
    Then $r_i$ and $s_i$ $(\delta'-k)$-criss by \cref{crissCross} for both~$i$.
    If $r_1$ and $r_2$ $\delta$-cross, then we are done.
    Hence we may assume that $r_1 \le_\delta r_2$, say.

    Assume first that $F = (O, \mathcal{V})$ is a graphic-$k$-daisy.
    See \cref{fig:NewerUniqueAsterisk} for an illustration of the setting.
    Since $r_1$ is displayed by $F$, there exist two subpaths $P_{A_1}, P_{B_1}$ of $O$ such that $V(O) = V(P_{A_1}) \sqcup V(P_{B_1})$, and $A_1 = V \rest P_{A_1}$ and $B_1 = V \rest P_{B_1}$.
    Let $P_{C_1}, P_{D_1}, P_{A_2}, P_{B_2}, P_{C_2}$ and $P_{D_2}$ be defined similarly.
    Let $e_1, f_1$ be the only edges of $O$ that are not edges of $P_{C_1} \cup P_{D_1}$. Since $(A_1, B_1)$ and $(C_1, D_1)$ cross, exactly one of them, say $e_1$, is an edge of $P_{A_1}$.
    Similarly, let $e_2, f_2$ be the only edges of $O$ that are not edges of $P_{C_2} \cup P_{D_2}$. Since $(A_2, B_2)$ and $(C_2, D_2)$ cross, exactly one of them, say $f_2$, is an edge of $P_{B_2}$.

    We first show that $f_2 \in E(P_{B_1})$.
    Since $r_1 \leq_{\delta} r_2$, we have $|A_1 \cap B_2| \leq |A_1 \setminus A_2| + |A_2 \cap B_2| \leq \delta + k$.
    Let $R_1, R_2$ be the two components of $P_{B_2} - \{f_2\}$.
    Each $V \rest R_j$ contains one of the proper $r_2s_2$-corners, so $|V \rest R_j| > \delta'-k$ by \cref{crissCross}~\ref{crissCross:crossImpliesCriss}.
    Suppose for a contradiction that $f_2 \notin E(P_{B_1})$, so $f_2$ lies in $P_{A_1}$ or is one of its boundary edges.
    If one of the endpoints of $P_{B_2}$ lies in $P_{A_1}$ then the corresponding $R_j$ is contained in $P_{A_1}$, otherwise $P_{A_1}$ is contained in $P_{B_2}$.
    In the former case, we have $|A_1 \cap B_2| \geq |V \rest R_j| > \delta'-k \geq \delta+k$.
    In the latter case, we have $|A_1 \cap B_2| \geq |A_1| \geq \delta' > \delta+k$.
    In both cases we obtain a contradiction so $f_2 \in E(P_{B_1})$.
    Symmetrically, we have $e_1 \in P_{A_2}$.
    In particular, $e_1 \neq f_2$.
    
    Then, $O - \{e_1, f_2\}$ has exactly two connected components, which are paths, call them~$Q_1$ and~$Q_2$.
    Let $s \coloneqq (U, W) \coloneqq (V \rest Q_1, V \rest Q_2)$, and observe that $s \in D(F)$. 

    \begin{figure}[ht]
        \centering
        \includegraphics[height=16\baselineskip]{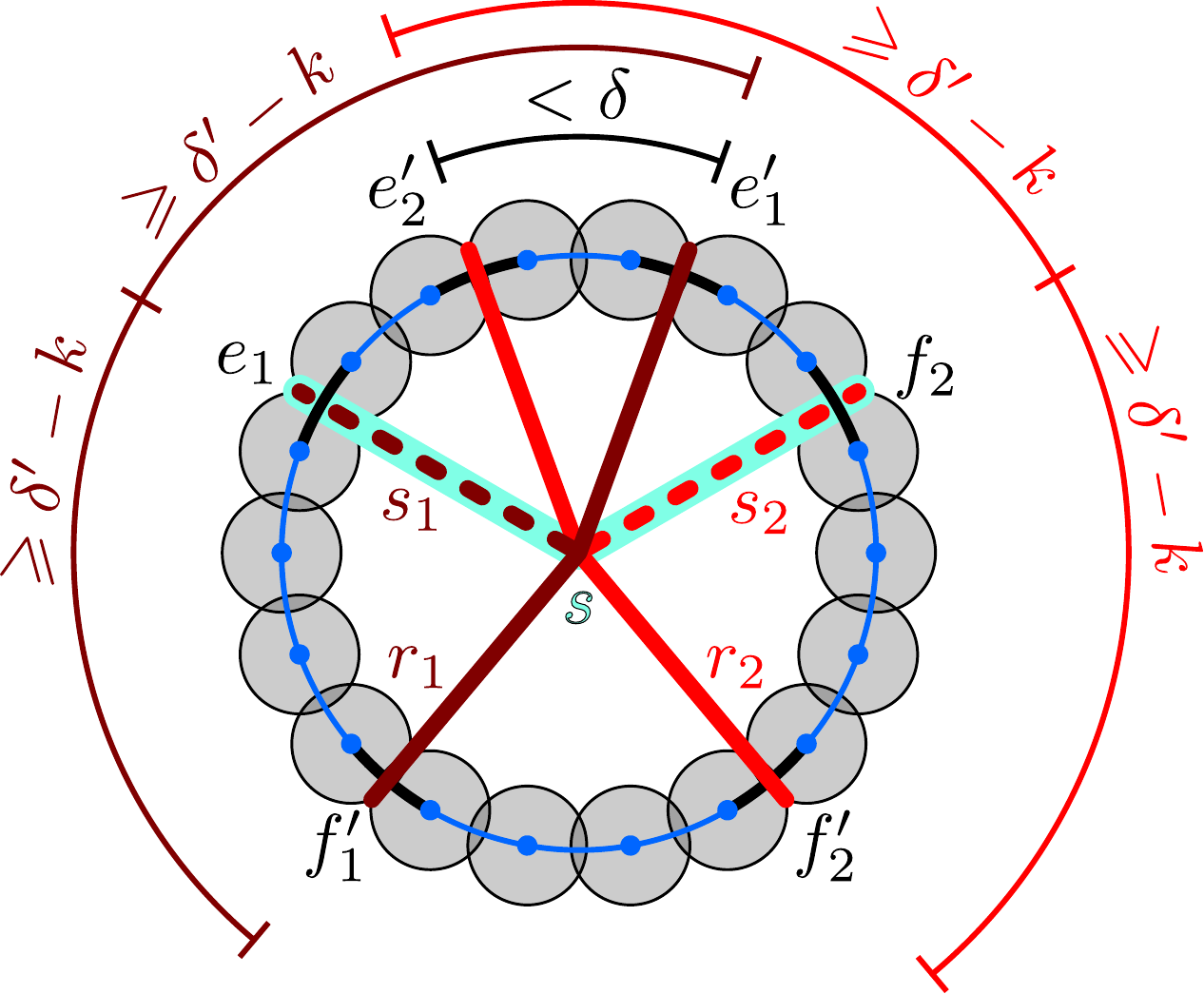}
        \caption{The situation in the proof of \cref{NewerUniqueAsterisk}. Only one half of each $s_i$ is drawn.
        Here, $s$ is the turquoise separation.}
        \label{fig:NewerUniqueAsterisk}
    \end{figure}

    \begin{claim}
        $s$ $\delta$-crosses both $r_1$ and~$r_2$.
    \end{claim}

    \begin{cproof}
        Since $f_2 \in E(P_{B_1})$ and $e_1 \in E(P_{A_1})$, the cut $\{e_1, f_2\}$ splits $P_{A_1}$ at $e_1$, just as the cut defining $s_1$ does.
        Thus, one of the proper $A_1C_1$-corner and the proper $A_1D_1$-corner is contained in the proper $A_1U$-corner and the other is contained in the proper $A_1W$-corner.
        Thus, the proper $A_1U$-corner and the proper $A_1W$-corner both have size greater than $\delta'-k$ by \cref{crissCross}~\ref{crissCross:crossImpliesCriss}.
        Similarly, since $e_1 \in P_{A_2}$ and $f_2 \in P_{B_2}$, one of the proper $B_2C_2$-corner and the proper $B_2D_2$-corner is contained in the proper $B_2U$-corner and the other is contained in the proper $B_2W$-corner.
        Thus, the proper $B_2U$-corner and the proper $B_2W$-corner both have size greater than $\delta'-k$.

        Since $r_1 \leq_{\delta} r_2$, it follows that \[|(B_1 \setminus A_1) \cap (U \setminus W)| \geq |(B_2 \setminus A_2) \cap (U \setminus W)| - |A_1 \setminus A_2| > \delta' - k - \delta \geq \delta\] by the above.
        Similarly, we also have $|(B_1 \setminus A_1) \cap (W \setminus U)| > \delta$.
        Thus, $s$ and $r_1$ $\delta$-criss so they $\delta$-cross by \cref{crissCross}~\ref{crissCross:crissImpliesCross}.
        Symmetrically, $s$ and $r_2$ $\delta$-cross.
    \end{cproof}
    
    Assume now that $F = (S, \mathcal{V})$ is a graphic-$k$-anemone with centre $c$.
    Note that $|(A_1 \cap B_2) \setminus (A_2 \cup B_1)| \leq |A_1 \setminus A_2| \leq \delta$ since $r_1 \leq_{\delta} r_2$.
    Let $X = V_c$, so the separator of every proper separation induced by $F$ is equal to $X$.
    Let $P \coloneqq \{p \in V(S) : V_p \subseteq (A_1 \cap C_1) \cup (B_2 \cap D_2)\}$ and let $s \eqqcolon (U, W) \eqqcolon (\bigcup_{p \in P}V_p, \bigcup_{p \notin P}V_p)$.
    Note that $s$ is displayed by $F$.

    \begin{claim}
        $s$ $\delta$-crosses both $r_1$ and~$r_2$.
    \end{claim}

    \begin{cproof}
        By \cref{crissCross} and since $\delta \le \delta'-\delta-2k$, it suffices to prove that $s$ $(\delta'-\delta-2k)$-crisses both $r_1$ and~$r_2$.
        We prove that $s$ $(\delta'-\delta-2k)$-crisses $r_1$; the proof for $r_2$ is symmetric.

        The separations $r_1$ and $s_1$ $\delta'$-cross, and both have their separator equal to $X$ so $|(A_1 \cap C_1) \sm X| > \delta'-k$ by \cref{crissCross}~\ref{crissCross:crossLargeCorner}.
        Let $v \in (A_1 \cap C_1) \sm X$.
        There exists a unique $p \in V(S)$ such that $v \in V_p$.
        Since $r_1$ and $s_1$ are displayed by $F$, we have $V_p \subseteq A_1 \cap C_1$, so $p \in P$, and therefore $v \in U$.
        Since $s$ is displayed by $F$, its separator is contained in $X$.
        Therefore, the $A_1U$-corner contains $(A_1 \cap C_1) \sm X$, so has size greater than $\delta'-k$.

        As above, we have $|(B_2 \cap D_2) \sm X| > \delta'-k$.
        Let $v \in (B_2 \cap D_2) \sm X$. As before, we have $v \in U$.
        Furthermore, $|(B_2 \sm X) \sm B_1| \leq |(A_1 \cap B_2) \sm X| \leq \delta$.
        Therefore, $|(B_1 \cap U) \sm X| \geq |(B_2 \cap U) \sm X| - |(B_2 \sm X) \sm B_1| \geq |(B_2 \cap D_2) \sm X| - \delta > \delta'-\delta-k$.
        Therefore, the $B_1U$-corner has size greater than $\delta'-\delta-k$.

        As above, we have $|(A_1 \cap D_1) \sm X| > \delta'-k$.
        Let $v \in (A_1 \cap D_1) \sm X$. Suppose further that $v \notin B_2$.
        There exists a unique $p \in V(S)$ such that $v \in V_p$. Since $v \notin A_1 \cap C_1$ and $v \notin B_2 \cap D_2$, we have $p \notin P$, so $v \notin U$ so $v \in W \sm U$.
        Therefore, the $A_1W$-corner has size at least $|(A_1 \cap D_1) \sm (X \cup B_2)| \geq |(A_1 \cap D_1) \sm X| - |(A_1 \cap B_2) \sm X| > \delta'-k-\delta$.

        As above, we have $|(B_2 \cap C_2) \sm X| > \delta'-k$.
        Let $v \in (B_2 \cap C_2) \sm X$. Suppose further that $v \notin A_1$.
        There exists a unique $p \in V(S)$ such that $v \in V_p$. Since $v \notin A_1 \cap C_1$ and $v \notin B_2 \cap D_2$, we have $p \notin P$, so $v \notin U$ so $v \in W \sm U$.
        Therefore, $|(B_1 \cap W) \sm X| \geq |(B_1 \cap B_2 \cap C_2) \sm X| = |(B_2 \cap C_2) \sm X| - |(B_2 \cap C_2) \sm (B_1 \cup X)| > \delta'-k - |B_2 \sm B_1| \geq \delta' - k - \delta$. 
    \end{cproof}
    
\end{proof}

\begin{definition}[Spread]
    Let $S$ be a set of separations of~$G$.
    The \defn{spread} of~$S$ is the maximum number $\gamma\in\N$ for which some two separations in~$S$ $\gamma$-cross, and we set it to $-\infty$ if no two separations in~$S$ cross.
    The \defn{spread} of a graphic-$k$-flower is the spread of the set of its displayed separations.
\end{definition}

\begin{definition}[Relevant]
    Assume \cref{set:FlowerFinding}.
    Let $\delta,\delta'\in\N$ such that $\delta \geq \gamma_{\ref{HugeSeparationsUniformlyCrossing}}(q, k)$ and $\delta'\geq\delta'_{\ref{NewerUniqueAsterisk}}(k,\delta)$.
    We say that a nontrivial $\gamma_\cA$-asterisk $A\in\cA$ is \defn{$\delta'$-relevant} if $\cA(D(F_A),\delta')\neq\emptyset$, that is, if $D(F_A)$ has spread $\ge\delta'$.
    In this event, we denote the unique nontrivial $\delta$-asterisk provided by \cref{NewerUniqueAsterisk} by~$\defnMath{A^\delta}$, where $\delta'$ will be clear from context.
    Note that $A^\delta\se D(F_A)$.
\end{definition}

\begin{lem}\label{NewerLargelyCrossingEntersAsterisk}
    Assume \cref{set:FlowerFinding}.
    Let $\delta,\delta'\in\N$ such that 
    \[
        \delta > \delta_{\ref{NewerLargelyCrossingEntersAsterisk}}(q,k,\gamma_\cA,\eps_F):= \max\{\,2k+3\gamma_\cA + 3\eps_F\,,\,\gamma_{\ref{HugeSeparationsUniformlyCrossing}}(q, k)\,\} \qquad \text{and} \qquad \delta'\geq\delta'_{\ref{NewerUniqueAsterisk}}(k,\delta).
    \]
    Let $A\in\cA$ have spread $\ge\zeta_{\ref{NewerLargelyCrossingEntersAsterisk}}(\delta',\eps_F):=\delta'+2\eps_F$.
    Then:
    \begin{enumerate}
        \item\label{NewerLargelyCrossingEntersAsterisk:relevant} $A$ is $\delta'$-relevant and $A^\delta$ exists;
        \item\label{NewerLargelyCrossingEntersAsterisk:deltaCross} whenever $d_1,d_2\in D(F_A)$ $\delta$-cross, we have $d_1,d_2\in A$;
        \item\label{NewerLargelyCrossingEntersAsterisk:Adelta} $A^\delta\se A$.
    \end{enumerate}
\end{lem}

\begin{proof}
    \cref{NewerLargelyCrossingEntersAsterisk:relevant}.
    Let $s'_1,s'_2\in A$ $(\delta'+2\eps_F)$-cross.
    Let $s_1,s_2\in D(F_A)$ be $\eps_F$-close to $s'_1,s'_2$, respectively, which we find with \ref{set:FlowerFinding:CloseFlower}.
    Then $s_1$ and $s_2$ $\delta'$-cross by \cref{lem:close-separations-still-cross}.
    Note that $\delta\ge \gamma_{\ref{HugeSeparationsUniformlyCrossing}}(q, k)$.
    Hence $A$ is $\delta'$-relevant, and $A^\delta$ exists by \cref{NewerUniqueAsterisk}.

    \cref{NewerLargelyCrossingEntersAsterisk:deltaCross}.
    Let $d_1,d_2\in D(F_A)$ $\delta$-cross.
    By symmetry, it suffices to show $d_1\in A$.
    Since $d_1$ and $d_2$ $\delta$-cross, $d_1$ is $(\delta-k)$-broad by \cref{LargeCrossedImpliesBroad}, with $\delta-k \geq \gamma_{\mathcal{A}} \geq \beta$ by \ref{set:FlowerFinding:Inequalities}.
    Furthermore, $\beta > q$ by \ref{set:FlowerFinding:Inequalities} and $G$ is $(q,k-1)$-unbreakable so $d_1$ is a $\beta$-broad $k$-separation of $G$ so $d_1 \in S_0$.
    If $d_1$ $(\gamma_\cA + \eps_F)$-crosses $s_1$ or $s_2$, then $d_1$ $\gamma_{\cA}$-crosses $s'_1$ or $s'_2$ by \cref{CloseVsCross}, so $d_1\in A$.
    So we may assume that $d_1$ is $(\gamma_\cA + \eps_F)$-nested with both $s_1$ and~$s_2$.
    Then $s_1,s_2\le_{\gamma_\cA + \eps_F} d_1$, say.
    Write $s_i = (A_i, B_i)$ and $d_i = (C_i, D_i)$ for $i \in [2]$.

    \begin{claim}\label{cl:exists-crossing-separation}
        There is a displayed separation $s = (U, W)$ of $F_A$ with the following properties (up to renaming $C_2$ and $D_2$): \begin{itemize}
            \item $(A_1 \cap B_2) \cup (D_1 \cap C_2) \subseteq U$, 
            \item $(A_1 \cap A_2) \cup (D_1 \cap D_2) \subseteq_{\gamma_{\cA} + \eps_F} W$.
        \end{itemize} 
    \end{claim}

    \begin{cproof}
        Suppose first that $F_A \eqqcolon (S, \cV)$ is a graphic-$k$-anemone with centre $c$. 
        Let $X \coloneqq V_c$, so the separator of every proper separation induced by $F_A$ is equal to $X$.
        Let \[P \coloneqq \{p \in V(S) \setminus \{c\} \mid V_p \subseteq (A_1 \cap B_2) \cup (D_1 \cap C_2)\}\] and let $s \coloneqq (\bigcup_{p \in P}V_p, \bigcup_{p \notin P}V_p)$.
        Since $s_1$ and $s_2$ $\delta'$-cross, by \cref{crissCross}~\ref{crissCross:crossLargeCorner} we have $|A_1 \cap B_2| > \delta' \geq\delta'_{\ref{NewerUniqueAsterisk}}(k,\delta) \geq k$, so there exists $v \in (A_1 \cap B_2) \setminus X$.
        Let $p \in V(S) \setminus \{c\}$ be the unique node such that $v \in V_p$. 
        Since $s_1$ and $s_2$ are displayed by $F_A$, we have $V_p \subseteq A_1 \cap B_2$, so $p \in P$ and $P \neq \emptyset$.
        Thus, $s$ is a proper displayed separation of $F_A$, so $X = \sep(s)$.
        Write $s \eqqcolon (U, W)$.
        As above, if $v \in ((A_1 \cap B_2) \cup (D_1 \cap C_2)) \setminus X$, the unique $p \in V(S)$ such that $v \in V_p$ satisfies $p \in P$, so $v \in U$. 
        Thus, $(A_1 \cap B_2) \cup (D_1 \cap C_2) \subseteq U$.
        Then, let $N \coloneqq |((A_1 \cap B_2) \cup (D_1 \cap C_2)) \cap ((A_1 \cap A_2) \cup (D_1 \cap D_2)) \setminus X|$. Using that $s_1 \le_{\gamma_{\cA} + \eps_F}d_1$, we have \begin{align*}
            N &= |((A_1 \cap A_2 \cap B_2) \cup (A_1 \cap B_2 \cap D_1 \cap D_2) \cup (D_1 \cap C_2 \cap A_1 \cap A_2) \cup (D_1 \cap C_2 \cap D_2)) \setminus X| \\
            &\leq |(A_2 \cap B_2) \setminus X| + |(A_1 \cap D_1) \setminus X| + |(C_2 \cap D_2) \setminus X| \\
            &\leq |(A_1 \cap B_1) \setminus X| + |(D_1 \setminus B_1) \setminus X| \\
            &\leq \gamma_{\cA} + \eps_F.
        \end{align*}
        Thus, to prove that $(A_1 \cap A_2) \cup (D_1 \cap D_2) \subseteq_{\gamma_{\cA} + \eps_F} W$, it suffices to prove that $((A_1 \cap A_2) \cup (D_1 \cap D_2)) \setminus ((A_1 \cap B_2) \cup (D_1 \cap C_2)) \subseteq W$, where we use $X=\sep(s) \subseteq W$.
        Let $v \in ((A_1 \cap A_2) \cup (D_1 \cap D_2)) \setminus ((A_1 \cap B_2) \cup (D_1 \cap C_2))$ and let $p \in V(S) \setminus \{c\}$ be the unique node such that $v \in V_p$. 
        Since $v \notin (A_1 \cap B_2) \cup (D_1 \cap C_2)$, we have $p \notin P$ so $v \in W$.
    
        Suppose now that $F_A \eqqcolon (O, \cV)$ is a graphic-$k$-daisy.
        Given two edges $e, e' \in E(O)$, $O - \{e, e'\}$ has exactly two connected components, which are paths, call them $Q_1$ and $Q_2$. We say that the separations $(V \rest Q_1, V \rest Q_2)$ and $(V \rest Q_2, V \rest Q_1)$ are \emph{displayed by the edge-cut $\{e, e'\}$.}
        There exist edges $e_1, e_2, e'_1, e'_2, f_1, f_2, f'_1, f'_2$ of $O$ such that each $s_i$ is displayed by the edge-cut $\{e_i, e'_i\}$ and each $d_i$ is displayed by the edge-cut $\{f_i, f'_i\}$.
        We arbitrarily fix one of the two orientations of the cycle $O$.
        For $i \in \{1, 2\}$, up to exchanging $A_i$ and $B_i$, we may assume that $A_i$ is the union of the bags of the nodes on the arc from $e'_i$ to $e_i$ in this orientation, and $B_i$ is the union of the bags of the nodes on the arc from $e_i$ to $e'_i$ in this orientation, see \cref{fig:NewerLargelyCrossingEntersAsterisk}.
        Since $s_1$ and $s_2$ cross, we may assume that the edges $e_2,e_1,e'_2,e'_1$ appear in this order along $O$ (following this orientation), as in \cref{fig:NewerLargelyCrossingEntersAsterisk}.
        Up to renaming $f_1$ and $f'_1$, we may assume that $e_2, f_1, f'_1$ appear in this order along $O$.
        Since $d_1$ and $d_2$ cross, up to renaming $f_2$ and $f'_2$, we may assume that $e_2, f_1, f_2, f'_1$ appear in this order along $O$, again see \cref{fig:NewerLargelyCrossingEntersAsterisk}.
        Let $Q_1, Q_2$ be the two connected components of $O - \{e_2, f_2\}$. By the above, up to renaming $Q_1$ and $Q_2$, we may assume that $f_1 \in E(Q_1)$ and $f'_1 \in E(Q_2)$.
        Let $s \coloneqq (V \rest Q_1, V \rest Q_2)$.
        Then, $s$ is a displayed separation of $F_A$.
        Write $s \eqqcolon (U, W)$.
        Since $f_1 \in E(Q_1)$ and $f'_1 \in E(Q_2)$, we have $D_1 \cap C_2 \subseteq U$ and $D_1 \cap D_2 \subseteq W$ (up to renaming $C_2$ and $D_2$), see \cref{fig:NewerLargelyCrossingEntersAsterisk}.
        We now argue that $e_1, f_2, e'_2$ appear in this order along $O$. 
        If not, $e_1, e'_2, f_2$ appear in this order along $O$. 
        Thus, either $e_1, e'_2, f_2, e_2$ appear in this order along $O$ or $e_1, e'_2, e_2, f_2$ appear in this order along $O$.
        By symmetry, suppose that $e_1, e'_2, f_2, e_2$ appear in this order along $O$.
        Then, $e'_2, f_2, f'_1, e_2$ appear in this order along $O$.
        Combined with the fact that $e_2, f_1, f'_1$ appear in this order along $O$, this implies that the $D_1D_2$-corner is contained in the $D_1A_2$-corner (see \cref{fig:NewerLargelyCrossingEntersAsterisk}).
        However, $d_1$ and $d_2$ $\delta$-cross so they $(\delta-k)$-criss by \cref{crissCross}~\ref{crissCross:crossImpliesCriss}.
        Thus, the $D_1D_2$-corner has size $> \delta - k \geq \gamma_{\cA} + \eps_F$, so the $D_1A_2$-corner also has size $> \gamma_{\cA} + \eps_F$, so $|D_1 \setminus B_2| > \gamma_{\cA} + \eps_F$, which contradicts $s_2\le_{\gamma_\cA + \eps_F} d_1$.
        Thus, $e_1, f_2, e'_2$ indeed appear in this order along $O$. 
        This immediately yields $A_1 \cap B_2 \subseteq U$ and $A_1 \cap A_2 \subseteq W$.
    \end{cproof}

    \begin{figure}[ht]
        \centering
        \includegraphics[height=16\baselineskip]{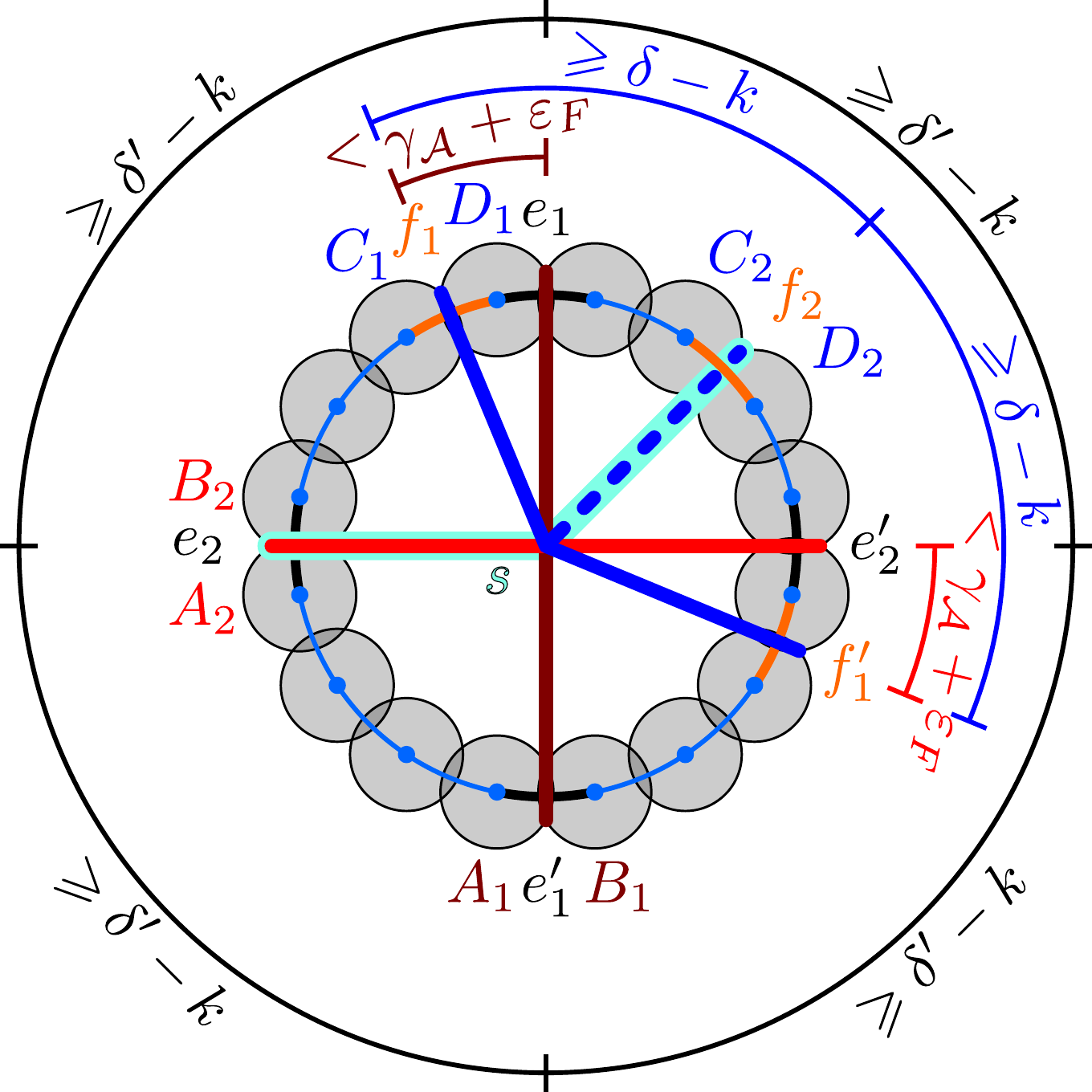}
        \caption{The situation in the proof of \cref{NewerLargelyCrossingEntersAsterisk}.
        Here, the cycle~$O$ is oriented clockwise, $s$ is the turquoise separation, and $\delta'\gg\delta\gg\gamma_\cA, k$.}
        \label{fig:NewerLargelyCrossingEntersAsterisk}
    \end{figure}

    Let~$s$ be the separation obtained from \cref{cl:exists-crossing-separation}.
    \begin{claim}
        All extended corners of the separations $s$ and $d_1$, and of the separations $s$ and $s_1$ have size greater than $(\delta-2\gamma_{\cA} - 2\eps_F)$.
    \end{claim}
    
    \begin{cproof}
        We encourage the reader to check this statement by looking at \cref{fig:NewerLargelyCrossingEntersAsterisk}, as the proof simply consists of computations.
        
        Using \cref{crissCross}~\ref{crissCross:crossLargeCorner} and that $s$ was obtained from \cref{cl:exists-crossing-separation} throughout, we have 
        \[
            |A_1 \cap U| \geq |A_1 \cap B_2| > \delta'
        \]
        since $s_1$ and $s_2$ $\delta'$-cross.
        We have 
        \[
            |A_1 \cap W| \geq |A_1 \cap A_2 \cap W| > \delta' - (\gamma_{\cA} + \eps_F)
        \]
        since $s_1$ and $s_2$ $\delta'$-cross.
        We have 
        \[
            |B_1 \cap W| \geq |B_1 \cap D_1 \cap D_2 \cap W| \geq |D_1 \cap D_2 \cap W| - |D_1 \sm B_1| > \delta - (\gamma_{\cA} + \eps_F) - (\gamma_{\cA} + \eps_F)
        \]
        since $d_1$ and $d_2$ $\delta$-cross and $s_1 \leq_{\gamma_{\cA} + \eps_F}d_1$.
        We have 
        \[
            |B_1 \cap U| \geq |B_1 \cap D_1 \cap C_2| \geq |D_1 \cap C_2| - |D_1 \sm B_1| > \delta - (\gamma_{\cA} + \eps_F)
        \]
        since $d_1$ and $d_2$ $\delta$-cross and $s_1 \leq_{\gamma_{\cA} + \eps_F}d_1$.
        We have 
        \[
            |C_1 \cap U| \geq |A_1 \cap B_2 \cap C_1| \geq |A_1 \cap B_2| - |A_1 \sm C_1| > \delta' - (\gamma_{\cA} + \eps_F)
        \] 
        since $s_1$ and $s_2$ $\delta'$-cross and $s_1 \leq_{\gamma_{\cA} + \eps_F}d_1$.
        We have 
        \[
            |C_1 \cap W| \geq |A_1 \cap A_2 \cap C_1 \cap W| \geq |A_1 \cap A_2 \cap W| - |A_1 \sm C_1| > \delta' - (\gamma_{\cA} + \eps_F) - (\gamma_{\cA} + \eps_F)
        \]
        since $s_1$ and $s_2$ $\delta'$-cross and $s_1 \leq_{\gamma_{\cA} + \eps_F}d_1$.
        We have 
        \[
            |D_1 \cap W| \geq |D_1 \cap D_2 \cap W| > \delta - (\gamma_{\cA} + \eps_F)
        \]
        since $d_1$ and $d_2$ $\delta$-cross.
        We have 
        \[
            |D_1 \cap U| \geq |D_1 \cap C_2| > \delta
        \]
        since $d_1$ and $d_2$ $\delta$-cross.
    \end{cproof}

    Thus, the separations $s$ and $d_1$ $(\delta-2k-2\gamma_{\cA} - 2\eps_F)$-criss, so $s$ has breadth at least $2(\delta-2k-2\gamma_\cA - 2\eps_F) \ge \gamma_\cA\ge\beta$, using that $\delta \geq 2k+3\gamma_{\cA} + 3\eps_F$ and~\ref{set:FlowerFinding:Inequalities}.
    Furthermore, $\beta > q$ and $G$ is $(q,k-1)$-unbreakable so $s\in S_0$.
    Since $s'_1$ is $\eps_F$-close to~$s_1$, we have that $s'_1$ $(\delta-2k-2\gamma_\cA-3\eps_F)$-crosses~$s$ by \cref{CloseVsCross}.
    Now $s\in S_0$ and $s'_1\in A$ with 
    \[
        \delta-2k-2\gamma_\cA-3\eps_F\ge\gamma_\cA
    \]
    yields $s\in A$.
    Hence $d_1\in A$ as well.

    \cref{NewerLargelyCrossingEntersAsterisk:Adelta} follows from \cref{NewerLargelyCrossingEntersAsterisk:relevant} and \cref{NewerLargelyCrossingEntersAsterisk:deltaCross}.
\end{proof}

\begin{cor}\label{highSpreadFlowersDistinct}
    Assume \cref{set:FlowerFinding}.
    Let $\delta>\delta_{\ref{NewerLargelyCrossingEntersAsterisk}}(q,k,\gamma_\cA,\eps_F)$ and $\delta'\ge \delta'_{\ref{NewerUniqueAsterisk}}(k,\delta)$.
    Let $A_1 \neq A_2\in\cA$ have spread $\ge\zeta_{\ref{NewerLargelyCrossingEntersAsterisk}}(\delta',\eps_F)$.
    Then $D(F_{A_1})\neq D(F_{A_2})$.
    In particular, $F_{A_1}$ and $F_{A_2}$ are not $\mathrm{id}_G$-isomorphic, and $F_{A_1}\neq F_{A_2}$.
\end{cor}
\begin{proof}
    Assume for a contradiction that $D(F_{A_1})=D(F_{A_2})$.
    By \cref{NewerLargelyCrossingEntersAsterisk}, both $A_i$ are $\delta'$-relevant, their $A^\delta_i$ exists, and $A^\delta_i\se A_i$.
    Since $D(F_{A_1})=D(F_{A_2})$, we find that $A^\delta_1=A^\delta_2$ by definition.
    Hence $A_1$ and $A_2$ intersect, contradicting that $A_1$ and $A_2$ are disjoint by definition.
\end{proof}

\begin{lem}
    \label{NewerLargeAsteriskNotCloseButBetter}
    Assume~\cref{set:FlowerFinding}.
    Let $\delta\in\N$ such that $\delta\ge\gamma_N$.
    Let $s,s'\in S_0$ $\delta$-cross.
    Then neither $s$ nor $s'$ is $(\delta-\gamma_N)$-close to any $t\in N$.
\end{lem}
\begin{proof}
    By symmetry, it suffices to focus on~$s$.
    Assume for a contradiction that $s$ is $(\delta-\gamma_N)$-close to some~$t\in N$.
    Then $s'$ $\gamma_N$-crosses~$t$ by \cref{CloseVsCross}.
    Recall that all separations in $S_0$ are narrow by \cref{FlowerFindingBroadSeparationsNice}.
    This contradicts the assumption \ref{set:FlowerFinding:NtotallyNested} that every separation in $N$ is narrow-$\gamma_N$-nested.
\end{proof}

\begin{lem}
    \label{NewerLargeAsteriskUniqueSplittingStar}
    Assume~\cref{set:FlowerFinding}.
    Let $\delta\in\N$ such that $\delta\ge\delta_{\ref{NewerLargeAsteriskUniqueSplittingStar}}(\gamma_N,\eps_N):=\gamma_N + 2\eps_N$.
    Let $A\se S_0$ be nontrivially $\delta$-asteriskal.
    For every $s\in A$ let
    \begin{equation}\label{NewerLargeAsteriskUniqueSplittingStar:eq:Orientation}
        O(s):=\{\,t\in N : t \le_{\eps_N} s\text{ or }t\le_{\eps_N} s^*\,\}
    \end{equation}
    and let $\sigma(s)$ consist of the $\le$-maximal elements of~$O(s)$.
    Then $O(s)$ is a consistent orientation of~$N$ and $\sigma(s)$ is the unique splitting star of~$N$ such that $s$ $\eps_N$-interlaces~$\sigma(s)$.
    Moreover, $O(s)=O(s')$ and $\sigma(s)=\sigma(s')$ for every two $s,s'\in A$.
\end{lem}
\begin{proof}
    Since $A$ is nontrivially $\delta$-asteriskal and $s \in A$, $s$ $\delta$-crosses some separation of $A$.
    By \cref{LargeCrossedImpliesBroad}, $s$ is $2(\delta-k)$-broad.
    Using that $\gamma_N \geq \beta + k \geq k+1$ by \ref{set:FlowerFinding:Inequalities}, it follows that $s$ has breadth at least $2(\gamma_N+2\eps_N-k) \geq 2\eps_N+1$.
    
    Since $\delta\ge\gamma_N$ and $\delta-\gamma_N\ge\eps_N$ by assumption, no $t\in N$ is $\eps_N$-close to~$s$ by \cref{NewerLargeAsteriskNotCloseButBetter}.
    We have $s \in A \subseteq S_0$ so $s$ is $\beta$-broad, hence $b$-broad by \ref{set:FlowerFinding:Inequalities}.
    Since $G$ is $(b, k)$-narrow, this implies that $s$ is narrow.
    By \ref{set:FlowerFinding:Inequalities} and \ref{set:FlowerFinding:NtotallyNested}, $s$ is $\eps_N$-nested with every separation in~$N$.
    By \cref{lem:gamma-close-or-gamma-interlaces},
    $O(s)$ is a consistent orientation of~$N$, and $\sigma(s)$ is a splitting star of~$N$.
    Note that $s$ $\eps_N$-interlaces~$\sigma(s)$.
    
    Since $O(s_1)=O(s_2)$ implies $\sigma(s_1)=\sigma(s_2)$,
    it remains to show the former for every two $s_1,s_2\in A$.
    For this, let $s_1\neq s_2\in A$.
    Since $A$ is $\delta$-crossing-connected,
    we may assume without loss of generality that $s_1$ and $s_2$ $\delta$-cross.
    Now assume for a contradiction that $O(s_1) \neq O(s_2)$.
    Then we can assume without loss of generality that there exists $t \in O(s_2) \sm O(s_1)$.
    Then, either $s_1 \leq_{\eps_N} t$ or $(s_1)^* \leq_{\eps_N} t$, and either $t \le_{\eps_N} s_2$ or $t\le_{\eps_N} (s_2)^*$.
    In all cases, since $\delta\ge 2\eps_N$ and with \cref{approxLeQuasiTransitive}, this contradicts our assumption that $s_1$ and~$s_2$ $\delta$-cross.
\end{proof}

\begin{proposition}\label{HighSpreadUniqueSplittingStar}
    Assume \cref{set:FlowerFinding}.
    Let
    \[
        \delta >\delta_{\ref{HighSpreadUniqueSplittingStar}}(q,k,\gamma_\cA,\eps_F,\gamma_N,\eps_N):=\max\,\{\,\delta_{\ref{NewerLargelyCrossingEntersAsterisk}}(q,k,\gamma_\cA,\eps_F), 
        \delta_{\ref{NewerLargeAsteriskUniqueSplittingStar}}(\gamma_N,\eps_N)\,\}
    \]
    and $\delta' \geq\delta'_{\ref{NewerUniqueAsterisk}}(k,\delta)$.
    Let $A\in\cA$ have spread $\ge \zeta_{\ref{NewerLargelyCrossingEntersAsterisk}}(\delta',\eps_F)$.
    Then $A$ is $\delta'$-relevant, $A^\delta$~exists and $A^\delta$ $\eps_N$-interlaces a unique splitting star of~$N$.
\end{proposition}
\begin{proof}
    We combine \cref{NewerLargelyCrossingEntersAsterisk} and \cref{NewerLargeAsteriskUniqueSplittingStar}.
\end{proof}

\subsection{Every splitting star is interlaced by at most one high-spread asterisk}

Here, we prove \cref{NewerDistinctSplittingStars}: if $A_1,A_2$ are two distinct asterisks of large enough spread, then $A_1^\delta$ and $A_2^\delta$ interlace distinct splitting stars.
For this, we first reduce to the case where $A_1^\delta$ $\gamma_N$-interlaces the petal-star of $F_{A_2}$ (see \cref{NewerDistinctSplittingStars:Claim}), which we then show is impossible.
Later, we will combine \cref{NewerDistinctSplittingStars} with \cref{HighSpreadUniqueSplittingStar} to obtain the flowers~$F_t$ in \cref{mainTutte}.

\begin{lem}\label{NewerDistinctSplittingStars:Claim}
    Assume \cref{set:FlowerFinding}.
    Let
    \[
        \delta>\delta_{\ref{NewerDistinctSplittingStars:Claim}}(q,k,\gamma_\cA,\eps_F,\gamma_N,\eps_N):=\max\,\{\,\delta_{\ref{HighSpreadUniqueSplittingStar}}(q,k,\gamma_\cA,\eps_F,\gamma_N,\eps_N),\;k+2\gamma_N+\eps_N+1\,\}
    \]
    and $\delta' \geq\delta'_{\ref{NewerUniqueAsterisk}}(k,\delta)$.
    Let $A_1\neq A_2\in\cA$ have spread $\ge \zeta_{\ref{NewerLargelyCrossingEntersAsterisk}}(\delta',\eps_F)$.
    So $A_1^\delta$ and $A_2^\delta$ exist by \cref{HighSpreadUniqueSplittingStar} and $\eps_N$-interlace unique splitting stars $\sigma_1$ and $\sigma_2$ of~$N$, respectively.
    Then either $A^\delta_1$ $\gamma_N$-interlaces $\pi(F_{A_2})$ or $\sigma_1\neq\sigma_2$.
\end{lem}

\begin{proof}
    Let $\pi_2$ denote the petal-star~$\pi(F_{A_2})$.
    Assume that $A^\delta_1$ does not $\gamma_N$-interlace~$\pi_2$.
    Let $r_1\in A^\delta_1$ witness this.
    Let $r_2\in A^\delta_2$.
    Recall that $A_i^\delta\se D(F_{A_i})$, so $r_i$ is a displayed separation of~$F_{A_i}$.
    
    Since every separation in $\pi_2$ is $\gamma_N$-nested with $r_1$ by our assumption \ref{set:FlowerFinding:S1totallyNested}, there is $p\in\pi_2$ such that $r_1\le_{\gamma_N} p$, say.
    We also have $p\le r_2$, say, as $p$ is a petal-separation of $F_{A_2}$ while $r_2$ is a displayed separation of~$F_{A_2}$.
    Hence $r_1\le_{\gamma_N} p \le r_2$.
    Since both $r_i$ are $2(\delta-k)$-broad by \cref{LargeCrossedImpliesBroad}, and since $2(\delta-k)\ge\gamma_N+\eps_N$, we have that $p$ is $\eps_N$-broad.
    So $p\in S_1$ is $\eps_N$-close to some $t\in N$ by~\ref{set:FlowerFinding:Napproximates}.
    Now $r_1\le_{\gamma_N+\eps_N} t\le_{\eps_N} r_2$ by \cref{approxLeQuasiTransitive}.

    \begin{claim}\label{NewerDistinctSplittingStars:ClaimClaim}
        We have $r_1 <_{\eps_N} t$ or $r^*_1 <_{\eps_N} t$.
    \end{claim}
    \begin{cproof}
        Assume for a contradiction that neither holds.
        Since $A^\delta_1$ contains $r_1$ and $\eps_N$-interlaces $\sigma_1$, and since $\sigma_1$ is a splitting star of~$N$ and $t \in N$, \cref{lem:almost-interlaces-then-nested-with-N} implies that either~$t<_{\eps_N} r_1$ or $t<_{\eps_N} r^*_1$. 

        \casen{Case $t<_{\eps_N} r_1$.}
        Then $r_1$ is $(\gamma_N+\eps_N)$-close to~$t$.
        Since $\delta-\gamma_N\ge\gamma_N+\eps_N$, this contradicts \cref{NewerLargeAsteriskNotCloseButBetter}.

        \casen{Case $t<_{\eps_N} r^*_1$.}
        Then $r_1 \le_{2\eps_N+\gamma_N} r^*_1$ by \cref{approxLeQuasiTransitive}. 
        Write $(U,W):=r_1$.
        Then $|U\sm W|\le 2\eps_N+\gamma_N$, so $r_1$ has breadth $\le 2\eps_N+\gamma_N$.
        But since $r_1$ is $\delta$-crossed, $r_1$ has breadth $\ge 2(\delta-k)$ by \cref{LargeCrossedImpliesBroad}.
        This yields a contradiction as $2(\delta-k)> 2\eps_N+\gamma_N$.
    \end{cproof}

    \begin{claim}
        We have $t<_{\eps_N} r_2$.
    \end{claim}
    \begin{cproof}
        As we already have $t\le_{\eps_N} r_2$, it remains to show $r_2\not\le_{\eps_N} t$.
        We cannot have $r_2\le_{\eps_N} t$, as then $r_2$ would be $\eps_N$-close to~$t$, contradicting \cref{NewerLargeAsteriskNotCloseButBetter} as 
        $\delta-\gamma_N\ge\eps_N$.
    \end{cproof}

    By the two claims, $O(r_1)$ and $O(r_2)$ (as in \cref{NewerLargeAsteriskUniqueSplittingStar}) contain different elements of $\{t,t^*\}$.
    Hence $\sigma_1\neq\sigma_2$.
\end{proof}

Let $\tau,\sigma$ be two stars of separations of the same graph~$G$.
We let $\tau\;\defnMath{\preceq}\;\sigma$ if for every $t\in\tau$ there is $s\in\sigma$ such that~$t\le s$, as in~\cite{ASS}.

\begin{lem}\label{NewerDistinctSplittingStars:TotallyNested}
    Assume \cref{set:FlowerFinding}.
    Let $\delta>\delta_{\ref{NewerDistinctSplittingStars:Claim}}$
    and $\delta' \geq\delta'_{\ref{NewerUniqueAsterisk}}(k,\delta)$.
    Let $A_1\neq A_2\in\cA$ both have spread $\ge \zeta_{\ref{NewerLargelyCrossingEntersAsterisk}}(\delta',\eps_F)$.
    So $A_1^\delta$ and $A_2^\delta$ exist by \cref{HighSpreadUniqueSplittingStar} and $\eps_N$-interlace unique splitting stars $\sigma_1$ and $\sigma_2$ of~$N$, respectively.
    Let~$S$ be the set of all narrow proper $k$-separations of~$G$.
    Assume that the common centre of $A_1$ is an $S$-totally-nested $k$-separator.
    Then $\sigma_1\neq \sigma_2$.
\end{lem}
\begin{proof}
    Let $F_i:=F_{A_i}$ for $i=1,2$.
    Both $A_i$ are $\delta'$-relevant by \cref{HighSpreadUniqueSplittingStar}.
    Hence there are $\delta'$-crossing displayed separations $r_i,r'_i\in D(F_i)$.
    By \cref{NewerUniqueAsterisk}, we have $r_i,r'_i\in A^\delta_i$.
    Moreover, by \cref{LargeCrossedImpliesBroad}, $r_i$ and $r'_i$ have breadth $\geq 2(\delta'-k) \geq 2\delta > 2\eps_N$.
    Let $X$ denote the common centre of~$A_1$.
    By assumption, $X$ is an $S$-totally-nested $k$-separator of~$G$.
    Every separation in $S_0$ is a $\beta$-broad $k$-separation of $G$, so $S_0 \subseteq S$ by \cref{FlowerFindingBroadSeparationsNice}.
    By definition of $X$, there is a separation in $S_0 \subseteq S$ with separator equal to~$X$, so there is a half-connected separation in $S$ with separator equal to~$X$.
    Thus, every half-connected separation in $S$ with separator equal to~$X$ is $S$-totally-nested by \cref{lem:equivalence-totally-nested-set}.
    Let $F$ be the graphic-$X$-anemone, and let $\pi_X \coloneqq \pi(F) \setminus \{(X,V(G))\}$ be the star of proper petal-separations of~$F$.
    Since there is a narrow separation with separator $X$, the neighbourhood of each component of $G-X$ is equal to $X$.
    Thus, every separation in $\pi_X$ is a narrow proper $k$-separation with separator $X$, so is $S$-totally-nested.
    Then $\pi_X\se N$ by \ref{set:FlowerFinding:NcontainsTotallyNested}.
    As both $r_1,r'_1$ have $X$ as their separator, both $r_1,r'_1$ interlace~$\pi_X$.

    Assume for a contradiction that $\sigma_1=\sigma_2$.
    Then both $r_2,r'_2$ $\eps_N$-interlace~$\sigma_1$.
    Both $r_1,r'_1$ $\eps_N$-interlace~$\sigma_1$, and both interlace~$\pi_X$.
    By the latter, we have $\pi_X\se O(r_1)=O(r'_1)$ (as in \cref{NewerLargeAsteriskUniqueSplittingStar:eq:Orientation}) by \cref{NewerLargeAsteriskUniqueSplittingStar}, which yields $\pi_X\preceq\sigma_1$.
    Therefore, also $r_2,r'_2$ $\eps_N$-interlace~$\pi_X$.
    Let $p \in \pi_X$ be arbitrary. Since $r_2$ $\eps_N$-interlaces $\pi_X$, we have $p <_{\eps_N} r_2$ or $p <_{\eps_N} (r_2)^*$, say the former. Then, $r_2 \not\le_{\eps_N} p$. Moreover, if $(r_2)^* \leq_{\eps_N} p$ then we would obtain $(r_2)^* \leq_{2\eps_N}r_2$ by \cref{approxLeQuasiTransitive}, which by \cref{lem:assymmetric-nested-implies-not-broad} would contradict that $r_2$ has breadth $> 2\eps_N$.
    Therefore, we have $r_2,r_2^*\not\le_{\eps_N} p$, and similarly for~$r'_2$.
    Hence, since the separations in $\pi_X$ are $S$-totally-nested, we even get that $r_2,r'_2$ interlace~$\pi_X$ (without~`$\eps_N$').
    Hence $r_2,r'_2\in D(F)$.
    By~\ref{set:FlowerFinding:XAnemone}, we have $F=F_1$.
    Thus $r_2,r'_2\in D(F_1)$.
    As $r_2,r'_2$ $\delta'$-cross, we have $r_2,r'_2\in A_1^\delta$ by \cref{NewerUniqueAsterisk}.
    Since $A_1^\delta\se A_1$ by \cref{NewerLargelyCrossingEntersAsterisk},
    we have $r_2,r'_2\in A_1$, contradicting that $A_1\neq A_2$.
\end{proof}

\begin{lem}\label{totallyNestedSeparatorEasyCrossing}
    Let~$s$ be a $k$-separation of~$G$, and let~$t$ be a narrow proper $k$-separation of~$G$ that crosses~$s$.
    Assume that~$t$ does not separate vertices of~$\sep(s)$.
    Then $\sep(t)=\sep(s)$.
\end{lem}

\begin{proof}
    Let $s =: (A, B)$ and $t =: (C, D)$.
    Assume for a contradiction that $\sep(t) \neq \sep(s)$.
    Since~$s$ and~$t$ are both $k$-separations, we have that $\sep(s) \setminus \sep(t)$ and $\sep(t) \setminus \sep(s)$ are both nonempty.
    Without loss of generality, assume that the $At$-link and the $sC$-link are nonempty.
    Since~$t$ does not separate vertices of~$\sep(s)$, the $sD$-link is empty.
    Since~$t$ is narrow and the $At$-link is nonempty, the $BD$-corner must be empty.
    Thus, the $AD$-corner is nonempty.
    Since~$t$ is narrow, this implies that the $Bt$-link is also empty, contradicting that~$s$ and~$t$ cross.
\end{proof}

\begin{proposition}\label{NewerDistinctSplittingStars}
    Assume \cref{set:FlowerFinding}.
    Let $\delta>\delta_{\ref{NewerDistinctSplittingStars:Claim}}$ and
    \[
        \delta' > \delta'_{\ref{NewerDistinctSplittingStars}}(q,k,b,\gamma_N,\delta):=\max\,\{\,\delta'_{\ref{NewerUniqueAsterisk}}(k,\delta),\enspace\beta_{\ref{NewerInterlacingPetalStar}}(q,k,b,\gamma_N)+k,\enspace\delta+2\eps_{\ref{NewerInterlacingPetalStar}}(q,k,b,\gamma_N)\,\}.
    \]
    Let $A_1\neq A_2\in\cA$ have spread $\ge \zeta_{\ref{NewerLargelyCrossingEntersAsterisk}}(\delta',\eps_F)$.
    So $A_1^\delta$ and $A_2^\delta$ exist by \cref{HighSpreadUniqueSplittingStar} and $\eps_N$-interlace unique splitting stars $\sigma_1$ and $\sigma_2$ of~$N$, respectively.
    Then $\sigma_1\neq\sigma_2$.
\end{proposition}

\begin{proof}
    Let $F_i:=F_{A_i}$ for $i=1,2$.
    Both $A_i$ are $\delta'$-relevant by \cref{HighSpreadUniqueSplittingStar}.
    Hence there are $\delta'$-crossing displayed separations $r_i,r'_i\in D(F_i)$.
    By \cref{NewerUniqueAsterisk}, we have $r_i,r'_i\in A^\delta_i$.
    Let $\pi_2:=\pi(F_2)$.
    By \cref{NewerDistinctSplittingStars:Claim}, we may assume that both $r_1$ and $r'_1$ $\gamma_N$-interlace $\pi_2$.
    Both $r_1,r'_1$ have breadth
    \[
    \overset{\ref{LargeCrossedImpliesBroad}}{\ge} 2(\delta' - k)\ge \beta_{\ref{NewerInterlacingPetalStar}}>b.
    \]
    Note that both $r_1,r'_1$ are narrow as $G$ is $(b,k)$-narrow.
    We apply \cref{NewerInterlacingPetalStar} to both $r_1,r'_1$, which leads to one of two outcomes.
    Let $S$ be the set of all narrow proper $k$-separations of~$G$.

    \casen{Outcome~1: the separator of $r_1$ or $r'_1$ is an $S$-totally-nested $k$-separator of~$G$.}
    Without loss of generality, the separator $X$ of $r_1$ is $S$-totally-nested.
    Since $r'_1$ crosses~$r_1$ and both $r_1,r'_1$ have order~$k$,
    the separator of $r'_1$ equals $X$ by \cref{totallyNestedSeparatorEasyCrossing}.
    Hence $X$ is the common centre of~$A_1$.
    Now \cref{NewerDistinctSplittingStars:TotallyNested} yields $\sigma_1\neq\sigma_2$.

    \casen{Outcome~2: $r_1,r'_1$ are $\eps_{\ref{NewerInterlacingPetalStar}}$-close to displayed separations $d,d'\in D(F_2)$, respectively.} 
    Since $r_1$ and $r'_1$ $\delta'$-cross, \cref{lem:close-separations-still-cross} yields that $d$ and $d'$ $(\delta'-2\eps_{\ref{NewerInterlacingPetalStar}})$-cross.
    As $\delta'-2\eps_{\ref{NewerInterlacingPetalStar}}\ge\delta$, we find that $d,d'\in D(F_2)$ $\delta$-cross.
    Hence $d,d'\in A_2$ by \cref{NewerLargelyCrossingEntersAsterisk}.
    Since $r_1$ and $r'_1$ $\delta'$-cross and $r'_1$ is $\eps_{\ref{NewerInterlacingPetalStar}}$-close to $d'$, we have that $r_1$ and $d'$ $(\delta'-\eps_{\ref{NewerInterlacingPetalStar}})$-cross.
    Using 
    \begin{equation}\label{deltaGeqGammaA}
        \delta'-\eps_{\ref{NewerInterlacingPetalStar}}\ge\delta\ge\delta_{\ref{NewerDistinctSplittingStars:Claim}}\ge\delta_{\ref{HighSpreadUniqueSplittingStar}}\ge\delta_{\ref{NewerLargelyCrossingEntersAsterisk}}\ge\gamma_\cA,
    \end{equation}
    we find that $r_1\in A_2$.
    But also $r_1\in A_1$ by \cref{NewerLargelyCrossingEntersAsterisk}.
    Hence $A_1=A_2$, a~contradiction.
\end{proof}

\subsection{Interlacing separations and torso connectivity}

We show that if a splitting star~$\sigma$ is $\eps_N$-interlaced by a sufficiently broad $k$-separation $s$ that is not too close to~$N$, then there is a unique high-spread asterisk~$A$ such that $A^\delta$ $\eps_N$-interlaces~$\sigma$.
Moreover, $s$ belongs to~$A$ and is $\eps_F$-close to a displayed separation of~$F_A$ (see \cref{NewerInterlacesCloseToDisplayed}).
If no such separation exists, we show that the torso of~$\sigma$ is $(q',k)$-unbreakable for a suitable~$q'$ (see \cref{highlyConnectedTorso}).

\begin{lem}\label{lem:separations-almost-interlace-splitting-stars}
    Let~$N$ be a symmetric nested set of proper separations of a graph~$G$, and let~$\sigma \neq \sigma'$ be two splitting stars of~$N$.
    Let~$s,s'$ be two separations of~$G$ that almost-$\eps$-interlace~$\sigma,\sigma'$, respectively.
    Then there is a $t \in \sigma'$ such that
    \begin{enumerate}
        \item $s \leq_\eps t$ or $s^* \leq_\eps t$, and
        \item $t \leq_\eps s'$ or $t \leq_\eps (s')^*$.
    \end{enumerate}
\end{lem}
\begin{proof}
    Since $\sigma \neq \sigma'$, we can use \cref{TDCfromNestedSetConstruct} and \cref{TDCfromNesteSetVerify} to find $r \in \sigma$ and $r' \in \sigma'$ such that $r^* \leq r'$. 
    Since~$s$ almost-$\eps$-interlaces~$\sigma$, we have $r \leq_\eps s$ or $r \leq_\eps s^*$.
    So, $s^* \leq_\eps r^*$ or $s \leq_\eps r^*$.
    With \cref{approxLeQuasiTransitive}, we have $s^* \leq_\eps r'$ or $s \leq_\eps r'$.
    Similarly, we get $r' \leq_\eps s'$ or $r' \leq_\eps (s')^*$.
    So, for $t:=r'$ the statement follows.
\end{proof}

\begin{keylemma}\label{NewerInterlacesCloseToDisplayed}
    Assume \cref{set:FlowerFinding}.
    Let $\delta>\delta_{\ref{NewerDistinctSplittingStars:Claim}}$, $\delta' >\delta'_{\ref{NewerDistinctSplittingStars}}(\ldots ,\delta)$,
    \begin{align*}
        \zeta \geq \; \zeta_{\ref{NewerInterlacesCloseToDisplayed}}(k,\delta',\eps_F) := {} & \zeta_{\ref{NewerLargelyCrossingEntersAsterisk}}(\delta'+2k+1,\eps_F), \text{ and }\\
        \eps_{\ref{NewerInterlacesCloseToDisplayed}}(k,\eps_F,\eps_N,\zeta) := {} & \,2\zeta+3\eps_F+\eps_N+3k.
    \end{align*}
    Let $\sigma$ be a splitting star of~$N$.
    Assume that $\sigma$ is $\eps_N$-interlaced by an $\eps_{\ref{NewerInterlacesCloseToDisplayed}}$-broad $k$-separation $s$ of~$G$.
    Assume that $s$ is not $\eps_{\ref{NewerInterlacesCloseToDisplayed}}$-close to any separation in~$N$.
    Then we have all of the following:
    \begin{enumerate}
        \item\label{NewerInterlacesCloseToDisplayed:Adelta} There exists a unique $A\in\cA$ such that $A$ has spread~$\ge\zeta$, $A$ is $\delta'$-relevant and $A^\delta$ $\eps_N$-interlaces~$\sigma$.
        Moreover, $F_A$ has spread~$\ge \zeta$.
        \item\label{NewerInterlacesCloseToDisplayed:InA} $s\in A$; in particular, $s$ is $\eps_F$-close to a displayed separation of $F_A$ by~\ref{set:FlowerFinding:CloseFlower}.
        \item\label{NewerInterlacesCloseToDisplayed:LargelyCrossed} $s$ is $\zeta$-crossed by some separation in $A\cap D(F_A)$.
    \end{enumerate}
\end{keylemma}
\begin{proof}
    First, we claim that $s$ is $\gamma_\cA$-crossed by some separation in~$S_0$.
    Assume for a contradiction that it is not.
    Note that $s\in S_0$ as
    $\eps_{\ref{NewerInterlacesCloseToDisplayed}}\ge \eps_N\ge\beta$ by \ref{set:FlowerFinding:Inequalities}.
    Then $\{s,s^*\}$ forms a trivial $\gamma_\cA$-asterisk of~$S_0$, so $s\in S_1$ by definition of~$S_1$.
    Since $N$ $\eps_N$-approximates~$S_1\ni s$ by \ref{set:FlowerFinding:Napproximates} and since $s$ has breadth~$\ge\eps_{\ref{NewerInterlacesCloseToDisplayed}}\ge\eps_N$, there is $t\in N$ that is $\eps_N$-close to~$s$.
    But since $\eps_N\le\eps_{\ref{NewerInterlacesCloseToDisplayed}}$, this contradicts our assumption that $s$ is not $\eps_{\ref{NewerInterlacesCloseToDisplayed}}$-close to any separation in~$N$.

    Since $s\in S_0$ is $\gamma_\cA$-crossed by some separation in~$S_0$, there is an $A\in\cA$ such that~$s\in A$.
    Let $d$ be a displayed separation of $F_A$ that is $\eps_F$-close to~$s$, which exists by~\ref{set:FlowerFinding:CloseFlower}.
    We apply \cref{FlowerClosePetalOrCrossed} with~$\xi:=\zeta+\eps_F+k$,
    which yields one of the following two outcomes:

    \casen{Outcome~1: $d$ is $2\xi$-close to some petal-separation $p$ of~$F_A$ or its inverse.}
    We will show that this outcome is impossible under our assumptions.
    If $p$ is $\eps_N$-broad, then $p$ is $\eps_N$-close to some $t\in N$ by~\ref{set:FlowerFinding:Napproximates}, and then $s$ is $(2\xi + \eps_N + \eps_F)$-close 
    to~$t$ or $t^*$, both of which belong to $N$ since $N$ is symmetric, by \cref{approxLeQuasiTransitive}, contradicting our assumptions as $2\xi+\eps_N+\eps_F\le \eps_{\ref{NewerInterlacesCloseToDisplayed}}$.
    Hence we may assume that $p$ has breadth~$<\eps_N$.
    Then $s$ has breadth~$<2\xi+\eps_N+\eps_F+k\le\eps_{\ref{NewerInterlacesCloseToDisplayed}}$ 
    by \cref{LowBreadthClose}, contradicting our assumption that $s$ has breadth $\ge\eps_{\ref{NewerInterlacesCloseToDisplayed}}$.

    \casen{Outcome~2: $F_A$ has a displayed separation $d'$ that $(\xi -k)$-crisses~$d$.}
    Then $d$ and $d'$ $(\xi-k)$-cross by \cref{crissCross}~\ref{crissCross:crissImpliesCross}.
    Since $\xi-k \geq \zeta$, it follows that $F_A$ has spread $\geq \zeta$.
    Note that $d,d' \in S_0$ since they are $2(\xi-k)$-broad (with $2(\xi-k) \geq 2\zeta_{\ref{NewerLargelyCrossingEntersAsterisk}}(\delta'+2k+1,\eps_F) \geq 2\delta' \geq 2\delta \geq 2\gamma_N \geq \beta$ by \ref{set:FlowerFinding:Inequalities}).
    As $s$ is $\eps_F$-close to~$d$, we have that $d'$ $(\xi-k-\eps_F)$-crosses~$s$ by \cref{CloseVsCross}.
    As 
    \begin{equation}\label{NewerInterlacesCloseToDisplayed:Inequalities}
        \xi-k-\eps_F
        \ge \zeta \ge \zeta_{\ref{NewerInterlacesCloseToDisplayed}}
        \ge \delta'
        >\delta'_{\ref{NewerDistinctSplittingStars}}        >\delta>\delta_{\ref{NewerDistinctSplittingStars:Claim}}\ge\delta_{\ref{HighSpreadUniqueSplittingStar}}\ge\delta_{\ref{NewerLargelyCrossingEntersAsterisk}}\ge\gamma_\cA
    \end{equation}
    and $s\in A$, we have $d'\in A$ (which shows~\cref{NewerInterlacesCloseToDisplayed:LargelyCrossed}), and so $A$ has spread~$\ge\zeta$.
    Hence $A$ is $\delta'$-relevant, $A^\delta$ exists, and $A^\delta$ $\eps_N$-interlaces a unique splitting star~$\sigma'$ of~$N$ by \cref{HighSpreadUniqueSplittingStar}
    (using $\delta>\delta_{\ref{HighSpreadUniqueSplittingStar}}$ from \cref{NewerInterlacesCloseToDisplayed:Inequalities}).
    Note that $d, d' \in D(F_A)$ $(\xi-k)$-cross with $\xi-k \geq \delta'$ as above, so both belong to $A^{\delta}$ by definition.

    We claim that $\sigma'=\sigma$.
    Assume for a contradiction that $\sigma'\neq\sigma$.
    By \cref{lem:separations-almost-interlace-splitting-stars}, there is $t\in N$ such that $s\le_{\eps_N} t\le_{\eps_N} d'$, say.
    Hence $s\le_{2\eps_N}d'$ by \cref{approxLeQuasiTransitive}.
    But we have seen above that $s$ and $d'$ $(\xi-k-\eps_F)$-cross, so they $\zeta$-cross and 
    \[
        \zeta
        \overset{\cref{NewerInterlacesCloseToDisplayed:Inequalities}}{>}
        \delta_{\ref{HighSpreadUniqueSplittingStar}}\ge\delta_{\ref{NewerLargeAsteriskUniqueSplittingStar}}\ge 2\eps_N,
    \]
    a contradiction.
    Therefore, $\sigma'=\sigma$.
    
    The uniqueness of $A$ in \cref{NewerInterlacesCloseToDisplayed:Adelta} is given by \cref{NewerDistinctSplittingStars}.
\end{proof}

\begin{definition}[Lift, induced separation]
Let $\sigma$ be a star of separations of a graph~$G$.
Let $H$ be the torso of~$\sigma$.
Let $s=(A,B)$ be a separation of~$H$, and let $\hat s=(\hat A,\hat B)$ be a separation of~$G$.
Assume that $A=\hat A\cap V(H)$, $B=\hat B\cap V(H)$ and $\sep(s)=\sep(\hat s)$.
Then we call $\hat s$ a \defn{lift} of~$s$ (from $H$ to~$G$).
Conversely, we say that $s$ is \defn{induced} by~$\hat s$ (on the torso~$H$).
\end{definition}

\begin{lemma}\label{highlyConnectedTorso}
    Let $\sigma$ be a star of $k$-separations of a $(q,k-1)$-unbreakable graph~$G$.
    Let $q',\eps\in\N$ such that $q'\ge\eps+k$.
    Assume that there is no $(\le k)$-separation $s$ of $G$ such that
    $s$ almost-interlaces~$\sigma$,
    $s$ has breadth~$>q'$,
    and $s$ is not $\eps$-close to any separation in~$\sigma\cup\sigma^*$.
    Then the torso of $\sigma$ is $(q',k)$-unbreakable.
\end{lemma}
\begin{proof}
    Let $H$ be the torso of~$\sigma$.
    Assume for a contradiction that $H$ fails to be $(q',k)$-unbreakable.
    Then $H$ has a $(\le k)$-separation $(A,B)$ of breadth~$>q'$.
    Let $(\hat A,\hat B)$ be a lift of $(A,B)$ to~$G$ such that $(\hat{A}, \hat{B})$ is a $(\le k)$-separation of $G$ that almost-interlaces~$\sigma$, which exists by \cite[2.6.4, its proof and the definition below]{TriAiC}.
    The lift $(\hat A,\hat B)$ has breadth~$>q'$ by construction.

    \begin{claim}
        No $r\in\sigma\cup\sigma^*$ is $\eps$-close to~$(\hat A,\hat B)$.
    \end{claim}

    \begin{cproof}
        Let $r \eqqcolon (U, W) \in \sigma$.
        Since $H$ is the torso of $\sigma$, we have $V(H) \subseteq W$.
        Then, \begin{align*}
            |\hat{A} \setminus U| &\geq |A \setminus U| \geq |A \cap W| - |U \cap W| \geq |A \cap V(H)| - |U \cap W| > q' - k \geq \eps, \text{ and} \\
            |\hat{B} \setminus U| &\geq |B \setminus U| \geq |B \cap W| - |U \cap W| \geq |B \cap V(H)| - |U \cap W| > q' - k \geq \eps.
        \end{align*}
        The first inequality rules out $(\hat{A}, \hat{B}) \approx_{\eps}(U, W)$ and the second rules out $(\hat{A}, \hat{B}) \approx_{\eps}(W, U)$.
    \end{cproof}
    
    Hence $s:=(\hat A,\hat B)$ contradicts our assumptions.
\end{proof}

\subsection{Closeness of petal-stars and splitting stars}

The next two key lemmas (\ref{NewerPetalSeparationsCloseToSplittingStar} and \ref{NewerSplittingStarCloseToPetalStar}) show that the petal-star of each constructed graphic-flower and its corresponding splitting star approximate each other.

Recall that when $(H,\cV)$ is an $H$-decomposition of a graph~$G$ and $X\se H$, then we write $V\rest X$ for $\bigcup_{x\in V(X)}V_x$.

\begin{lem}\label{NewerTLemma}
    Let $\delta\in\N$.
    Let $(O,\cV)$ be a cycle-decomposition of a graph~$G$.
    Let $P\se O$ be a path and $p\in P$ a vertex.
    Assume that $|V_p|>\delta$.
    Then we have at least one of the following:
    \begin{enumerate}
        \item\label{NewerTLemma:LargeCross} There is an edge $xy$ on~$P$ (with $x$ before~$y$) such that both $V\rest Px$ and $V\rest yP$ have size $>\delta$.
        \item\label{NewerTLemma:DoubleSmall} Both $V\rest P\mathring{p}$ and $V\rest \mathring{p}P$ have size $\le\delta$.
    \end{enumerate}
\end{lem}
\begin{proof}
    Assume that \cref{NewerTLemma:DoubleSmall} fails.
    By symmetry, we may assume without loss of generality that $|V\rest\mathring{p}P|>\delta$.
    We also have $|V\rest Pp|>\delta$ as $|V_p|>\delta$ by assumption.
    Let $x:=p$ and let $y$ be the successor of~$p$ on the path~$P$.
    Then $xy$ satisfies~\cref{NewerTLemma:LargeCross}.
\end{proof}

\begin{keylemma}\label{NewerPetalSeparationsCloseToSplittingStar}
    Assume \cref{set:FlowerFinding}.
    Let $\delta>\delta_{\ref{HighSpreadUniqueSplittingStar}}$ and $\delta'>\delta'_{\ref{NewerUniqueAsterisk}}(k,\delta)$.
    Let $\eps_{\ref{NewerPetalSeparationsCloseToSplittingStar}}(k,\eps_N,\delta'):=2(k+\eps_N+\delta')$.
    Let $A\in\cA$ have spread
    $\ge\zeta_{\ref{NewerInterlacesCloseToDisplayed}}(k,\delta',\eps_F)$.
    So $A^\delta$ exists by \cref{HighSpreadUniqueSplittingStar} and $\eps_N$-interlaces a unique splitting star $\sigma$ of~$N$.
    Let $\pi$ denote the petal-star of~$F_A$.
    Then $\sigma$ $\eps_{\ref{NewerPetalSeparationsCloseToSplittingStar}}$-approximates~$\pi$.
\end{keylemma}
\begin{proof}
    Choose $s_1, s_2 \in A$ that $\zeta_{\ref{NewerInterlacesCloseToDisplayed}}(k,\delta',\eps_F)$-cross, that is, $(\delta'+2k+1+2\eps_F)$-cross.
    By \ref{set:FlowerFinding:CloseFlower}, there exist $s'_1, s'_2 \in D(F_A)$ such that $s_i$ and $s'_i$ are $\eps_F$-close for $i \in \{1, 2\}$.
    Then, $s'_1$ and $s'_2$ $(\delta'+2k+1)$-cross by \cref{lem:close-separations-still-cross}.
    If~$F_A$ is a graphic-$k$-anemone, let $F\eqqcolon(O,\cV)$ be the $k$-anemone given by \cref{AnemoneAsDaisy} that displays both $s'_1$ and $s'_2$.
    Otherwise, $F_A$ is a graphic-$k$-daisy, and we let $(O,\cV):=F:=F_A$.
    In both cases, we have $D(F) \subseteq D(F_A)$ and $s'_1, s'_2 \in D(F)$. Moreover, $s'_1, s'_2 \in A^{\delta}$ by \cref{NewerUniqueAsterisk}.
    
    Let $p\in\pi$ be $\eps_{\ref{NewerPetalSeparationsCloseToSplittingStar}}$-broad.
    Let $x\in O$ be the node that induces the petal-separation~$p$.
    Let $E^\delta$ be the set of all edges of $O$ that are used by 2-edge-cuts of $O$ to define displayed separations that are contained in~$A^\delta$.
    Note that $|E^\delta|\ge 4$ as $A^\delta$ is nontrivial.
    Let $P$ denote the component of $O-E^\delta$ that contains~$x$, and note that $P$ is a path.
    Let $Q:=O-P$, which is a path with at least two edges as $|E^\delta| \ge 4$.

    \begin{claim}\label{NewerPetalSeparationsCloseToSplittingStar:edgeBigParts}
        There are two edges $uv,u'v'$ on $Q$ (with $u$ before~$v$ before~$u'$ before~$v'$) such that $V\rest Qu$, $V \rest vQu'$ and $V\rest v'Q$ all have size~$>\delta'+k$.
    \end{claim}
    \begin{cproof}
        By our choice of $F$, the separations $s'_1, s'_2 \in D(F) \cap A^{\delta}$ $(\delta'+2k+1)$-cross.
        The corresponding edge-cuts of $O$ therefore give four distinct edges $e_1, e_2, e_3, e_4 \in E^{\delta}$ whose deletion splits $O$ into four consecutive paths $P_1,\ldots,P_4$ such that $V\rest P_i$ has size~$\ge (\delta'+2k+1)-k > \delta'+k$ by \cref{crissCross}.
        
        Since all four edges $e_1,\ldots,e_4$ are contained in~$E^\delta$, it follows that~$P_1$, say, includes the entire path~$P$.
        The other three paths avoid $P$ and hence are included in~$Q$.
        Then we let $uv$ be the edge on $O$ in between $P_2$ and~$P_3$ with $u$ coming before $v$ on~$Q$, and let $u'v'$ be the edge on $O$ in between $P_3$ and~$P_4$ with $u'$ coming before $v'$ on~$Q$.
    \end{cproof}
    
    Consider the displayed separation $s:=(V\rest P,V\rest Q)$ of~$F$.

    \begin{claim}\label{NewerPetalSeparationsCloseToSplittingStar:notDeltaCrossed}
        No displayed separation of $F$ $\delta'$-crosses~$s$.
    \end{claim}
    \begin{cproof}
        Assume for a contradiction that some $d\in D(F)$ $\delta'$-crosses~$s$.
        Since $s, d \in D(F) \subseteq D(F_A)$ and $s, d$ $\delta'$-cross, \cref{NewerUniqueAsterisk} implies that $s, d \in A^{\delta}$.
        Consider the 2-edge-cuts of $O$ that induce $s$ and~$d$, respectively.
        These must cross.
        For this, the 2-edge-cut that induces $d$ must use an edge of~$P$.
        But since $d\in A^\delta$, this edge also has to be contained in~$E^\delta$, a contradiction.
    \end{cproof}

    Since $x\in P$ we have $p\le s$.
    Since $p$ has breadth~$>\delta'+k$, we have $|V_x|> \delta'+k$.
    We apply \cref{NewerTLemma} (with $\delta'+k$), which has two possible outcomes.

    \casen{Case~1: \cref{NewerTLemma} yields outcome~\cref{NewerTLemma:LargeCross}.}
    We claim that this outcome is impossible.
    Assume for a contradiction that we get~it, that is, there is an edge $yz$ on $P$ (with $y$ before $z$ on $P$) such that both $V\rest Py$ and $V\rest zP$ have size~$>\delta'+k$.
    Let $uv\in E(Q)$ be as provided by \cref{NewerPetalSeparationsCloseToSplittingStar:edgeBigParts}.
    Let $d\in D(F)$ be induced by the 2-edge-cut $\{yz,uv\}$ of~$O$.
    Then $d$ $\delta'$-crisses~$s$, and so $d$ $\delta'$-crosses~$s$ by \cref{crissCross}.
    This contradicts \cref{NewerPetalSeparationsCloseToSplittingStar:notDeltaCrossed}.

    \casen{Case~2: \cref{NewerTLemma} yields outcome~\cref{NewerTLemma:DoubleSmall}}.
    That is to say that $V\rest P\mathring{x}$ and $V\rest\mathring{x}P$ have size~$\le\delta'+k$.
    Hence $p\le s$ and $s\le_{2(\delta'+k)}p$, so
    \begin{equation}\label{eq:NewerPetalSeparationsCloseToSplittingStar:pCloseTos}
        p\approx_{2(\delta'+k)} s.
    \end{equation}

    As $p\in S_1$ has breadth~$\ge\eps_{\ref{NewerPetalSeparationsCloseToSplittingStar}}\ge\eps_N$, by \ref{set:FlowerFinding:Napproximates} there is $t\in N$ such that $t$ is $\eps_N$-close to~$p$.
    If $t\in\sigma$, then we are done as $\eps_N\le\eps_{\ref{NewerPetalSeparationsCloseToSplittingStar}}$.
    Hence we may assume that $t\notin\sigma$.

    We claim that there are $d_1,d_2\in A^\delta$ such that $s=d_1\wedge d_2$. 
    By \cref{NewerPetalSeparationsCloseToSplittingStar:edgeBigParts} there are two edges $uv,u'v'$ of~$Q$ (with $u$ before~$v$ before~$u'$ before~$v'$) such that all of $V \rest Qu$, $V \rest vQu'$ and $V \rest v'Q$ have size~$> \delta' + k$.
    Let $yz,y'z'$ be the two edges that together induce the separation~$s$, and are such that $y,z,y',z',u,v,u',v'$ appear in this cyclic order along~$O$.
    Then let~$d_1$ be induced by $\{yz,uv\}$ and let~$d_2$ be induced by~$\{y'z',u'v'\}$ such that $s=d_1\wedge d_2$.
    Note that they $\delta'$-criss, hence they $\delta'$-cross by \cref{crissCross}.
    Hence, $d_1,d_2 \in A^\delta$.
    
    Since $p\le s$, we have $p\le d_1,d_2$ and, moreover, $t\le_{\eps_N} d_1,d_2$.
    With the notations of \cref{NewerLargeAsteriskUniqueSplittingStar}, we have $t \in O(d_1)$, so there exists $r \in \sigma(d_1)$ such that $t \leq r$. Since $d_1 \in A^{\delta}$, \cref{NewerLargeAsteriskUniqueSplittingStar} implies that $\sigma = \sigma(d_1)$, so $r \in \sigma$.
    Since $A^\delta$ $\eps_N$-interlaces~$\sigma$, we have $r<_{\eps_N}d_i$ or $r<_{\eps_N}d^*_i$ for $i=1,2$.
    We claim that $r<_{\eps_N}d_1,d_2$, and assume for a contradiction that $r<_{\eps_N}d^*_1$, say.
    Then
    \[
        p\le_{\eps_N} t\le r<_{\eps_N} d^*_1
    \]
    which implies $p\le_{2\eps_N}d^*_1$ by \cref{approxLeQuasiTransitive}.
    But also $p\le d_1$, which is equivalent to $d^*_1\le p^*$, and combines with the previous sentence to $p\le_{2\eps_N}p^*$.
    Thus \cref{lem:assymmetric-nested-implies-not-broad} yields that $p$ has breadth~$\le 2\eps_N$, contradicting that $p$ has breadth~$\ge\eps_{\ref{NewerPetalSeparationsCloseToSplittingStar}}>2\eps_N$.
    Therefore, we have $r<_{\eps_N}d_1,d_2$ as claimed.
    Hence $r\le_{2\eps_N}d_1\wedge d_2=s$ by \cref{approxLeQuasiInfimum}.
    Therefore,
    \[
        p\le_{\eps_N} t \le r\le_{2\eps_N} s\overset{\cref{eq:NewerPetalSeparationsCloseToSplittingStar:pCloseTos}}{\le}_{2(\delta'+k)} p.
    \]
    Hence $p$ is $2(\delta'+k+\eps_N)$-close to $r\in\sigma$ by \cref{approxLeQuasiTransitive}.
\end{proof}

\begin{lem}\label{fishBroad}
    Let $(A,B)\le r\wedge s$, where $r$ and $s$ are $k$-separations and~$(A,B)$ is a separation of a graph~$G$.
    Assume that $r$ and $s$ $\gamma$-cross.
    Then $|B\sm A|\ge 3(\gamma-k)$.
\end{lem}
\begin{proof}
    By \cref{crissCross}, $r$ and~$s$ $(\gamma-k)$-criss, so all four corners have size~$\ge \gamma-k$.
    Three of these corners are included in~$B\sm A$.
\end{proof}

\begin{lem}
    \label{lem:petal-separation-of-graphic-anemone-narrow}
    Let~$G$ be a $(b,k)$-narrow graph and $X \subseteq V(G)$ of size $|X|=k$.
    Assume that the graphic-$X$-anemone $\cA$ has a $b$-broad displayed separation.
    Then every petal-separation of~$\cA$ is narrow.
\end{lem}
\begin{proof}
    Let~$s$ be a $b$-broad displayed separation of the graphic-$X$-anemone~$\cA$.
    Note that $X=\sep(s)$.
    Since~$G$ is $(b,k)$-narrow, the separation~$s$ is narrow.
    Hence, every component~$K$ of $G-X$ has $N(K) = X$.
    In particular, every petal-separation of~$\cA$ is narrow.
\end{proof}

\begin{keylemma}\label{NewerSplittingStarCloseToPetalStar}
    Assume \cref{set:FlowerFinding}.
    Let $\delta>\delta_{\ref{HighSpreadUniqueSplittingStar}}$, $\delta'>\max\{\delta'_{\ref{NewerUniqueAsterisk}}(k,\delta), \,\delta+2\eps_{\ref{NewerInterlacingPetalStar}}(q,k,b,\gamma_N)\}$ and
    \[
    \eps_{\ref{NewerSplittingStarCloseToPetalStar}}(q,k,b,\gamma_N,\eps_N,\delta'):=\max\{\,3\eps_{\ref{NewerInterlacingPetalStar}}(q,k,b,\gamma_N)+2\gamma_N+2k,\,\eps_{\ref{NewerPetalSeparationsCloseToSplittingStar}}+\gamma_N+1,\,\beta_{\ref{NewerInterlacingPetalStar}}(q,k,b,\gamma_N)\,\}.
    \]
    Let $A\in\cA$ have spread $\ge\zeta_{\ref{NewerInterlacesCloseToDisplayed}}(k,\delta',\eps_F)$.
    So $A$ is $\delta'$-relevant, $A^\delta$ exists by \cref{HighSpreadUniqueSplittingStar} and $\eps_N$-interlaces a unique splitting star $\sigma$ of~$N$.
    Let $\pi$ denote the petal-star of~$F_A$.
    Then $\pi$ $\eps_{\ref{NewerSplittingStarCloseToPetalStar}}$-approximates~$\sigma$.
\end{keylemma}
\begin{proof}
    Let $s \eqqcolon (U, W)\in\sigma$.
    We may assume that $s$ is $\eps_{\ref{NewerSplittingStarCloseToPetalStar}}$-broad, as otherwise there is nothing to show.
    Since $s \in \sigma \subseteq N$, $s$ is a $k$-separation of $G$ so $s$ is $\gamma_N$-nested with all the separations in~$\pi$ by \ref{set:FlowerFinding:S1totallyNested}.
    
    Hence we are in one of the following two cases:

    \casen{Case~1: $s$ almost-$\gamma_N$-interlaces~$\pi$.}
    We assumed earlier that $s$ is $\eps_{\ref{NewerSplittingStarCloseToPetalStar}}$-broad, and we have $\eps_{\ref{NewerSplittingStarCloseToPetalStar}}\ge \beta_{\ref{NewerInterlacingPetalStar}}$.
    We apply \cref{NewerInterlacingPetalStar}, which has two possible outcomes:

    \casen{Case~1a: $s$ is $\eps_{\ref{NewerInterlacingPetalStar}}$-close to some displayed separation~$d$ of~$F_A$.}
    Since $s$ is narrow-$\gamma_N$-nested by \ref{set:FlowerFinding:NtotallyNested}, and $d$ is $\eps_{\ref{NewerInterlacingPetalStar}}$-close to~$s$, \cref{CloseVsCross} yields that $d$ cannot be $(\eps_{\ref{NewerInterlacingPetalStar}}+\gamma_N)$-crossed by any narrow proper $k$-separation.
    Even more, since $\gamma_N \geq \beta + k \geq b$ by \ref{set:FlowerFinding:Inequalities}, and using \cref{crissCross}, we have that~$d$ cannot be $(\eps_{\ref{NewerInterlacingPetalStar}}+\gamma_N)$-crossed by \emph{any} $k$-separation.
    So by \cref{FlowerClosePetalOrCrossed}, $d$ is $2(\eps_{\ref{NewerInterlacingPetalStar}}+\gamma_N+k)$-close to some $p\in\pi$ or its inverse.
    Hence $s$ is $(3\eps_{\ref{NewerInterlacingPetalStar}}+2\gamma_N+2k)$-close to some~$p \in \pi$ or its inverse by \cref{approxLeQuasiTransitive}.

    Suppose for a contradiction that $s$ is $(3\eps_{\ref{NewerInterlacingPetalStar}}+2\gamma_N+2k)$-close to some $p^* \eqqcolon (Q, P) \in \pi^*$.
    Since $A$ is $\delta'$-relevant, there exist two $\delta'$-crossing separations $r_1, r_2 \in A^{\delta} \subseteq D(F_A)$.
    Since $A^{\delta}$ $\eps_N$-interlaces $\sigma$ and $s \in \sigma$, we can orient each $r_i$ so that $s \leq_{\eps_N} r_i$.
    Write $r_1 \eqqcolon (A_1, B_1)$ and $r_2 \eqqcolon (A_2, B_2)$.
    Since $P$ meets at most one of their four corners, $Q$ contains two corners that are disjoint from $A_1 \cap A_2$.
    Each has size $> \delta'-k$ by \cref{crissCross}~\ref{crissCross:crossImpliesCriss}.
    However, $|U \setminus (A_1 \cap A_2)| \leq |U \setminus A_1| + |U \setminus A_2| \leq 2\eps_N$.
    Thus, using that $\delta' \geq \delta + 2\eps_{\ref{NewerInterlacingPetalStar}} \geq \gamma_N+2\eps_N+2\eps_{\ref{NewerInterlacingPetalStar}}\geq \gamma_N + \eps_N + 2k + 2\eps_{\ref{NewerInterlacingPetalStar}}$ by \ref{set:FlowerFinding:Inequalities}, it follows that
    \[
        3\eps_{\ref{NewerInterlacingPetalStar}}+2\gamma_N+2k \geq |Q \setminus U| \geq |Q \setminus (A_1 \cap A_2)| - |U \setminus (A_1 \cap A_2)| > 2(\delta'-k-\eps_N) \geq 3\eps_{\ref{NewerInterlacingPetalStar}}+2\gamma_N+2k~,
        \] a contradiction.
    Thus, $s$ is $(3\eps_{\ref{NewerInterlacingPetalStar}}+2\gamma_N+2k)$-close, hence $\eps_{\ref{NewerSplittingStarCloseToPetalStar}}$-close, to some~$p \in \pi$.

    \casen{Case~1b: $X:=\sep (s)$ is $S$-totally-nested, where $S$ is the set of all narrow proper $k$-separations of~$G$, and $s\in D(F)$ for the graphic-$X$-anemone~$F$. Moreover, $s$ is not half-connected.}
    Let $\pi_X \coloneqq \pi(F) \setminus \{(X,V(G))\}$ be the star of proper petal-separations of~$F$.
    Using that $s$ is $\eps_{\ref{NewerSplittingStarCloseToPetalStar}}$-broad with $\eps_{\ref{NewerSplittingStarCloseToPetalStar}} > \gamma_N \geq \beta \geq q$ by \ref{set:FlowerFinding:Inequalities} and that $G$ is $(q,k-1)$-unbreakable, we have that $|X| = k$ and $F$ has a $b$-broad displayed separation, namely $s$. 
    By \cref{lem:petal-separation-of-graphic-anemone-narrow}, all the separations in $\pi_X$ are narrow.
    Thus, $\pi_X \cup \pi_X^*$ is the set of all narrow proper half-connected $k$-separations with separator~$X$, and $\pi_X \neq \emptyset$ as $s$ is a narrow proper $k$-separation with separator $X$.
    Since $X$ is $S$-totally-nested, so are the separations in~$\pi_X$ by \cref{lem:equivalence-totally-nested-set}.
    Then the separations in~$\pi_X$ are contained in~$N$ by \ref{set:FlowerFinding:NcontainsTotallyNested}. 
    
    We now show that $\pi_X\preceq\sigma$.
    Let $p \in \pi_X$.
    As $s\in D(F)$, we have that $s$ almost-interlaces~$\pi_X$, so $p \leq s$ or $p \leq s^*$.
    If $p \leq s$ we are done, so assume $p \leq s^*$.
    Since $\sigma$ is a splitting star of~$N$ and $p \in \pi_X \subseteq N$, there exists $t \in \sigma$ such that $t \not<p$ and $t \not<p^*$.
    Since $N$ is nested, we have $p \leq t$ or $p^* \leq t$.
    If $p \leq t$ we are done, so assume $p^* \leq t$.
    If $t=s$ then $s \leq p^* \leq t=s$ so $s = p^*$. 
    This is impossible since $p$ is half-connected and $s$ is not.
    Thus $t \neq s$, so,
    using that $\sigma$ is a star, we have $s \leq p^* \leq t \leq s^*$, which would contradict that $s$ is proper.
    Therefore, we indeed have $\pi_X\preceq\sigma$.
    
    By assumption, $A^\delta$ exists.
    Let $r_1,r_2\in A^\delta$ $\delta$-cross.
    Thus, $r_1, r_2$ have breadth $\geq 2(\delta - k)$ by \cref{LargeCrossedImpliesBroad}, with $\delta>\delta_{\ref{HighSpreadUniqueSplittingStar}} \geq \delta_{\ref{NewerLargelyCrossingEntersAsterisk}} \geq \gamma_{\cA} + k$, so they have breadth $\geq \gamma_{\cA} \geq \beta$ by~\ref{set:FlowerFinding:Inequalities}. 
    Then both $r_i$ are narrow and have order exactly~$k$, by \cref{FlowerFindingBroadSeparationsNice}.
    Hence $r_1, r_2 \in S$.
    Both $r_i$ $\eps_N$-interlace~$\sigma$ by assumption.
    Since $\pi_X\preceq\sigma$, both $r_i$ $\eps_N$-interlace~$\pi_X$ as well.
    We claim that both $r_i$ interlace~$\pi_X$.
    Since the $r_i$ $\eps_N$-interlace $\pi_X$,
    neither $r_i$ lies in~$\pi_X$, so by symmetry it suffices to show that $r_1$ almost-interlaces~$\pi_X$.
    For this, let $\ell\in\pi_X$.
    If $\ell\le r_1$ or $\ell\le r^*_1$, then we are done, so suppose that this fails.
    As the elements of $\pi_X$ are $S$-totally-nested and both $r_i$ lie in~$S$, we then have $r_1\le\ell$ or $r^*_1\le\ell$.
    Since $r_1$ $\eps_N$-interlaces~$\pi_X$, we also have $\ell<_{\eps_N}r_1$ or $\ell<_{\eps_N}r^*_1$.
    The two cases $r_1\le\ell<_{\eps_N}r_1$ and $r_1^*\le\ell<_{\eps_N}r_1^*$ are clearly impossible.
    If $r_1\le\ell<_{\eps_N} r_1^*$, then $r_1$ has breadth $\leq \eps_N$ by \cref{lem:assymmetric-nested-implies-not-broad}, contradicting that $r_1$ is $2(\delta-k)$-broad with $2(\delta-k) > 2(\delta_{\ref{HighSpreadUniqueSplittingStar}}-k) \geq 2(\delta_{\ref{NewerLargeAsteriskUniqueSplittingStar}}-k) \geq 2(2\eps_N-k) \geq 2\eps_N \geq \eps_N$ (using $\eps_N \geq k$ by \ref{set:FlowerFinding:Inequalities}).
    Similarly, the case $r_1^*\le\ell<_{\eps_N} r_1$ yields a contradiction.
    Therefore, $r_1$ almost-interlaces~$\pi_X$.
    This completes the proof of the claim that both $r_i$ interlace~$\pi_X$.
    
    Then $r_1$ and $r_2$ must cross with centre equal to~$X$ and $r_i\in D(F)$.    
    Since $A^\delta\se A$ by \cref{NewerLargelyCrossingEntersAsterisk}, we have $r_1,r_2\in A$, and so $X$ is the common centre of~$A$.
    By~\ref{set:FlowerFinding:XAnemone}, we have $F=F_A$.
    So $s\in D(F)=D(F_A)$.
    Then the premise of Case~1a is satisfied, and so we conclude the argument as in Case~1a.

    \casen{Case~2: there is $p\in\pi$ such that $s\le_{\gamma_N}p$ or $s^*\le_{\gamma_N}p$.}
    We first claim that $s^*\le_{\gamma_N}p$ is impossible, and assume for a contradiction that $s^*\le_{\gamma_N}p$.
    Then $s^*\le_{\eps_N}p$ since $\eps_N\ge\gamma_N$ by \ref{set:FlowerFinding:Inequalities}.
    Let $r\in A^\delta$.
    Then $r\in A^\delta\se A$ by \cref{NewerLargelyCrossingEntersAsterisk}.
    As $r$ $\eps_N$-interlaces~$\sigma$, we have $s<_{\eps_N} r$ or $s<_{\eps_N} r^*$, so $s\in O(r)$ with $O(r)$ from \cref{NewerLargeAsteriskUniqueSplittingStar}~\cref{NewerLargeAsteriskUniqueSplittingStar:eq:Orientation}.
    As $r\in A^\delta\se D(F_A)$ by \cref{NewerLargelyCrossingEntersAsterisk}, we also have $p\le r$ or $p\le r^*$.
    Combining this with $s^*\le_{\eps_N}p$ yields $s^*\le_{\eps_N} r$ or $s^*\le_{\eps_N} r^*$.
    Hence $s^*\in O(r)$.
    But $s,s^*\in O(r)$ contradicts that $O(r)$ is an orientation of~$N$ by \cref{NewerLargeAsteriskUniqueSplittingStar}.

    Therefore, $s\le_{\gamma_N}p$.
    We claim that $p$ is $\eps_{\ref{NewerPetalSeparationsCloseToSplittingStar}}$-broad.
    Let $(U,W):=p$ and $(C,D):=s$.
    Since $s\le_{\gamma_N}p$ and since $s$ is $\eps_{\ref{NewerSplittingStarCloseToPetalStar}}$-broad,
    we have 
    \[
        |U\sm W| \geq |(U \sm W) \cap (C \sm D)| = |C \sm D| - |W \sm D| \ge \eps_{\ref{NewerSplittingStarCloseToPetalStar}}-\gamma_N\ge \eps_{\ref{NewerPetalSeparationsCloseToSplittingStar}}~.
    \]
    As $A$ is $\delta'$-relevant by \cref{HighSpreadUniqueSplittingStar}, there are $r_1,r_2\in A^\delta$ that $\delta'$-cross.
    Hence both $r_i$ have breadth~$\ge 2(\delta'-k)$ by \cref{LargeCrossedImpliesBroad}.
    Since $p\in\pi$, we have $p\le r_1\wedge r_2$, say, so $|W\sm U|\ge 3(\delta'-k)$ by \cref{fishBroad}.
    Using $\delta'>\delta'_{\ref{NewerUniqueAsterisk}}=2\delta+2k$, we find that
    \begin{align*}
        |W\sm U| & {}\ge 3(\delta'-k)\\
        & {} \geq 2\delta' + (\delta'_{\ref{NewerUniqueAsterisk}}-3k)\\
        &{} \ge 2\delta' + 2\delta -k\\
        & {} \ge 2\delta' + 2\delta_{\ref{HighSpreadUniqueSplittingStar}} - k\\
        & {} \ge 2\delta' + \delta_{\ref{NewerLargelyCrossingEntersAsterisk}} + \delta_{\ref{NewerLargeAsteriskUniqueSplittingStar}} - k \\
        & {} \ge 2\delta' + k + 2\eps_N+\gamma_N\\
        & {} \overset{\text{\ref{set:FlowerFinding:Inequalities}}}{\ge} 2\delta' + 2k+2\eps_N=\eps_{\ref{NewerPetalSeparationsCloseToSplittingStar}}.
    \end{align*}
    Therefore, $p$ is $\eps_{\ref{NewerPetalSeparationsCloseToSplittingStar}}$-broad, as desired.
    
    By \cref{NewerPetalSeparationsCloseToSplittingStar}, there is $t\in\sigma$ that is $\eps_{\ref{NewerPetalSeparationsCloseToSplittingStar}}$-close to~$p$.
    If $s=t$ then we are done, so assume $s \neq t$.
    Now $s\le_{\gamma_N} p \le_{\eps_{\ref{NewerPetalSeparationsCloseToSplittingStar}}} t$ yields $s\le_{\gamma_N+\eps_{\ref{NewerPetalSeparationsCloseToSplittingStar}}} t$ by \cref{approxLeQuasiTransitive}.
    Since $s,t\in\sigma$ and $s \neq t$, by the definition of star we also have $s\le t^*$.
    As $s\le_{\gamma_N+\eps_{\ref{NewerPetalSeparationsCloseToSplittingStar}}} t$ is equivalent to $t^* \le_{\gamma_N+\eps_{\ref{NewerPetalSeparationsCloseToSplittingStar}}} s^*$, we obtain $s\le_{\gamma_N+\eps_{\ref{NewerPetalSeparationsCloseToSplittingStar}}} s^*$.
    Hence $s$ has breadth~$\le\gamma_N + \eps_{\ref{NewerPetalSeparationsCloseToSplittingStar}}$ by \cref{lem:assymmetric-nested-implies-not-broad}.
    But we assumed at the beginning of the proof that $s$ has breadth~$\ge\eps_{\ref{NewerSplittingStarCloseToPetalStar}}>\gamma_N+\eps_{\ref{NewerPetalSeparationsCloseToSplittingStar}}$, a contradiction.
\end{proof}

\subsection{The torso structure theorem}

Now, we are able to prove the main results for this section (\cref{keylemma:FlowerFinding} and \cref{keylemma:FlowerCriterion}).
This can be done by combining the key lemmas from the previous sections.

\begin{lem}\label{interlacingSeparationNotCloseSplittingStar}
    Let $N$ be a nested symmetric set of separations of a graph~$G$.
    Let $\sigma$ be a splitting star of~$N$.
    Let $s$ be a $k$-separation of~$G$ that $\eps$-interlaces~$\sigma$,
    but such that $s$ is not $\eps'$-close to any separation in $\sigma\cup\sigma^*$.
    Assume that $\eps'\ge\eps$ and that $s$ has breadth at least~$\beta_{\ref{interlacingSeparationNotCloseSplittingStar}}(\eps,\eps'):=\eps'+\eps+1$.
    Then $s$ is not $\eps'$-close to any separation in~$N$ either.
\end{lem}
\begin{proof}
    Assume for a contradiction that $s$ is $\eps'$-close to some $t\in N$.
    Since $\sigma$ is a splitting star of~$N$, there is $r\in\sigma$ such that $t\le r$ or $t^*\le r$.
    Since $s$ $\eps$-interlaces~$\sigma$, we have $r<_{\eps} s$ or $r<_{\eps} s^*$.
    We now consider all four cases, and derive a contradiction in each of them.

    \casen{Case $t\le r$ and $r<_{\eps} s$.}
    Here we have $s\le_{\eps'} t\le r <_{\eps}s$, so $s$ and $r$ are $\eps'$-close (using $\eps'\ge\eps$), a contradiction.

    \casen{Case $t\le r$ and $r<_{\eps}s^*$.}
    Here we have $s\le_{\eps'} t\le r <_{\eps}s^*$, so $s\le_{\eps'+\eps}s^*$ by \cref{approxLeQuasiTransitive}.
    By \cref{lem:assymmetric-nested-implies-not-broad}, $s$ has breadth~$\leq \eps'+\eps$.
    This contradicts that $s$ is a $\beta_{\ref{interlacingSeparationNotCloseSplittingStar}}$-broad separation.

    \casen{Case $t^*\le r$ and $r<_\eps s$.}
    From the assumption $t\le_{\eps'} s$ we obtain $s^*\le_{\eps'} t^*$.
    Then \[s^*\le_{\eps'} t^*\le r <_\eps s\]
    and we derive a contradiction in analogy to the previous case.

    \casen{Case $t^*\le r$ and $r<_\eps s^*$.}
    We derive a contradiction as in the first case.
\end{proof}

\begin{theorem}\label{keylemma:FlowerFinding}
    Assume \cref{set:FlowerFinding}.
    Let $\delta>\delta_{\ref{NewerDistinctSplittingStars:Claim}}$, $\delta' >\delta'_{\ref{NewerDistinctSplittingStars}}(\ldots,\delta)$, $\zeta \geq \zeta_{\ref{NewerInterlacesCloseToDisplayed}}(k,\delta',\eps_F)$, $\beta_{\ref{keylemma:FlowerFinding}}:=\beta_{\ref{interlacingSeparationNotCloseSplittingStar}}(\eps_N,\eps_{\ref{NewerInterlacesCloseToDisplayed}})$, and $q'_{\ref{keylemma:FlowerFinding}} :=\max\,\{\,q,\;\eps_{\ref{NewerInterlacesCloseToDisplayed}}+k,\;\beta_{\ref{keylemma:FlowerFinding}}\,\}$.
    Let $\sigma$ be a splitting star of~$N$.
    Then we have one of the following outcomes:
    \begin{defenum}
        \item\label{keylemma:FlowerFinding:FlowerStructure}
        There is a unique $A\in\cA$ of spread~$\geq\zeta$ such that $A^\delta$ exists and $\eps_N$-interlaces~$\sigma$.\\
        The graphic-$k$-flower $F_A$ has spread~$\ge\zeta$.\\
        The splitting star $\sigma$ and the petal-star of~$F_A$ $\eps_{\ref{NewerSplittingStarCloseToPetalStar}}$-approximate each other.
        \item\label{keylemma:FlowerFinding:HighlyConnected} The torso of $\sigma$ is $(q'_{\ref{keylemma:FlowerFinding}},k)$-unbreakable.
    \end{defenum}
\end{theorem}

\begin{theorem}\label{keylemma:FlowerCriterion}
    Assume the premises of \cref{keylemma:FlowerFinding}.
    Let $S(\sigma)$ be the set of all $\beta_{\ref{keylemma:FlowerFinding}}$-broad $k$-separations of~$G$ that $\eps_N$-interlace~$\sigma$.
    \begin{enumerate}
        \item\label{keylemma:FlowerCriterion:FlowerType} If some $s\in S(\sigma)$ is not $\eps_{\ref{NewerInterlacesCloseToDisplayed}}$-close to any separation in~$\sigma\cup\sigma^*$, then
        \begin{itemize}
            \item $\sigma$ has type~\ref{keylemma:FlowerFinding:FlowerStructure};
            \item $s\in A$ and $s$ is $\eps_F$-close to some displayed separation of~$F_A$; and
            \item $s$ is $\zeta$-crossed by some $s'\in D(F_A)\cap A$.
        \end{itemize}
        \item\label{keylemma:FlowerCriterion:HighConType} 
        Otherwise, $\sigma$ has type~\ref{keylemma:FlowerFinding:HighlyConnected}.
    \end{enumerate}
\end{theorem}

\begin{proof}[Proof of \cref{keylemma:FlowerFinding,keylemma:FlowerCriterion}]
    \cref{keylemma:FlowerCriterion:FlowerType}. Assume first that there is $s\in S(\sigma)$ that is not $\eps_{\ref{NewerInterlacesCloseToDisplayed}}$-close to any separation in~$\sigma\cup\sigma^*$.
    Observe that $\eps_{\ref{NewerInterlacesCloseToDisplayed}} \geq \eps_N$.
    Since $\sigma$ is a splitting star of~$N$, \cref{interlacingSeparationNotCloseSplittingStar} yields that $s$ is not $\eps_{\ref{NewerInterlacesCloseToDisplayed}}$-close to any separation in~$N$.
    Then we are done by \cref{NewerInterlacesCloseToDisplayed}, \cref{NewerPetalSeparationsCloseToSplittingStar} and \cref{NewerSplittingStarCloseToPetalStar}.

    \cref{keylemma:FlowerCriterion:HighConType}.
    Assume otherwise that every separation in $S(\sigma)$ is $\eps_{\ref{NewerInterlacesCloseToDisplayed}}$-close to some separation in~$\sigma \cup \sigma^*$.
    Note that $q'_{\ref{keylemma:FlowerFinding}}\ge \eps_{\ref{NewerInterlacesCloseToDisplayed}}+k$.
    To show that the torso of $\sigma$ is $(q'_{\ref{keylemma:FlowerFinding}},k)$-unbreakable, by \cref{highlyConnectedTorso} it suffices to show that there is no $(\le k)$-separation $s$ of~$G$ that almost-interlaces~$\sigma$, has breadth~$>q'_{\ref{keylemma:FlowerFinding}}$, and is not $\eps_{\ref{NewerInterlacesCloseToDisplayed}}$-close to any separation in~$\sigma\cup\sigma^*$.
    Assume for a contradiction that there is such a $(\le k)$-separation~$s$.

    Since $s$ has breadth~$>q'_{\ref{keylemma:FlowerFinding}}\ge q$ and $G$ is $(q,k-1)$-unbreakable, it follows that $s$ has order~$k$.
    The separation $s$ has breadth~$\ge\beta_{\ref{keylemma:FlowerFinding}}$ as $q'_{\ref{keylemma:FlowerFinding}}\ge\beta_{\ref{keylemma:FlowerFinding}}$.

    We claim that $s$ $\eps_N$-interlaces~$\sigma$.
    For this, let $r\in\sigma$.
    As $s$ almost-interlaces~$\sigma$ by assumption, we have $r\le s$, say.
    Hence $r\le_{\eps_N} s$.
    We have to show that $s\not\le_{\eps_N}r$.
    If we had $s\le_{\eps_N}r$, then $s$ would be $\eps_N$-close to~$r\in\sigma$.
    As $\eps_N\le\eps_{\ref{NewerInterlacesCloseToDisplayed}}$, this would contradict our assumption that $s$ is not $\eps_{\ref{NewerInterlacesCloseToDisplayed}}$-close to any separation in~$\sigma\cup\sigma^*$.
    Hence $s$ $\eps_N$-interlaces~$\sigma$.

    We have $s\in S(\sigma)$ by the above claims.
    We chose $s$ so that additionally $s$ is not $\eps_{\ref{NewerInterlacesCloseToDisplayed}}$-close to any separation in~$\sigma\cup\sigma^*$.
    This now contradicts the assumption of~\cref{keylemma:FlowerCriterion:HighConType}.
    Hence the torso of~$\sigma$ is $(q'_{\ref{keylemma:FlowerFinding}},k)$-unbreakable after all.
\end{proof}

\section{Decompositions of narrow graphs}
\label{sec:main-theorem-for-narrow}

In this section, we prove \cref{thm:main-intro} for narrow graphs, see \cref{thm:cactusForNarrow}.
To do so, in \cref{subsec:Tutte-for-narrow} we first establish the corresponding Tutte-type statement (\cref{mainTutte}) for narrow graphs, see \cref{mainForNarrow}. In \cref{subsec:cactus-for-narrow} we then refine the resulting tree-decomposition using the corresponding flowers to obtain the desired cactus-decomposition.

\subsection{An approximate Tutte-type decomposition for narrow graphs}\label{subsec:Tutte-for-narrow}

Before proving \cref{mainForNarrow}, we just need one last result.

\begin{lem}\label{totallyNestedCriterion}
    Let $q,k\in\N$ and let $G$ be a $(q,k-1)$-unbreakable graph.
    Let $s$ be a narrow proper $k$-separation of~$G$ that is $\gamma$-nested with every narrow proper $k$-separation of~$G$.
    Assume that $\gamma\ge q+k$.
    Then $s$ is narrow-$\gamma$-nested.
\end{lem}
\begin{proof}
    Let $r$ be an arbitrary narrow proper $(< k)$-separation of~$G$.
    We have to show that $s$ is $\gamma$-nested with~$r$.
    Since $G$ is $(q,k-1)$-unbreakable, $r$ has breadth~$\le q$.
    So some corner of $r$ and $s$ has size~$\le q$.
    In particular, $r$ and $s$ do not $q$-criss.
    Hence $r$ and $s$ are $(q+k)$-nested by the contrapositive of \cref{crissCross}~\ref{crissCross:crossImpliesCriss}.
\end{proof}

\begin{definition}[$\psi$-isomorphic for sets of graph-decompositions]
    Let $G$ be a graph, let $\mathcal{G}_1, \mathcal{G}_2$ be two sets of graph-decompositions of $G$ and let $\psi$ be an automorphism of $G$.
    We say that $\mathcal{G}_1$ and $\mathcal{G}_2$ are \defn{$\psi$-isomorphic} if there exists a bijection $\Psi\colon \cG_1\to\cG_2$ such that $\mathcal{H}$ and $\Psi(\mathcal{H})$ are $\psi$-isomorphic (as in \cref{defn:psi-isomorphic-graph-decomposition}) for every $\mathcal{H}\in\cG_1$.
\end{definition}

\begin{definition}[Canonical for sets of graph-decompositions]\label{def:cano-set-graph-decomp}
    Let $G$ be a graph, and let $\cG$ be a set of graph-decompositions of~$G$.
    We say that $\cG$ is \defn{canonical} if, for every automorphism $\psi$ of~$G$, $\cG$ is $\psi$-isomorphic to itself, i.e.\ there is a permutation $\Psi$ of~$\cG$ such that $\cH$ and $\Psi(\cH)$ are $\psi$-isomorphic (as in \cref{defn:psi-isomorphic-graph-decomposition}) for every $\mathcal{H}\in\cG$.
\end{definition}

We can now prove \cref{mainTutte} under the additional assumption that the graph to be decomposed is $(b, k)$-narrow.
For a tree-decomposition $\cT = (T,\cV)$ and a node $t$ of~$T$, the \defn{$t$-star} is the set of all displayed separations of $\cT$ that are induced by cuts $(X, Y)$ of $T$ such that there is a unique $X$--$Y$ edge in $T$, which is additionally incident to $t \in Y$.
We say that a graphic-$k$-flower $F$ \defn{$\eps$-refines} $(T,\cV)$ at~$t$ if the petal-star of~$F$ and the $t$-star of~$(T,\cV)$ $\eps$-approximate each other.

\begin{theorem}\label{mainForNarrow}
    For every $q,k,b\in\N$ there are $\gamma,\kappa,\rho,\bar\eps\in\N$ such that for all $\zeta\in\N$ there are $q',\eps \in \N$, such that the following holds.
    Let $G$ be a $(q,k-1)$-unbreakable $(b,k)$-narrow graph.
    There exists a set $\cF$ of graphic-$k$-flowers of~$G$ of spread $\geq \zeta$ and a symmetric nested set~$N$ of narrow-$\gamma$-nested proper $k$-separations of~$G$ with corresponding tree-decomposition $(T,\cV)$ of $G$ and a bipartition $X\sqcup Y= V(T)$ such that:
    \begin{enumerate}
        \item\label{mainForNarrowTorsoQKC} Every node $t\in X$ has a $(q',k)$-unbreakable torso.
        \item\label{mainForNarrowTorsoFlower} Every node $t\in Y$ comes with a graphic-$k$-flower $F_t\in\cF$ of $G$ such that $F_t$ $\rho$-refines $(T,\cV)$ at~$t$. 
        Furthermore, if $F_t$ is a graphic-$k$-daisy, then $F_t$ has the clover-property.
        \item\label{mainForNarrowDisplayedInterlace}
        Every displayed separation of some~$F_t$ (for $t \in Y$) that is $\rho$-crossed by another displayed separation of~$F_t$ almost-$\rho$-interlaces the splitting star of~$(T,\cV)$ corresponding to~$t$.
        \item\label{mainForNarrowSeparationAnalysis} Every $\eps$-broad $k$-separation of $G$ either is $\eps$-close to a separation in~$N$ or is $\eps$-close to some displayed separation of some~$F_t$.
        \item\label{mainForNarrowCanonical} 
        $\cF$ is canonical.
        Both the tree-decomposition $(T,\cV)$ and the set $N$ are $\kappa$-canonical.
        Given $N$ and $(T, \cV)$, the bipartition $X\sqcup Y=V(T)$ is canonical and the assignment $Y\ni t\mapsto F_t\in\cF$ is a canonical bijection.
        \item\label{mainForNarrowFlowersDiffer}
        $D(F)\neq D(F')$ for every two distinct graphic-$k$-flowers $F,F'\in\cF$.
        \item\label{mainForNarrowCloseToNestedSeparation} For every $\bar\eps$-broad $k$-separation of~$G$ there is an $\bar\eps$-close $k$-separation of~$G$ nested with~$N$.
    \end{enumerate}
\end{theorem}

\begin{proof}
    Let $q,k,b,\zeta\in\N$ be given, and let $G$ be a $(q,k-1)$-unbreakable $(b,k)$-narrow graph.
    We first establish \cref{set:FlowerFinding}.
    To this end, let
    \begin{align*}
        \beta & {}:=\max\,\{\,\beta_{\ref{keylem:flower-approximation}}(q,k),\enspace b,\enspace q+1\,\}\ge 2q+k\\
        \eps_F & {}:=\eps_{\ref{keylem:flower-approximation}}(q,k) \\
        \gamma_\cA & {}:=\max\{\,\gamma_{\ref{HugeSeparationsUniformlyCrossing}}(q,k),\enspace\gamma_{\ref{keylem:flower-approximation}}(q,k)+1,\enspace \gamma_{\ref{keylem:petal-separations-totally-nested}}(q,k,\eps_F),\enspace\beta\,\}.
    \end{align*}

    The choice of $\beta$ and $\gamma_\cA$ establishes the first two points of~\ref{set:FlowerFinding:Inequalities}.
    Let $S_0$ be the set of all $\beta$-broad $k$-separations of~$G$.
    Since $\beta \geq b$ and $G$ is $(b, k)$-narrow, all separations in $S_0$ are narrow.
    By \cref{HugeSeparationsUniformlyCrossing} and using $\beta\ge q$ as well as $\gamma_\cA\ge \gamma_{\ref{HugeSeparationsUniformlyCrossing}}$, we find that the set $S_0$ is uniformly $\gamma_\cA$-crossing.

    \casen{Flowers.}
    Let $\cA\triv$ and $\cA$ consist of all trivial and nontrivial $\gamma_\cA$-asterisks of~$S_0$, respectively.
    For each $A\in\cA$, let $F_A:=F_{\ref{keylem:flower-approximation}}(A)$ be the $q$-lean graphic-$k$-flower found by \cref{keylem:flower-approximation}.
    By \cref{keylem:flower-approximation}~\cref{flower-canonical}, these graphic-$k$-flowers are canonical.
    When the common centre $X$ of~$A$ has size~$k$, we are in a special case that we considered in \cref{subsec:anemone}.
    In that special case, we defined $F_{\ref{keylem:flower-approximation}}(A)=F_A$ to be the graphic-$X$-anemone, which establishes~\ref{set:FlowerFinding:XAnemone}.
    In general, every separation in~$A$ is $\eps_F$-close to some displayed separation of~$F_A$, by \cref{keylem:flower-approximation}~\ref{item:keylem:flower-approximation-1} and \cref{closeVSapprox}, which establishes~\ref{set:FlowerFinding:CloseFlower}.
    We will later find $\cF$ as a subset of $\{\,F_A:A\in\cA\,\}$.

    \casen{Petal-separations.} Next, we let
    \begin{align*}
        \gamma_\pi & {}:= \max\,\{\,\gamma'_{\ref{keylem:petal-separations-totally-nested}}(q,k),\enspace \gamma_\cA,\enspace \eps_{\ref{keylem:flower-approximation}}(q, k),\enspace \beta+k\,\}\\
        \gamma_N & {}:=\gamma'_{\ref{keylem:max-nested-set}}(k,b,\gamma_\pi)\ge\gamma_\pi \geq \beta+k\ge q+k.
    \end{align*}
    The choice of $\gamma_N$ ensures that it meets the requirements of~\ref{set:FlowerFinding:Inequalities}.
    Let $A \in \cA$ be arbitrary.
    Then, $A$ is a nontrivial $\gamma_{\cA}$-asteriskal set of $k$-separations of~$G$, and~$F_A$ is a $q$-lean graphic-$k$-flower of~$G$ that $\eps_F$-approximates~$A$ by \cref{keylem:flower-approximation}~\ref{item:keylem:flower-approximation-1}.
    By \cref{keylem:petal-separations-totally-nested}, every $k$-separation of~$G$ is $\gamma_{\pi}$-compatible with~$F_A$.
    By \cref{keylem:flower-approximation}~\ref{item:keylem:flower-approximation-2}, since $\gamma_{\pi} \geq \eps_{\ref{keylem:flower-approximation}}(q, k)$, we obtain that 
    \begin{equation}\label{eq:narrowMainPetalsNested}
        \text{every $k$-separation of $G$ is $\gamma_{\pi}$-nested with all the petal-separations of $F_A$,}
    \end{equation}
     which establishes~\ref{set:FlowerFinding:S1totallyNested} since $\gamma_N \geq \gamma_{\pi}$.

    \casen{Nested sets of separations.}
    Let~$M$ be the set of all narrow-$\gamma_\pi$-nested $k$-separations of~$G$.
    Let~$M'$ be the set of all half-connected narrow-$\gamma_N$-nested $\beta_M$-broad $k$-separations of~$G$ that are $\eps_M$-close to some separation in~$M$ for
    \begin{align*}
        \beta_M & {}:=\beta_{\ref{keylem:max-nested-set}}(q,k,b)\\
        \eps_M &{}:= \delta_{\ref{keylem:max-nested-set}}(k,b,\gamma_\pi).
    \end{align*}
    Let~$N' \subseteq M'$ be a maximal nested subset.
    Let $S$ be the set of all narrow proper $k$-separations of~$G$.
    Note that $N'\se S$, as every separation in~$N'$ is $\beta_M$-broad where $\beta_M\ge b$ and $G$ is $(b,k)$-narrow.
    We obtain $N$ from~$N'$ by adding all the $S$-totally-nested $k$-separations of~$G$.
    Then $N$ is a nested symmetric set of proper $k$-separations of $G$ and satisfies \ref{set:FlowerFinding:NcontainsTotallyNested} by definition.
    
    Every separation in~$N$ is either contained in $N'\se M'$ and so is narrow-$\gamma_N$-nested by definition of~$M'$, or is $S$-totally-nested and so is narrow-$\gamma_N$-nested by \cref{totallyNestedCriterion} using $\gamma_N\ge q+k$.
    Hence all separations in~$N$ are narrow-$\gamma_N$-nested, 
    so $N$ satisfies~\ref{set:FlowerFinding:NtotallyNested}.

    By \cref{keylem:max-nested-set},~$N'$ $\eps_N$-approximates $M$ for
    \[
        \eps_N :=\max\{\,\eps_{\ref{keylem:max-nested-set}}(q,k,b,\gamma_\pi),\enspace\gamma_N\,\} \geq q+1,b,\gamma_{\pi}.
    \]
    Note that $\eps_N\ge\gamma_N$, so \ref{set:FlowerFinding:Inequalities} is fully established now.
    
    Let $S_1:=\bigcup\,\{\,\pi(F_A)\cup (\pi(F_A))^*:A\in\cA\,\}\cup\bigcup\cA\triv$.

    \begin{claim}\label{mainForNarrow:claim:SoneTotallyNested}
        $N$ $\eps_N$-approximates~$S_1$.
    \end{claim}
    
    \begin{cproof}
        For this, let $s\in S_1$ be $\eps_N$-broad.
        Since $\eps_N>q$ and $G$ is $(q,k-1)$-unbreakable, $s$~has order~$k$.
        Since $G$ is $(b,k)$-narrow and $\eps_N\ge b$, we have that $s$ is narrow.
        
        Assume first that $s\in\pi(F_A)\cup (\pi(F_A))^*$ for some~$A\in\cA$.
        By \cref{eq:narrowMainPetalsNested}, we have that $s$ is $\gamma_\pi$-nested with every (narrow proper) $k$-separation of~$G$.
        Hence $s \in M$ by \cref{totallyNestedCriterion} and $\gamma_\pi\ge q+k$.
        Then some separation in $N'\se N$ is $\eps_N$-close to~$s$.
        
        We assume next that $s\in A$ for some trivial~$A\in\cA\triv$.
        We show that $s \in M$.
        If not, $s$ is $\gamma_\pi$-crossed by some narrow proper $(\leq k)$-separation $r$ of~$G$.
        As $\gamma_\pi\ge\gamma_\cA$ and $A$ is trivial, we cannot have~$r\in S_0$.
        To obtain a contradiction, we will show that~$r\in S_0$.
        Since $s$ and $r$ $\gamma_\pi$-cross, \cref{LargeCrossedImpliesBroad} yields that $r$ has breadth $\ge 2(\gamma_\pi-k)\ge\beta$.
        Since $G$ is $(q,k-1)$-unbreakable with $\beta > q$, $r$ has order $k$, so $r \in S_0$, a contradiction.
        Thus, $s \in M$, and as above some separation in $N'\se N$ is $\eps_N$-close to~$s$.
    \end{cproof}    
    
    \cref{mainForNarrow:claim:SoneTotallyNested} establishes~\ref{set:FlowerFinding:Napproximates}.
    We now have completely established~\cref{set:FlowerFinding}.

    \casen{Tree-decomposition.}
    In order to apply \cref{keylemma:FlowerFinding}, we let
    \begin{align*}
        \delta & {}:=\delta_{\ref{NewerDistinctSplittingStars:Claim}}(q,k,\gamma_\cA,\eps_F,\gamma_N,\eps_N)+1~,\\
        \delta' & {} := \delta'_{\ref{NewerDistinctSplittingStars}}(q,k,b,\gamma_N,\delta)+1~,\\
        \rho & {} := \eps_{\ref{NewerSplittingStarCloseToPetalStar}}(q,k,b,\gamma_N,\eps_N,\delta')~,\\
        \zeta' & {} := \max\{ \zeta\, ,\, \zeta_{\ref{NewerInterlacesCloseToDisplayed}}(k,\delta',\eps_F)\}~,\\
        q' & {} := q'_{\ref{keylemma:FlowerFinding}}(q, k,\beta,\eps_F,\eps_N,\delta',\zeta')~,\\
        \eps' & {} := \eps_{\ref{NewerInterlacesCloseToDisplayed}}(k,\eps_F,\eps_N,\zeta')\ge\eps_F,\eps_N~,\\
        \eps & {} := \max\,\{\,\beta_{\ref{keylemma:FlowerFinding}}(\eps_N,\eps')\, ,\enspace 2\gamma_N+1\,\}\ge\eps_F,\eps_N, \eps'~.
    \end{align*}
    Let $(T,\cV)$ be the tree-decomposition of~$G$ that is uniquely determined by~$N$, see \cref{TDCfromNesteSetVerify}.
    Setting $\gamma \coloneqq \gamma_N$, we have that the separations in $N$ are indeed narrow-$\gamma$-nested.
    Let $\sigma\mapsto t(\sigma)$ be the natural bijection between the set of all splitting stars of~$N$ and~$V(T)$.

    Let $Y_N$ be the set of all splitting stars of~$N$ that have type~\ref{keylemma:FlowerFinding:FlowerStructure} for the parameter $\zeta'$, and let $X_N$ be the set of all splitting stars of~$N$ that are not in~$Y_N$.
    Let $X,Y$ be the images of $X_N,Y_N$ under the map $\sigma\mapsto t(\sigma)$.
    Given~$N$ and~$(T,\cV)$, the bipartition $X\sqcup Y=V(T)$ is canonical.

    \casen{Definition of $\cF$.}
    Let $\cF$ consist of all graphic-$k$-flowers $F_A$ such that both~$A$ and~$F_A$ have spread~$\geq \zeta'$.

    \casen{Verification.} First, observe that all the parameters are independent of $\zeta$, except $\zeta', q', \varepsilon'$ and $\varepsilon$.

    \cref{mainForNarrowTorsoQKC}.
    For every $t=t(\sigma)\in X$, the torso of $\sigma$ (which coincides with the torso of~$V_t$ by \cref{TDCfromNesteSetVerify}) is $(q',k)$-unbreakable by~\ref{keylemma:FlowerFinding:HighlyConnected}.

    \cref{mainForNarrowTorsoFlower}.
    For each $t=t(\sigma)\in Y$, set $F_t:=F_A$ for the unique $A\in\cA$ provided for~$\sigma$ by~\ref{keylemma:FlowerFinding:FlowerStructure}.
    So, in particular, $F_t$ has spread~$\ge\zeta' \geq \zeta$, and $F_t \in \cF$.
    Then $F_t$ $\rho$-refines $(T,\cV)$ at~$t$ by~\ref{keylemma:FlowerFinding:FlowerStructure}.
    $F_t$~was constructed using \cref{keylem:flower-approximation}, so the ``Furthermore'' part follows from \cref{keylem:flower-approximation}~\ref{item:keylem:flower-approximation-3}.

    \cref{mainForNarrowDisplayedInterlace}.
    Let $t=t(\sigma)\in Y$, let $A \in \cA$ be the asterisk provided for~$\sigma$ by~\ref{keylemma:FlowerFinding:FlowerStructure}. Recall that $F_t=F_A$.
    Note that 
    \[
        \rho
        = \eps_{\ref{NewerSplittingStarCloseToPetalStar}}(q,k,b,\gamma_N,\eps_N,\delta')
        \geq \eps_{\ref{NewerPetalSeparationsCloseToSplittingStar}}(k,\eps_N,\delta')
        = 2(k+\eps_N+\delta')
        \geq \delta',\eps_N~.
    \]
    If $s\in D(F_t)$ is $\rho$-crossed by another separation in~$D(F_t)$, then $   s\in\bigcup\cA(D(F_A),\delta')\subseteq A^\delta$ by \cref{NewerUniqueAsterisk} and $\rho \geq\delta'$.
    Since $A^\delta$ $\eps_N$-interlaces~$\sigma$ by \ref{keylemma:FlowerFinding:FlowerStructure}, it follows that~$s$ almost-$\rho$-interlaces~$\sigma$.
    
    \cref{mainForNarrowSeparationAnalysis}.
    For this, let $s$ be an arbitrary $\eps$-broad $k$-separation of~$G$.
    In order to apply \cref{lem:gamma-close-or-gamma-interlaces}, we need to show two things: that $s$ is $(2\gamma_N+1)$-broad and that $s$ is $\gamma_N$-nested with every separation in~$N$.
    The former is immediate since $s$ is $\eps$-broad with $\eps \geq 2\gamma_N+1$.
    For the latter, it suffices to show that $s$ is narrow, as the separations in $N$ are narrow-$\gamma_N$-nested.
    To show that $s$ is narrow, it suffices to show that $\eps\ge b$ as $G$ is $(b,k)$-narrow.
    And indeed, $\eps\ge\eps_N\ge b$.
    
    By \cref{lem:gamma-close-or-gamma-interlaces}, $s$ either is $\gamma_N$-close to some separation in~$N$ (in which case we are done as $\gamma_N\le\eps$), or there is a unique splitting star~$\sigma$ of~$N$ such that $s$ $\gamma_N$-interlaces~$\sigma$ (which we assume from now on).
    As $\eps_N\ge\gamma_N$, we in particular have that $s$ $\eps_N$-interlaces~$\sigma$ (unless $s$ is $\eps_N$-close to some separation in $\sigma \cup \sigma^* \se N$, in which case we are done as $\eps_N\le\eps$).
    So we have $s\in S(\sigma)$ with $S(\sigma)$ as in~\cref{keylemma:FlowerCriterion}.
    If $s$ is $\eps'$-close to some separation in $\sigma\cup\sigma^*$, then we are done as $\sigma\cup\sigma^*\se N$ and $\eps'\le\eps$.
    Otherwise, \cref{keylemma:FlowerCriterion}~\ref{keylemma:FlowerCriterion:FlowerType} applies (with $\zeta'$ in place of~$\zeta$),
    by which $s$ is $\eps_F$-close (hence $\eps$-close) to some displayed separation of~$F_t$ for $t=t(\sigma)$.

    \cref{mainForNarrowCanonical}.
    To show that $\cF$ is canonical, let $\psi$ be any automorphism of~$G$.
    Then $\psi$ induces a permutation of~$\cA$.
    By definition, for each $F\in\cF$, there exists $A\in\cA$ of spread~$\ge\zeta'$ such that~$F=F_A$. This $A$ is unique by \cref{highSpreadFlowersDistinct}.
    Let $\Psi$ be the permutation of $\cF$ that maps each $F_A$ to $\Psi(F_A):=F_{\psi(A)}$.
    Then $F_A$ is $\psi$-isomorphic to $\Psi(F_A)$, as the graphic-$k$-flowers in $\cF$ are canonical by \cref{keylem:flower-approximation}~\ref{flower-canonical}.

    Given~$N$ and~$(T,\cV)$, we already argued that the bipartition $X\sqcup Y=V(T)$ is canonical. Moreover, the assignment $Y\ni t\mapsto F_t\in\cF$ is defined avoiding arbitrary decisions, hence it is also canonical.
    The assignment is injective by \ref{keylemma:FlowerFinding:FlowerStructure} and \cref{NewerDistinctSplittingStars}.
    To see surjectivity, let $F_A\in\cF$.
    Recall that, by definition of~$\cF$, the asterisk~$A$ has spread~$\geq \zeta' \geq \zeta_{\ref{NewerInterlacesCloseToDisplayed}}(k,\delta',\eps_F)
    = \zeta_{\ref{NewerLargelyCrossingEntersAsterisk}}(\delta'+2k+1,\eps_F) = \delta'+2k+1+2\eps_F \geq \zeta_{\ref{NewerLargelyCrossingEntersAsterisk}}(\delta',\eps_F)$.
    So by \cref{HighSpreadUniqueSplittingStar}, 
    $A^\delta$ exists and $\eps_N$-interlaces a unique splitting star~$\sigma$ of~$N$.
    By \cref{NewerPetalSeparationsCloseToSplittingStar} and \cref{NewerSplittingStarCloseToPetalStar}, $\sigma$ and $\pi(F_A)$
    $\rho$-approximate each other.
    By definition of~$\cF$, both~$A$ and~$F_A$ have spread~$\geq\zeta'$, and \cref{NewerDistinctSplittingStars} ensures that no other $B\in\cA$ of spread~$\geq\zeta'$ has~$B^\delta$ $\eps_N$-interlacing~$\sigma$.
    Together with the $\rho$-approximation proved above, this shows that~$\sigma$ has type~\ref{keylemma:FlowerFinding:FlowerStructure} for the parameter~$\zeta'$.
    Thus $\sigma\in Y_N$ and $F_{t(\sigma)}=F_A$.
    It follows that the assignment $Y\ni t\mapsto F_t\in\cF$ is surjective, and hence bijective.
    
    Let $\kappa \coloneqq 2\eps_N$.
    By \cref{canonical:max-nested-set} we have that $N'$ is $\kappa$-canonical, and $S$ is canonical so the set of all $S$-totally-nested $k$-separations of $G$ too.
    Hence $N$ is $\kappa$-canonical, so $(T,\cV)$ as well by definition.

    \ref{mainForNarrowFlowersDiffer}.
    This is immediate from \cref{highSpreadFlowersDistinct}.

    \ref{mainForNarrowCloseToNestedSeparation}.
    We let~$\bar\eps := \max\{\eps'_{\ref{keylem:max-nested-set}}(q,k,b,\gamma_\pi), \delta'_{\ref{keylem:max-nested-set}}(k,b,\gamma_\pi) + b\}$.
    Let~$s$ be an $\bar\eps$-broad $k$-separation of~$G$.
    By the `moreover'-part of \cref{keylem:max-nested-set}, there is a $\delta'_{\ref{keylem:max-nested-set}}(k,b,\gamma_\pi)$-close $k$-separation~$s'$ of~$G$ that is nested with~$N'$.
    Observe that~$s'$ is $b$-broad and hence narrow since~$G$ is $(b,k)$-narrow.
    So,~$s'$ is nested with every narrow-nested $k$-separation of~$G$ and hence nested with~$N$.
\end{proof}

\subsection{The approximate cactus-decomposition for narrow graphs}\label{subsec:cactus-for-narrow}

To combine the tree-decomposition and flowers from \cref{mainForNarrow} into a cactus-decomposition, we first adjust each flower to the corresponding splitting star. We do this by perturbing its displayed separations to make them nested with the tree-decomposition, completing each attached branch to a clique, and then reconstructing a flower. The next lemmas ensure that these perturbations preserve crossing-connectedness and that each resulting clique is contained in a single bag.

\begin{lemma}\label{lem:close-separations-crossing-connected}
    Let~$G$ be a graph and~$\gamma,\eps\in\N$ with $\gamma \geq 2\eps$.
    Let~$S$ be a symmetric $\gamma$-crossing-connected set of separations of~$G$.
    For each separation~$s \in S$ let~$s'$ be an $\eps$-close separation of~$G$ such that $(s^*)' = (s')^*$ for every $s \in S$
    Denote $S':=\{s' : s \in S\}$.
    Then~$S'$ is $(\gamma-2\eps)$-crossing-connected.
\end{lemma}
\begin{proof}
    It suffices to show that~$s'$ and~$t'$ $(\gamma-2\eps)$-cross whenever~$s$ and~$t$ $\gamma$-cross.
    This follows from \cref{CloseVsCross}.
\end{proof}

\begin{lem}\label{lem:subpath}
    Let $G$ be a graph and let $F = (O, \cV)$ be a graphic-$k$-daisy of~$G$ that has the clover-property.
    Let $X \coloneqq \bigcap_{o \in V(O)} V_o$ be the centre of $F$.
    Let~$K$ be a clique of~$G$ with $K \not \subseteq X$.
    Then, the set of vertices $x \in V(O)$ such that $K \subseteq V_x$ induces a nonempty subpath of $O$.
    In particular, there is a vertex~$x \in V(O)$ whose bag~$V_x$ contains~$K$.
\end{lem}

\begin{proof}
    Since $F$ has the clover-property, $F$ has four adhesion-sets $V_{e_0},V_{e_1},V_{e_2},V_{e_3}$ with $V_{e_i} \cap V_{e_j} = X$ for all $i \neq j$.

    Fix an orientation of~$O$, so every edge of $O$ has a tail and a head.
    For every $v \in V(G)$, let $O_v$ be the set of all nodes $x \in V(O)$ such that $v \in V_x$.
    Let $u \in K \sm X$.
    Without loss of generality, we may assume that the edges $e_0, e_1, e_2$ and $e_3$ appear in that cyclic order along $O$.
    Since $u \notin X$ we have that~$u$ belongs to at most one of the sets $V_{e_0},V_{e_1},V_{e_2},V_{e_3}$.
    By symmetry, we may assume that all nodes in~$O_u$ lie between the head of~$e_0$ and the tail of~$e_2$ (both included) in~$O$.
    Since~$K$ is complete, any two vertices of~$K$ must appear in a common bag.
    Thus, for every vertex $v \in K \sm X$, all nodes in $O_v$ lie between the head of~$e_0$ and the tail of~$e_2$ (both included) in~$O$, or between the head of~$e_1$ and the tail of~$e_3$ (both included) in $O$, or between the head of~$e_3$ and the tail of~$e_1$ (both included) in~$O$.
    Observe that the last two possibilities are incompatible: there do not exist $v, v' \in K \sm X$ such that all nodes in~$O_v$ lie between the head of~$e_1$ and the tail of~$e_3$ (both included) in~$O$ and all nodes in~$O_{v'}$ lie between the head of~$e_3$ and the tail of~$e_1$ (both included) in~$O$.
    By symmetry, say that the last option does not happen for any $v \in K$.
    Then, the tail~$y$ of~$e_0$ satisfies $V_y \cap (K \sm X) = \emptyset$.

    For every vertex $u \in K \sm X$, $O_u$ induces a connected subgraph of $O-y$, hence a subpath of $O-y$, and for all $u \neq v \in K \sm X$ we have $O_u \cap O_v \neq \emptyset$.
    By the Helly property of subpaths of a path, the set $\bigcap_{v \in K \sm X}O_v$ induces a nonempty subpath of $O-y$.
    Since $X = \bigcap_{o \in V(O)}V_o$, it follows that the set of all vertices $x \in V(O)$ such that $K \subseteq V_x$ induces a nonempty subpath of $O-y$.
    This concludes the proof.
\end{proof}

A cactus-decomposition $\mathcal{C}$ of a graph $G$ is \defn{clique-faithful} if every clique of $G$ is contained in a bag of a node of $C$ that is not a $P$-node.

\begin{theorem}\label{thm:cactusForNarrow}
    For every $q,k,b \in \N$ there exists $\eta = \eta_{\ref{thm:cactusForNarrow}}(q, k, b) \in \N$ such that every $(q,k-1)$-unbreakable $(b,k)$-narrow graph~$G$ has a clique-faithful cactus-decomposition~$\cC$ such that every displayed separation has order~$\leq k$ and the set of displayed separations of $\cC$ $\eta$-approximates the set of all $k$-separations of~$G$.
\end{theorem}

\begin{proof}
    Let $q,k,b \in \N$ be arbitrary.
    Let $\gamma,\kappa,\rho,\bar\eps \in \N$ be the constants guaranteed by \cref{mainForNarrow}.
    We define
    \begin{align*}
        \delta &:= \max\{ \max\{\gamma_{\ref{HugeSeparationsUniformlyCrossing}}(q,k),\, \gamma_{\ref{keylem:flower-approximation}}(q,k)\} + 2\bar\eps+1,\, \max\{b, \,\beta_{\ref{keylem:flower-approximation}}(q,k), \,\gamma+\rho\} + \bar\eps + k+1\}~,\\
        \zeta&:=\delta'_{\ref{NewerUniqueAsterisk}}(k,\delta)~.
    \end{align*}
    
    Then by \cref{mainForNarrow} (applied to these choices of $q,k,b,\zeta$) there exist constants~$q',\eps \in \N$ and a set~$\cF$ of graphic-$k$-flowers of~$G$ of spread $\geq \zeta$ and a symmetric nested set~$N$ of narrow-$\gamma$-nested proper $k$-separations of~$G$ with corresponding tree-decomposition $(T,\cV)$ of~$G$ and a bipartition $X \sqcup Y = V(T)$ such that \cref{mainForNarrow}~\ref{mainForNarrowTorsoQKC}--\ref{mainForNarrowCloseToNestedSeparation} hold.
    For every node $t \in V(T)$, we denote by~$\sigma_t$ the splitting star of~$N$ corresponding to $t$.

    For a flower~$F_t$ (with $t \in Y$) and $\alpha \in \N$,
    we define~$S_{t,\alpha}$ to be the set of displayed separations of~$F_t$ that are $\alpha$-crossed by some displayed separation of~$F_t$.
    Moreover, we let~$P_{t,\alpha} := D(F_t) \sm S_{t,\alpha}$ be the set of displayed separations of~$F_t$ that are $\alpha$-nested with~$D(F_t)$.
    For every $\alpha\in\N$ with $\alpha\geq\gamma_{\ref{HugeSeparationsUniformlyCrossing}}(q,k)$, the set~$S_{t,\alpha}$ is uniformly $\alpha$-crossing by \cref{lem:gamma-crossed-uniform}.
    For every such~$\alpha$,~$\cA(D(F_t),\alpha)$ denotes the set of all $\alpha$-asterisks of~$S_{t,\alpha}$, which is well-defined (see \cref{defn:cA-gamma-asterisks}) and partitions~$S_{t,\alpha}$.

    \begin{claim}\label{claim:cactus-asteriskal-sets}
        There is a unique nontrivial $A_t\in\cA(D(F_t),\delta)$ such that $S_{t,\zeta}\subseteq A_t$.
    \end{claim}
    \begin{cproof}
        Since~$F_t$ has spread $\geq \zeta$, we have $\cA(D(F_t),\zeta) \neq \emptyset$.
        So we can apply \cref{NewerUniqueAsterisk} (since $\delta \geq \gamma_{\ref{HugeSeparationsUniformlyCrossing}}(q,k)$ and $\zeta \geq \delta'_{\ref{NewerUniqueAsterisk}}(k,\delta)$) and get that there is a unique $A_t \in \cA(D(F_t),\delta)$ with $S_{t,\zeta} = \bigcup \cA(D(F_t),\zeta) \subseteq A_t$.
        Since $S_{t,\zeta}$ contains two $\zeta$-crossing separations, and $\zeta\geq\delta$, the $\delta$-asterisk~$A_t$ is nontrivial.
    \end{cproof}

    Since $A_t$ is a nontrivial $\delta$-asterisk, every separation in~$A_t$ is $\delta$-crossed by another separation in~$A_t$.
    Therefore every separation in~$A_t$ is $2(\delta-k)$-broad by \cref{LargeCrossedImpliesBroad}, hence $b$-broad since $\delta \geq b+k$.
    Similarly, every separation in~$A_t$ has breadth~$> q$ since $\delta \geq \beta_{\ref{keylem:flower-approximation}}(q,k)+k+1 \geq q+k+1$. 
    Since~$G$ is $(q,k-1)$-unbreakable, every separation in~$A_t$ is a $k$-separation, hence is narrow since~$G$ is $(b, k)$-narrow.
    
    Additionally, every separation in~$A_t$ is $\bar\eps$-broad since $\delta \geq \bar\eps + k$.
    By \cref{mainForNarrow}~\ref{mainForNarrowCloseToNestedSeparation}, every separation~$s$ in~$A_t$ is $\bar\eps$-close to some $k$-separation~$s'$ that is nested with~$N$.
    For each $s \in A_t$ we fix one such~$s'$, chosen so that $(s^*)' = (s')^*$ for all $s \in A_t$.
    We let~$A_t' = \{s' : s \in A_t\}$.
    Since~$A_t$ is a $\delta$-asterisk, it is $\delta$-crossing-connected, so~$A_t'$ is $(\delta-2\bar\eps)$-crossing-connected by \cref{lem:close-separations-crossing-connected}, with $\delta-2\bar\eps > \gamma_{\ref{HugeSeparationsUniformlyCrossing}}(q,k), \gamma_{\ref{keylem:flower-approximation}}(q,k)$.
    Moreover, note that every separation in~$A_t'$ is $(2\delta-2k-\bar\eps)$-broad by \cref{LowBreadthClose}.

    \begin{claim}
    \label{claim:cactusForNarrow-almost-interlacing}
        Let $s \in A_t$.
        Then~$s'$ almost-interlaces~$\sigma_t$.
    \end{claim}
    \begin{cproof}
        First, we argue that~$s$ almost-$\rho$-interlaces~$\sigma_t$.
        Since~$A_t$ is a nontrivial $\delta$-asterisk, there is a separation in~$A_t$ that $\delta$-crosses~$s$.
        Since $\delta \geq\rho$, \cref{mainForNarrow}~\ref{mainForNarrowDisplayedInterlace} implies that~$s$ almost-$\rho$-interlaces~$\sigma_t$.
        
        Suppose that~$s'$ does not almost-interlace~$\sigma_t$.
        Since~$s'$ is nested with~$N$, there is a separation $n \in \sigma_t$ with $s' \leq n$ or $(s')^* \leq n$, say the former.
        Since~$s$ almost-$\rho$-interlaces~$\sigma_t$, we have $n\leq_\rho s$ or $n\leq_\rho s^*$.
        In the latter case, we have $s' \leq n \leq_\rho s^* \approx_{\bar\eps} (s')^*$, so $s' \leq_{\rho + \bar\eps}(s')^*$, which implies that~$s'$ has breadth~$\leq \rho+\bar\eps$ by \cref{lem:assymmetric-nested-implies-not-broad}, contradicting that~$s' \in A_t'$ has breadth $\geq 2\delta-2k-\bar\eps > \rho+\bar\eps$.
        In the former case, we have $s \approx_{\bar\eps} s' \leq n \leq_\rho s$.
        In particular, $s \approx_{\rho+\bar\eps} n$.
        Recall that~$n \in N$ is narrow-$\gamma$-nested and that~$s$ is narrow and proper.
        By \cref{obs:totally-nested-close-separation},~$s$ is narrow-$(\gamma+\rho+\bar\eps)$-nested.
        With $\delta \geq \gamma+\rho+\bar\eps$, $s$ is narrow-$\delta$-nested.
        However, since $s \in A_t$, there is some $k$-separation~$r \in A_t$ that $\delta$-crosses~$s$.
        Recall that~$r$ is (proper and) narrow.
        This contradicts that~$s$ is narrow-$\delta$-nested.
    \end{cproof}

    For $t \in Y$, we let~$\tilde G_t$ be the graph obtained from~$G$ by making $G[A]$ into a clique for every $(A,B) \in \sigma_t$.
    So,~$\tilde G_t$ is again $(q,k-1)$-unbreakable and $(b,k)$-narrow.
    Note that every separation in~$A_t'$ is a separation of~$\tilde G_t$ since it almost-interlaces~$\sigma_t$ by \cref{claim:cactusForNarrow-almost-interlacing}.
    Recall that every separation in~$A_t'$ is $(2\delta-2k-\bar\eps)$-broad, and hence $q$-broad by our choice of~$\delta$.
    Thus, $A'_t$ is uniformly $\gamma_{\ref{HugeSeparationsUniformlyCrossing}}(q,k)$-crossing by \cref{HugeSeparationsUniformlyCrossing}, hence uniformly $(\delta-2\bar\eps)$-crossing since $\delta-2\bar\eps \geq \gamma_{\ref{HugeSeparationsUniformlyCrossing}}(q,k)$.
    Recall that the set~$A_t'$ is symmetric.
    Moreover, since~$A_t$ is nontrivial, there are $\delta$-crossing separations~$r,s$ in~$A_t$.
    By \cref{lem:close-separations-still-cross}, the separations~$r',s'$ in~$A_t'$ $(\delta-2\bar\eps)$-cross, so~$A_t'$ is nontrivial.
    Therefore, $A'_t$ is nontrivial $(\delta-2\bar\eps)$-asteriskal (both in $G$ and in $\tilde{G}_t$).
    
    Now we apply \cref{keylem:flower-approximation} to~$\tilde G_t$ and~$A'_t$ with parameter $\delta-2\bar\eps > \gamma_{\ref{keylem:flower-approximation}}(q, k)$.
    We can do it since the separations in~$A'_t$ are $(2(\delta-k)-\bar\eps)$-broad with $2(\delta-k)-\bar\eps \geq \beta_{\ref{keylem:flower-approximation}}(q, k),b$ and hence also narrow.
    This gives a graphic-$k$-flower~$\tilde F_t$ such that \cref{keylem:flower-approximation}~\ref{item:keylem:flower-approximation-1}--\ref{flower-canonical} hold.
    In particular, every separation in~$A'_t$ is $\eps_{\ref{keylem:flower-approximation}}$-close to some displayed separation of~$\tilde F_t$ by \cref{keylem:flower-approximation}~\ref{item:keylem:flower-approximation-1}.
    We choose \[\eta' := \eps + \max\{\,\bar\eps+\eps_{\ref{keylem:flower-approximation}} \,,\, 2(\zeta+k)+\rho\,\}.\]

    \begin{claim}
    \label{claim:cactus-approx}
        The set $N \cup \bigcup_{t \in Y}D(\tilde F_t)$ $\eta'$-approximates the set of all $k$-separations of~$G$.
    \end{claim}
    \begin{proof}
        Let~$s$ be an $\eta'$-broad $k$-separation of~$G$.
        We need to show that~$s$ is $\eta'$-close to a separation in~$N$ or to a displayed separation of some~$\tilde F_t$.
        By \cref{mainForNarrow}~\ref{mainForNarrowSeparationAnalysis},~$s$ is $\eps$-close to a separation in~$N$ or $\eps$-close to a displayed separation of some~$F_t$.
        We are done in the former case (since $\eta' \geq \eps$), so assume the latter.
        
        So,~$s$ is $\eps$-close to some displayed separation~$d$ of some~$F_t$ (with $t \in Y$).
        Recall that $D(F_t)=S_{t,\zeta} \sqcup P_{t,\zeta}$.
        If $d \in S_{t,\zeta} \subseteq A_t$, then there is a separation $d' \in A'_t$ that is $\bar\eps$-close to~$d$.
        By \cref{keylem:flower-approximation}~\ref{item:keylem:flower-approximation-1}, there is a displayed separation~$d''$ of~$\tilde F_t$ that is $\eps_{\ref{keylem:flower-approximation}}$-close to~$d'$.
        Thus, by \cref{approxLeQuasiTransitive} and with $\eta' \geq \eps + \bar\eps + \eps_{\ref{keylem:flower-approximation}}$,~$s$ is $\eta'$-close to~$d''$ and we are done.
        Otherwise, $d$ is in~$P_{t,\zeta}$.
        So,~$d$ is not $\zeta$-crossed (and hence not $\zeta$-crissed) by a displayed separation of~$F_t$.
        Then~$d$ is $2(\zeta+k)$-close to a petal-separation~$p$ of~$F_t$ or its inverse by \cref{FlowerClosePetalOrCrossed}.
        Then, $p$ or~$p^*$ is $\eps + 2(\zeta+k)$-close to~$s$ by \cref{approxLeQuasiTransitive}, so both are $\rho$-broad by \cref{LowBreadthClose}.
        Thus,~$p$ is $\rho$-close to a separation of~$N$ by \cref{mainForNarrow}~\ref{mainForNarrowTorsoFlower} (so its inverse too because $N$ is symmetric).
        Therefore, by \cref{approxLeQuasiTransitive} and with $\eta' \geq \eps + 2(\zeta+k)+\rho$,~$s$ is $\eta'$-close to a separation in~$N$ and we are done.
    \end{proof}

    We are now ready to build the desired cactus-decomposition~$\cC = (C, \cW)$ of~$G$.
    \begin{claim}
    \label{claim:cactus-XY-edges}
        For every $t \in Y$ and every separation~$n=(A,B) \in \sigma_t$, there is a bag of $\tilde F_t$ that contains $A$.
    \end{claim}
    \begin{cproof}
        Let~$t \in Y$ and~$n=(A,B) \in \sigma_t$.
        Recall that $\tilde G_t[A]$ is complete and that $\tilde F_t$ is a graphic-$k$-flower of $\tilde G_t$.
        First, we consider the case where~$\tilde F_t$ is a graphic-$k$-daisy, in which case~$\tilde F_t$ has the clover-property by~\cref{keylem:flower-approximation}~\ref{item:keylem:flower-approximation-3}.
        Then, some bag of $\tilde F_t$ contains $A$: this is obvious if $A$ is contained in the centre of~$\tilde F_t$, and follows from \cref{lem:subpath} otherwise.
        Now, assume that~$\tilde F_t$ is a graphic-$k$-anemone.
        Since graphic-$k$-anemones are star-decompositions (hence tree-decompositions), it follows from the Helly property of subtrees of a tree that some bag of~$\tilde F_t$ contains~$A$.
    \end{cproof}
    
    Without loss of generality, we can assume that no two nodes of~$Y$ are adjacent in~$T$ (otherwise, we can subdivide this edge~$yy'$ of~$T$, add the subdivision-node to~$X$ and define the bag of the subdivision-node to be the intersection $V_y \cap V_{y'}$).
    For each $t\in Y$, we write $\tilde F_t=(H_t,(V'_x:x\in V(H_t)))$.
    We construct the graph~$C$ from the disjoint union of the forest~$T[X]$ and the graphs~$H_t$ for $t \in Y$.
    For every $t \in Y$ and every edge $tt'$ of $T$ (hence with $t' \in X$), let $n = (A, B) \in \sigma_t$ be the separation induced by $tt'$ in $(T, \cV)$.
    By \cref{claim:cactus-XY-edges}, there is a node $x\in V(H_t)$ with $A\subseteq V'_x$.
    Pick such a node arbitrarily and add the edge~$xt'$ to~$C$.
    We define~$C$ to be the decomposition graph of~$\cC$.
    Now, we define the bags of~$\cC$.
    For a node $t \in X$, we define its bag to be~$W_t \coloneqq V_t$.
    For a node~$x$ of a decomposition-graph of some~$\tilde F_t$ with $t\in Y$, with bag $V'_x$ in $\tilde F_t$, we define its bag to be $W_x \coloneqq V'_x \cap V_t$.
    For each graphic-$k$-anemone $\tilde{F}_t$, we declare the centre of $\tilde{F}_t$ as a $P$-node. There are no other $P$-nodes.
    This concludes the definition of~$\cC = (C, \cW)$.

    \begin{claim}\label{cl:is-cactus}
        $\cC$ is a cactus-decomposition of~$G$.
    \end{claim}
    \begin{cproof}
        Observe that~$C$ is connected and that the only cycles in $C$ are the decomposition-cycles of the graphic-$k$-daisies.
        Hence~$C$ is a strong cactus.
        It remains to verify~\ref{itm:graph-decomposition-1} and~\ref{itm:graph-decomposition-2}.

        \ref{itm:graph-decomposition-1}.
        For each $t\in Y$, restricting the bags of~$\tilde F_t$ to~$V_t$ gives a graph-decomposition of~$G[V_t]$, since~$\tilde F_t$ is a graph-decomposition of~$\tilde G_t \supseteq G$.
        Thus every vertex and every edge of~$G$ is contained in some~$W_z$ (with $z \in V(C)$), since this holds for the original tree-decomposition~$(T,\cV)$.

        \ref{itm:graph-decomposition-2}.
        Let $v \in V(G)$ be arbitrary.
        The nodes of~$T$ whose bag~$V_t$ contains~$v$ induce a subtree~$T_v$ of~$T$.
        Each node $t\in Y\cap V(T_v)$ is replaced by the nonempty connected subgraph of~$H_t$ induced by the nodes of~$H_t$ whose bag contains~$v$.
        If an edge~$tt'$ of~$T_v$ (with $t\in Y$ and $t' \in X$) is replaced in $C$ by~$xt'$ for some $x \in V(H_t)$, then $v\in V_t\cap V_{t'}\subseteq A\subseteq V'_x$, where $(A,B)\in\sigma_t$ is the separation induced by~$tt'$.
        Thus both~$W_x$ and~$W_{t'}$ contain~$v$.
        Every edge of~$T_v$ between two nodes of~$X$ is also an edge of~$C$.
        Therefore the nodes of~$C$ whose bag contains~$v$ induce a connected subgraph of $C$.
    \end{cproof}

    \begin{claim}\label{cl:small-adhesion}
        Every displayed separation of~$\cC$ has order~$\leq k$.
    \end{claim}
    \begin{proof}
        Let~$(A,B)$ be a displayed separation of $\cC$ induced by a cut~$(U,W)$ of~$C$, and let~$F$ be the set of edges of~$C$ with one end in~$U$ and the other in~$W$.
        We distinguish the three types of displayed separations.

        \casen{Case 1: $|F|\leq1$.}
        If $F=\emptyset$, then $A\cap B=\emptyset$, so the assertion holds.
        Assume now that $|F|=1$.
        If the unique edge of~$F$ arises from an edge of~$T$, then $(A,B)\in N$ and hence has order~$k$.
        Otherwise, it lies in the star~$H_t$ of some graphic-$k$-anemone~$\tilde F_t$, so $A\cap B$ is contained in the central bag of~$\tilde F_t$ and hence has size~$\leq k$.

        \casen{Case 2: $|F|=2$ and both edges of~$F$ lie on a common cycle of~$C$.}
        This cycle is~$H_t$ for some graphic-$k$-daisy~$\tilde F_t$.
        By~\ref{itm:graph-decomposition-2} and since $W_x\subseteq V'_x$ for every $x\in V(H_t)$, we have
        \[
            A\cap B
            =\bigcup_{xy\in F}(W_x\cap W_y)
            \subseteq\bigcup_{xy\in F}(V'_x\cap V'_y).
        \]
        The set on the right is the separator of the separation displayed by these two edges in~$\tilde F_t$, and hence has size~$\leq k$.

        \casen{Case 3: $|F|\geq 2$ and all edges are incident with a $P$-node of~$C$.}
        Let~$c$ be this $P$-node.
        By definition of the $P$-nodes, $c$ is the centre of the star $H_t$ for some graphic-$k$-anemone $\tilde{F}_t$.
        By~\ref{itm:graph-decomposition-2}, we have
        \[
            A\cap B
            =\bigcup_{cx\in F}(W_c\cap W_x)
            \subseteq W_c
            \subseteq V'_c.
        \]
        Since~$V'_c$ is the central bag of a graphic-$k$-anemone, it has size~$\leq k$.
    \end{proof}

    \begin{claim}
        \label{claim:cactusForNarrow-extend}
        For every $t \in Y$ and every separation $d = (A, B) \in D(\tilde{F}_t)$, there is a separation $d' = (A', B')$ displayed by $\cC$ such that $A' \subseteq A$ and $B' \subseteq B$.
        In particular, $d$ and $d'$ are $k$-close.
    \end{claim}
    \begin{cproof}
        Let $t \in Y$ and $d = (A, B) \in D(\tilde{F}_t)$.
        Write $\tilde F_t = (H, \cV')$.
        If $A=\emptyset$ or $B=\emptyset$, we may take $d'=d$.
        Hence assume that both sides are nonempty.
        By the construction of $C$, every component of $C-E(H)$ contains exactly one node of $H$.
        For each $x \in V(H)$, let $K_x$ denote the component containing $x$.

        First, note that for every $x \in V(H)$, we have $\bigcup_{z \in V(K_x)}W_z \subseteq V'_x$.
        To see it, observe first that $W_x = V'_x \cap V_t \subseteq V'_x$.
        Then, every other node of $K_x$ comes from a component of $T-t$ corresponding to a separation $n = (A_n, B_n) \in \sigma_t$ whose chosen attachment node in $H$ is $x$.
        Each bag of $(C, \cW)$ coming from this component of $T-t$ is contained in a bag of $(T, \cV)$ from that component, and hence in $A_n$.
        Since $x$ was chosen as the attachment node in $H$ for this component, we have $A_n \subseteq V'_x$, as desired.

        Let $(U, W)$ be a cut of $H$ displaying $d$.
        Consider the cut $(U', W') \coloneqq \left(\bigcup_{x \in U}V(K_x), \bigcup_{x \in W}V(K_x)\right)$ of $C$.
        The edges between $U'$ and $W'$ in $C$ are exactly the edges between $U$ and $W$ in $H$.
        If $H$ is a cycle, these are two edges on the same cycle, and if $H$ is a star, these are edges incident to the centre of the star, which is a $P$-node.
        In both cases, the separation $d' = (A', B')$ induced by $(U', W')$ is displayed by $\cC$.
        
        By the previous observation, we have $A' = \bigcup_{x \in U} \bigcup_{z \in V(K_x)} W_z \subseteq \bigcup_{x \in U} V'_x = A$ and similarly $B' \subseteq B$.
        As $d$ and $d'$ have order at most $k$, it follows that they are $k$-close.
    \end{cproof}
    
    \begin{claim}
        \label{claim:cactusForNarrow-approximation}
        Let $\eta:=\eta'+k$.
        Then the set of displayed separations of~$\cC$ $\eta$-approximates the set of all $k$-separations of~$G$.
    \end{claim}
    \begin{cproof}
        Let~$s$ be an $\eta$-broad $k$-separation of~$G$.
        By \cref{claim:cactus-approx},~$s$ is $\eta'$-close to some separation in $N \cup \bigcup_{t \in Y}D(\tilde F_t)$.
        Note that every separation in~$N$ is displayed by~$\cC$.
        So, if~$s$ is $\eta'$-close to a separation in~$N$, it is $\eta$-close to a displayed separation of~$\cC$.
        Otherwise, $s$ is $\eta'$-close to some displayed separation $d=(A,B) \in D(\tilde F_t)$ for some $t \in Y$.
        Then, $d$ is $k$-close to a separation displayed by $\cC$ by \cref{claim:cactusForNarrow-extend}.
        Therefore,~$s$ is $(\eta'+k)$-close to a displayed separation of~$\cC$.
    \end{cproof}

    \begin{claim}
        \label{claim:cactusForNarrow-cliques}
        $\mathcal{C}$ is clique-faithful.
    \end{claim}
    \begin{cproof}
        Let~$K$ be a clique of~$G$.
        By the Helly property of subtrees of a tree, there is a node~$t$ of $T$ such that~$K \subseteq V_t$.
        If $t \in X$, then we are done since $V_t = W_t$ and $t$ is not a $P$-node.
        Otherwise, we have $t \in Y$.
        Recall that~$\tilde F_t$ is a graphic-$k$-flower of~$\tilde G_t \supseteq G$.
        First assume that~$\tilde F_t = (O, \cV')$ is a graphic-$k$-daisy. 
        Then~$\tilde F_t$ has the clover-property.
        Therefore, there is a node $x$ of $O$ such that $K \subseteq V'_x$: this is obvious if $K$ is contained in the centre of~$\tilde F_t$, and follows from \cref{lem:subpath} otherwise.
        Moreover, $x$ is not a $P$-node since~$\tilde F_t$ is a graphic-$k$-daisy.
        Since $K \subseteq V_t$, we have $K \subseteq V'_x \cap V_t = W_x$.
        Suppose now that~$\tilde F_t = (S, \cV')$ is a graphic-$k$-anemone, and let $c$ be the centre of $S$.
        Since~$\tilde F_t$ is a tree-decomposition, there is a node $x$ of $S$ such that~$K \subseteq V'_x$.
        Furthermore, since $V'_c$ is contained in all other bags of $\tilde F_t$, we can choose $x$ distinct from $c$, so that $x$ is not a $P$-node.
        Since $K \subseteq V_t$, we have $K \subseteq V'_x \cap V_t = W_x$.
    \end{cproof}

    \cref{cl:is-cactus,cl:small-adhesion,claim:cactusForNarrow-approximation,claim:cactusForNarrow-cliques} together conclude the proof of the theorem.
\end{proof}

\section{The approximate cactus-decomposition}\label{sec:RemoveNarrow}

In this section, we prove the main result, \cref{thm:main-intro}, in full generality.
Our strategy for this is to show that we can reduce to the case where $G$ is $(\beta, k)$-narrow. 
We do this in \cref{subsec:reduce-narrow}, and then prove \cref{thm:main-intro} in \cref{subsec:the-proof}.

\subsection{Reducing to a narrow graph}\label{subsec:reduce-narrow}

Let~$s$ be a separation of a graph~$G$, and let~$K$ be a component of $G-\sep(s)$.
We say that the separation $(V(K) \cup N(K), V(G-K))$ is a \defn{componental separation at~$s$}.
Note that $s^*$ might not be componental.
A separation~$s'$ is \defn{$(\beta,k)$-irrelevant} if it is a componental $(<k)$-separation at some $\beta$-broad $k$-separation of~$G$.
Recall that~$G$ is \defn{$(\beta,k)$-narrow} if every $\beta$-broad $k$-separation is narrow.

Our goal in this section is to show that we can easily reduce from the case where $G$ is $(q,k-1)$-unbreakable to the case where $G$ is additionally $(\beta, k)$-narrow.
To do so, we show that the $(\beta, k)$-irrelevant separations form a star and that the torso obtained by cutting at all of them is $(q,k-1)$-unbreakable and $(\beta,k)$-narrow.

\begin{restatable}{keylemma}{irrelevantstardec}\label{keylem:irrelevant-star-dec}
    Let~$G$ be a $(q,k-1)$-unbreakable graph and let $\beta > \beta_{\ref{cor:irrelevant-separations-disjoint}}(q,k)$.
    Then the set~$\sigma$ of all $(\beta,k)$-irrelevant separations forms a star.
    Moreover, the torso of~$\sigma$ is $(q,k-1)$-unbreakable and $(\beta,k)$-narrow.
\end{restatable}

The following lemma motivates the study of $(\beta,k)$-irrelevant separations.

\begin{lem}\label{lem:narrow-equivalent-no-irrelevant}
    Let~$G$ be a graph.
    Then~$G$ is $(\beta,k)$-narrow if and only if it does not have a $(\beta,k)$-irrelevant separation.
\end{lem}
\begin{proof}
    If~$s'$ is a $(\beta,k)$-irrelevant separation, then it is a componental $(<k)$-separation at some $\beta$-broad $k$-separation~$s$ of~$G$.
    Then~$s$ is not narrow, so~$G$ is not $(\beta,k)$-narrow.

    On the other hand, if~$G$ is not $(\beta,k)$-narrow, then there is a $\beta$-broad $k$-separation~$s$ of~$G$ that is not narrow.
    Hence, there is a componental $(<k)$-separation~$s'$ at~$s$.
    Then~$s'$ is $(\beta,k)$-irrelevant.
\end{proof}

We now prove some results about $\beta$-broad $k$-separations and $(\beta, k)$-irrelevant separations in $(q,k-1)$-unbreakable graphs.

\begin{lem}\label{lem:irrelevant-small}
    Let $G$ be a $(q,k-1)$-unbreakable graph, let $\beta > q$ and let $(C, D)$ be a $(\beta, k)$-irrelevant separation of $G$.
    Then, $|C \sm D| \leq q$.
\end{lem}

\begin{proof}
    Since $(C, D)$ is a $(\beta, k)$-irrelevant separation of $G$, there exists a $\beta$-broad $k$-separation $s = (A, B)$ of $G$ such that $(C, D)$ is a componental $(<k)$-separation at $s$.
    After possibly swapping the names $A$ and $B$, we may assume without loss of generality that $C \subseteq A$, so $D \supseteq B$.
    Hence $|D \sm C| \geq |B \sm A| \geq \beta > q$.
    Since $(C, D)$ is a $(<k)$-separation and $G$ is $(q,k-1)$-unbreakable, we then have $|C \sm D| \leq q$.
\end{proof}

\begin{lem}
\label{lem:irrelevant-separations-disjoint}
    Let~$s$ and~$t$ be two $\beta$-broad $k$-separations of a $(q,k-1)$-unbreakable graph~$G$.
    Let~$K$ be a component of $G - \sep(s)$ of size~$h$.
    Assume $\beta > 2(q+h)+k$.
    Then~$V(K)$ is disjoint from~$\sep(t)$.
\end{lem}
\begin{proof}
    Let~$s=:(A,B)$ and $t=:(C,D)$.
    Assume for a contradiction that~$V(K)$ is not disjoint from~$\sep(t)$.
    Without loss of generality,~$K$ is a component of~$G[A \sm B]$ and intersects the $At$-link.
    We define $\hat s:=(\hat A,\hat B):=(A \sm V(K),B\cup V(K))$. 

    By \cref{lem:two-opposite-corners}, there are two opposite corner-separations of~$s$ and~$t$ of order~$k$ and breadth $\geq (\beta-k)/2$, say the $AC$-corner-separation and the $BD$-corner-separation.
    Now, consider the two separations~$\hat s$ and~$t$.
    Then the $\hat AC$-corner has size~$\geq (\beta-k)/2 - h$.
    Hence, the corresponding corner-separation has breadth~$\geq (\beta-k)/2 - h > q$.
    Since~$K$ intersects the $At$-link, the $\hat AC$-corner-separation has order~$<k$.
    This contradicts that~$G$ is $(q,k-1)$-unbreakable.
\end{proof}

\begin{cor}
\label{cor:irrelevant-separations-disjoint}
    Let~$G$ be a $(q,k-1)$-unbreakable graph.
    Let $\beta > \beta_{\ref{cor:irrelevant-separations-disjoint}}(q,k) \coloneqq 4q+k$.
    Let~$s'=(A',B')$ be a $(\beta,k)$-irrelevant separation, and let~$t$ be a $\beta$-broad $k$-separation of~$G$.
    Then~$A' \sm B'$ is disjoint from~$\sep(t)$.
\end{cor}
\begin{proof}
    Since~$s'$ is $(\beta,k)$-irrelevant, there is a $\beta$-broad $k$-separation~$s$ such that~$s'$ is a componental $(<k)$-separation at~$s$.
    That is, there is a component~$K$ of $G-\sep(s)$ such that $A' \sm B' = V(K)$.
    Note that $h:=|V(K)| \leq q$ since~$G$ is $(q,k-1)$-unbreakable.
    Hence we have $\beta > 4q+k \geq 2(q+h)+k$ and we can apply \cref{lem:irrelevant-separations-disjoint} to the separations~$s$ and~$t$.
    Therefore, $V(K)=A' \sm B'$ is disjoint from~$\sep(t)$.
\end{proof}

\begin{lem}\label{irrelevantSeparationsFormStar}
    Let~$G$ be a $(q,k-1)$-unbreakable graph and let $\beta > \beta_{\ref{cor:irrelevant-separations-disjoint}}(q,k)$.
    Then the set~$\sigma$ of all $(\beta,k)$-irrelevant separations of~$G$ forms a star.
\end{lem}
\begin{proof}
    Let~$s'\neq t'$ be two $(\beta,k)$-irrelevant separations of~$G$.
    That is, there are $\beta$-broad $k$-separations $s \coloneqq (A, B)$ and~$t \coloneqq (C, D)$ of~$G$, such that~$s'$ and~$t'$ are componental $(<k)$-separations at~$s$ and~$t$, respectively.
    Let 
    \begin{align*}
        &(A',B'):=s'=(V(K)\cup N(K), V(G-K))\\
        \text{and }& (C',D'):=t'=(V(L)\cup N(L), V(G-L))
    \end{align*}
    where~$K$ is a component of $G-\sep(s)$ and~$L$ is a component of $G-\sep(t)$.

    First, we show that~$V(L)$ is disjoint from~$V(K) \cup N(K)$.
    Consider the two separations~$s'$ and~$t$.
    By \cref{cor:irrelevant-separations-disjoint},~$V(K)$ is disjoint from~$\sep(t)$.
    Since~$K$ is connected, we have without loss of generality $V(K) \subseteq C \sm D$.
    Similarly, we have $V(L) \subseteq A \sm B$, say.
    Hence $K$ and $L$ are components of $G-(\sep(s)\cup\sep(t))$. 
    If $K$ and $L$ intersect, then $K=L$, and so $s' = t'$, a contradiction.
    Therefore, $K$ and $L$ are disjoint.
    \cref{cor:irrelevant-separations-disjoint} yields $V(L) \cap N(K) \subseteq V(L) \cap \sep(s) = \emptyset$.

    Since~$V(L)$ is disjoint from~$V(K) \cup N(K)$, we have $V(L) \subseteq B' \sm A'$ and $C'=V(L) \cup N(L) \subseteq B'$.
    By symmetry, we also have $A' \subseteq D'$.
    Therefore, $s'=(A',B') \leq (D',C')=(t')^*$.
    Hence,~$\sigma$ is a star.
\end{proof}

\begin{lem}\label{lem:few-vtces-attach-to-small-sep}
    Let $G$ be a $(q,k-1)$-unbreakable graph such that $|V(G)| \geq \nu_{\ref{lem:few-vtces-attach-to-small-sep}}(q,k):= 3q+k$.
    Let $X$ be a set of less than $k$ vertices of~$G$.
    Then, the union of all connected components of $G-X$ of size~$\leq q$ has size~$\leq q$.
\end{lem}

\begin{proof}
    By contradiction, suppose that more than $q$ vertices belong to a connected component of $G-X$ of size at most $q$.
    Denote these connected components by $C_1, \ldots, C_c$.
    Let $i \in [c]$ be minimum such that $|C_1 \cup \dots \cup C_i| > q$ (which exists since $|C_1 \cup \dots \cup C_c| > q$).
    Then, $|C_1 \cup \dots \cup C_i| \leq 2q$ by minimality of~$i$.
    Thus, 
    \[
        \big(\,V(C_1 \cup \dots \cup C_i) \cup X,\enspace V(G - (C_1 \cup \dots \cup C_i))\,\big)
    \]
    is a $(<k)$-separation of $G$ which is $(>q)$-broad since $|V(G)| \geq 3q+k$.
    This contradicts that $G$ is $(q,k-1)$-unbreakable.
\end{proof}

The next lemma will be instrumental in the proofs of \cref{keylem:irrelevant-star-dec}, \cref{thm:main-intro} and later \cref{mainGeneral}. Roughly speaking, it provides a way to ``lift'' separations of the torso of a star of $(\beta, k)$-irrelevant separations to separations of the entire graph.

\begin{definition}[Lift]
Let $\sigma$ be a star of separations of a graph~$G$.
Let $Z$ be the torso of~$\sigma$.
Let $(C, \cV)$ be a clique-faithful cactus-decomposition of $Z$. 
We define a cactus-decomposition $(\hat{C}, \hat{\cV})$ of $G$ as follows.
Let $s\in\sigma$ be arbitrary.
Since $\sigma$ is a star, we have $\sep(s) \subseteq V(Z)$.
Since $Z$ is the torso of $\sigma$, $\sep(s)$ is a clique in $Z$.
Thus, there exists a node $x \in V(C)$ that is not a $P$-node such that $\sep(s) \subseteq V_x$.
Choose an arbitrary such node $x(s)\in V(C)$.
We obtain $\hat C$ from $C$ by adding a new leaf $t(s)$ for each $s\in\sigma$, joining each leaf $t(s)$ to the node $x(s)$ we chose for~it.
Then we obtain $\hat \cV$ by extending $\cV$ to the set of new leaves $t(s)$ by letting $\hat V_{t(s)}:=A$ for $(A,B):=s$.
We declare that the $P$-nodes in $(\hat C,\hat\cV)$ are the same as in $(C, \cV)$.
Then $(\hat C,\hat\cV)$ is a cactus-decomposition of~$G$.

We say that $(\hat{C}, \hat{\cV})$ is a \defn{lift} of $(C, \mathcal{V})$ from $Z$ to $G$.
Let $S$ be the set of separations displayed by $(C, \cV)$ and let $\hat{S}$ be the set of separations displayed in $(\hat C,\hat\cV)$ by the cuts of $\hat C$ that do not separate any new leaf from its neighbour.
We say that $\hat{S}$ is the \defn{lift} of $S$ from $Z$ to $G$ corresponding to $(\hat{C}, \hat{\cV})$.

Given a cut $(U, W)$ of $C$, there is a unique cut $(\hat{U}, \hat{W})$ of $\hat{C}$ that does not separate any new leaf from its neighbour and such that $\hat{U} \cap V(C) = U$ and $\hat{W} \cap V(C) = W$. We call $(\hat{U}, \hat{W})$ the \defn{lift} of $(U, W)$ from $C$ to $\hat{C}$.
\end{definition}

This definition is justified by the following result.

\begin{lem}\label{lem:lifting}
    Let $\sigma$ be a star of separations of a graph $G$.
    Let $Z$ be the torso of $\sigma$.
    Let $(C, \cV)$ be a clique-faithful cactus-decomposition of $Z$.
    Let $(\hat{C}, \hat{\cV})$ be a lift of $(C, \cV)$ from $Z$ to $G$.
    Let $S$ be the set of separations displayed by $(C, \cV)$ and let $\hat{S}$ be the lift of $S$ from $Z$ to $G$ corresponding to $(\hat{C}, \hat{\cV})$. The following properties hold: \begin{enumerate}
        \item\label{item:induced-sep} 
        The set of separations displayed by $(\hat{C}, \hat{\cV})$ equals $\hat{S} \cup \sigma \cup \sigma^*$.
        \item\label{item:projection}\label{item:exists-lift}
        Let $(U, W)$ be a cut of $C$, let $(\hat{U}, \hat{W})$ be the lift of $(U, W)$ from $C$ to $\hat{C}$.
        Let $s$ be the separation induced by $(U, W)$ in $(C, \cV)$,
        and let $\hat{s}$ be the separation induced by $(\hat{U}, \hat{W})$ in $(\hat{C}, \hat{\cV})$. 
        Then $\hat s$ is a lift of~$s$.
        \item\label{item:same-torso} If~$C$ is a tree and $x \in V(C)$, then the torso at $x$ in $(\hat{C}, \hat{\cV})$ is equal to the torso at $x$ in $(C, \cV)$.
    \end{enumerate}
\end{lem}

\begin{proof}
    \cref{item:induced-sep}. By construction of $(\hat{C}, \hat{\cV})$, every separation in $\hat{S} \cup \sigma \cup \sigma^*$ is displayed by $(\hat{C}, \hat{\cV})$.
    Conversely, let $\hat{s}$ be a separation displayed by $(\hat{C}, \hat{\cV})$, and let $(\hat{U}, \hat{W})$ be a corresponding cut of $\hat{C}$.
    If $(\hat{U}, \hat{W})$ does not separate any new leaf from its neighbour then $\hat{s} \in \hat{S}$.
    Otherwise, since no cycle separates a leaf from its neighbour and since the neighbours of the new leaves are not $P$-nodes, $(\hat{U}, \hat{W})$ is the cut between a new leaf $t(s)$ for some $s \in \sigma$ and the rest of $\hat{C}$. In that case, we have $\hat{s} = s$ or $\hat{s} = s^*$.
    
    \cref{item:projection}. This is immediate since $\hat{C}$ is built from $C$ by adding leaves attached to nodes of $C$ but without changing the bags of the nodes in $C$, and since $(C, \cV)$ is a cactus-decomposition of $Z$.

    \cref{item:same-torso}. 
    Since~$C$ is a tree and~$\hat C$ is obtained by adding leaves, both decomposition graphs are trees.
    Let $\hat{W}$ be the torso at $x$ in $(\hat{C}, \hat{\cV})$ and $W$ be the torso at $x$ in $(C, \cV)$.
    Since $(\hat{C}, \hat{\cV})$ is a lift of $(C, \cV)$ from $Z$ to $G$, we have $V(\hat{W}) = \hat{V}_x = V_x = V(W)$, so it suffices to argue that $W$ and $\hat{W}$ have the same edges.
    
    First we show \casen{$E(W) \subseteq E(\hat W)$}.
    For this, let $uv \in E(W)$. Then, either $uv \in E(Z)$ or there exists a node $x' \neq x \in V(C)$ such that $\{u, v\} \subseteq V_{x'}$. In the second case, we also have $\{u, v\} \subseteq \hat{V}_{x'}$ so $uv \in E(\hat{W})$.
    In the first case, either $uv \in E(G)$ or there exists $(A, B) \in \sigma$ such that $\{u, v\} \subseteq A \cap B$. If $uv \in E(G)$, we have $uv \in E(\hat{W})$. If there exists $(A, B) \in \sigma$ such that $\{u, v\} \subseteq A \cap B$, then there exists a node $x' \neq x \in V(\hat{C})$ such that $\{u, v\} \subseteq \hat{V}_{x'}$ so $uv \in E(\hat{W})$.
    
    Now, we show \casen{$E(\hat W) \subseteq E(W)$}.
    For this, let $uv \in E(\hat{W})$. Then, either $uv \in E(G)$ or there exists a node $x' \neq x \in V(\hat{C})$ such that $\{u, v\} \subseteq \hat{V}_{x'}$. If $uv \in E(G)$ then $uv \in E(Z)$ so $uv \in E(W)$.
    Thus, we can assume that there exists a node $x' \neq x \in V(\hat{C})$ such that $\{u, v\} \subseteq \hat{V}_{x'}$.
    If $x'$ was already a node of $C$ then we have $\hat{V}_{x'} = V_{x'}$ so there exists a node $x' \neq x \in V(C)$ such that $\{u, v\} \subseteq V_{x'}$, and thus $uv \in E(W)$.
    Otherwise, there exists $(A, B) \in \sigma$ such that $\hat{V}_{x'} = A$, and thus $\{u, v\} \subseteq A \cap B$ so $uv \in E(Z)$ so $uv \in E(W)$.
    
    This finishes the proof that $W = \hat{W}$.
\end{proof}

We now prove \cref{keylem:irrelevant-star-dec}, which we restate for convenience.

\irrelevantstardec*

\begin{proof}[Proof of \cref{keylem:irrelevant-star-dec}]
    We have that $\sigma$ is a star by \cref{irrelevantSeparationsFormStar}.

    Let $Z$ be the torso of $\sigma$.
    By \cref{lem:lifting}~\ref{item:exists-lift}, every $(>q)$-broad separation of~$Z$ lifts to a $(>q)$-broad separation of~$G$ of the same order.
    Since~$G$ is $(q,k-1)$-unbreakable,~$Z$ is $(q,k-1)$-unbreakable too.
    It remains to show that~$Z$ is $(\beta,k)$-narrow.
    Otherwise, by \cref{lem:narrow-equivalent-no-irrelevant}, there would exist a $(\beta,k)$-irrelevant separation~$r' \eqqcolon (C, D)$ of~$Z$.
    That is, there is a $\beta$-broad $k$-separation~$r$ of~$Z$ such that~$r'$ is a componental $(<k)$-separation at~$r$.
    Consider a lift~$\hat r$ of~$r$ to~$G$ as in \cref{lem:lifting}~\ref{item:exists-lift}.
    Then~$\hat r$ is a $\beta$-broad $k$-separation of~$G$.
    Since $r'$ is $(\beta, k)$-irrelevant, $C \sm D$ is the vertex set of a connected component of $Z - \sep(r)$, and $N(C \sm D) = \sep(r') \subsetneq \sep(r) = \sep(\hat{r})$. Thus, all vertices in $C \sm D$ are in the same component $K$ in $G - \sep(\hat{r})$, and $N(K) = \sep(r') \subsetneq \sep(\hat{r})$.
    Now, consider the separation $\hat{s} \coloneqq (V(K) \cup \sep(r'), V(G - K))$.
    Therefore, $\hat{s}$ is a componental $(<k)$-separation at~$\hat{r}$, so $\hat{s} \in \sigma$. 
    Since $Z$ is the torso of $\sigma$, this implies $V(K) \cap V(Z) = \emptyset$, so $C \sm D = \emptyset$, contradicting that $r'$ is proper.
\end{proof}

\subsection{Lifting the cactus-decomposition}\label{subsec:the-proof}

Having established \cref{keylem:irrelevant-star-dec}, we now turn to the proof of \cref{thm:main-intro}.

First, we prove that broad separations of $G$ project to broad separations in the torso of a star $\sigma$ of irrelevant separations.
For all integers $n,k \in \N$, we set \defn{$\binom{n}{\leq k}$} $\coloneqq \sum_{i=0}^k \binom{n}{i}$, where we use the convention $\binom{n}{i}=0$ for $i>n$.

\begin{lem}\label{lem:projection-still-broad}
    Let $G$ be a $(q,k-1)$-unbreakable graph and let $\beta > \beta_{\ref{cor:irrelevant-separations-disjoint}}(q,k)$. Let $\alpha \in \Nbb$ and let $\alpha' \geq \alpha'_{\ref{lem:projection-still-broad}}(q, k, \alpha) \coloneqq q \cdot \binom{\alpha}{\leq k-1} + \alpha$.
    Let $\sigma$ be a star of $(\beta, k)$-irrelevant separations of $G$ and let $Z$ denote the torso of $\sigma$.
    Let $s \coloneqq (A, B)$ be a $\beta$-broad $k$-separation of $G$ with $|A| > \alpha'$.
    Then, $|A \cap V(Z)| > \alpha$.
\end{lem}

\begin{proof}
    First, since $s$ is $\beta$-broad and $\beta > \beta_{\ref{cor:irrelevant-separations-disjoint}}(q, k) = 4q+k$, we have $|V(G)| \geq 3q+k\ge \nu_{\ref{lem:few-vtces-attach-to-small-sep}}(q,k)$.

    By contradiction, suppose that $|A \cap V(Z)| \leq \alpha$.
    Let $v \in A \sm V(Z)$.
    There exists a unique $t \eqqcolon (C, D) \in \sigma$ such that $v \in C \sm D$.
    Since $(C, D)$ is $(\beta, k)$-irrelevant, we have $|\sep(t)| < k$ and $G[C \sm D]$ is a connected component of $G - \sep(t)$.
    By \cref{lem:irrelevant-small}, we have $|C \sm D| \leq q$.
    By \cref{cor:irrelevant-separations-disjoint}, $C \sm D$ is disjoint from $\sep(s)$, so $C \sm D \subseteq A \sm B$ since $v \in A$.

    Since $t$ is $(\beta, k)$-irrelevant, we have $\sep(t) = N(C \sm D)$, so $\sep(t) \subseteq A$, and $\sep(t) \subseteq V(Z)$ since $\sigma$ is a star.
    Then, $\sep(t)$ is a set of less than $k$ vertices in $A \cap V(Z)$ such that $C \sm D$ is a connected component of $G - \sep(t)$ with $\le q$ vertices.
    As $|A\cap V(Z)|\le\alpha$ by assumption,
    there are at most $\binom{\alpha}{\leq k-1}$ sets $X$ of less than $k$ vertices in $A \cap V(Z)$, and for each such set $X$, at most $q$ vertices belong to a connected component of $G-X$ of size at most $q$ by \cref{lem:few-vtces-attach-to-small-sep}.
    Thus, $|A \sm V(Z)| \leq q \cdot \binom{\alpha}{\leq k-1}$.
    Therefore, $|A| \leq q \cdot \binom{\alpha}{\leq k-1} + \alpha \leq \alpha'$, a contradiction.
\end{proof}

\begin{lem}
    \label{lem:lifts-are-still-nested}
    Let~$q,k,\beta,\gamma,\delta \in\N$ with $\beta > q$ and $\delta:=\delta_{\ref{lem:lifts-are-still-nested}}(q,k,\gamma):=\max\{\gamma + q \binom{k+\gamma}{\leq k-1}, \nu_{\ref{lem:few-vtces-attach-to-small-sep}}(q,k) - 1\}$.
    Let~$G$ be a $(q,k-1)$-unbreakable graph, let~$\sigma$ be a star of $(\beta,k)$-irrelevant separations of~$G$ and let $Z$ denote the torso of $\sigma$.
    Let~$s,t$ be two $(\leq k)$-separations of~$Z$, and let~$\hat s, \hat t$ be lifts of~$s,t$ to~$G$, respectively.
    If $s \leq_\gamma t$ then $\hat s \leq_\delta \hat t$.
\end{lem}
\begin{proof}
    If $|V(G)| < \nu_{\ref{lem:few-vtces-attach-to-small-sep}}(q,k)$ then $\hat{s} \leq_\delta \hat{t}$, so suppose $|V(G)|\ge \nu_{\ref{lem:few-vtces-attach-to-small-sep}}(q,k)$.
    This allows us to apply \cref{lem:few-vtces-attach-to-small-sep} later.

    Let $\hat{s} \eqqcolon (\hat A, \hat B)$ and $\hat{t} \eqqcolon (\hat C, \hat D)$ and note that $s=(\hat A \cap V(Z), \hat B \cap V(Z))$ and $t=(\hat C \cap V(Z), \hat D \cap V(Z))$ by definition of lift.
    By symmetry, it suffices to show that $|\hat A \sm \hat C| \leq \delta$.
    Let $v \in (\hat A \sm \hat C) \sm V(Z)$.
    There exists a unique $(E, F) \in \sigma$ such that $v \in E \sm F$.
    Since $(E \sm F) \cap V(Z) = \emptyset$ and $\hat A \cap \hat B \subseteq V(Z)$ and $\hat C \cap \hat D \subseteq V(Z)$ by definition of lift, it follows that $E\sm F$ is disjoint from both $\hat A \cap \hat B$ and $\hat C \cap \hat D$.
    Then, $v \in (\hat A \sm \hat B) \cap (\hat D \sm \hat C)$. 
    Since $G[E \sm F]$ is connected, this implies $E \sm F \subseteq (\hat A \sm \hat B) \cap (\hat D \sm \hat C)$.
    Then, $E \cap F = N(E \sm F) \subseteq \hat A \cap \hat D$.
    Furthermore, we have $E \cap F \subseteq V(Z)$, so $E \cap F$ is a set of less than~$k$ vertices in $\hat A \cap \hat D \cap V(Z)$.
    However, $|\hat A \cap \hat D \cap V(Z)| \leq |\hat A \cap \hat B \cap V(Z)| + |(\hat D \sm \hat B) \cap V(Z)| \leq k + \gamma$ since $s \leq_\gamma t$ implies $|(\hat D \sm \hat B)\cap V(Z)| \leq \gamma$.
    Furthermore, $G[E \sm F]$ is a connected component of $G - (E \cap F)$ of size at most~$q$ by \cref{lem:irrelevant-small}.
    Thus, $v$ is one of the at most $q$ vertices that belong to a connected component of $G - (E \cap F)$ of size at most $q$ by \cref{lem:few-vtces-attach-to-small-sep}.
    This yields $|(\hat A \sm \hat C) \sm V(Z)| \leq q \binom{k+\gamma}{\leq k-1}$.
    Thus, $|\hat A \sm \hat C| = |(\hat A \sm \hat C)\cap V(Z)| + |(\hat A \sm \hat C)\sm V(Z)| \leq \gamma + q \binom{k+\gamma}{\leq k-1} \leq \delta$.
\end{proof}

This implies that two separations of $G$ which project to two close separations in the torso of~$\sigma$ must be close in $G$. 

\begin{lem}\label{lem:projections-close-then-close}
    Let $G$ be a $(q,k-1)$-unbreakable graph and $\beta > q$.
    Let $\gamma \in \Nbb$ and let $\delta \geq \delta_{\ref{lem:projections-close-then-close}}(q, k, \gamma) \coloneqq \delta_{\ref{lem:lifts-are-still-nested}}(q,k,\gamma)$.
    Let $\sigma$ be a star of $(\beta, k)$-irrelevant separations of $G$ and let $Z$ denote the torso of $\sigma$.
    Let $s, t$ be two $\gamma$-close $(\leq k)$-separations of $Z$ and let $\hat{s}$ (resp.\ $\hat{t}$) be a lift of $s$ (resp.\ $t$).
    Then, $\hat{s}$ and $\hat{t}$ are $\delta$-close.
\end{lem}
\begin{proof}
    Follows from \cref{lem:lifts-are-still-nested}.
\end{proof}

We now have all the tools to prove the main result of this paper, which we restate for convenience.

\cactusgeneral*

\begin{proof}
    Given $q, k \in \N$, set \begin{align*}
        \beta &\coloneqq \beta_{\ref{cor:irrelevant-separations-disjoint}}(q,k)+1, \\
        \eta &\coloneqq \eta_{\ref{thm:cactusForNarrow}}(q, k, \beta), \\
        \varepsilon &\coloneqq \max\{\alpha'_{\ref{lem:projection-still-broad}}(q, k, \eta+k)+1, \beta\}, \text{ and } \\
        d &\coloneqq \max\{\varepsilon+k, 2\delta_{\ref{lem:projections-close-then-close}}(q, k, \eta)\}.
    \end{align*} 
    Let $G$ be a $(q,k-1)$-unbreakable graph.
    Let $\sigma$ be the set of all $(\beta, k)$-irrelevant separations.
    By \cref{keylem:irrelevant-star-dec}, $\sigma$ forms a star and its torso $Z$ is $(q,k-1)$-unbreakable and $(\beta, k)$-narrow.
    By \cref{thm:cactusForNarrow}, $Z$ has a clique-faithful cactus-decomposition $\check{\cC} = (\check{C}, \check{\cV})$ such that every displayed separation has order~$\leq k$ and the set of displayed separations of $\check{\cC}$ $\eta$-approximates the set of all $k$-separations of~$Z$.
    Let $\check{S}$ be the set of separations displayed by $(\check{C}, \check{\cV})$.
    Let $(C, \cV)$ be a lift of $(\check{C}, \check{\cV})$ from $Z$ to $G$, and let $S$ be the lift of $\check{S}$ from $Z$ to $G$ corresponding to $(C, \cV)$.
    Then, $(C, \cV)$ is a cactus-decomposition of $G$ such that every displayed separation has order $\leq k$.

    It remains to prove that every $(\leq k)$-separation of~$G$ has edit-distance at most~$d$ to a separation in~$S$.
    Let $s$ be a $(\leq k)$-separation of $G$.
    If $s$ is not $\varepsilon$-broad then $s$ has edit-distance at most~$d$ to $(V(G),\emptyset)$ or $(\emptyset,V(G))$, both of which belong to~$S$.
    We may therefore assume that $s$ is $\varepsilon$-broad.
    Since $\varepsilon\ge\beta>q$ and~$G$ is $(q,k-1)$-unbreakable, $s$ has order $k$.
    It follows from \cref{cor:irrelevant-separations-disjoint} that $\sep(s) \subseteq V(Z)$, so $s$ induces a $k$-separation $\check s$ on~$Z$ (i.e.\ $s$ is a lift of $\check s$ from~$Z$ to~$G$).
    By \cref{lem:projection-still-broad}, since $\beta > \beta_{\ref{cor:irrelevant-separations-disjoint}}(q,k)$ and $\varepsilon \geq \max\{\alpha'_{\ref{lem:projection-still-broad}}(q, k, \eta+k)+1, \beta\}$, it follows that $\check{s}$ is an $\eta$-broad $k$-separation of $Z$.
    By \cref{thm:cactusForNarrow}, $\check{s}$ is $\eta$-close to a displayed separation $\check{s}'$ of $(\check{C}, \check{\cV})$. 
    By \cref{lem:lifting}~\ref{item:exists-lift}, there exists a lift $s' \in S$ of $\check{s}'$.
    By \cref{lem:projections-close-then-close}, since $\beta = \beta_{\ref{cor:irrelevant-separations-disjoint}}(q,k)+1 > q$, $s$ and $s'$ are $\delta_{\ref{lem:projections-close-then-close}}(q,k,\eta)$-close.
    By \cref{closeVSapprox}, their edit-distance is therefore at most $2\delta_{\ref{lem:projections-close-then-close}}(q,k,\eta)\le d$.
\end{proof}

\section{An approximate Tutte-type decomposition}\label{sec:TutteNarrow}

Finally, we extend \cref{mainForNarrow} to all $(q, k-1)$-unbreakable graphs by dropping the assumption that $G$ is $(b, k)$-narrow.
This establishes \cref{mainTutte}.

\subsection{Lifting flowers}

We start by proving some results to lift graphic-$k$-flowers.

\begin{lem}\label{lem:subpath-torso}
    Let $G$ be a graph, let $\sigma$ be a star of separations of $G$ and let $Z$ denote the torso of $\sigma$.
    Let $F = (O, \cV)$ be a graphic-$k$-daisy of $Z$ that has the clover-property.
    Let $X \coloneqq \bigcap_{o \in V(O)} V_o$ be the centre of $F$.
    Let $s = (C, D) \in \sigma$ be such that $\sep(s) \not \subseteq X$. 
    Then, the set of vertices $x \in V(O)$ such that $\sep(s) \subseteq V_x$ induces a nonempty subpath of $O$.
\end{lem}

\begin{proof}
    Since $\sigma$ is a star, we have $\sep(s) \subseteq V(Z)$.
    Since $Z$ is the torso of $\sigma$, we have that $Z[\sep(s)]$ is a complete graph.
    The lemma then follows from \cref{lem:subpath}.
\end{proof}

The next lemma (\cref{lem:lifting-flowers}) is very close to \cref{lem:lifting}. Essentially, it says that graphic-$k$-flowers of the torso of $\sigma$ can be lifted to graphic-$k$-flowers of $G$.
When lifting a graphic-$k$-daisy, each component outside the torso whose neighbourhood lies in the centre can be attached to any bag. Taking all possible assignments yields a canonical family of lifts whose members differ only in the placement of these small components.

\begin{definition}[$(q, k)$-coherent]\label{defnQKcoherent}
    Let $q, k\in\N$.
    Let $G$ be a $(q,k-1)$-unbreakable graph.
    Let $\mathcal{F}$ be a set of graphic-$k$-flowers of $G$, all with the same centre $X$.
    Let $R$ be the set of all vertices of $G$ that belong to some connected component $K$ of $G-X$ of size $|V(K)|\le q$ such that $|N_G(K)| < k$.
    For every $F \in \mathcal{F}$, let $\check{F}$ be the graph-decomposition of $G-R$ obtained from $F$ by removing the vertices of $R$ from the bags that contain them.
    We say that $\mathcal{F}$ is \defn{$(q, k)$-coherent} if $\check{F} = \check{F'}$ for all $F, F' \in \mathcal{F}$.
\end{definition}

\begin{remark}
    Assume the context of \cref{defnQKcoherent} with $|V(G)| \geq \nu_{\ref{lem:few-vtces-attach-to-small-sep}}(q,k)$.
    In \cref{defnQKcoherent}, we have that \cref{lem:few-vtces-attach-to-small-sep} implies $|R| \leq qk$.
    Thus, one should think of a $(q, k)$-coherent set $\mathcal{F}$ of graphic-$k$-flowers as a way to represent a graphic-$k$-flower which is essentially canonical, up to a bounded number of vertices, which additionally have a tame behaviour.
\end{remark}

\begin{lem}\label{lem:lifting-flowers}
    Let $G$ be a $(q,k-1)$-unbreakable graph and let $\beta > \beta_{\ref{cor:irrelevant-separations-disjoint}}(q,k)$.
    Let $\sigma$ denote the star of all $(\beta, k)$-irrelevant separations of $G$ (see \cref{keylem:irrelevant-star-dec}).
    Let $Z$ denote the torso of $\sigma$.
    Let $F$ be a graphic-$k$-flower of~$Z$, such that $F$ has the clover-property if it is a graphic-$k$-daisy.
    Then, there exists a nonempty $(q, k)$-coherent set $\hat{\mathcal{F}}$ of graphic-$k$-flowers of $G$ with the following properties: \begin{enumerate}
        \item\label{item:flower-projection}
        For every $\hat{F} \in \hat{\cF}$, every $\hat s\in D(\hat F)$ induces a separation $s$ of $Z$. If $s$ is proper then $s \in D(F)$.
        Moreover, if $\hat{s}$ is a petal-separation of $\hat{F}$ and $s$ is proper, then
        $s$ is a petal-separation of $F$.
        \item\label{item:flower-exists-lift} For every $\hat{F} \in \hat{\cF}$, every separation $s \in D(F)$ has a lift $\hat s\in D(\hat{F})$.
        Moreover, if $s$ is a petal-separation of~$F$, then we can choose $\hat{s}$ as a petal-separation of $\hat{F}$.
        \item\label{item:flower-lift-canonical} $\hat{\cF} = \hat{\cF}(F)$ is uniquely determined by $F$ (given~$G$).
        Furthermore, if $\psi$ is an automorphism of~$G$ and $F'$ is a graphic-$k$-flower of $Z$ such that $F$ and $F'$ are $\psi\rest Z$-isomorphic, then $\hat{\cF}(F)$ and $\hat{\cF}(F')$ are $\psi$-isomorphic.
    \end{enumerate}
\end{lem}

\begin{proof}
    \casen{Suppose first that $F = (S, \cV)$ is a graphic-$k$-anemone.}
    Let $c$ be the centre of $S$.
    We define a graphic-$k$-anemone $\hat{F} = (\hat{S}, \hat{\cV})$ of $G$ as follows.    
    Start with $(\hat{S}, \hat{\cV}) = (S, \cV)$ and let $s \coloneqq (C, D) \in \sigma$.
    Since $\sigma$ is a star, we have $\sep(s) \subseteq V(Z)$.
    Since $Z$ is the torso of $\sigma$, $Z[\sep(s)]$ is a complete graph.
    If $\sep(s) \subseteq V_c$, we add a new leaf $t(s)$ to $\hat{S}$ with $c$ as its only neighbour, and set $V_{t(s)} = C \cup V_c$.
    Otherwise, there exists a unique leaf $\ell(s) \in V(S)$ such that $\sep(s) \subseteq V_{\ell(s)}$, and we add $C \sm D$ to $\hat{V}_{\ell(s)}$.
    After iterating this process over all separations $s \in \sigma$, we eventually obtain a graphic-$k$-anemone $\hat{F} = (\hat{S}, \hat{\cV})$ of $G$.
    We set $\hat{\cF} \coloneqq \{\hat{F}\}$.
    The set~$\hat{\mathcal{F}}$ is trivially a nonempty $(q, k)$-coherent set of graphic-$k$-flowers of~$G$.

    \cref{item:flower-projection}. 
    Let $\hat{s} = (\hat{A}, \hat{B}) \in D(\hat{F})$.
    There is a cut $(\hat U,\hat W)$ of $\hat S$ such that $\hat A=\hat \cV\rest \hat U$ and $\hat B=\hat \cV\rest \hat W$.
    If $\hat U=\emptyset$ or $\hat W=\emptyset$, then $\hat s$ induces $(\emptyset,V(Z))$ or $(V(Z),\emptyset)$, respectively, so the assertion holds.
    Hence assume that $\hat U$ and $\hat W$ are nonempty.
    Let $U:=\hat U\cap V(S)$ and $W:=\hat W\cap V(S)$.
    Let $A:=\cV\rest U$, $B:=\cV\rest W$ and $s:=(A,B)$.
    We distinguish two cases.

    In the first case, both $U$ and $W$ are nonempty.
    Then $\hat s$ induces $(A,B)$, which is displayed by $F$.
    Moreover, if $\hat s$ is a petal-separation, then $\hat U$ consists of a single node, and so does $U$, hence $s$ is a petal-separation as well.

    In the second case, $U$ or $W$ is empty, say $U=\emptyset$.
    Then $\hat U=\{\,t(r):r\in\sigma'\,\}$ for some subset $\sigma'\se\sigma$.
    In this case, $\hat s$ induces $(V_c,V(Z))$, which is a petal-separation of~$F$.

    \cref{item:flower-exists-lift}. The proof is straightforward, as in \cref{lem:lifting}~\ref{item:exists-lift}.

    \cref{item:flower-lift-canonical}. 
    This is immediate since the construction of $\hat{F}$ involves no arbitrary choices.

    \medskip

    \casen{Suppose now that $F = (O, \cV)$ is a graphic-$k$-daisy.}
    Let $O_1$ be the 1-subdivision of the cycle~$O$. 
    For every edge $e \in E(O)$, denote by $d(e)$ the corresponding vertex of $O_1$.
    We extend $\cV$ to all vertices $d(e) \in V(O_1) \sm V(O)$ by letting $V_{d(e)} := V_e$, where $V_e$ is the adhesion-set at~$e$.
    Then, $(O_1, \cV)$ is a graphic-$k$-daisy of $Z$.
    Let $X \coloneqq \bigcap_{x \in V(O)}V_x$, and note that $X = \bigcap_{x \in V(O_1)}V_x$, so the centre remains unchanged.
    Since $F$ has the clover-property, we have $|X| < k$.
    
    We now describe an algorithm to build a graphic-$k$-daisy $\hat{F} = (O_1, \hat{\cV})$ of $G$.
    Start with $(O_1, \hat{\cV}) = (O_1, \cV)$ and let $s \coloneqq (C, D) \in \sigma$.
    Since $\sigma$ is a star, we have $\sep(s) \subseteq V(Z)$.
    Since $Z$ is the torso of $\sigma$, $Z[\sep(s)]$ is a complete graph.
    If $\sep(s) \subseteq X$, let $x \in V(O_1)$ be an arbitrary vertex, and add $C \sm D$ to $\hat{V}_x$.
    Otherwise, by \cref{lem:subpath-torso}, the set of vertices $x \in V(O)$ such that $\sep(s) \subseteq V_x$ induces a subpath of $O$. Thus, the set of vertices $x \in V(O_1)$ such that $\sep(s) \subseteq V_x$ induces a path in $O_1$ with an odd number of vertices.
    Let $y$ be the vertex at the middle of that path, and add $C \sm D$ to $\hat{V}_y$.
    After iterating this process over all separations $s \in \sigma$, we eventually obtain a cycle-decomposition $\hat{F} = (O_1, \hat{\cV})$ of $G$.
    The construction of $\hat{F}$ involves some arbitrary choices. We then define $\hat{\cF}$ as the (nonempty) set of all graphic-$k$-flowers $\hat{F}$ of $G$ that can be built by this algorithm. Note that every $\hat{F} \in \hat{\cF}$ has centre $X$. 
    Let $R$ be the set of all vertices of $G$ that belong to some connected component $K$ of $G-X$ of size $|V(K)|\le q$ such that $|N_G(K)| < k$. If $s \coloneqq (C, D) \in \sigma$ is such that $\sep(s) \subseteq X$, then $G[C \sm D]$ is a connected component of $G-X$ of size at most $q$ by \cref{lem:irrelevant-small}, and $|N_G(C \sm D)| \leq |X| < k$, so $C \sm D \subseteq R$. Thus, we only make arbitrary choices regarding vertices in $R$. This shows that $\hat{\cF}$ is a nonempty $(q, k)$-coherent set of graphic-$k$-flowers of $G$.

    The proofs of \cref{item:flower-projection} and \cref{item:flower-exists-lift} are similar to the proofs for the case where $F$ is a graphic-$k$-anemone, so we do not reproduce them.

    \cref{item:flower-lift-canonical}. Since $\hat{\cF}$ contains all the graphic-$k$-flowers that can be built by the above algorithm, $\hat{\cF} = \hat{\cF}(F)$ is uniquely determined by $F$ (given~$G$).
    Let $\psi$ be an automorphism of~$G$ and $F' \eqqcolon (O', \cV')$ be a graphic-$k$-flower of $Z$ such that $F$ and $F'$ are $\psi\rest Z$-isomorphic.
    Let $\Phi\colon V(O) \to V(O')$ be an isomorphism such that $V'_{\Phi(o)} = \psi(V_o)$ for every $o \in V(O)$. 
    Then, $\Phi$ naturally extends to $\Phi\colon V(O_1) \to V(O'_1)$.
    Note that the centre $X'$ of $F'$ satisfies $X' = \psi(X)$.
    Let $\hat{F} \eqqcolon (O_1, \hat{\cV}) \in \hat{\cF}(F)$.
    We show that the graphic-$k$-flower $\hat{F'} = (O'_1, \hat{\cV'})$ of $\psi(G)$ with $\hat{V}'_{\Phi(o)} = \psi(\hat{V}_o)$ for every $o \in V(O_1)$ belongs to $\hat{\cF}(F')$.
    Let $s \coloneqq (C, D) \in \sigma$. 
    Since $\sigma$ is canonical, we have that $s' \coloneqq (\psi(C), \psi(D))$ is a $(\beta, k)$-irrelevant separation of~$\psi(G)$.
    
    If $\sep(s) \subseteq X$, then $C \sm D$ was added to an arbitrary bag $\hat{V}_o$ by the algorithm. In that case, $\sep(s') = \psi(\sep(s)) \subseteq \psi(X) = X'$, so the algorithm can add $\psi(C) \sm \psi(D)$ to the bag $\hat{V}'_{\Phi(o)}$.

    If $\sep(s) \not\subseteq X$, then the set $P_s$ of vertices $x \in V(O_1)$ such that $\sep(s) \subseteq V_x$ induces a path in $O_1$ with an odd number of vertices, with middle vertex $y$ say. 
    Then, we also have $\sep(s') = \psi(\sep(s)) \not\subseteq \psi(X) = X'$, and the set $P'_{s'}$ of vertices $x' \in V(O'_1)$ such that $\sep(s') \subseteq V'_{x'}$ induces a path in $O'_1$ with an odd number of vertices, with middle vertex $y'$ say.
    However, $\Phi\colon V(O_1) \to V(O'_1)$ satisfies $V'_{\Phi(o)} = \psi(V_o)$ for every $o \in V(O)$, so $P'_{s'} = \Phi(P_s)$, and thus $y' = \Phi(y)$.
    In that case, the algorithm adds $C \sm D$ to $\hat{V}_y$ to build $\hat{F}$ from $F$, and adds $\psi(C) \sm \psi(D)$ to $\hat{V}'_{y'}$ to build $\hat{F'}$ from $F'$.
    Thus, the algorithm for $F'$ can build the graphic-$k$-flower $\hat{F'} \coloneqq (O'_1, \hat{\cV'})$ such that $\hat{V}'_{\Phi(o)} = \psi(\hat{V}_o)$ for every $o \in V(O_1)$.
    Then, $\hat{F} \in \hat{\cF}(F)$ is $\psi$-isomorphic to $\hat{F'} \in \hat{\cF}(F')$.
    It then follows immediately that $\hat{\cF}(F)$ and $\hat{\cF}(F')$ are $\psi$-isomorphic.
\end{proof}

We refer to the set $\hat\cF$ in \cref{lem:lifting-flowers} as the \defn{set of lifts} of~$F$.
This is justified by \cref{lem:lifting-flowers}~\cref{item:flower-lift-canonical}.

\subsection{Lifting the approximate Tutte-type decomposition}

We now have all the tools to prove the Tutte-type version of our results. We need one last definition to explain in what way the flowers in \cref{mainTutte} are canonical.
    
\begin{definition}[$(q,k)$-canonical]\label{qkcanonical}
    Let $q, k\in\N$.
    Let $G$ be a $(q,k-1)$-unbreakable graph, and let $\cG$ be a collection of sets of graphic-$k$-flowers of~$G$.
    We say that $\cG$ is \defn{$(q, k)$-canonical} if the following two conditions are satisfied:
    \begin{itemize}
        \item Every set $\cF\in\cG$ of graphic-$k$-flowers is $(q,k)$-coherent.
        \item For every automorphism $\psi$ of~$G$, there is a permutation $\Phi$ of~$\cG$ such that $\cF$ and $\Phi(\cF)$ are $\psi$-isomorphic for every $\cF \in \cG$.
    \end{itemize}
\end{definition}
    
The following is a strengthening of \cref{mainTutte}.

\begin{theorem}\label{mainGeneral}
    For every $q,k\in\N$ there are $\gamma,\kappa,\rho,q',\eps\in\N$ such that the following holds.
    Let~$G$ be a $(q,k-1)$-unbreakable graph.
    There exists a collection $\cG$ of nonempty sets of graphic-$k$-flowers of $G$ and a symmetric nested set $M$ of totally-$\gamma$-nested proper $(\leq k)$-separations of~$G$ with corresponding tree-decomposition $(T,\cV)$ of~$G$ and a bipartition $X\sqcup Y= V(T)$ such that:
    \begin{enumerate}
        \item\label{mainOutroTorsoQKC} Every node $t\in X$ has a $(q',k)$-unbreakable torso.
        \item\label{mainOutroTorsoFlower} Every node $t\in Y$ comes with a set $\mathcal{F}_t \in \mathcal{G}$ of graphic-$k$-flowers of $G$ such that every $F \in \mathcal{F}_t$ $\rho$-refines $(T,\cV)$ at~$t$.
        \item\label{mainOutroSeparationAnalysis} For every $\eps$-broad $k$-separation $s$ of $G$, either $s$ is $\eps$-close to a separation in~$M$ or there is $t \in Y$ such that for every $F \in \mathcal{F}_t$, $s$ is $\eps$-close to some displayed separation of~$F$.
        \item\label{mainOutroCanonical} 
        The set $\mathcal{G}$ is $(q, k)$-canonical.
        Both the tree-decomposition $(T,\cV)$ and the set $M$ are $\kappa$-canonical. 
        Given~$M$ and~$(T,\cV)$, the bipartition $X\sqcup Y= V(T)$ is canonical and the assignment $Y\ni t \mapsto \mathcal{F}_t \in \mathcal{G}$ is a canonical bijection.
    \end{enumerate}
\end{theorem}

\begin{proof}
    Let $b:=\beta_{\ref{cor:irrelevant-separations-disjoint}}(q,k)+1$.
    Let $\gamma_0,\kappa_0,\rho_0,q'_0,\eps_0\in\N$ be the constants obtained by applying \cref{mainForNarrow} for $q, k, b$.
    Let
    \begin{align*}
        \gamma_1 &\coloneqq \max\{\,\gamma_0,\enspace b\,\}, \\
        \eps_1 &\coloneqq \max\{\,q+1,\enspace \eps_0\,\}, \\
        \gamma &\coloneqq \max\{\,b+k,\enspace \alpha'_{\ref{lem:projection-still-broad}}(q, k, \gamma_1+k)+k\,\}, \\
        \kappa &\coloneqq \delta_{\ref{lem:projections-close-then-close}}(q, k, \kappa_0+k)+1, \\
        \rho &\coloneqq \max\{\,b,\, \alpha'_{\ref{lem:projection-still-broad}}(q, k, \rho_0+k)+1,\enspace \delta_{\ref{lem:projections-close-then-close}}(q, k, \rho_0)\,\}, \\
        q' &\coloneqq \max\{\,\nu_{\ref{lem:few-vtces-attach-to-small-sep}}(q,k),\enspace q'_0\,\}, \text{ and } \\
        \eps &\coloneqq \max\{\,\nu_{\ref{lem:few-vtces-attach-to-small-sep}}(q,k),\enspace \alpha'_{\ref{lem:projection-still-broad}}(q, k, \varepsilon_1+k)+1,\enspace \delta_{\ref{lem:projections-close-then-close}}(q, k, \eps_0)\,\}.
    \end{align*}

    Let $G$ be a $(q,k-1)$-unbreakable graph.
    If $|V(G)| < \nu_{\ref{lem:few-vtces-attach-to-small-sep}}(q,k)$, then we set $\mathcal{G} = \emptyset$ and $M = \emptyset$. The corresponding tree-decomposition consists of a single node $t$, with $V_t = V(G)$. We declare that $t \in X$. It is then straightforward to verify that the desired properties hold.

    Suppose now that $|V(G)| \geq \nu_{\ref{lem:few-vtces-attach-to-small-sep}}(q,k)$.
    Let~$\sigma$ be the set of all $(b,k)$-irrelevant separations of~$G$.
    By \cref{keylem:irrelevant-star-dec},~$\sigma$ is a star such that the torso~$Z$ of~$\sigma$ is $(q,k-1)$-unbreakable and $(b,k)$-narrow.
    Thus,~$Z$ satisfies the premises for \cref{mainForNarrow}.
    Hence, there exist a set $\check{\cF}$ of graphic-$k$-flowers of~$Z$, a symmetric nested set~$\check{N}$ of narrow-$\gamma_0$-nested proper $k$-separations of~$Z$ with corresponding tree-decomposition~$(\check{T},\check{\cV})$ of~$Z$, a bipartition $\check{X} \sqcup \check{Y} = V(\check{T})$, and an assignment $\check{Y}\ni \check{t}\mapsto \check{F}_{\check{t}}\in\check{\cF}$, which satisfy the outcome of \cref{mainForNarrow}.

    For every $\check{F}\in\check{\cF}$, let $\cF(\check{F})$ be the nonempty set of lifts of $\check{F}$ given by \cref{lem:lifting-flowers} (this set is uniquely determined by $\check{F}$ by \cref{lem:lifting-flowers}~\ref{item:flower-lift-canonical}).
    Let $\cG:=\{\cF(\check{F}):\check{F}\in\check{\cF}\}$.

    Let $(T,\cV)$ be a lift of $(\check T,\check\cV)$ from $Z$ to $G$. 
    Let $N$ be the set of separations induced in $(T,\cV)$ by the oriented edges of $\check T\subseteq T$, and put $M:=N\cup\sigma\cup\sigma^*$.

    First, we show that the separations in $M$ are totally-$\gamma$-nested.
    By contradiction, suppose that some separation $s \coloneqq (A, B) \in M$  is $\gamma$-crossed by a proper $(\leq k)$-separation $s' \coloneqq (C, D)$ of $G$.
    Then, $s$ is $2(\gamma-k)$-broad (by \cref{LargeCrossedImpliesBroad}) with $2(\gamma-k) \geq b > \beta_{\ref{cor:irrelevant-separations-disjoint}}(q,k) \geq q$, so $s \in N$ by \cref{lem:irrelevant-small}.
    We have that $s$ and $s'$ $(\gamma-k)$-criss by \cref{crissCross}.
    Thus, the $AC$-, \mbox{$AD$-,} $BC$- and $BD$-corner separations of $G$ are $(> \gamma - k)$-broad.
    Since $G$ is $(q,k-1)$-unbreakable and $\gamma - k \geq b > q$, they all have order at least $k$.
    By submodularity (\cref{lem:submodularity}), they all have order exactly~$k$.
    By applying \cref{lem:projection-still-broad} to each of these four corner separations, while using $\gamma - k \geq b > \beta_{\ref{cor:irrelevant-separations-disjoint}}(q,k)$ and $\gamma - k \geq \alpha'_{\ref{lem:projection-still-broad}}(q, k, \gamma_1+k)$, we get 
    \[
        |A \cap C \cap V(Z)|,\enspace |A \cap D \cap V(Z)|,\enspace |B \cap C \cap V(Z)|,\enspace |B \cap D \cap V(Z)|\enspace > \gamma_1+k.
    \]
    Therefore, the induced separations 
    \[
        \check{s} \coloneqq (A \cap V(Z), B \cap V(Z))\qquad \text{and} \qquad \check{s}' \coloneqq (C \cap V(Z), D \cap V(Z))
    \]
    $\gamma_1$-criss, so they $\gamma_1$-cross by \cref{crissCross}.
    Thus, the $(\leq k)$-separation $\check{s}'$ of $Z$ is $\gamma_1$-broad, hence $b$-broad, and thus narrow since $Z$ is $(b, k)$-narrow.
    Therefore, $\check{s}'$ is a narrow proper $(\leq k)$-separation of $Z$ which $\gamma_0$-crosses $\check{s}$. 
    However, by \cref{lem:lifting}~\ref{item:induced-sep} and \cref{lem:lifting}~\ref{item:projection}, we have $\check{s} \in \check{N}$, so $\check{s}$ is narrow-$\gamma_0$-nested, which is a contradiction.

    Since $(T, \cV)$ is a lift of $(\check{T}, \check{\cV})$ from $Z$ to $G$, $T$ can be obtained from $\check{T}$ by adding leaves attached to nodes of~$\check{T}$. 
    We extend the bipartition $(\check{X}, \check{Y})$ of $V(\check{T})$ to a bipartition $(X, Y)$ of $V(T)$ by declaring that all the newly added leaves are in $X$. 

    \medskip

    \cref{mainOutroTorsoQKC}. Let $t \in X$. 
    Suppose first that $t$ is one of the leaves that were added to go from $\check{T}$ to $T$. 
    Since $(T, \cV)$ is a lift of $(\check{T}, \check{\cV})$ from $Z$ to $G$, there exists $(C, D) \in \sigma$ such that $V_t = C$.
    Thus, $|V_t| = |C| \leq q+k$ by \cref{lem:irrelevant-small}, so the torso at $t$ is indeed $(q',k)$-unbreakable since $q' \geq \nu_{\ref{lem:few-vtces-attach-to-small-sep}}(q,k) \geq q+k$.
    Otherwise, $t$ was an original node of $\check{T}$ and $t \in \check{X}$.
    By \cref{lem:lifting}~\ref{item:same-torso}, the torso at $t$ in $(T, \cV)$ is equal to the torso at $t$ in $(\check{T}, \check{\cV})$. 
    We are then done since the latter is $(q',k)$-unbreakable by \cref{mainForNarrow}~\ref{mainForNarrowTorsoQKC}.

    \medskip

    \cref{mainOutroTorsoFlower}. Let $t \in Y = \check{Y}$ and let $\check{F}_t \in \check{\cF}$ be the corresponding graphic-$k$-flower of $Z$.
    Let $\cF_t$ be the set of lifts of $\check{F}_t$ given by \cref{lem:lifting-flowers} (this set is unique by \cref{lem:lifting-flowers}~\ref{item:flower-lift-canonical}).
    By definition of $\cG$, we have $\cF_t \in \cG$.
    Let $F \in \cF_t$ be arbitrary. 
    We show that $F$ $\rho$-refines $(T, \cV)$ at~$t$.
    
    First, let $s$ be a $\rho$-broad petal-separation of $F$. 
    By \cref{lem:projection-still-broad}, since $s$ is $\rho$-broad with $\rho \geq b> \beta_{\ref{cor:irrelevant-separations-disjoint}}(q,k)$ and $\rho > \alpha'_{\ref{lem:projection-still-broad}}(q, k, \rho_0+k)$, the separation $\check{s}$ that $s$ induces on $Z$ is $\rho_0$-broad.
    By \cref{lem:lifting-flowers}~\ref{item:flower-projection}, $\check{s}$ is a petal-separation of $\check{F}_t$.    
    By \cref{mainForNarrow}~\ref{mainForNarrowTorsoFlower}, there exists a separation $\check{s}'$ of the $t$-star of~$(\check{T},\check{\cV})$ which is $\rho_0$-close to $\check{s}$.
    By \cref{lem:lifting}~\ref{item:exists-lift}, there is a lift $s'$ of $\check{s}'$ in the $t$-star of $(T, \cV)$.
    By \cref{lem:projections-close-then-close}, and since $\rho \geq \delta_{\ref{lem:projections-close-then-close}}(q, k, \rho_0)$, we have that 
    $s'$ is $\rho$-close to $s$.

    Conversely, let $s'$ be a $\rho$-broad separation of the $t$-star of $(T, \cV)$.
    By \cref{lem:irrelevant-small} and since $\rho \geq b > q$, it follows that $s' \notin \sigma \cup \sigma^*$.
    By \cref{lem:lifting}~\ref{item:induced-sep}, we have $s' \in N$.
    Then, by \cref{lem:lifting}~\ref{item:projection}, the separation $\check{s}'$ that $s'$ induces on~$Z$ belongs to the $t$-star of $(\check{T}, \check{\cV})$.
    By \cref{lem:projection-still-broad}, since $s'$ is $\rho$-broad with $\rho \geq b> \beta_{\ref{cor:irrelevant-separations-disjoint}}(q,k)$ and $\rho \geq \alpha'_{\ref{lem:projection-still-broad}}(q, k, \rho_0+k)$, 
    we get that $\check{s}'$ is $\rho_0$-broad.
    By \cref{mainForNarrow}~\ref{mainForNarrowTorsoFlower}, there exists a petal-separation $\check{s}$ of $\check{F}_t$ which is $\rho_0$-close to $\check{s}'$.
    By \cref{lem:lifting-flowers}~\ref{item:flower-exists-lift}, there is a petal-separation $s$ of $F$ that is a lift of~$\check{s}$.
    By \cref{lem:projections-close-then-close}, and since $\rho \geq \delta_{\ref{lem:projections-close-then-close}}(q, k, \rho_0)$, we have that $s$ is $\rho$-close to~$s'$.

    \medskip

    \cref{mainOutroSeparationAnalysis}. Let $s \coloneqq (A, B)$ be an $\eps$-broad $k$-separation of $G$.
    Then, 
    \[
        |A|, |B| > \alpha'_{\ref{lem:projection-still-broad}}(q, k, \varepsilon_1+k).
    \]
    Let $\check{s}$ be the separation that $s$ induces on~$Z$.
    By \cref{lem:projection-still-broad}, $\check{s}$ is an $\eps_1$-broad $(\leq k)$-separation of~$Z$. 
    Since $\check{s}$ is $\varepsilon_1$-broad with $\varepsilon_1 > q$, and $Z$ is $(q,k-1)$-unbreakable, $\check{s}$ is a $k$-separation of~$Z$.
    Since $\eps_1 \geq \eps_0$, by \cref{mainForNarrow}~\ref{mainForNarrowSeparationAnalysis}, $\check{s}$ is either $\eps_0$-close to a separation in $\check{N}$ or $\varepsilon_0$-close to some displayed separation of some $\check{F}_t$. 
    Call this separation $\check{s}'$. 
    By \cref{lem:lifting}~\ref{item:exists-lift} or \cref{lem:lifting-flowers}~\ref{item:flower-exists-lift}, either there exists a lift $s' \in N$ of $\check{s}'$, or every $F \in \cF_t$ has a displayed separation $s'$ that is a lift of~$\check{s}'$.
    By \cref{lem:projections-close-then-close}, since $\eps \geq \delta_{\ref{lem:projections-close-then-close}}(q, k, \eps_0)$, every such $s'$ is $\varepsilon$-close to $s$.

    \medskip

    \cref{mainOutroCanonical}. 
    Every $\cF(\check{F}) \in \cG$ stems from \cref{lem:lifting-flowers}, and so is $(q, k)$-coherent.
    Let $\psi$ be an automorphism of~$G$.
    Since $\sigma$ is canonical, $\psi\rest Z$ is an automorphism of~$Z$.
    Since $\check{\cF}$ is canonical by \cref{mainForNarrow}~\ref{mainForNarrowCanonical},
    there is a permutation $\check{\Phi}$ of $\check{\cF}$ such that $\check{F}$ and $\check{\Phi}(\check{F})$ are $\psi \rest Z$-isomorphic for every $\check{F}\in\check{\cF}$.
    By the proof of \cref{highSpreadFlowersDistinct}, distinct $\check F,\check F'\in\check{\cF}$ have different sets of proper displayed separations.
    By \cref{lem:lifting-flowers}~\ref{item:flower-projection} and~\ref{item:flower-exists-lift}, this set is recovered from any lift by projecting its displayed separations and retaining the proper ones.
    Hence the nonempty families $\cF(\check F)$ and $\cF(\check F')$ are disjoint, and thus distinct.
    In particular, $\check F\mapsto\cF(\check F)$ is a bijection $\check\cF\to\cG$.
    Hence $\check{\Phi}$ induces a permutation $\Phi$ of~$\cG$, namely the one which maps $\cF(\check{F})$ to $\cF(\check{\Phi}(\check{F}))$ for every $\check{F}\in\check{\cF}$.
    By \cref{lem:lifting-flowers}~\ref{item:flower-lift-canonical}, the sets $\cF(\check{F})$ and $\cF(\check{\Phi}(\check{F}))=\Phi(\cF(\check{F}))$ are $\psi$-isomorphic for every $\check{F}\in\check{\cF}$.
    Hence $\cG$ is $(q,k)$-canonical.

    Since~$\sigma$ is canonically determined by~$G$, we recover $N=M\sm(\sigma\cup\sigma^*)$ from~$M$, and~$\check N$ by projecting~$N$ to~$Z$.
    Removing the added leaves corresponding to~$\sigma$ recovers~$(\check T,\check\cV)$.
    Hence, by \cref{mainForNarrow}~\ref{mainForNarrowCanonical},
    the bipartition~$\check X\sqcup\check Y=V(\check T)$ and the assignment~$\check Y\ni t\mapsto\check F_t\in\check{\cF}$ are canonical given~$M$ and~$(T,\cV)$.
    Since the added leaves are determined by~$\sigma$, assigning them to~$X$ preserves canonicity of the bipartition.
    Moreover, assigning to each~$t\in Y$ the family $\cF_t=\cF(\check F_t)$ of lifts of~$\check F_t$ is canonical given~$M$ and~$(T,\cV)$ by \cref{lem:lifting-flowers}~\ref{item:flower-lift-canonical}.
    It is bijective, because we have seen above that $\check\cF\to\cG,\check F\mapsto\cF(\check F)$ is a bijection.

    Finally, we show that $M$ and $(T,\cV)$ are $\kappa$-canonical.
    Let $\psi$ be an automorphism of~$G$.
    Since $\sigma$ is canonical, $\psi\rest Z$ is an automorphism of~$Z$.
    The set $\check{N}$ is $\kappa_0$-canonical by \cref{mainForNarrow}~\ref{mainForNarrowCanonical}, so the sets $\psi(\check{N})$ and $\check{N}$ $\kappa_0$-approximate each other.
    We show that $\psi(N)$ $\kappa$-approximates $N$.
    Let $s' \in N$. 
    By \cref{lem:lifting}~\ref{item:projection}, $s'$ is a lift of some $\check{s}' \in \check{N}$. 
    If $\check{s}'$ is not $\kappa_0$-broad, then $\check{s}'$ is $(\kappa_0+k)$-close to the separation $(V(Z), \emptyset)$, say.
    Then, $s'$ and $(V(G), \emptyset)$ are $\delta_{\ref{lem:projections-close-then-close}}(q, k, \kappa_0+k)$-close by \cref{lem:projections-close-then-close}, so $s'$ is not $\kappa$-broad and we have nothing to show.
    If $\check{s}'$ is $\kappa_0$-broad, there exists $\check{s} \in \check{N}$ such that $\psi(\check{s})$ and $\check{s}'$ are $\kappa_0$-close.
    There exists a lift $s$ of $\check{s}$ in $N$ by \cref{lem:lifting}~\ref{item:exists-lift}. 
    Then $\psi(s)$ is a lift of $\psi(\check{s})$, so $\psi(s)$ and $s'$ are $\kappa$-close by \cref{lem:projections-close-then-close}.
    This proves that $\psi(N)$ $\kappa$-approximates $N$. 
    By symmetry, $N$ also $\kappa$-approximates $\psi(N)$, so $N$ is $\kappa$-canonical.
    Since $\sigma$ is canonical and $M = N \cup \sigma \cup \sigma^*$, it follows that $M$ is $\kappa$-canonical.
    Since $(T, \cV)$ is entirely determined by $M$, it is also $\kappa$-canonical.
\end{proof}

\begin{proof}[Proof of \cref{mainTutte}]
    \cref{mainGeneral} implies \cref{mainTutte} by picking a graphic-$k$-flower $F_t\in\cF_t$ for each $t\in Y$.
\end{proof}

\section{Acknowledgements.}
In 2012, Fabian Hundertmark mentioned to the second author the idea to investigate 3-dimensional corner diagrams, yet it is only now that we succeeded in applying this idea in a meaningful way.
We are grateful to Nathan Bowler for spotting a masked error in an early version of our notion of totally-$\eps$-nested.
This led us to distinguish between the two notions totally-$\eps$-nested and narrow-$\eps$-nested.
We thank Vera Traub for asking us at Bertinoro whether our initial decompositions, as in \cref{mainTutte}, could be described in terms of cactus-decompositions.
We thank Joseph Devine for discussions in the initial phase of the project.
We thank Sandra Albrechtsen for the cactus illustration marking the final frontier of our extended abstract.
\smallskip

\section{AI methodology}\label{sec:AIstatement}
We developed the main structural ideas and original proof arguments without AI assistance. A complete draft containing the main results and their proofs was ready in April 2026 without any use of AI. During subsequent revisions, we used ChatGPT (OpenAI) and Claude (Anthropic) to examine and revise mathematical arguments, revise expository text, and proofread the manuscript. 

The main technical revision concerned the parameter dependencies in the proof of \cref{NewerInterlacingPetalStar}. AI identified incompatible breadth requirements in the application of an auxiliary lemma.
The subsequent AI-assisted revision involved modifying the shifting argument in \cref{lem:shift-to-interlace,lem:wiggle-lemma} and strengthening \cref{AlphaInterlacesCloseToInterlaces} so that the required breadth of the star no longer depended on the approximation parameter.

Minor AI-assisted revisions concerned the lifting of flowers in \cref{lem:lifting-flowers} and its use in the canonicity argument or simplifying the argument in \cref{highlyConnectedTorso}. AI also pointed out some small inaccuracies, such as not accounting for the separator size in an approximation bound, and some places where additional arguments were required, such as in \cref{lem:petal-separation-as-meet}.

AI was also used to generate a draft proof of \cref{lem:unbreakable-minor}. 
It proposed proof arguments, identified supporting results in the literature, and drafted the supplementary material in \cref{app:unbreakable-minors}. We checked this material and added it during the revision process.

AI was used to draft and revise portions of the introduction, the detailed proof overview, and explanatory passages throughout the manuscript. This assistance mostly consisted in refining the wording. It was also used to identify typographical errors, inconsistencies in notation, and other imprecisions.

We evaluated the suggestions and reviewed all AI-generated material incorporated into the manuscript.
We checked the mathematical arguments and supporting references and take responsibility for the final text.

\bibliographystyle{alphaurl}
\bibliography{Biblio.bib}

\appendix
\section{Unbreakable minors for iterated decompositions}
\label{app:unbreakable-minors}

In the recursive application of \cref{thm:main-intro} described
in the introduction, we retain the cyclic pieces and continue
within the bags whose nodes lie on no cycle of the decomposition
cactus. This appendix provides the auxiliary lemma used to
represent these bags by unbreakable minors of the original graph.
Working with minors preserves membership in any minor-closed
class, whereas passing to a torso by completing adhesion sets
into cliques need not do so.

A set $C\subseteq V(G)$ is \defn{$(q,k)$-unbreakable in $G$} if every
separation $(A,B)$ of $G$ of order at most $k$ satisfies $\min\bigl\{|C\cap(A\setminus B)|,\ |C\cap(B\setminus A)|\bigr\}\le q.$
A \emph{$C$-rooted minor} of $G$ is a graph $H$ with $C\subseteq V(H)$ that admits a minor model in $G$ in which the branch set corresponding to each $c\in C$ contains $c$.

The following lemma allows us to retain an unbreakable set as a set of distinct roots in an unbreakable minor.

\begin{lemma}[Unbreakable minor]
\label{lem:unbreakable-minor}
There is a computable function
$F\colon\mathbb{N}_0^2\to\mathbb{N}_0$ such that, whenever
$C\subseteq V(G)$ is $(q,k)$-unbreakable in a finite graph $G$, there
is a $C$-rooted minor $H$ of $G$ that is
$(F(q,k),k)$-unbreakable.
\end{lemma}

To justify the recursive application described in the introduction, let $d$ be the approximation bound for a cactus-decomposition given by \cref{thm:main-intro}, and let $C$ be one of its bags whose node lies on no cycle.
Every displayed separation has a side containing $C$, as does every improper separation.
Consequently, every separation at edit-distance at most $d$ from such a separation has at most $d$ vertices of $C$ on one of its proper sides.
Thus $C$ is $(d,k)$-unbreakable in $G$.
By \cref{lem:unbreakable-minor}, there is a $C$-rooted minor $H$ of $G$ that is $(F(d,k),k)$-unbreakable.
We may therefore apply \cref{thm:main-intro} to $H$ with $k$ replaced by $k+1$.

We first prove the existence of the bound using the labelled Graph Minor Theorem \cite[Theorem~1.7]{RS10}, together with replacement and minimality arguments. We then establish its computability.

This background-material appendix was generated by AI and checked by the authors. 

\subsection{Bounded replacements and unbreakability}

Fix $k\in\mathbb{N}_0$. For a graph $J$ with a specified set $R$ of
labelled roots and an ordered partition $(L,Z,U)$ of $R$, let
\[
 \lambda_{J,R}(L,Z,U)
 =
 \min\left\{
 |X\cap Y|:
 \begin{array}{l}
 (X,Y)\text{ is a separation of }J,\\
 R\cap(X\setminus Y)=L,\quad R\cap X\cap Y=Z,\\
 R\cap(Y\setminus X)=U
 \end{array}
 \right\},
\]
where the minimum of the empty set is $\infty$. The \emph{truncated
separation profile} of $(J,R)$ is the function
\[
 \pi^k_{J,R}(L,Z,U)
   =\min\{k+1,\lambda_{J,R}(L,Z,U)\}.
\]
Profiles on different root sets are compared using the root labels.

\begin{lemma}[Bounded rooted replacement]
\label{lem:bounded-rooted-replacement}
For every $r,k\in\mathbb{N}_0$ there is $b(r,k)$ such that every graph
$J$ with at most $r$ labelled roots $R$ has an $R$-rooted minor $J'$
satisfying
\[
 |V(J')|\le b(r,k)
 \qquad\text{and}\qquad
 \pi^k_{J',R}=\pi^k_{J,R}.
\]
The function $b$ can be chosen nondecreasing in its first argument.
\end{lemma}

\begin{proof}
Fix a set of $s$ root labels. Graphs with these roots are
well-quasi-ordered by the rooted minor relation, by the labelled
Graph Minor Theorem \cite[Theorem~1.7]{RS10}. To obtain this rooted
version, give the $s$ roots distinct labels $1,\ldots,s$, and give
every other vertex label $0$. Order this finite label set by equality.
In the resulting minor model, the branch set of each root must contain
the correspondingly labelled root, as required.

There are at most $(k+2)^{3^s}$ truncated separation profiles on these
roots. For each profile, its minimal realizations under the rooted
minor relation form a finite set, up to rooted isomorphism, since
distinct minimal realizations are incomparable. Moreover, every
realization contains a minimal realization of the same profile as a
rooted minor: successive proper rooted minors strictly decrease
$|V|+|E|$. Taking the maximum order of these minimal realizations,
over all profiles and all $s\le r$, gives the desired $b(r,k)$.
\end{proof}

We next record the replacement property of the profiles. Suppose that
$Q=J\cup K$, where $V(J)\cap V(K)=S\subseteq R\subseteq V(J)$, and
that $J'$ has the same truncated profile as $J$ on $R$. Form
$Q'=J'\cup K$ by identifying the corresponding vertices of $S$,
with all other vertices of $J'$ taken to be new. Let $D$ be a set of
distinguished vertices of $Q$ such that $D\cap V(J)\subseteq R$,
and use the corresponding vertices as distinguished vertices of $Q'$.
Then
\begin{equation}
\label{eq:profile-replacement}
 \pi^k_{Q,D}=\pi^k_{Q',D}.
\end{equation}
Indeed, restrict a separation of $Q$ of order at most $k$ to $J$
and $K$. Profile equality supplies a separation of $J'$ of no greater
order with the same assignment of the roots in $R$ to the two
exclusive sides and the separator. It agrees with the separation of
$K$ on $S$, so the two separations can be pasted. If their separator
orders are $a$ and $b$, and $z$ vertices of $S$ lie in the separator,
the pasted separation has order $a+b-z$. Thus its order does not
increase, and the assignment of $D$ is unchanged. Interchanging $J$
and $J'$ proves the reverse implication, giving
\eqref{eq:profile-replacement}. If $J'$ is an $R$-rooted minor of $J$,
the replacement also gives a $D$-rooted minor of $Q$.

\begin{proof}[Proof of the existence assertion in
Lemma~\ref{lem:unbreakable-minor}]
Set $F(q,k)=b(q+k,k)$. Among all $C$-rooted minors of $G$ having the
same truncated separation profile on $C$ as $G$, choose $H$ with the
fewest vertices. Such a minor exists, since $G$ itself is eligible.
Profile equality implies that $C$ remains $(q,k)$-unbreakable in $H$.

Let $(A,B)$ be a separation of $H$ of order at most $k$. By symmetry,
we may assume that $|C\cap(A\setminus B)|\le q$. Put
\[
 S=A\cap B
 \qquad\text{and}\qquad
 R=S\cup(C\cap A).
\]
Then $|R|\le k+q$. By Lemma~\ref{lem:bounded-rooted-replacement},
$H[A]$ has an $R$-rooted minor $J'$ with the same truncated profile
on $R$ and with at most $b(q+k,k)$ vertices. Replace $H[A]$ by $J'$
and retain $H[B]$, identifying the corresponding vertices of $S$.
The resulting graph $H'$ is a $C$-rooted minor of $H$. By
\eqref{eq:profile-replacement}, it has the same truncated profile on
$C$ as $H$, and hence as $G$. The minimality of $H$ now gives
\[
 |V(H)|
 \le |V(H')|
 = |V(H)|-|A|+|V(J')|,
\]
so $|A|\le |V(J')|\le F(q,k)$. In particular,
\[
 \min\{|A\setminus B|,\ |B\setminus A|\}\le F(q,k),
\]
as required.
\end{proof}

\subsection{Computability of the bound}
\label{app:unbreakable-minor-computability}

It remains to show that $b$, and hence $F$, can be chosen computable.
Well-quasi-ordering alone does not make the preceding choice of $b$
effective. We use the following consequence of the irrelevant-vertex
results for rooted topological minors.

\begin{lemma}[Bounded treewidth reduction]
\label{lem:profile-treewidth-reduction}
There is a computable function $w(r,k)$ such that every graph $J$
with $r$ labelled roots $R$ has an induced subgraph $J_0$ containing
$R$ with
\[
 \operatorname{tw}(J_0)\le w(r,k)
 \qquad\text{and}\qquad
 \pi^k_{J_0,R}=\pi^k_{J,R}.
\]
\end{lemma}

\begin{proof}
Put $\delta=r+2(k+1)$. The rooted topological $\delta$-folio of $J$
records its rooted topological minors $P$ with
$|E(P)|+\operatorname{is}(P)\le\delta$, where
$\operatorname{is}(P)$ is the number of isolated vertices of $P$;
see \cite[Definition~4.4]{FLPSZ20}. Rooted topological models map each
root to the vertex with the same label.

We first show that this folio, together with $J[R]$, determines
$\pi^k_{J,R}$. Consider an ordered partition $(L,Z,U)$ of $R$.
If $|Z|>k$, its profile value is $k+1$. An edge between $L$ and $U$
also makes this value $k+1$, since no separation with the prescribed
root assignment exists. Otherwise, let $\tau$ be the minimum size
of a set $T\subseteq V(J)\setminus R$ separating $L$ from $U$ in
$J-Z$. Then
\[
 \lambda_{J,R}(L,Z,U)=|Z|+\tau.
\]
By the vertex-capacitated form of Menger's theorem, $\tau$ is the
maximum number of $L$--$U$ paths with pairwise disjoint sets of
nonroot vertices. The paths may share their endpoints. Trimming
each path between two consecutive root vertices at which it passes
from $L$ to $U$, we may assume that no path has an internal root.

For $1\le t\le k+1$, the existence of $t$ such paths is determined
by the $\delta$-folio. Indeed, each path has an internal vertex,
since there is no edge between $L$ and $U$. Retain one internal
vertex on each path and suppress the other internal vertices.
Include all $r$ roots as labelled vertices, adding unused roots as
isolated vertices. The resulting rooted graph has $2t$ edges and
at most $r$ isolated vertices, so belongs to the $\delta$-folio.
Conversely, a rooted topological model of such a graph supplies
the required paths. Including every root prevents a root from
occurring internally on any of these paths. Thus the folio
determines $\min\{k+1,\tau\}$ and hence the truncated profile.

The irrelevant-vertex results of Fomin, Lokshtanov, Panolan, Saurabh
and Zehavi imply a computable threshold $w=w(r,k)$ such that, if
$\operatorname{tw}(J)>w$, some nonroot vertex can be deleted while
preserving the rooted $\delta$-folio. Specifically, combine
\cite[Lemmas~3.1, 4.1 and~6.1]{FLPSZ20}, as in the proof of
Lemma~8.1 there. These results give a bounded-treewidth, large-clique,
or large-flat-wall alternative, with a nonroot irrelevant vertex
in the latter two cases; they in fact preserve the extended folio.
Their hypothesis on the number of roots holds since
$r\le\delta\le16\delta^2$. All these references are to the full
version linked in the bibliography.

Repeatedly delete such vertices until the treewidth is at most $w$.
Each deletion preserves both the folio and the graph induced by
$R$, and hence the truncated separation profile. This yields $J_0$.
\end{proof}

\begin{proof}[Proof of computability in
Lemma~\ref{lem:unbreakable-minor}]
Fix $r,k$ and a set of $r$ root labels, and compute $w(r,k)$ from
Lemma~\ref{lem:profile-treewidth-reduction}. For each possible
truncated separation profile $\pi$, we compute a finite list
$\mathcal{L}_{\pi}$ of rooted graphs with profile $\pi$.

Start with the empty list. At each step, decide whether there is a
graph of treewidth at most $w(r,k)$ with these roots and profile
$\pi$ that contains no member of the current list as a rooted minor.
This is a decidable monadic second-order satisfiability question on
graphs of bounded treewidth; see \cite[Theorem~2.3]{AGK08}. Indeed,
for each root partition and each $a\le k$, the existence of a
separation of order at most $a$ with that root assignment is
MSO-definable. Boolean combinations of these conditions specify
the exact truncated profile. Containment of a fixed rooted minor
is also MSO-definable, using disjoint connected branch sets with
the prescribed roots in their corresponding branch sets. The root
labels can be represented by existentially quantified distinct
vertices. This is the bounded-treewidth MSO enumeration method
used in \cite[Section~3]{AGK08}.

If the answer is yes, find a witness by enumerating finite rooted
graphs and add it to $\mathcal{L}_{\pi}$. Repeat until the answer
is no. The procedure terminates: otherwise the successive witnesses
would form an infinite sequence in which no earlier member is a
rooted minor of a later member, contrary to labelled minor
well-quasi-ordering. This also handles profiles having no
realization, for which the list remains empty.

At termination, every graph of treewidth at most $w(r,k)$ with
profile $\pi$ contains a member of $\mathcal{L}_{\pi}$ as a rooted
minor. Lemma~\ref{lem:profile-treewidth-reduction} extends this
conclusion to every graph with profile $\pi$: first pass to the
induced subgraph of bounded treewidth and then take the rooted
minor from the list.

Perform this computation for every profile and every root count
$s\le r$, and take the maximum order of the representatives in
the resulting lists. This gives a computable choice of $b(r,k)$
in Lemma~\ref{lem:bounded-rooted-replacement}, nondecreasing in $r$.
Using this choice in the existence proof gives the computable
function $F(q,k)=b(q+k,k)$ asserted in
Lemma~\ref{lem:unbreakable-minor}.
\end{proof}

The minor in Lemma~\ref{lem:unbreakable-minor} can also be found by a
uniform algorithm. Enumerate the $C$-rooted minors of $G$, compute
their truncated separation profiles on $C$, and choose one with
the fewest vertices among those preserving the profile of $G$.
The existence proof shows that this minor has the required
unbreakability. No running-time bound is claimed.
\end{document}